\documentclass[11pt,twoside]{article}
\usepackage{fancyhdr}
\usepackage[colorlinks,citecolor=blue,urlcolor=blue,linkcolor=blue,bookmarks=false]{hyperref}
\usepackage{amsfonts,epsfig,graphicx}
\usepackage{afterpage}
\usepackage{amsmath,amssymb,amsthm} 
\usepackage{fullpage}
\usepackage{epsf} 
\usepackage{graphics} 
\usepackage{amsfonts,amsmath}
\usepackage[sort,numbers]{natbib} 
\usepackage{psfrag,xspace}
\usepackage{color,etoolbox}

\newtheorem{theorem}{Theorem} 
\newtheorem{lemma}{Lemma}
\newtheorem{proposition}{Proposition}
\newtheorem{corollary}{Corollary}
\theoremstyle{definition}
\newtheorem{definition}{Definition}

\usepackage[utf8]{inputenc} 
\usepackage[T1]{fontenc}    
\usepackage{hyperref}       
\usepackage{url}            
\usepackage{booktabs}       
\usepackage{amsfonts}       
\usepackage{nicefrac}       
\usepackage{microtype}      

\usepackage{microtype}
\usepackage{graphicx}
\usepackage{float}
\usepackage[export]{adjustbox}
\usepackage{subcaption}
\usepackage{booktabs}
\usepackage{xcolor}
\usepackage{algorithm}
\usepackage{algorithmic}
\usepackage{enumerate}
\usepackage[shortlabels]{enumitem}
\usepackage{bm}
\usepackage{mathtools}

\DeclareFontFamily{U}{mathx}{\hyphenchar\font45}
\DeclareFontShape{U}{mathx}{m}{n}{<-> mathx10}{}
\DeclareSymbolFont{mathx}{U}{mathx}{m}{n}
\DeclareMathAccent{\wb}{0}{mathx}{"73}

\newcommand{\Reals}{\mathbb{R}}

\newcommand{\abs}[1]{\left \lvert #1 \right \rvert}

\newcommand{\set}[1]{\left\{#1\right\}}

\newcommand{\1}{\mathbf{1}}

\DeclareMathOperator*{\argmin}{argmin}

\newcommand{\Rd}{\Reals^d}

\DeclarePairedDelimiterX{\norm}[1]{\lVert}{\rVert}{#1}
\DeclarePairedDelimiterX{\seminorm}[1]{\lvert}{\rvert}{#1}

\newcommand{\Leb}{L}
\newcommand{\mc}[1]{\mathcal{#1}}
\newcommand{\mbb}[1]{\mathbb{#1}}

\newcommand{\Pbb}{\mathbb{P}}

\newcommand{\Ebb}{\mathbb{E}}
\newcommand{\Qbb}{\mathbb{Q}}

\newcommand{\wt}[1]{\widetilde{#1}}
\newcommand{\wh}[1]{\widehat{#1}}

\newcommand{\dist}{\mathrm{dist}}
\newcommand{\vol}{\mathrm{vol}}
\newcommand{\cut}{\mathrm{cut}}
\newcommand{\diam}{\mathrm{diam}}

\makeatletter
\long\def\@makecaption#1#2{
        \vskip 0.8ex
        \setbox\@tempboxa\hbox{\small {\bf #1:} #2}
        \parindent 1.5em  
        \dimen0=\hsize
        \advance\dimen0 by -3em
        \ifdim \wd\@tempboxa >\dimen0
                \hbox to \hsize{
                        \parindent 0em
                        \hfil 
                        \parbox{\dimen0}{\def\baselinestretch{0.96}\small
                                {\bf #1.} #2
                                } 
                        \hfil}
        \else \hbox to \hsize{\hfil \box\@tempboxa \hfil}
        \fi
        }
\makeatother

\begin{document}

\begin{center} {\Large{\bf{Statistical Guarantees for Local Spectral Clustering on Random Neighborhood Graphs}}}

\vspace*{.3cm}

{\large{
\begin{center}
Alden Green~~~~~ Sivaraman Balakrishnan~~~~~ Ryan J. Tibshirani\\
\vspace{.2cm}
\end{center}

\begin{tabular}{c}
Department of Statistics and Data Science \\
Carnegie Mellon University
\end{tabular}

\vspace*{.2in}

\begin{tabular}{c}
\texttt{\{ajgreen,sbalakri,ryantibs\}@stat.cmu.edu}
\end{tabular}
}}

\vspace*{.2in}

\today
\vspace*{.2in}
\end{center}

\begin{abstract}
	We study the Personalized PageRank (PPR) algorithm, a local spectral method for clustering, which extracts clusters using locally-biased random walks around a given seed node.  In contrast to previous work, we adopt a classical statistical learning setup, where we obtain samples from an unknown nonparametric distribution, and aim to identify sufficiently salient clusters.  We introduce a trio of population-level functionals---the \emph{normalized cut}, \emph{conductance}, and \emph{local spread}, analogous to graph-based functionals of the same name---and prove that PPR, run on a neighborhood graph, recovers clusters with small population normalized cut and large conductance and local spread. We apply our general theory to establish that PPR identifies connected regions of high density (density clusters) that satisfy a set of natural geometric conditions. We also show a converse result, that PPR can fail to recover geometrically poorly-conditioned density clusters, even asymptotically. Finally, we provide empirical support for our theory.
\end{abstract}

\section{Introduction}
In this paper, we consider the problem of clustering: splitting a given data set into groups that satisfy some notion of within-group similarity and between-group difference.  Our particular focus is on spectral clustering methods, a family of powerful nonparametric clustering algorithms. Roughly speaking, a spectral algorithm first constructs a geometric graph $G$, where vertices correspond to samples, and edges correspond to proximities between samples. The algorithm then estimates a feature embedding based on a (suitable) Laplacian matrix of $G$, and applies a simple clustering technique (like $k$-means clustering) in the embedded feature space.

When applied to geometric graphs built from a large number of samples, global spectral clustering methods can be computationally cumbersome and insensitive to the local geometry of the underlying distribution \citep{leskovec2010,mahoney2012}.  This has led to increased interest in \emph{local} spectral clustering algorithms, which leverage locally-biased spectra computed using random walks around some user-specified seed node.  A popular local clustering algorithm is the Personalized PageRank (PPR) algorithm, first introduced by \citet{haveliwala2003}, and then further developed by several others \citep{spielman2011,spielman2014,andersen2006,mahoney2012,zhu2013}.

Local spectral clustering techniques have been practically very successful \citep{leskovec2010,andersen2012,gleich2012,mahoney2012,wu2012}, which has led many authors to develop supporting theory \citep{spielman2013,andersen2009,gharan2012,zhu2013} that gives worst-case guarantees on traditional graph-theoretic notions of cluster quality (such as normalized cut and conductance). In contrast, in this paper we adopt a classical statistical viewpoint, and examine what the output of local clustering on a data set reveals about the underlying density $f$ of the samples. We establish conditions on $f$ under which PPR, when appropriately tuned and initialized inside a candidate cluster $\mc{C} \subseteq \Rd$, will approximately recover this candidate cluster. We pay special attention to the case where $\mc{C}$ is a \emph{density cluster} of $f$---defined as a connected component of the upper level set $\{x \in \Rd : f(x) \geq \lambda\}$ for some $\lambda > 0$---and show precisely how PPR accounts for both geometry and density in estimating a cluster.

Before giving a more detailed overview of our main results, we formally define PPR on a neighborhood graph, review some of the aforementioned worst-case guarantees, and introduce the population-level functionals that govern the behavior of local clustering in our statistical context. 

\subsection{PPR Clustering}
\label{subsec:ppr_clustering}
We start by reviewing the PPR clustering algorithm. Let $G = (V,E)$ be an undirected, unweighted, and connected graph. We denote by $A \in \Reals^{n \times n}$ the adjacency matrix of $G$, with entries $A_{uv} = 1$ if $(u,v) \in E$ and $0$ otherwise.  We also denote by $D$ the diagonal degree matrix, with entries $D_{uu} := \sum_{v \in V} A_{uv}$, and by $I$ the $n \times n$ identity matrix. The PPR vector $p_v = p(v,\alpha;G)$ is defined with respect to a given seed node $v \in V$ and a teleportation parameter $\alpha \in [0,1]$, as the solution of the following linear system:
\begin{equation}
\label{eqn:ppr_vector}
p_v = \alpha e_{v} + (1 - \alpha) p_v W,
\end{equation}
where $W = (I + D^{-1}A)/2$ is the lazy random walk matrix over $G$ and $e_{v}$ is the indicator vector for node $v$ (that has a 1 in position $v$ and 0 elsewhere). 

In practice, exactly solving the system of equations~\eqref{eqn:ppr_vector} to compute the PPR vector may be too computationally expensive. To address this limitation, \citet{andersen2006} introduced the \emph{$\varepsilon$-approximate} PPR vector (aPPR), which we will denote by \smash{$p_v^{(\varepsilon)}$}. We refer the curious reader to \citet{andersen2006} for a formal algorithmic definition of the aPPR vector, and limit ourselves to highlighting a few salient points: the aPPR vector can be computed in order $\mc{O}(1/(\varepsilon\alpha))$ time, while satisfying the following uniform error bound: 
\begin{equation}
\label{eqn:appr_error}
p_v(u) - \varepsilon D_{uu}\leq p_v^{(\varepsilon)}(u) \leq p_v(u), 
\quad \text{for all $u \in V$}.
\end{equation}

Once $p_v$ or \smash{$p_v^{(\varepsilon)}$} is computed, the cluster estimate \smash{$\wh{C}$} is chosen by taking a particular sweep cut. For a given level $\beta > 0$, the \emph{$\beta$-sweep cut} of $p_v = (p_v(u))_{u \in V}$ is 
\begin{equation}
\label{eqn:sweep_cuts}
S_{\beta,v} := \set{u \in V: \frac{p_v(u)}{D_{uu}} > \beta}.
\end{equation}
To define \smash{$\wh{C}$}, one computes $S_{\beta,v}$ over all \smash{$\beta \in (L, U)$} (where the range $(L,U)$ is user-specified), and then outputs the cluster estimate \smash{$\wh{C} = S_{\beta^*,v}$} with minimum normalized cut. For a set $C \subseteq V$ with complement $C^c = V \!\setminus\! C$, the \emph{cut} and \emph{volume} are respectively,
\begin{equation}
\label{eqn:cut_volume}
\cut(C;G) := \sum_{u \in C} \sum_{v \in C^c}
\1\{(u,v) \in E\},~~ \vol(C; G) := \sum_{u \in C}  \sum_{v \in V} \1\{(u,v) \in E\},
\end{equation}
and the \emph{normalized cut} of $C$ is
\begin{equation}
\label{eqn:normalized_cut}
\Phi(C; G) := \frac{\cut(C;G)}{\min \set{\vol(C; G), \vol(C^c; G)}}.
\end{equation} 

\subsection{Worst-Case Guarantees for PPR Clustering}
As mentioned, most analyses of local clustering have focused on worst-case guarantees, defined with respect to functionals of an a priori fixed graph $G = (V,E)$.  For instance, \cite{andersen2006} analyze the normalized cut of the cluster estimate $\wh{C}$ output by PPR, showing that when PPR is appropriately seeded within a candidate cluster $C \subseteq V$, the normalized cut $\Phi(\wh{C};G)$ is upper bounded by (a constant times)  $\sqrt{\Phi(C;G)}$. \cite{zhu2013} build on this: they introduce a second functional, the \emph{conductance} $\Psi(G)$, defined as 
\begin{equation}
\label{eqn:conductance}
\Psi(G) := \min_{S \subseteq V} \Phi(S;G),
\end{equation}
and show that if $\Phi(C;G)$ is much smaller than $\Psi(G[C])^2$---where $G[C] = (C,E \cap (C \times C))$ is the subgraph of $G$ induced by $C$--- then (in addition to having a small normalized cut) the cluster estimate $\wh{C}$ approximately recovers $C$. Our own analysis builds on that of~\cite{zhu2013}, and we present a more detailed summary of their results in Section~\ref{sec:ub_symmetric_set_difference}. For now, we merely reiterate that the conclusions of \citet{andersen2006,zhu2013} cannot be straightforwardly applied to our setting, where the input data are random samples $\{x_1,\ldots,x_n\}$ drawn from a distribution $\Pbb$, the graph $G$ is a random neighborhood graph formed from the samples, and the candidate cluster is a set $\mc{C} \subseteq \Rd$.\footnote{Throughout, we use calligraphic notation to refer to subsets of $\Rd$.}

\subsection{PPR on a Neighborhood Graph}
\label{subsec:ppr_neighborhood_graph}
We now formally describe the statistical setting in which we operate, as well as the method we will study: PPR on a neighborhood graph. Let $X = \{x_1,\ldots, x_n\}$ be samples drawn i.i.d.\ from a distribution $\Pbb$ on $\Rd$. We will assume throughout that $\Pbb$ has a density $f$ with respect to the Lebesgue measure $\nu$ on $\Rd$ . For a radius $r > 0$, we define $G_{n,r}=(V,E)$ to be the \emph{$r$-neighborhood graph} of $X$, an unweighted, undirected graph with vertices $V=X$, and an edge $(x_i,x_j) \in E$ if and only if $i \neq j$ and $\|x_i - x_j\| \leq r$, where $\|\cdot\|$ is the Euclidean norm. Once the neighborhood graph $G_{n,r}$ is formed, the PPR vector $p_v$ is then computed over $G_{n,r}$, with a resulting cluster estimate \smash{$\wh{C} \subseteq X$}. The precise algorithm is summarized in Algorithm~\ref{alg:ppr}.   

\begin{algorithm}
	\caption{PPR on a neighborhood graph}
	\label{alg:ppr}	
	{\bfseries Input:} data $X=\{x_1,\ldots,x_n\}$, radius $r > 0$, teleportation
	parameter $\alpha \in [0,1]$, seed $v \in X$, sweep cut range $(L,U)$. \\     
	{\bfseries Output:} cluster estimate $\wh{C} \subseteq V$.
	\begin{algorithmic}[1]
		\STATE Form the neighborhood graph $G_{n,r}$.
		\STATE Compute the PPR vector $p_v=p(v, \alpha; G_{n,r})$ as in
		\eqref{eqn:ppr_vector}.  
		\STATE Compute sweep cuts $S_{\beta}$ as in \eqref{eqn:sweep_cuts}, for each $\beta \in (L,U)$.\footnotemark 
		\STATE Return \smash{$\wh{C} = S_{\beta^*}$}, where  
		$$
		\beta^* = \argmin_{\beta \in (L,U)}~ \Phi(S_{\beta}; G_{n,r}).
		$$
	\end{algorithmic}
\end{algorithm}
\footnotetext{Technically speaking, for each $\beta \in (L,U) \cap \{p_v(u)/D_{uu}: u \in V\}$.}

\subsection{Cluster Accuracy}
We need a metric to assess the accuracy with which $\wh{C}$ estimates the candidate cluster $\mc{C}$. One commonly used metric is the misclassification error, i.e., the size of the symmetric set difference between \smash{$\wh{C}$} and the empirical cluster $\mc{C}[X] = \mc{C} \cap X$ \citep{korostelev1993,polonik1995,rigollet2009}. We will consider a related metric, the volume of the symmetric set difference, which weights misclassified points according to their degree in $G_{n,r}$. To keep things simple, for a given set $S \subseteq X$ we write $\vol_{n,r}(S) := \vol(S;G_{n,r})$. 

\begin{definition}
	\label{def:volume_symmetric_set_difference}
	For an estimator \smash{$\wh{C} \subseteq X$} and a set $\mc{C} \subseteq \Rd$, their symmetric set difference is  
	\begin{equation*}
	\wh{C} \vartriangle \mc{C}[X] :=
	\bigl(\wh{C} \setminus \mc{C}[X]\bigr) \cup
	\bigl(\mc{C}[X] \setminus \wh{C}\bigr).
	\end{equation*}
	Furthermore, we denote the volume of the symmetric set difference by 
	$$
	\Delta(\wh{C}, \mc{C}[X]) := \vol_{n,r}(\wh{C} \vartriangle \mc{C}[X]). 
	$$
\end{definition}

\subsection{Population Normalized Cut, Conductance, and Local Spread}
Next we define three population-level functionals of $\mc{C}$---the normalized cut $\Phi_{\Pbb,r}(\mc{C})$, conductance $\Psi_{\Pbb,r}(\mc{C})$, and the local spread $s_{\Pbb,r}(\mc{C})$---which we will use to upper bound the volume of the symmetric set difference \smash{$\Delta(\wh{C},\mc{C}[X])$} with high probability.

Let the population-level \emph{cut} of $\mc{C}$ be the expectation (up to a rescaling) of $\cut_{n,r}(\mc{C}[X]) := \cut(\mc{C}[X]; G_{n,r})$,  and likewise let the population-level \emph{volume} of $\mc{C}$ be the expectation (up to a rescaling) of $\vol_{n,r}(\mc{C}[X]) := \vol(\mc{C}[X]; G_{n,r})$; i.e., let
\begin{equation*}
\mathrm{cut}_{\Pbb,r}(\mc{C}) := \int_{\mc{C}} \int_{\mc{C}^c} \1\{\|x - y\| \leq r\} \,d\Pbb(y) \,d\Pbb(x),~~ \mathrm{vol}_{\Pbb,r}(\mc{C}) := \int_{\mc{C}} \int_{\Rd} \1\{\|x - y\| \leq r\} \,d\Pbb(y) \,d\Pbb(x),
\end{equation*}
where $\mc{C}^c := \Rd \!\setminus\! \mc{C}$. Also let $\deg_{\Pbb,r}(x) := \int_{\Rd} \1\{\|y - x\| \leq r\} \,d\Pbb(y)$ be the expected degree of $x$ in $G_{n,r}$. 
\begin{definition}[Population normalized cut]
	For a set $\mc{C} \subseteq \Rd$, distribution $\Pbb$, and radius $r > 0$, the \emph{population normalized cut} is
	\begin{equation}
	\label{eqn:population_normalized_cut}
	\Phi_{\Pbb,r}(\mc{C}) := \frac{\mathrm{cut}_{\Pbb,r}(\mc{C})}{\min\{\mathrm{vol}_{\Pbb,r}(\mc{C}), \mathrm{vol}_{\Pbb,r}(\mc{C}^c)\}}.
	\end{equation}
\end{definition}

Let \smash{$\wt{\Pbb}(\cdot) = \Pbb(\cdot|x \in \mc{C})$} be the conditional distribution of $x$, i.e., let \smash{$\wt{\Pbb}(\mc{S}) = \wt{\Pbb}(\mc{S} \cap \mc{C})/\wt{\Pbb}(\mc{C})$} for measurable sets $\mc{S} \subseteq \Rd$.

\begin{definition}[Population conductance]
	For a set $\mc{C} \subseteq \Rd$, distribution $\Pbb$ and radius $r > 0$, the \emph{population conductance} is
	\begin{equation}
	\label{eqn:population_conductance}
	\Psi_{\Pbb,r}(\mc{C}) = \inf_{\mc{S} \subseteq \mc{C}} \Phi_{\wt{\Pbb},r}(\mc{S}).
	\end{equation}
\end{definition}

\begin{definition}[Population local spread]
	For a set $\mc{C} \subseteq \Rd$, distribution $\Pbb$ and radius $r > 0$, the \emph{population local spread} is
	\begin{equation}
	\label{eqn:local_spread}
	s_{\Pbb,r}(\mc{C}) := \min_{x \in \mc{C}} \biggl\{\frac{\bigl(\deg_{\wt{\Pbb},r}(x)\bigr)^2}{\vol_{\wt{\Pbb},r}(\mc{C})} \biggr\},
	\end{equation}
\end{definition}

It is quite natural that $\Phi_{\Pbb,r}(\mc{C})$ and $\Psi_{\Pbb,r}(\mc{C})$ should help quantify the role geometry plays in local spectral clustering. Indeed, these functionals are the population-level analogues of the empirical quantities $\Phi_{n,r}(\mc{C}[X]) := \Phi(\mc{C}[X];G_{n,r})$ and $\Psi_{n,r}(\mc{C}[X]) := \Psi(G_{n,r}\bigl[\mc{C}[X]\bigr])$, and as we have already mentioned, these empirical quantities can be used to upper bound the volume of the symmetric set difference. For this reason, similar population-level functionals are used by \cite{shi2009,schiebinger2015,garciatrillos19} in the analysis of \emph{global} spectral clustering in a statistical context. We will comment more on the relationship between these works and our own results in Section~\ref{subsec:related_work}. 

The role played by $s_{\Pbb,r}(\mc{C})$ is somewhat less obvious. For now, we mention only that it plays an essential part in obtaining tight bounds on the mixing time of a particular random walk that is closely related to the PPR vector, and defer further discussion until later in Section~\ref{sec:ub_symmetric_set_difference}.

\subsection{Main Results}
We now informally state our two main upper bounds, regarding the recovery of a generic cluster $\mc{C}$, and a density cluster $\mc{C}_{\lambda}$. Theorem~\ref{thm:volume_ssd_ub_informal} informally summarizes the first of our main results (formally stated in Theorem~\ref{thm:volume_ssd_ub}) regarding the recovery of a generic cluster $\mc{C}$.
\begin{theorem}[Informal]
	\label{thm:volume_ssd_ub_informal}
	If $\mc{C} \subseteq \Rd$ and $\Pbb$ satisfy appropriate regularity conditions, and Algorithm~\ref{alg:ppr} is initialized properly with respect to $\mc{C}$, then for all sufficiently large $n$, with high probability it holds that
	\begin{equation*}
	\frac{\Delta(\wh{C},\mc{C}[X])}{\vol_{n,r}(\mc{C}[X])} \leq c \cdot\Phi_{\Pbb,r}(\mc{C}) \cdot  \biggl(\frac{\log(1/s_{\Pbb,r}(\mc{C}))}{\Psi_{\Pbb,r}(\mc{C})}\biggr)^2.
	\end{equation*}
\end{theorem}

(Above, and throughout, $c$ stands for a universal constant that may change from line to line.) Put more succinctly, we find that \smash{$\Delta(\wh{C},\mc{C}[X])$} is small when $\Phi_{\Pbb,r}(\mc{C})$ is small relative to $(\Psi_{\Pbb,r}(\mc{C})/\log(1/s_{\Pbb,r}(\mc{C})))^2$. To the best of our knowledge, this gives the first population-level guarantees for local clustering in the nonparametric statistical context.

Next, Theorem~\ref{thm:density_cluster_volume_ssd_ub_informal} informally summarizes the second of our main results (formally stated in Theorem~\ref{thm:density_cluster_volume_ssd_ub}) regarding the recovery of a $\lambda$-density cluster $\mc{C}_{\lambda}$ by PPR. For reasons that we explain later in Section~\ref{sec:ppr_density_cluster}, our cluster recovery statement will actually be with respect to the $\sigma$-thickened set $\mc{C}_{\lambda,\sigma} := \{x \in \Rd: \mathrm{dist}(x,\mc{C}_{\lambda}) < \sigma\}$, for a given $\sigma > 0$. The upper bound we establish is a function of various parameters that measure the conditioning of both the density cluster $\mc{C}_{\lambda,\sigma}$ and density $f$ for recovery by PPR. We assume that $\mc{C}_{\lambda,\sigma}$ is the image of a convex set $\mc{K}$ of finite diameter $\mathrm{diam}(\mc{K}) \leq \rho < \infty$ under a Lipschitz, measure-preserving mapping $g$, with Lipschitz constant $M$. We also assume that $f$ is bounded away from $0$ and $\infty$ on $\mc{C}_{\lambda,\sigma}$:
\begin{equation*}
0 < \lambda_{\sigma} \leq f(x) \leq \Lambda_{\sigma} < \infty~~\text{for all $x \in \mc{C}_{\lambda,\sigma}$},
\end{equation*}
and additionally satisfies the following low-noise condition:
\begin{equation*}
\inf_{y \in \mc{C}_{\lambda,\sigma}} f(y) - f(x) \geq  \theta \cdot \dist(x, \mc{C}_{\lambda,\sigma})^{\gamma}~~\text{for all $x$ such that $0 < \dist(x,\mc{C}_{\lambda,\sigma}) \leq \sigma$.}
\end{equation*}
(Here $\dist(x,\mc{C}) := \inf_{y \in \mc{C}} \|y - x\|$.) 
\begin{theorem}[Informal]
	\label{thm:density_cluster_volume_ssd_ub_informal}
	If $\mc{C}_{\lambda} \subseteq \Rd$ is a $\lambda$-density cluster of a distribution $\Pbb$, which satisfies appropriate regularity conditions, and Algorithm~\ref{alg:ppr} is initialized properly with respect to $\mc{C}_{\lambda,\sigma}$, then for all sufficiently large $n$, with high probability it holds that
	\begin{equation*}
	\frac{\Delta(\wh{C},\mc{C}_{\lambda,\sigma}[X])}{\vol_{n,r}(\mc{C}_{\lambda,\sigma})} \leq c \cdot d^4 \cdot \frac{M^2\rho^2}{\sigma r} \cdot \frac{\Lambda_{\sigma}^2 \lambda (\lambda - \theta \frac{r^{\gamma}}{\gamma + 1})}{\lambda_{\sigma}^4} \cdot \log^2\biggl(\frac{\Lambda_{\sigma}^{2/d} M\rho}{\lambda_{\sigma}^{2/d}2r}\biggr).
	\end{equation*}
\end{theorem}

The above result reveals the separate roles played by geometry and density in the ability of PPR to recover a density cluster. Here $M$, $\rho$, and $\sigma$ capture whether $\mc{C}_{\lambda,\sigma}$ is geometrically well-conditioned (short and fat) or poorly-conditioned (long and thin) for recovery by PPR. Likewise, the parameters $\lambda_{\sigma}, \Lambda_{\sigma}, \gamma$, and $\theta$ measure whether the density $f$ is well-conditioned (approximately uniform over the density cluster, and having thin tails outside of it) or poorly conditioned (vice versa). Theorem~\ref{thm:density_cluster_volume_ssd_ub_informal} says that if the thickened density cluster $\mc{C}_{\lambda,\sigma}$ is geometrically well-conditioned---meaning, $M^2\rho^2/(\sigma r) \approx 1$---and the density $f$ is well-conditioned near $\mc{C}_{\lambda,\sigma}$---meaning, $\Lambda_{\sigma} \approx \lambda \approx \lambda_{\sigma}$ and $\lambda - \theta r^{\gamma}/(\gamma + 1)$ is much less than $\lambda_{\sigma}$---then PPR will approximately recover $\mc{C}_{\lambda,\sigma}$.

\subsection{Related Work}
\label{subsec:related_work}
We now summarize some related work (in addition to the background material already given above), regarding the theory of spectral clustering, and of density cluster recovery.

\subsubsection{Spectral Clustering} 
In the stochastic block model (SBM), arguably one of the simplest models of network formation, edges between nodes independently occur with probability based on a latent community membership. In the SBM, the ability of spectral algorithms to perform clustering---or community detection---is well-understood, dating back to \citet{mcsherry2001} who gives conditions under which the entire community structure can be recovered. In more recent work, \citet{rohe2011} upper bound the fraction of nodes misclassified by a spectral algorithm for the high-dimensional (large number of blocks) SBM, and \citet{lei2015} extend these results to the sparse (low average degree) regime. Relatedly, \citet{clauset08,balakrishnan2011,li2018} analyze the misclassification rate when the block model exhibits some hierarchical structure. The framework we consider, in which nodes correspond to data points sampled from an underlying density, and edges between nodes are formed based on geometric proximity, is quite different than the SBM, and therefore so is our analysis.

In general, the study of spectral algorithms on neighborhood graphs has been focused on establishing asymptotic convergence of eigenvalues and eigenvectors of certain sample objects to the eigenvalues and eigenfunctions of corresponding limiting operators. \citet{koltchinskii2000} establish convergence of spectral projections of the adjacency matrix to a limiting integral operator, with similar results obtained using simplified proofs in \citet{rosasco10}. \citet{vonluxburg2008} studies convergence of eigenvectors of the Laplacian matrix for a neighborhood graph of fixed radius. \citet{belkin07} and \citet{garciatrillos18} extend these results to the regime where the radius $r \to 0$ as $n \to \infty$. 

These results are of fundamental importance. However, they remain silent on the following natural question: do the spectra of these continuum operators induce a partition of the sample space which is ``good'' in some sense? \citet{shi2009,schiebinger2015,garciatrillos19,hoffmann2019} address this question, showing that spectral algorithms will recover the latent labels in certain well-conditioned nonparametric mixture models. These works are probably the most similar to our own: the conditioning of these mixture models depends on population-level functionals resembling the population normalized cut and conductance introduced above, and the resulting bounds on the error of spectral clustering are comparable to those we establish in Theorem~\ref{thm:volume_ssd_ub}. However, these results focus on global rather than local methods, and impose global rather than local conditions on $\Pbb$. Moreover, they do not explicitly consider recovery of density clusters, which is an important concern of our work. We comment further on the relationship between our results and these works after Theorem~\ref{thm:volume_ssd_ub}.

\subsubsection{Density Clustering} 
For a given threshold $\lambda \in (0,\infty)$, we denote by \smash{$\mbb{C}_f(\lambda)$} the connected components of the density upper level set $\{x \in \Rd: f(x) \geq \lambda\}$. In the density clustering problem, initiated by~\cite{hartigan1975}, the goal is to recover $\mbb{C}_{f}(\lambda)$. By now, density clustering, and the related problem of level-set estimation, have been thoroughly studied. For instance, \citet{polonik1995,rigollet2009, rinaldo2010, steinwart2015} study density clustering under the symmetric set difference metric, \citet{tsybakov1997,singh2009,jiang2017} describe minimax optimal level-set and cluster estimators under Hausdorff loss, and \citet{hartigan1981,chaudhuri2010,kpotufe11,balakrishnan2013,steinwart2017,wang2019} consider estimation of the cluster tree $\{\mbb{C}_f(\lambda): \lambda \in (0,\infty)\}$.

We emphasize that our goal is not to improve on these results, nor is it to offer a better algorithm for density clustering. Indeed, seen as a density clustering algorithm, PPR has none of the optimality guarantees found in the aforementioned works. Rather, we hope to better understand the implications of our general theory by applying it within an already well-studied framework. We should also note that since we study a local algorithm, our interest will be in a local version of the density clustering problem, where the goal is to recover a single density cluster $\mc{C}_{\lambda} \in \mbb{C}_f(\lambda)$. 

\subsection{Organization}
We now outline the rest of the paper.
\begin{itemize}
	\item In Section~\ref{sec:ub_symmetric_set_difference}, we derive bounds on the error of PPR as a function of sample normalized cut, conductance, and local spread. We then show that under certain conditions the sample normalized cut, conductance, and local spread are close to their population-level counterparts, with high probability for sufficient number of samples. As a result, we obtain an upper bound on \smash{$\Delta(\wh{C},\mc{C}[X])/\vol_{n,r}(\mc{C}[X])$} purely in terms of these population-level functionals (Theorem~\ref{thm:density_cluster_volume_ssd_ub}).
	\item In Section~\ref{sec:ppr_density_cluster}, we focus on the special case where the candidate cluster $\mc{C} = \mc{C}_{\lambda}$ is a \emph{$\lambda$-density cluster}---that is, a connected component of the upper level set $\{x: f(x) \geq \lambda\}$. We derive bounds on the population normalized cut, conductance, and local spread of the density cluster, which depend on $\lambda$ as well as some other natural parameters. This leads to a bound on the symmetric set difference between \smash{$\wh{C}$} and the $\lambda$-density cluster (Theorem~\ref{thm:density_cluster_volume_ssd_ub}).
	\item In Section~\ref{sec:lower_bound}, we prove a negative result: we give a hard distribution $\Pbb$ with corresponding density cluster $\mc{C}_{\lambda}$ for which the symmetric set difference between \smash{$\wh{C}$} and the $\lambda$-density cluster is provably large.
	\item In Section~\ref{sec:experiments} we empirically investigate some of our conclusions, before ending with some discussion in Section~\ref{sec:discussion}.
\end{itemize}

\section{Recovery of a Generic Cluster with PPR}
\label{sec:ub_symmetric_set_difference}

In the main result (Theorem~\ref{thm:volume_ssd_ub}) of this section, we give a high probability upper bound on the volume of the symmetric set difference \smash{$\Delta(\wh{C}, \mc{C}[X])$}, in terms of the population normalized cut $\Phi_{\Pbb,r}(\mc{C})$, conductance $\Psi_{\Pbb,r}(\mc{C})$, and local spread $s_{\Pbb,r}(\mc{C})$. We build to this theorem slowly, giving new structural results in two distinct directions. First, we build on some previous work (mentioned in the introduction) to relate \smash{$\Delta(\wh{C}, \mc{C}[X])$} to the sample normalized cut $\Phi_{n,r}(\mc{C}[X])$, conductance $\Psi_{n,r}(\mc{C}[X])$, and local spread $s_{n,r}(\mc{C}[X]) := s(G_{n,r}[\mc{C}[X]])$. Second, we argue that when $n$ is large, each of these graph functionals can be bounded by their population-level analogues with high probability.

\subsection{The Fixed Graph Case}
\label{subsec:ppr_cluster_recovery_fixed_graph}
When PPR is run on a fixed graph $G = (V,E)$ with the goal of recovering a candidate cluster $C \subseteq V$, \cite{zhu2013} provide the sharpest known bounds on the volume of the symmetric set difference between the cluster estimate \smash{$\wh{C}$} and candidate cluster $C$. Since these results will play a major part in our analysis, in Lemma~\ref{lem:zhu} we restate them for the convenience of the reader.\footnote{Lemma~\ref{lem:zhu} improves on Lemma~3.4 of~\cite{zhu2013} by some constant factors, and for completeness we prove Lemma~\ref{lem:zhu} in the Appendix. Nevertheless, to be clear the essential idea of Lemma~\ref{lem:zhu} is no different than that of~\cite{zhu2013}, and we do not claim any novelty.}

In their most general form, the results of~\citet{zhu2013} depend on the mixing time of a lazy random walk over the induced subgraph $G[C]$. The \emph{mixing time} of a lazy random walk over a graph $G$ is
\begin{equation}
\label{eqn:mixing_time}
\tau_{\infty}(G) := \min\set{ t: \frac{{\pi}(u) - {q}_{v}^{(t)}(u)}
	{{\pi}(u)} \leq \frac{1}{4}, \; \text{for all $u,v \in V$}};
\end{equation}
here $q_v^{(t)} := e_v W^t$ is the distribution of a lazy random walk over $G$ initialized at node $v$ and run for $t$ steps, and $\pi := \lim_{t \to \infty} q_v^{(t)}$ is the limiting distribution of $q_v^{(t)}$.

\begin{lemma}[Lemma~3.4 of \cite{zhu2013}]
	\label{lem:zhu}
	For a set $C \subseteq V$, suppose that
	\begin{equation}
	\label{eqn:zhu_condition}
	\alpha \leq \min\Bigl\{\frac{1}{45}, \frac{1}{2\tau_{\infty}(G[C])}\Bigr\},~~ \beta \leq \frac{1}{5\vol(C;G)}.
	\end{equation}
	Then there exists a set $C^g \subseteq C$ with $\vol(C^g;G) \geq \frac{1}{2}\vol(C;G)$ such that for any $v \in C^g$, the sweep cut $S_{\beta,v}$ satisfies
	\begin{equation}
	\label{eqn:zhu_ub}
	\vol(S_{\beta,v} \vartriangle C;G) \leq 6\frac{\Phi(C;G)}{\alpha \beta}.
	\end{equation}
\end{lemma}
The upper bound in~\eqref{eqn:zhu_ub} does not obviously depend on the conductance $\Psi(G[C])$. However, as \cite{zhu2013} point out, letting $\pi_{\min}(G) := \min_{u \in V}\{\pi(u)\}$, it follows from Cheeger's inequality~\citep{chung1997} that 
\begin{equation}
\label{eqn:mixing_time_cheeger}
\tau_{\infty}(G) \leq \frac{\log(1/\pi_{\min}(G))}{\Psi(G)^2}.
\end{equation}
Therefore, setting (for instance) \smash{$\alpha = \frac{\Psi(G[C])^2}{2\log(1/\pi_{\min}(G))}$} and \smash{$\wh{C} = S_{\beta_0,v}$} for \smash{$\beta_0 = \frac{1}{5 \vol(C;G)}$}, we obtain from~\eqref{eqn:zhu_ub} that 
\begin{equation}
\label{eqn:zhu_ub2}
\frac{\vol(C \vartriangle \wh{C};G)}{\vol(C; G)} \leq 60\frac{\Phi(C;G) \log\bigl( 1/\pi_{\min}(G[C])\bigr)}{\Psi(G[C])^2}.
\end{equation}

\subsection{Improved Bounds on Mixing Time} 
\label{subsec:mixing_time}
Having reviewed the conclusions of~\cite{zhu2013}, we return now to our own setting, where the data is not a fixed graph $G$ but instead random samples $\{x_1,\ldots,x_n\}$, and our goal is to recover a candidate cluster $\mc{C} \subseteq \Rd$. Ideally, we would like to apply~\eqref{eqn:zhu_ub2} with $C = \mc{C}[X]$ and $G = G_{n,r}$, replace $\Phi_{n,r}(\mc{C}[X])$ and $\Psi_{n,r}(\mc{C}[X])$ by $\Phi_{\Pbb,r}(\mc{C})$ and $\Psi_{\Pbb,r}(\mc{C})$ inside~\eqref{eqn:zhu_ub2}, and thereby obtain an upper bound on \smash{$\Delta(\wh{C};\mc{C}[X])$} that depends only on $\Pbb$ and $\mc{C}$. Unfortunately, however, there is a catch: when the graph is $G = G_{n,r}$ and the candidate cluster is $C = \mc{C}[X]$, as $n \to \infty$ the sample normalized cut $\Phi_{n,r}(\mc{C}[X])$ and conductance $\Psi_{n,r}(\mc{C}[X])$ each converge to their population-level analogues, but $\pi_{\min}(G_{n,r}[\mc{C}[X]]) \asymp 1/n$.\footnote{For sequences $(a_n)$ and $(b_n)$, we say $a_n \asymp b_n$ if there exists a constant $c \geq 1$ such that $a_n/c \leq b_n \leq c a_n$ for all $n \in \mbb{N}$.} Therefore the right hand side of~\eqref{eqn:zhu_ub2} diverges at a $\log n$ rate, which turns~\eqref{eqn:zhu_ub2} into a vacuous upper bound whenever the number of samples is sufficiently large.

To address this, in Proposition~\ref{prop:pointwise_mixing_time} we improve the upper bound on the mixing time in~\eqref{eqn:mixing_time_cheeger}. Specifically, in~\eqref{eqn:pointwise_mixing_time} the ``start penalty'' of $\log(1/\pi_{\min}(G))$ is replaced by $\log(1/s(G))$, where $s(G)$ is the graph local spread, defined as
\begin{equation*}
s(G) := d_{\min}(G) \cdot \pi_{\min}(G),
\end{equation*}
for $d_{\min}(G) = \min_{u \in V}\bigl\{\deg(u;G)\bigr\}$, and likewise $d_{\max}(G) = \max_{u \in V}\bigl\{\deg(u;G)\bigr\}$. Notice that $s(G) \geq \pi_{\min}(G)$.
\begin{proposition}
	\label{prop:pointwise_mixing_time}
	Assume $d_{\max}(G)/d_{\min}(G)^2 \leq 1/16$. Then,
	\begin{equation}
	\label{eqn:pointwise_mixing_time}
	\tau_{\infty}(G) \leq \frac{17}{\ln(2)} \biggl(\frac{\ln\bigl(32/s(G)\bigr)}{\Psi(G)}\biggr)^2.
	\end{equation}
\end{proposition}

While Proposition~\ref{prop:pointwise_mixing_time} can be applied to \emph{any} graph $G$ (provided that the ratio of maximum degree to squared minimum degree is at most $1/16$), it is particularly useful for geometric graphs: when $G = G_{n,r}[\mc{C}[X]]$ for a fixed radius $r > 0$, we have $d_{\min}(G_{n,r}[\mc{C}[X]]) \asymp n$, and thus $s(G_{n,r}[\mc{C}[X]]) \asymp 1$. We give a precise upper bound on $s(G_{n,r}[\mc{C}[X]])$ in Proposition~\ref{prop:sample_to_population_1}, which does not grow with $n$, and in combination with Proposition~\ref{prop:pointwise_mixing_time} this allows us to remove the unwanted $\log n$ factor from the upper bound in~\eqref{eqn:zhu_ub2}. 

The local spread $s(G)$ plays an intuitive role in the analysis of mixing time. Indeed, in any graph $G$ sufficiently small sets are expanders---that is, if a set $R \subseteq V$ has cardinality less than the minimum degree, the normalized cut $\Phi(R;G)$ will be much larger than the conductance $\Psi(G)$. As a consequence, a random walk over $G$ will rapidly mix over all small sets $R$, and in our analysis of the mixing time we may therefore ``pretend'' that the random walk was given a warm start over a larger set $S$. The local spread $s(G)$ simply delineates small sets $R$ from larger sets $S$. Of course, the proof of Proposition~\ref{prop:pointwise_mixing_time} requires a much more intricate analysis, and---as with the proofs of all results in this paper---it is deferred to the appendix.

\subsection{Sample-to-Population Results}
\label{subsec:sample_to_population}
In Propositions~\ref{prop:sample_to_population_1} and~\ref{prop:sample_to_population_2}, we establish high probability bounds on the sample normalized cut, conductance, and local spread in terms of their population-level analogues. To establish these bounds, we impose the following regularity conditions on \smash{$\wt{\Pbb}$} and $\mc{C}$.
\begin{enumerate}[label=(A\arabic*)]
	\item 
	\label{asmp:bounded_density} 
	The distribution \smash{$\wt{\Pbb}$} has a density \smash{$\wt{f}: \mc{C} \to (0,\infty)$} with respect to Lebesgue measure. There exist $0 < f_{\min} \leq f_{\max} < \infty$ for which
	\begin{equation*}
	(\forall x \in \mc{C})~~ f_{\min} \leq \wt{f}(x) \leq f_{\max}.
	\end{equation*}
	For convenience, we will assume $f_{\min} \leq 1$ and $f_{\max} \geq 1$.
	\item 
	\label{asmp:domain} 
	The candidate cluster $\mc{C} \subseteq \Rd$ is a bounded, connected, open set, and for $d \geq 2$, it has a Lipschitz boundary $\partial \mc{C}$, meaning it is locally the graph of a Lipschitz function (e.g., see Definition 9.57 of~\cite{leoni2017}).
\end{enumerate}

In what follows, we use $b_1,b_2,\ldots$ and $B_1,B_2,\ldots$ to refer to positive constants that may depend on $\Pbb$, $\mc{C}$, and $r$, but do not depend on $n$ or $\delta$. We explicitly keep track of all constants in our proofs.

\begin{proposition}
	\label{prop:sample_to_population_1}
	Fix $\delta \in (0,1/3)$. Suppose \smash{$\wt{\Pbb}$} and $\mc{C}$ satisfy~\ref{asmp:bounded_density} and~\ref{asmp:domain}. Then each of the following statements hold.
	\begin{itemize}
		\item With probability at least $1 - 3\exp\{-b_1\delta^2n\}$,
		\begin{equation}
		\label{eqn:sample_to_population_normalized_cut}
		\Phi_{n,r}(\mc{C}[X]) \leq (1 + 3\delta) \Phi_{\Pbb,r}(\mc{C}).
		\end{equation}
		\item For any $n \in \mbb{N}$ for which 
		\begin{equation}
		\label{eqn:sample_to_population_local_spread_sample_complexity}
		\frac{1}{n} \leq \delta \cdot \frac{2\Pbb(\mc{C})}{3},
		\end{equation}
		the following inequality holds with probability at least \smash{$1 - (n + 2)\exp\{-b_2\delta^2n\}$}:
		\begin{equation}
		\label{eqn:sample_to_population_local_spread}
		s_{n,r}(\mc{C}[X]) \geq (1 - 4\delta) s_{\Pbb,r}(\mc{C}).
		\end{equation}
	\end{itemize}
\end{proposition}
Let $p_d := 1/2$ if $d = 1$, $p_d := 3/4$ if $d = 2$, and otherwise $p_d := 1/d$ for $d \geq 3$. 
\begin{proposition}
	\label{prop:sample_to_population_2}
	Fix $\delta \in (0,1/2)$. Suppose \smash{$\wt{\Pbb}$} and $\mc{C}$ satisfy~\ref{asmp:bounded_density} and~\ref{asmp:domain}. For any $n \in \mbb{N}$ satisfying
	\begin{equation}
	\label{eqn:sample_to_population_conductance_sample_complexity}
	B_1 \frac{(\log n)^{p_d}}{\min\{n^{1/2},n^{1/d}\}} \leq \delta,
	\end{equation}
	the following inequality holds with probability at least $1 - B_2/n - (n + 1)\exp\{-b_3n\}$:
	\begin{equation}
	\label{eqn:sample_to_population_conductance}
	\Psi_{n,r}(\mc{C}[X]) \geq (1 - 2\delta) \Psi_{\Pbb,r}(\mc{C}).
	\end{equation}
\end{proposition}

A note on the proof techniques: the upper bound in~\eqref{eqn:sample_to_population_normalized_cut} follows by applying Bernstein's inequality to control the deviations of $\cut_{n,r}(\mc{C}[X])$, $\vol_{n,r}(\mc{C}[X])$, and $\vol_{n,r}(\mc{C}^c[X])$ around their expectations (noting that each of these is an order-$2$ U-statistic). To prove the lower bound~\eqref{eqn:sample_to_population_local_spread}, we require a union bound to control the minimum degree $d_{\min}(G_{n,r}[\mc{C}[X]])$, but otherwise the proof is similarly straightforward. 

On the other hand, the proof of~\eqref{eqn:sample_to_population_conductance} is considerably more complicated. Our proof relies on the recent results of \citet{garciatrillos16b}, who upper bound the  $\Leb^{\infty}$-optimal transport distance between the empirical measure $\Pbb_n$ and $\Pbb$. For further details, we refer to Appendix~\ref{subsec:sample_to_population_conductance}, where we prove Proposition~\ref{prop:sample_to_population_2}, as well as~\citet{garciatrillos16}, who establish the asymptotic convergence of the sample conductance as $n \to \infty$ and $r \to 0$.

\subsection{Cluster Recovery}
\label{subsec:cluster_recovery}
As is typical in the local clustering literature, our algorithmic results will be stated with respect to specific ranges of each of the user-specified
parameters. In particular, for $\delta \in (0,1/4)$ and a candidate cluster $\mc{C} \in \Rd$, we require that some of the tuning parameters of Algorithm~\ref{alg:ppr} be chosen within specific ranges, 

\begin{equation}
\label{eqn:initialization}
\begin{aligned}
& \alpha \in \Bigl[(1 - 4\delta)^2, (1 - 2\delta)^2\Bigr) \cdot
\frac{\alpha_{\Pbb,r}(\mc{C},\delta)}{2} \\
& (L,U) \subseteq \Bigl(\frac{1}{5(1 + 2\delta)},\frac{1}{5(1 + \delta)}\Bigr) \cdot 
\frac{1}{n(n - 1)\vol_{\Pbb,r}(\mc{C})},
\end{aligned}  
\end{equation}
where
\begin{equation}
\label{eqn:alpha_initialization}
\alpha_{\Pbb,r}(\mc{C},\delta) := \frac{\ln(2)}{17} \cdot \frac{\Psi_{\Pbb,r}^2(\mc{C})}{\ln^2\Bigl(\frac{32}{(1 - 4\delta)s_{\Pbb,r}(\mc{C})}\Bigr)}.
\end{equation}

\begin{definition}
	When the input parameters to Algorithm \ref{alg:ppr} satisfy \eqref{eqn:initialization} for some $\mc{C} \subseteq \Rd$ and $\delta \in (0,1/4)$, we say the algorithm is $\delta$-\emph{well-initialized} with respect to $\mc{C}$.
\end{definition}

Of course, in practice it is not feasible to set tuning parameters based on the underlying (unknown) distribution $\Pbb$ and candidate cluster $\mc{C}$. Typically, one runs PPR over some range of tuning parameter values and selects the cluster which has the smallest normalized cut. 

By combining Lemma~\ref{lem:zhu} and Propositions~\ref{prop:pointwise_mixing_time}-\ref{prop:sample_to_population_2}, we obtain an upper bound on $\Delta(\wh{C},\mc{C}[X])$ that depends solely on the distribution $\Pbb$ and candidate cluster $\mc{C}$. To ease presentation, we introduce a \emph{condition number}, defined for a given $\mc{C} \subseteq \Rd$ and $\delta \in (0,1/4)$ as
\begin{equation}
\label{eqn:condition_number}
\kappa_{\Pbb,r}(\mc{C},\delta) := \frac{(1 + 3\delta)(1+2\delta)}{(1 - 4\delta)^2(1 - \delta)} \cdot \frac{\Phi_{\Pbb,r}(\mc{C})}{\alpha_{\Pbb,r}(\mc{C},\delta)}.
\end{equation}

\begin{theorem}
	\label{thm:volume_ssd_ub} 
	Fix $\delta \in (0,1/4)$. Suppose $\wt{\Pbb}$ and $\mc{C}$ satisfy~\ref{asmp:bounded_density} and~\ref{asmp:domain}. Then for any $n \in \mbb{N}$ which satisfies~\eqref{eqn:sample_to_population_local_spread_sample_complexity}, \eqref{eqn:sample_to_population_conductance_sample_complexity}, and
	\begin{equation}
	\label{eqn:volume_ssd_ub_sample_complexity}
	\frac{(1 + \delta)}{(1 - \delta)^4} \cdot B_3 \leq n,
	\end{equation} 
	the following statement holds with probability at least $1 - B_2/n - 4\exp\{-b_1\delta^2n\} - (2n + 2)\exp\{-b_2\delta^2n\} - (n + 1)\exp\{-b_3n\}$: there exists a set $\mc{C}[X]^g \subseteq \mc{C}[X]$ of large volume,
	$$
	\vol_{n,r}(\mc{C}[X]^g) \geq \frac{1}{2}\vol_{n,r}(\mc{C}[X]),
	$$ 
	such that if Algorithm~\ref{alg:ppr} is $\delta$-well-initialized with respect to $\mc{C}[X]$, and run with any seed node $v \in \mc{C}[X]^g$, then the PPR estimated cluster $\wh{C}$ satisfies
	\begin{equation}
	\label{eqn:volume_ssd_ub}
	\frac{\Delta(\wh{C};\mc{C}[X])}{\vol_{n,r}(\mc{C}[X])} \leq 60 \cdot \kappa_{\Pbb,r}(\mc{C},\delta).
	\end{equation}
\end{theorem}

We now make some remarks.
\begin{itemize}
	\item It is useful to compare Theorem~\ref{thm:volume_ssd_ub} with what is already known regarding \emph{global} spectral clustering in the context of nonparametric statistics. \cite{schiebinger2015} consider the following variant of spectral clustering: first embed the data $X$ into $\Reals^k$ using the bottom $k$ eigenvectors of the degree-normalized Laplacian $I - D^{-1/2}AD^{-1/2}$, and then partition the embedded data into estimated clusters \smash{$\wh{C}_1,\ldots,\wh{C}_k$} using $k$-means clustering. They derive error bounds on the misclassification error that depend on a difficulty function $\varphi(\Pbb)$. In our context, where the goal is to successfully distinguish $\mc{C}$ and $\mc{C}^c$, thus where $k = 2$, this difficulty function is roughly
	\begin{equation}
	\label{eqn:schiebinger}
	\varphi(\Pbb) \approx \sqrt{\Phi_{\Pbb,r}(\mc{C})} \cdot \max\biggl\{ \frac{1}{\Psi_{\Pbb,r}(\mc{C})^2}; \frac{1}{\Psi_{\Pbb,r}(\mc{C}^c)^2}\biggr\}.
	\end{equation}
	We point out two ways in which~\eqref{eqn:volume_ssd_ub} is a tighter bound than~\eqref{eqn:schiebinger}. First,~\eqref{eqn:schiebinger} depends on $\Psi_{\Pbb,r}(\mc{C}^c)$ in addition to $\Psi_{\Pbb,r}(\mc{C})$, and is thus a useful bound only if $\mc{C}^c$ and $\mc{C}$ are both internally well-connected. In contrast~\eqref{eqn:volume_ssd_ub} depends only on $\Psi_{\Pbb,r}(\mc{C})$, and is thus a useful bound if $\mc{C}$ has small conductance, regardless of the conductance of $\mc{C}^c$. This is intuitive: PPR is a local rather than global algorithm, and as such the analysis requires only local rather than global conditions. Second,~\eqref{eqn:schiebinger} depends on \smash{$\sqrt{\Phi_{\Pbb,r}(\mc{C})}$} rather than $\Phi_{\Pbb,r}(\mc{C})$, and since $\Phi_{\Pbb,r}(\mc{C}) \leq 1$ this results in a weaker bound.~\citet{schiebinger2015} provide experiments suggesting that the linear, rather than square-root, dependence is correct, and we theoretically confirm this in the local clustering setup. Of course, on the other hand~\eqref{eqn:volume_ssd_ub} depends on $\log^2(1/s_{\Pbb,r}(\mc{C}))$, which is due to the locally-biased nature of the PPR algorithm, and does not appear in~\eqref{eqn:schiebinger}.
	
	\item Although Theorem~\ref{thm:volume_ssd_ub} is stated with respect to the exact PPR vector $p_v$, for a sufficiently small choice of $\varepsilon$ the 
	application of \eqref{eqn:appr_error} within the proof of Theorem
	\ref{thm:volume_ssd_ub} leads to an analogous result which holds for the aPPR vector \smash{$p_v^{(\varepsilon)}$}. We formally state and prove this fact in Appendix~\ref{apdx:appr_misclassification_error}.
\end{itemize}

\section{Recovery of a Density Cluster with PPR}
\label{sec:ppr_density_cluster}
We now apply the general theory established in the last section to the special case where $\mc{C} = \mc{C}_{\lambda}$ is a $\lambda$-density cluster---that is, a connected component of the upper level set $\{x \in \Rd: f(x) \geq \lambda\}$. In Section~\ref{sec:lower_bound}, we also derive a lower bound, giving a ``hard problem'' for which PPR will provably fail to recover a density cluster. Together, these results can be summarized as follows: PPR recovers a density cluster $\mc{C}_{\lambda}$ if and only if both $\mc{C}_{\lambda}$ and $f$ are well-conditioned, meaning that $\mc{C}_{\lambda}$ is not too long and thin, and that $f$ is approximately uniform inside $\mc{C}_{\lambda}$ while satisfying a low-noise condition near its boundary.

\subsection{Recovery of Well-Conditioned Density Clusters}
\label{subsec:recovery_well-conditioned_density_clusters}

All results on density clustering assume the density $f$ satisfies some regularity conditions. A basic requirement is the need to avoid clusters which contain arbitrarily thin bridges or spikes, or more generally clusters which can be disconnected by removing a subset of (Lebesgue) measure $0$, and thus may not be resolved by any finite number of samples. To rule out such problematic clusters, we follow the approach of~\cite{chaudhuri2010}, who assume the density is lower bounded on a thickened version of $\mc{C}_{\lambda}$, defined as $\mc{C}_{\lambda,\sigma} := \{x \in \Rd: \dist(x,\mc{C}) < \sigma\}$ for a given $\sigma > 0$. Regardless of the dimension of $\mc{C}_{\lambda}$, the set $\mc{C}_{\lambda,\sigma}$ is full dimensional. Under typical uniform continuity conditions, the requirement that the density be lower bounded over $\mc{C}_{\lambda,\sigma}$ will be satisfied. Such continuity conditions can be weakened (see for instance \citet{rinaldo2010,steinwart2015}) but we do not pursue the matter further.

In summary, our goal is to obtain upper bounds on \smash{$\Delta(\wh{C},\mc{C}_{\lambda,\sigma}[X])$}, for some fixed $\lambda$ and $\sigma > 0$. We have already derived upper bounds on the symmetric set difference of \smash{$\wh{C}$} and a generic cluster $\mc{C}$ that depend on some population-level functionals of $\mc{C}$. What remains is to analyze these population-level functionals in the specific case where the candidate cluster is $\mc{C}_{\lambda,\sigma}$. To carry out this analysis, we will need to impose some conditions, and for the rest of this section we will assume the following.

\begin{enumerate}[label=(A\arabic*)]
	\setcounter{enumi}{2}
	\item
	\label{asmp:lambda_bounded_density}
	\emph{Bounded density within cluster:} There exist constants
	$0<\lambda_{\sigma}< \Lambda_{\sigma}<\infty$ such that 
	$$
	\lambda_{\sigma} \leq \inf_{x \in \mc{C}_{\lambda,\sigma}} f(x) \leq \sup_{x \in \mc{C}_{\lambda,\sigma}} f(x)
	\leq \Lambda_{\sigma}.
	$$
	
	\item 
	\label{asmp:low_noise_density}
	\emph{Low-noise density:} There exist $\theta \in (0,\infty)$ and $\gamma \in
	[0,1]$ such that for any $x \in \Rd$ with $0 < \dist(x, \mc{C}_{\lambda,\sigma}) \leq \sigma$,     
	$$
	\inf_{y \in \mc{C}_{\lambda,\sigma}} f(y) - f(x) \geq  \theta \cdot \dist(x, \mc{C}_{\lambda,\sigma})^{\gamma}.  
	$$
	Roughly, this assumption ensures that the density decays sufficiently quickly
	as we move away from the target cluster $\mc{C}_{\lambda,\sigma}$, and is a standard assumption in the level-set estimation literature (see for instance \citet{singh2009}).
	
	\item
	\label{asmp:embedding}
	\emph{Lipschitz embedding:}
	There exists a differentiable function $g: \Rd \to \Rd$, $\rho \in (0,\infty)$ and $M \in [1,\infty)$ such that
	\begin{enumerate}
		\item $\mc{C}_{\lambda,\sigma} = g(\mc{K})$, for a convex set $\mc{K} \subseteq \Rd$ with $\mathrm{diam}(\mc{K}) = \sup_{x,y \in \mc{K}}\|x - y\| \leq \rho < \infty$;
		\item $\det(\nabla g (x)) = 1$ for all $x \in \mc{K}$, where $\nabla g(x)$ is the Jacobian of $g$ evaluated at $x;$ and 
		\item for some $M \geq 1$,   
		$$
		\|g(x) - g(y)\| \leq M \|x - y\| ~	\text{for all $x,y \in \mc{K}$}. 
		$$
	\end{enumerate}
	Succinctly, we assume that $\mc{C}_{\lambda,\sigma}$ is the image of a convex set with finite
	diameter under a measure preserving, Lipschitz transformation. 
\end{enumerate}

\emph{For convenience only}, we will also make the following assumption.
\begin{enumerate}[label=(A\arabic*)]
	\setcounter{enumi}{5}
	\item
	\label{asmp:bounded_volume}
	\emph{Bounded volume:}
	The volume of $\mc{C}_{\lambda,\sigma}$ is no more than half the total volume of $\Rd$:
	$$
	\vol_{\Pbb,r}(\mc{C}_{\lambda,\sigma}) \leq \vol_{\Pbb,r}(\mc{C}_{\lambda,\sigma}^c). 
	$$
	This assumption implies that the normalized cut of $\mc{C}_{\lambda,\sigma}$ will be equal to the ratio of $\cut_{\Pbb,r}(\mc{C}_{\lambda,\sigma})$ to $\vol_{\Pbb,r}(\mc{C}_{\lambda,\sigma})$.
\end{enumerate}

\subsubsection{Normalized Cut, Conductance, and Local Spread of a Density Cluster} In Lemma~\ref{lem:density_cluster_local_spread}, Proposition~\ref{prop:density_cluster_normalized_cut}, and Proposition~\ref{prop:density_cluster_conductance}, we give bounds on the population local spread, normalized cut, and conductance of $\mc{C}_{\lambda,\sigma}$. These bounds depend on the various geometric parameters just introduced. 
\begin{lemma}
	\label{lem:density_cluster_local_spread}
	Assume $\mc{C}_{\lambda,\sigma}$ satisfies Assumptions~\ref{asmp:lambda_bounded_density} and~\ref{asmp:embedding} for some $\lambda_{\sigma},\Lambda_{\sigma},\rho$ and $M$. Then,
	\begin{equation}
	\label{eqn:density_cluster_local_spread}
	s_{\Pbb,r}(\mc{C}_{\lambda,\sigma}) \geq \frac{1}{4} \cdot \frac{\lambda_{\sigma}^2}{\Lambda_{\sigma}^2} \cdot \biggl(\frac{2r}{\rho}\biggr)^{d} \cdot \biggl(1 - \frac{r}{\sigma} \sqrt{\frac{d + 2}{2\pi}}\biggr).
	\end{equation}
\end{lemma}

\begin{proposition}
	\label{prop:density_cluster_normalized_cut}
	Assume $\mc{C}_{\lambda,\sigma}$ satisfies Assumptions~\ref{asmp:lambda_bounded_density},~\ref{asmp:low_noise_density} and~\ref{asmp:bounded_volume} for some $\lambda_{\sigma}, \Lambda_{\sigma}, \theta$, and $\gamma$, and additionally that $0 < r \leq \frac{\sigma}{4d}$. Then,
	\begin{equation}
	\label{eqn:density_cluster_normalized_cut}
	\Phi_{\Pbb,r}(\mc{C}_{\lambda,\sigma}) \leq \frac{16}{9} \cdot \frac{dr}{\sigma} \cdot \frac{\lambda\Bigl(\lambda_{\sigma} - \theta \frac{r^{\gamma}}{\gamma + 1}\Bigr)}{\lambda_{\sigma}^2}.
	\end{equation}
\end{proposition}

\begin{proposition}
	\label{prop:density_cluster_conductance}
	Assume $\mc{C}_{\lambda,\sigma}$ satisfies Assumptions~\ref{asmp:lambda_bounded_density} and~\ref{asmp:embedding} for some $\lambda_{\sigma}, \Lambda_{\sigma}, \rho$ and $M$. Then,
	\begin{equation}
	\label{eqn:density_cluster_conductance}
	\Psi_{\Pbb,r}(\mc{C}_{\lambda,\sigma}) \geq \Bigl(1 - \frac{r}{4\rho M}\Bigr) \cdot \Bigl(1 - \frac{r}{\sigma}\sqrt{\frac{d + 2}{2\pi}}\Bigr)^2 \cdot \frac{\sqrt{2\pi}}{36} \cdot \frac{r}{\rho M \sqrt{d + 2}} \cdot \frac{\lambda_{\sigma}^2}{\Lambda_{\sigma}^2}.
	\end{equation}
\end{proposition}

Some remarks are in order.
\begin{itemize}
	\item We prove Proposition~\ref{prop:density_cluster_normalized_cut} by separately upper bounding $\cut_{\Pbb,r}(\mc{C}_{\lambda,\sigma})$ and lower bounding the volume $\vol_{\Pbb,r}(\mc{C}_{\lambda,\sigma})$. Of these two bounds, the trickier to prove is the upper bound on the cut, which involves carefully estimating the probability mass of thin tubes around the boundary of $\mc{C}_{\lambda,\sigma}$. 
	\item Proposition~\ref{prop:density_cluster_conductance} is proved in a completely different way. The proof relies heavily on bounds on the isoperimetric ratio of convex sets (as derived by e.g., \cite{lovasz1990} or \cite{dyer1991b}), and thus the embedding assumption \ref{asmp:embedding} and Lipschitz parameter $M$
	play an important role in proving the upper bound in Proposition~\ref{prop:density_cluster_conductance}. 
	\item There is some interdependence between $M$ and $\sigma,\rho$, which might lead one to hope that \ref{asmp:embedding} is
	non-essential. However, it is not possible to eliminate condition \ref{asmp:embedding} without incurring an additional factor of at least
	$(\rho/\sigma)^d$ in \eqref{eqn:density_cluster_conductance}, achieved, for
	instance, when $\mc{C}_{\lambda,\sigma}$ is a dumbbell-like set consisting of two balls of diameter $\rho$ linked by a cylinder of radius $\sigma$. In contrast,~\eqref{eqn:density_cluster_conductance} depends polynomially on $d$, and many reasonably shaped sets---such as star-shaped sets as well as half-moon shapes of the type we consider in Section \ref{sec:experiments}---satisfy \ref{asmp:embedding} for reasonably small values of $M$ \citep{abbasi-yadkori2016a, abbasi-yadkori2016}.
\end{itemize}

Applying these results along with Theorem~\ref{thm:volume_ssd_ub}, we obtain an upper bound on \smash{$\Delta(\wh{C},\mc{C}_{\lambda,\sigma}[X])$}. In what follows, $C_{1,\delta},C_{2,\delta},\ldots$ are constants which may depend on $\delta$, but not on $n$, $\Pbb$ or $\mc{C}_{\lambda,\sigma}$, and which we keep track of in our proofs.

\begin{theorem}
	\label{thm:density_cluster_volume_ssd_ub}
	Let $C_{\lambda} \subseteq \Rd$ and $\delta \in (0,1/4)$. Suppose that $C_{\lambda,\sigma}$ satisfies~\ref{asmp:domain}-\ref{asmp:bounded_volume} for some $\lambda_{\sigma}, \Lambda_{\sigma}, \theta, \gamma, \rho$ and $M$, that $0 < r \leq \sigma/4d$, and that the sample size $n$ satisfies the same conditions as in Theorem~\ref{thm:volume_ssd_ub}. Then with probability at least $1 - B_2/n - 4\exp\{-b_1\delta^2n\} - (2n + 2)\exp\{-b_2\delta^2n\} - (n + 1)\exp\{-b_3n\}$, the following statement holds: there exists a set $\mc{C}_{\lambda,\sigma}[X]^g \subseteq \mc{C}_{\lambda,\sigma}[X]$ of large volume, 
	$$
	\vol_{n,r}(\mc{C}_{\lambda,\sigma}[X]^g) \geq \frac{1}{2} \vol_{n,r}(\mc{C}_{\lambda,\sigma}[X]).
	$$ 
	such that if Algorithm~\ref{alg:ppr} is $\delta$-well-initialized with respect to $\mc{C}_{\lambda,\sigma}$, and run with any seed node $v \in \mc{C}_{\lambda,\sigma}[X]^g$, then the PPR estimated cluster \smash{$\wh{C}$} satisfies
	\begin{equation}
	\label{eqn:density_cluster_volume_ssd_ub}
	\frac{\Delta(\wh{C};\mc{C}_{\lambda,\sigma}[X])}{\vol_{n,r}(\mc{C}[X])} \leq C_{1,\delta} \cdot d^3(d + 2) \cdot \frac{M^2\rho^2}{\sigma r} \cdot \frac{\Lambda_{\sigma}^2 \lambda (\lambda - \theta \frac{r^{\gamma}}{\gamma + 1})}{\lambda_{\sigma}^4} \cdot \log^2\biggl(C_{2,\delta}^{1/d} \frac{\Lambda_{\sigma}^{2/d} M\rho}{\lambda_{\sigma}^{2/d}2r}\biggr)
	\end{equation}
\end{theorem}

Several further remarks are as follows.
\begin{itemize}
	\item Observe that while the diameter $\rho$ is absent from our upper bound on normalized cut in Proposition \ref{prop:density_cluster_normalized_cut}, it enters the ultimate bound in Theorem~\ref{thm:density_cluster_volume_ssd_ub} through the conductance. This reflects (what may be regarded as) established wisdom regarding spectral partitioning algorithms more generally \citep{guattery1995, hein2010}, but newly applied to the density clustering setting: if the diameter $\rho$ is large, then PPR	may fail to recover $\mc{C}_{\lambda,\sigma}[X]$ even when $\mc{C}_{\lambda}$ is sufficiently well-conditioned to ensure that $\mc{C}_{\lambda,\sigma}[X]$ has a small normalized cut in $G_{n,r}$. This will be supported by simulations in Section~\ref{subsec:empirical_behavior_ppr}.   
	\item Several modifications of global spectral clustering have been proposed with the intent of making such procedures essentially independent of the shape of the density cluster $\mc{C}_{\lambda}$. For instance,~\citet{arias-castro2009,pelletier2011} introduce a cleaning step to remove low-degree vertices, whereas~\cite{little2020} use a weighted geometric graph, where the weights are computed with respect to a density-dependent distance. The resulting procedures come with stronger density cluster recovery guarantees. However, the key ingredient in such procedures is the explicitly density-dependent part of the algorithm, and spectral clustering functions as more of a post-processing step. These methods are as such very different in spirit to PPR, which is a bona fide (local) spectral clustering algorithm. 
	\item As mentioned in the discussion after Theorem~\ref{thm:volume_ssd_ub}, the population normalized cut and conductance also play a leading role in the analysis of global spectral clustering algorithms. It therefore seems likely that similar bounds to~\eqref{eqn:density_cluster_volume_ssd_ub} would apply to the output of global spectral clustering methods as well, but formalizing this is outside the scope of our work.
	\item The symmetric set difference does not measure whether \smash{$\wh{C}$} can (perfectly) distinguish any two distinct clusters $\mc{C}_{\lambda},\mc{C}_{\lambda}' \in \mbb{C}_f(\lambda)$. In Appendix~\ref{apdx:appr_misclassification_error}, we show that the PPR estimate $\wh{C}$ can in fact distinguish two distinct clusters $\mc{C}_{\lambda}$ and $\mc{C}_{\lambda}'$, but the result holds only under relatively restrictive conditions.
\end{itemize}

\section{Negative Result}
\label{sec:lower_bound}
\begin{figure}[tb]
	\centering
	\includegraphics[width=0.495\textwidth]{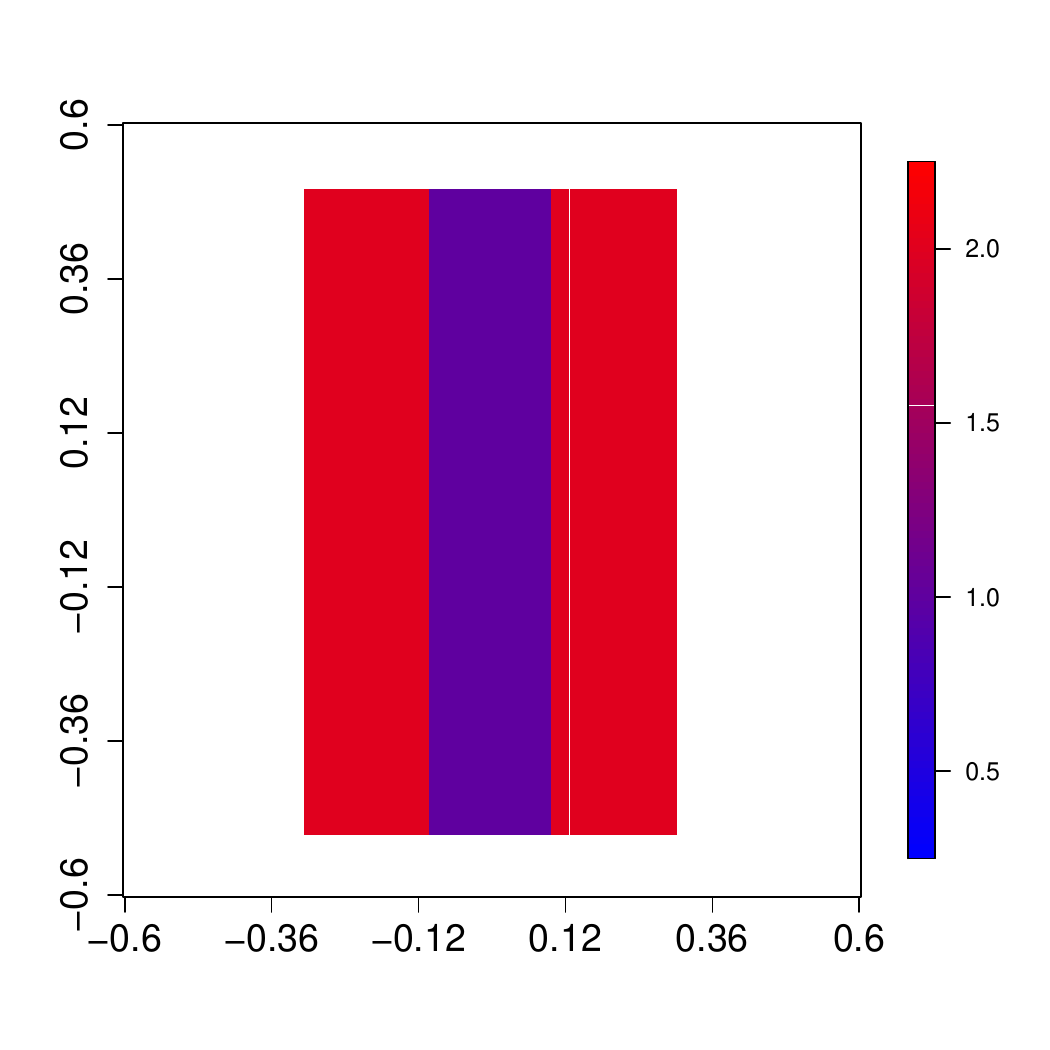}
	\includegraphics[width=0.495\textwidth]{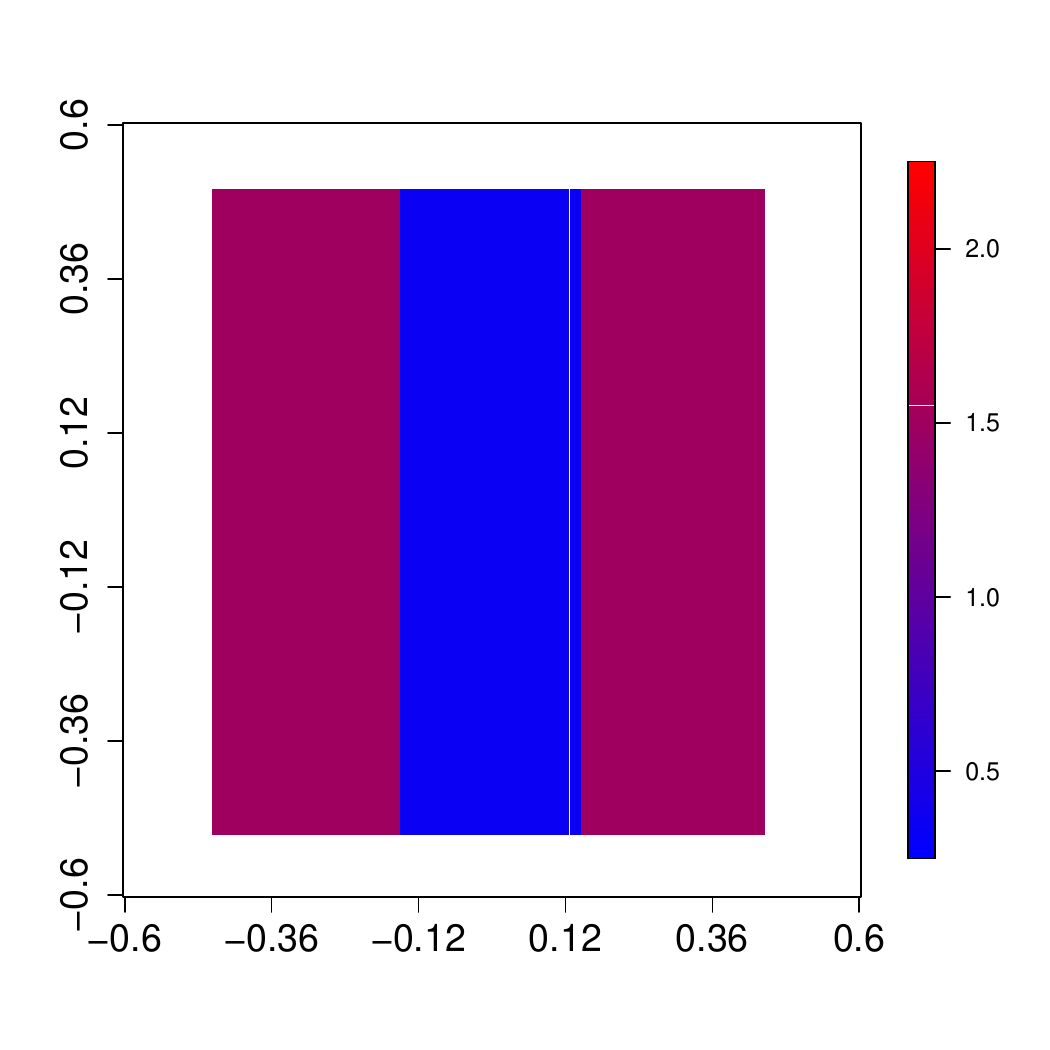}
	\caption{\small The density $f$ in \eqref{eqn:lb_density}, for
		$\rho=1$, and two different choices of $\epsilon$ and $\sigma$. Left:
		$\epsilon = 0.3$ and $\sigma = 0.1$; right: $\epsilon = 0.2$ and 
		$\sigma = 0.2$.} 
	\label{fig:hard_case}
\end{figure}

We now exhibit a hard case for density clustering using PPR, that is, a distribution $\Pbb$ for which PPR is unlikely to recover a density cluster. Let \smash{$\mc{C}^{(0)}, \mc{C}^{(1)}, \mc{C}^{(2)}$} be rectangles in $\Reals^2$,    
$$
\mc{C}^{(0)} = \biggl[-\frac{\sigma}{2}, \frac{\sigma}{2}\biggr] \times 
\biggl[-\frac{\rho}{2}, \frac{\rho}{2}\biggr], \quad 
\mc{C}^{(1)} = \mc{C}^{(0)} - (\sigma,0), \quad
\mc{C}^{(2)} = \mc{C}^{(0)} + (\sigma,0),
$$
where $0 < \sigma < \rho$, and let $\Pbb$ be the mixture distribution over \smash{$\mc{X} = \mc{C}^{(0)} \cup \mc{C}^{(1)} \cup \mc{C}^{(2)}$} given by   
$$
\Pbb = \frac{1 - \epsilon}{2} \Pbb_1 + \frac{1 - \epsilon}{2} \Pbb_2 +
\frac{\epsilon}{2} \Pbb_0, 
$$
where $\Pbb_k$ is the uniform distribution over $\mc{C}^{(k)}$ for $k = 0,1,2$.  
The density function $f$ of $\Pbb$ is simply
\begin{equation}
\label{eqn:lb_density}
f(x) = \frac{1}{\rho\sigma}\left(\frac{1 - \epsilon}{2}\1(x \in
\mc{C}^{(1)}) + \frac{1 - \epsilon}{2}\1(x \in \mc{C}^{(2)}) +
\epsilon\1(x \in \mc{C}^{(0)})  \right), 
\end{equation}
so that for any $\epsilon < \lambda < (1 - \epsilon)/2$, we have \smash{$\mbb{C}_f(\lambda) = \set{\mc{C}^{(1)}, \mc{C}^{(2)}}$}. Figure~\ref{fig:hard_case} visualizes the density $f$ for two different choices of $\epsilon, \sigma, \rho$.  

\subsection{Lower Bound on Symmetric Set Difference}
As the following theorem demonstrates, even when Algorithm~\ref{alg:ppr} is reasonably initialized, if the density cluster \smash{$\mc{C}^{(1)}$} is sufficiently geometrically ill-conditioned (in words, tall and thin) the cluster estimator $\wh{C}$ will fail to recover \smash{$\mc{C}^{(1)}$}. Let
\begin{equation}
\label{eqn:lower_set}
\mc{L} = \set{(x_1,x_2) \in \mc{X}: x_2 < 0}.
\end{equation}

In the following Theorem, $B_{1,\delta}$ and $B_{2,\delta}$ are constants which may depend on $\delta,\Pbb,\mc{C}_{\lambda,\sigma}$ and $r$, but not on $n$.
\begin{theorem}
	\label{thm:ppr_lb}
	Fix $\delta \in (0,1/7)$. Assume the neighborhood graph radius $r < \sigma/4$, that
	\begin{equation}
	\label{eqn:ppr_lb_condition}
	\max\biggl\{ B_{1,\delta} \cdot \frac{r}{\rho}, B_{2,\delta} \cdot \frac{1}{n} \biggr\} < \frac{1}{18}~~\text{and}~~n \geq 8 \frac{(1 + \delta)}{(1 - \delta)},
	\end{equation} 
	and that Algorithm~\ref{alg:ppr} is initialized using inputs $\alpha = 36 \cdot \Phi_{n,r}(\mc{L}[X])$, and $(L,U) = (0,1)$.  Then the following statement holds with probability at least $1 - (B_4 + 2n + 10)\exp\{-n\delta^2 b_4\}$: there exists a set $\mc{C}[X]^g$ of large volume, 
	$$
	\vol_{n,r}(\mc{C}[X]^g \cap \mc{C}^{(1)}[X]) \geq \frac{1}{8} \vol_{n,r}(\mc{C}^{(1)}[X];G_{n,r}),
	$$
	such that for any seed node $v \in \mc{C}[X]^g$, the PPR estimated cluster
	\smash{$\wh{C}$} satisfies    
	\begin{equation}
	\label{eqn:ppr_lb}
	\frac{\sigma \rho}{r^2 n^2} \cdot \vol_{n,r}(\wh{C} \vartriangle \mc{C}^{(1)}[X]) \geq \frac{1 - \delta}{2} -  C_{3,\delta} \cdot \frac{\sqrt{\sigma/\rho}}{\epsilon^2} \cdot \sqrt{ \log\left(C_{4,\delta} \cdot \frac{\rho \sigma}{\epsilon^2 r^2}\right) \frac{\sigma}{r}},   
	\end{equation} 
\end{theorem}

We make a couple of remarks.
\begin{itemize}
	\item Theorem~\ref{thm:ppr_lb} is stated with respect to a particular hard case, where the density clusters are rectangular subsets of $\Reals^2$.  We chose this
	setting to make the theorem simple to state, and our results are generalizable
	to $\Rd$ and to non-rectangular clusters. Technically, the rectangles $\mc{C}^{(0)},\mc{C}^{(1)},\mc{C}^{(2)}$ are not
	$\sigma$-expansions due to their sharp corners. To fix this, one can   
	simply modify these sets to have appropriately rounded corners, and our lower
	bound arguments do not need to change significantly, subject to some
	additional bookkeeping.  Thus we ignore this technicality in our subsequent
	discussion. 
	
	\item Although we state our lower bound with respect to PPR run on a neighborhood graph, the conclusion is likely to hold for a much broader class of spectral clustering algorithms. In the proof of Theorem~\ref{thm:ppr_lb}, we rely heavily on the fact that when $\epsilon^2$ is sufficiently greater than $\sigma/\rho$, the normalized cut of $\mc{C}^{(1)}$ will be much larger than that of $\mc{L}$. In this case, not merely PPR but any algorithm that approximates the minimum normalized cut is unlikely to recover $\mc{C}^{(1)}$. In particular, local spectral clustering methods that are based on truncated random walks \citep{spielman2013}, global spectral clustering algorithms \citep{shi00}, and $p$-Laplacian based spectral embeddings \citep{hein2010} all have provable upper bounds on the normalized cut of cluster they output, and thus we expect that they would all fail to estimate $\mc{C}^{(1)}$.
\end{itemize}

\subsection{Comparison Between Upper and Lower Bounds}
\label{subsec:comparison_upper_lower_bounds}
To better digest the implications of Theorem~\ref{thm:ppr_lb}, we translate the
results of our upper bound in Theorem~\ref{thm:density_cluster_volume_ssd_ub} to the density $f$ given in \eqref{eqn:lb_density}. Observe that $\mc{C}^{(1)}$ satisfies each of the Assumptions~\ref{asmp:lambda_bounded_density}--\ref{asmp:bounded_volume}:

\begin{enumerate}[label=(A\arabic*)]
	\setcounter{enumi}{6}
	\item The density $f(x) = \frac{1 - \epsilon}{2 \rho \sigma}$ for all $x \in	\mc{C}^{(1)}$.  
	\item The density $f(x) \leq \frac{\epsilon}{\rho\sigma}$ for all $x \in \Reals^2$ such that $0 < \dist(x,\mc{C}^{(1)}) \leq \sigma$. Therefore for all such $x$, 
	$$
	\inf_{x' \in \mc{C}^{(1)}} f(x') - f(x)  \geq \left\{\frac{1 - \epsilon}{2} - \epsilon \right\} \frac{1}{\rho \sigma},
	$$
	which meets the decay requirement with exponent $\gamma=0$.
	\item The set $\mc{C}^{(1)}$ is itself convex, and has diameter $\sqrt{\rho^2 + \sigma^2}$.
	\item By symmetry, \smash{$\vol_{\Pbb,r}(\mc{C}^{(1)}) = \vol_{\Pbb,r}(\mc{C}^{(2)})$}, and therefore \smash{$\vol_{\Pbb,r}(\mc{C}^{(1)}) \leq \frac{1}{2}\vol_{\Pbb,r}(\Rd)$}.   
\end{enumerate}

If the user-specified parameters are initialized according to~\eqref{eqn:initialization}, we may apply Theorem~\ref{thm:density_cluster_volume_ssd_ub}. This implies that there exists a set $\mc{C}^{(1)}[X]^g \subseteq \mc{C}^{(1)}[X]$ with \smash{$\vol_{n,r}(\mc{C}^{(1)}[X]^g) \geq \frac{1}{2}\vol_{n,r}(\mc{C}^{(1)}[X])$} such that for any seed node $v \in \mc{C}^{(1)}[X]^g$, and for large enough $n$, the PPR estimated cluster \smash{$\wh{C}$} satisfies with high probability
\begin{equation*}
\frac{\vol_{n,r}(\wh{C} \vartriangle \mc{C}^{(1)}[X])}{\vol_{n,r}(\mc{C}^{(1)}[X])} \leq 64 C_{1,\delta} \cdot \frac{\rho^2 + \sigma^2}{\sigma r} \cdot \frac{\epsilon}{1 - \epsilon} \cdot \log^2\biggl(\sqrt{C_{2,\delta}} \frac{\rho}{2r}\biggr)
\end{equation*}
To facilitate comparisons between our upper and lower bounds set $r = \sigma/8$.  Then the following statements each hold with high probability. 
\begin{itemize}
	\item If the user-specified parameters satisfy~\eqref{eqn:initialization}, and for some $a \geq 0$,
	\begin{equation*}
	\frac{\epsilon}{1 - \epsilon} \leq \frac{a}{512 C_{1,\delta}} \frac{\sigma^2}{(\rho^2 + \sigma^2) \log^2(\rho/\sigma \sqrt{C_{2,\delta}})},
	\end{equation*}
	then \smash{$\Delta(\wh{C}, \mc{C}^{(1)}[X]) \leq a \cdot
		\vol_{n,r}(\mc{C}^{(1)}[X])$}.
	
	\item The population-level volume $\vol_{\Pbb,r}(\mc{C}^{(1)}) \leq (1 - \epsilon)/2 \cdot \pi r^2/(\rho\sigma)$, and
	\begin{equation*}
	\vol_{n,r}(\mc{C}^{(1)}[X]) \leq (1 + \delta) \cdot n(n - 1) \vol_{\Pbb,r}(\mc{C}^{(1)}). 
	\end{equation*}
	Therefore, if the user-specified parameters are as in Theorem~\ref{thm:ppr_lb}, and
	\begin{equation*}
	\epsilon \geq \sqrt{8 C_{3,\delta}} \left({\frac{\sigma}{\rho}} \log \left(64C_{4,\delta} \cdot \frac{\rho}
	{\epsilon^2 \sigma}\right)\right)^{1/4},
	\end{equation*}
	then \smash{$\Delta(\wh{C}, \mc{C}^{(1)}[X]) \geq \frac{1}{20} \vol_{n,r}(\mc{C}^{(1)}[X])$}. 
\end{itemize}
Ignoring constants and log factors, we can summarize the above conclusions as follows: if $\epsilon$ is much less than $(\sigma/\rho)^2$, then PPR will approximately recover the density cluster $\mc{C}^{(1)}$, whereas if $\epsilon$ is much greater than $(\sigma/\rho)^{1/4}$ then PPR will fail to recover $\mc{C}^{(1)}$, even if it is reasonably initialized with a seed node $v \in \mc{C}^{(1)}$. Jointly, these upper and lower bounds give a relatively precise characterization of what it means for a density cluster to be well- or poorly-conditioned for recovery using PPR.\footnote{It is worth pointing out that the above conclusions are reliant on specific (albeit reasonable) ranges and choices of input parameters, which in some instances differ between the upper and lower bounds. We suspect that our lower bound continues to hold even when choosing input parameters as dictated by our upper bound, but do not pursue the details.}

Of course, it is not hard to show that in the example under consideration, classical plug-in density cluster estimators can consistently recover $\mc{C}^{(1)}$, even if $\epsilon$ is large compared to $\sigma/\rho$. That PPR has trouble recovering density clusters here (where standard plug-in approaches do not) is not meant to be a knock on PPR. Rather, it simply reflects that while classical density clustering approaches are specifically designed to identify high-density regions regardless of their geometry, PPR relies on geometry as well as density when forming the output cluster. 

\section{Experiments}
\label{sec:experiments}
We provide numerical experiments to investigate the tightness of our theoretical results in Section~\ref{sec:ppr_density_cluster}, and compare the performance of PPR with a density clustering algorithm on the ``two moons'' dataset. We defer details of the experimental settings to Appendix~\ref{apdx:experimental_details}.   

\subsection{Validating Theoretical Bounds}
\begin{figure}[htb]
	\centering
	\includegraphics[width=0.24\textwidth]{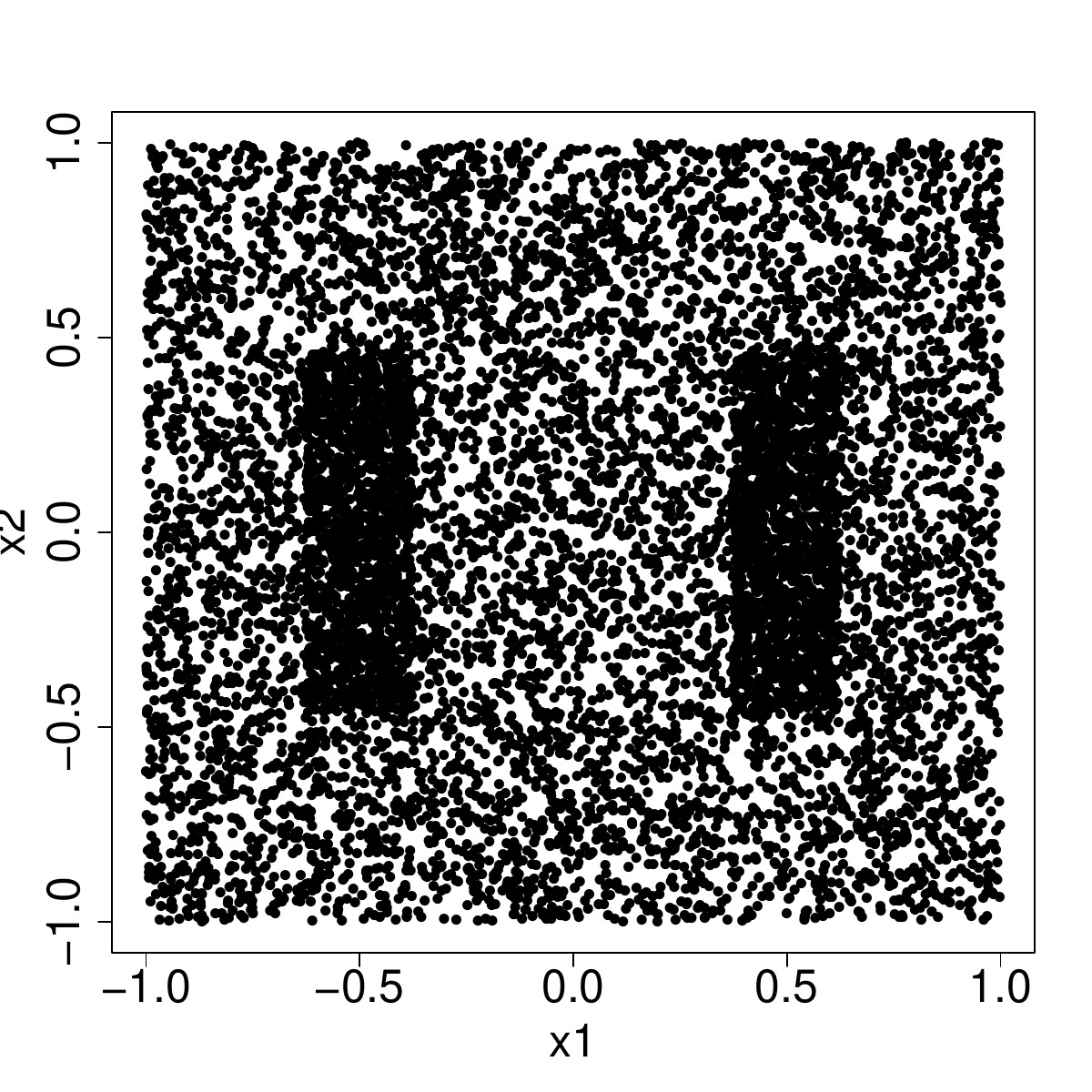}
	\includegraphics[width=0.24\textwidth]{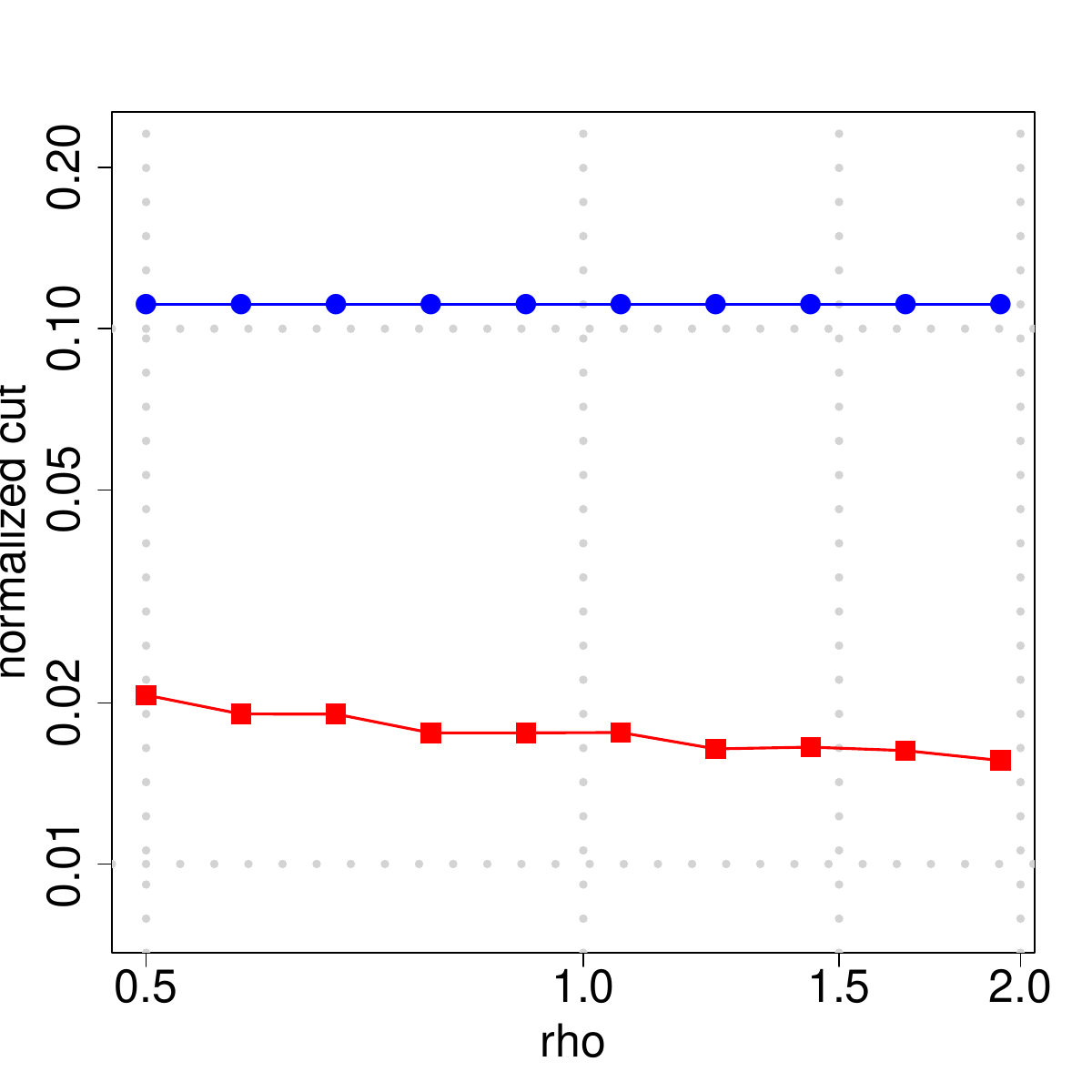}
	\includegraphics[width=0.24\textwidth]{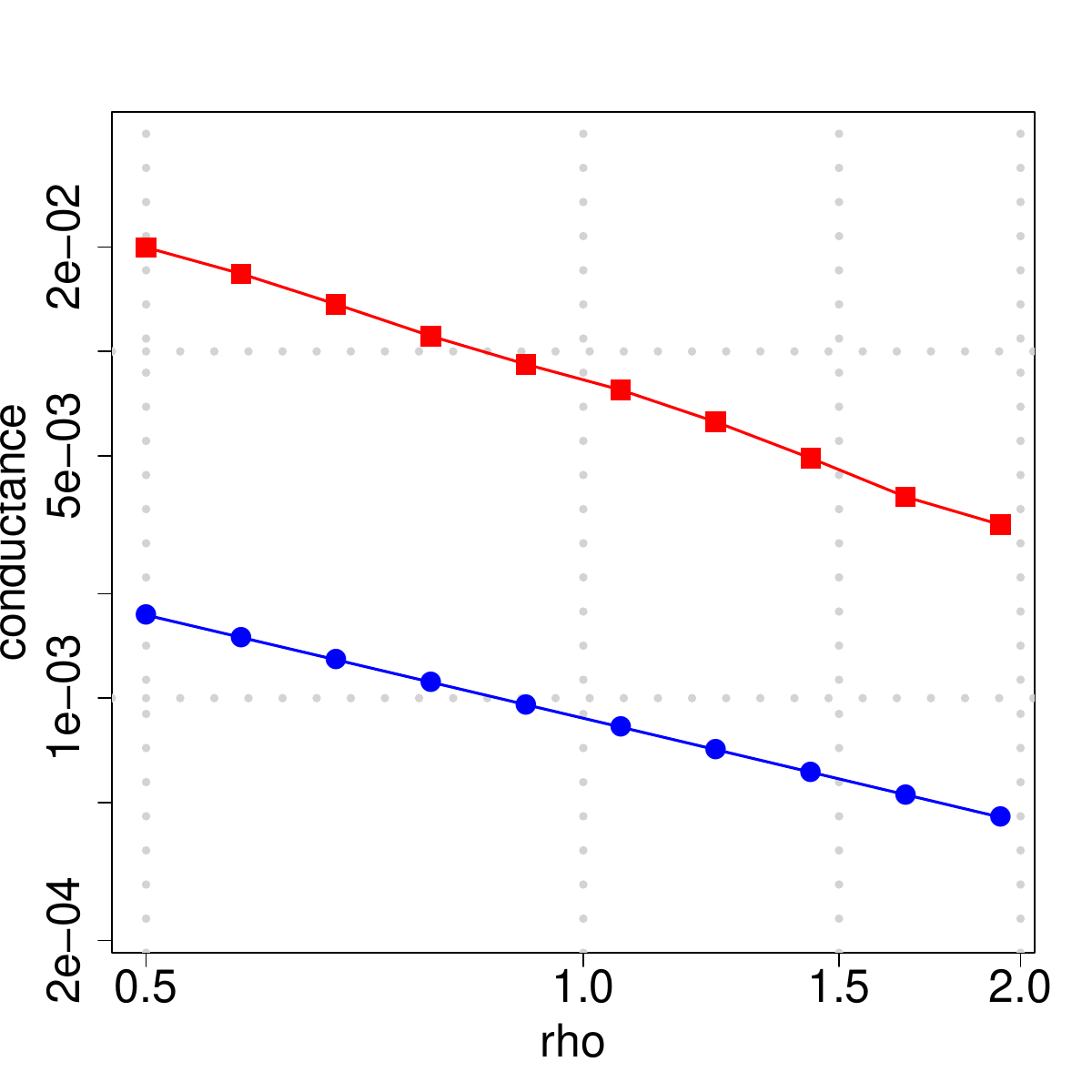}
	\includegraphics[width=0.24\textwidth]{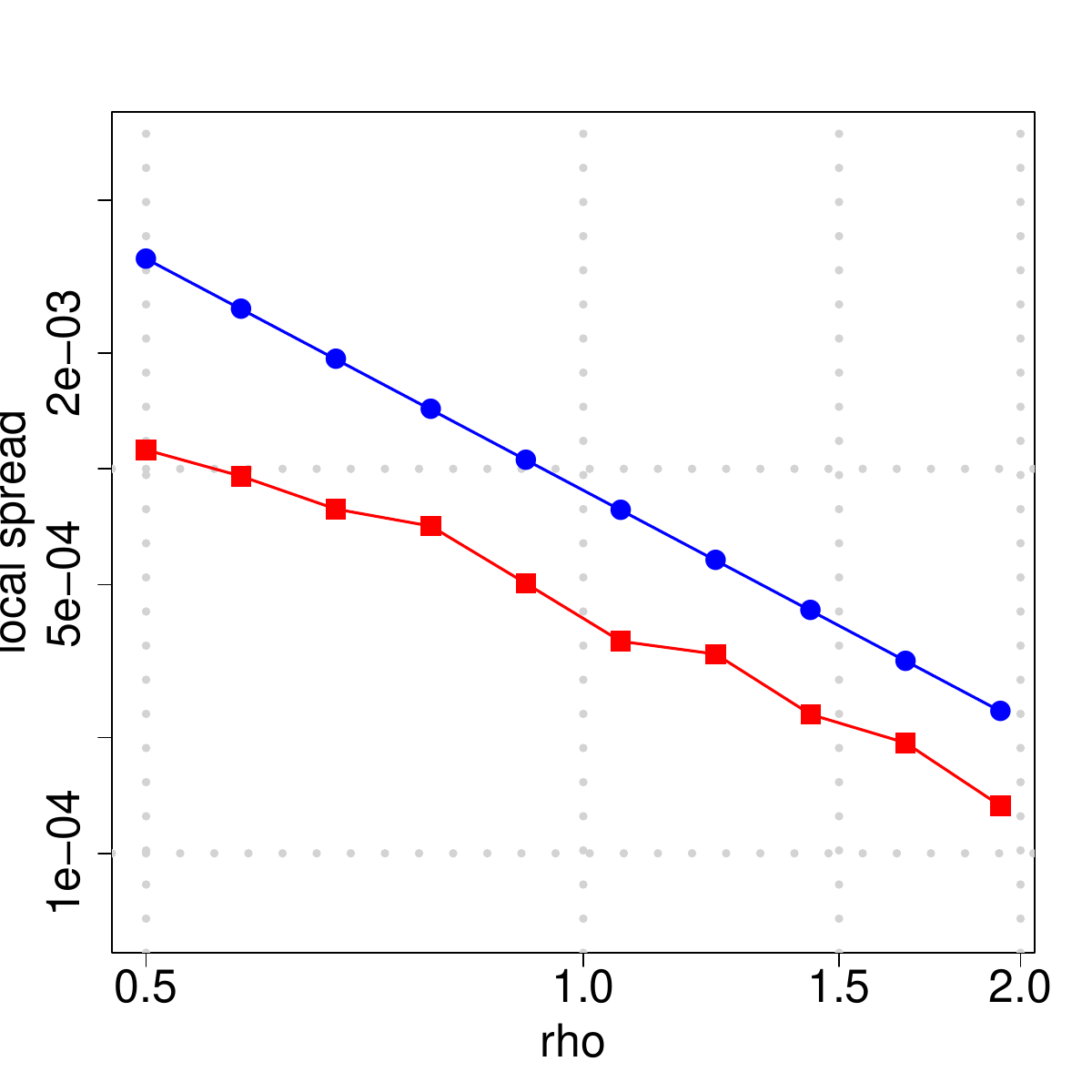}
	\includegraphics[width=0.24\textwidth]{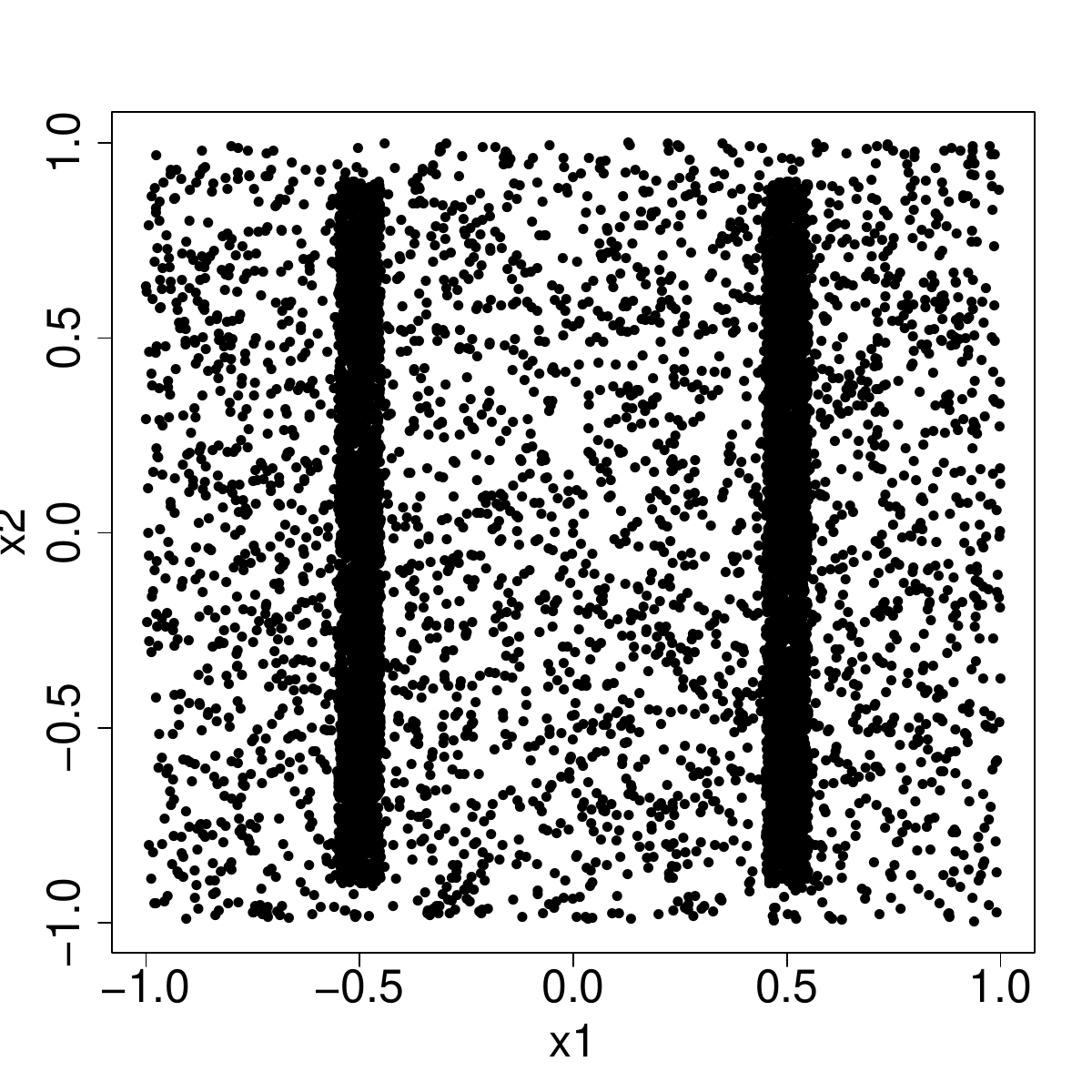}
	\includegraphics[width=0.24\textwidth]{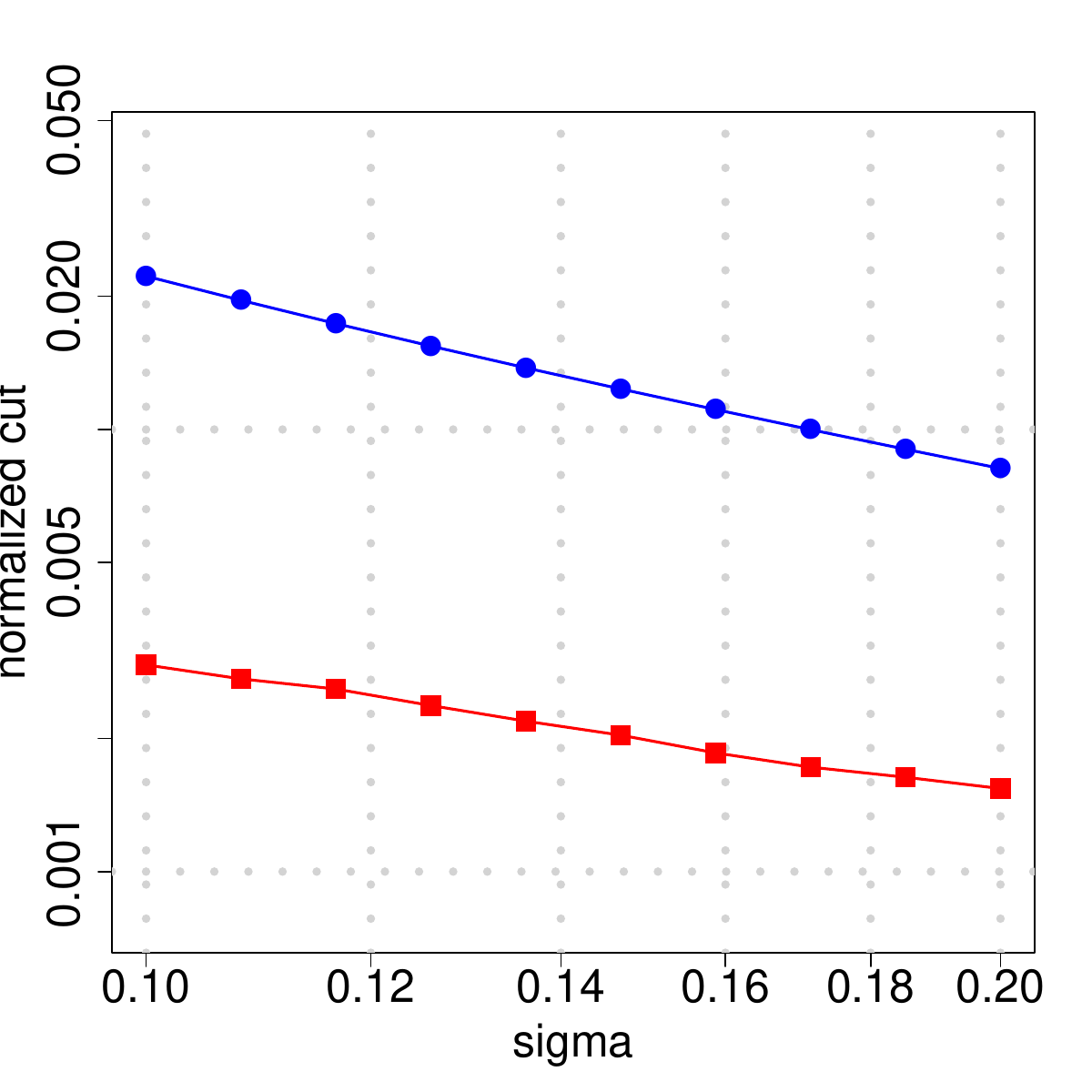}
	\includegraphics[width=0.24\textwidth]{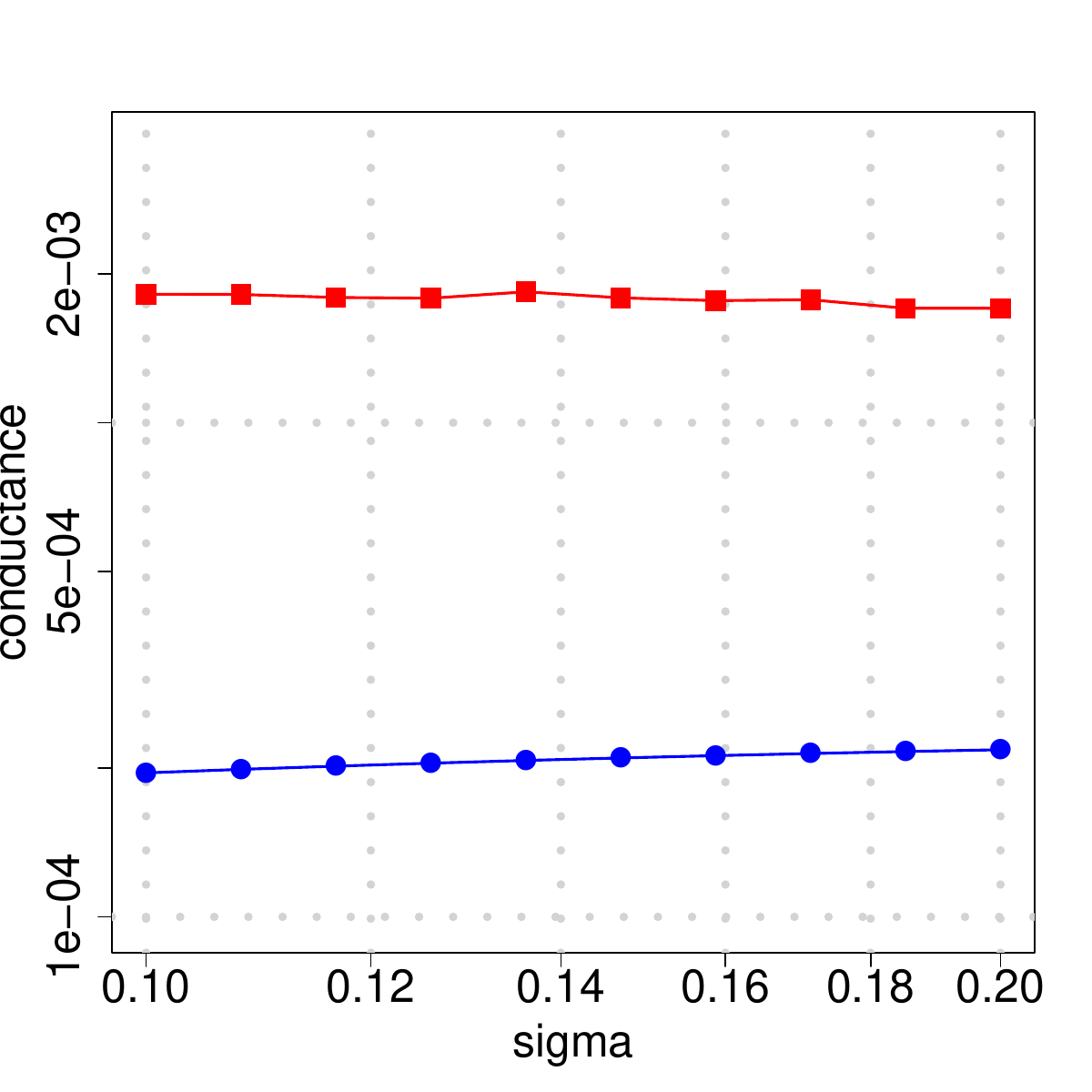}
	\includegraphics[width=0.24\textwidth]{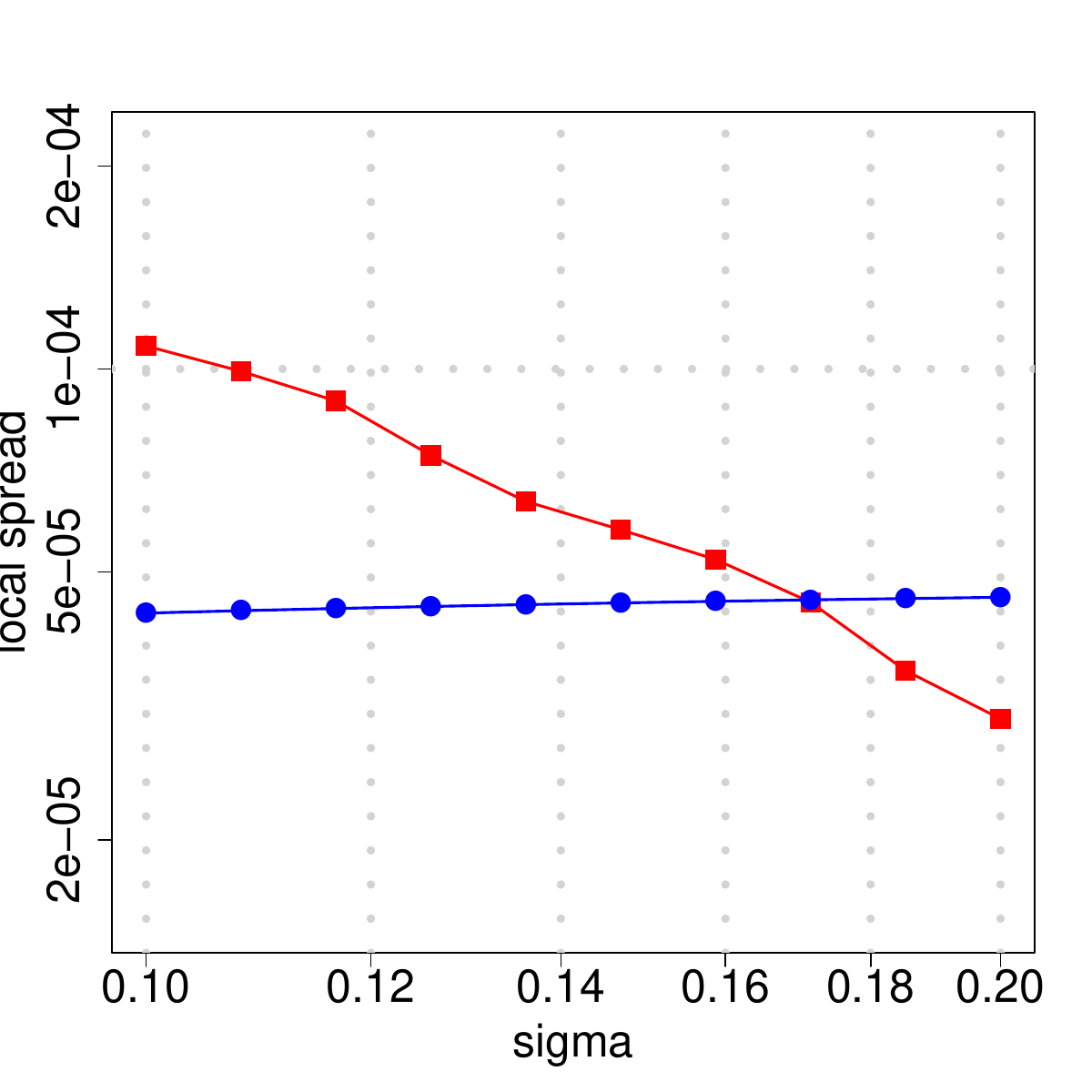}
	\includegraphics[width=0.24\textwidth]{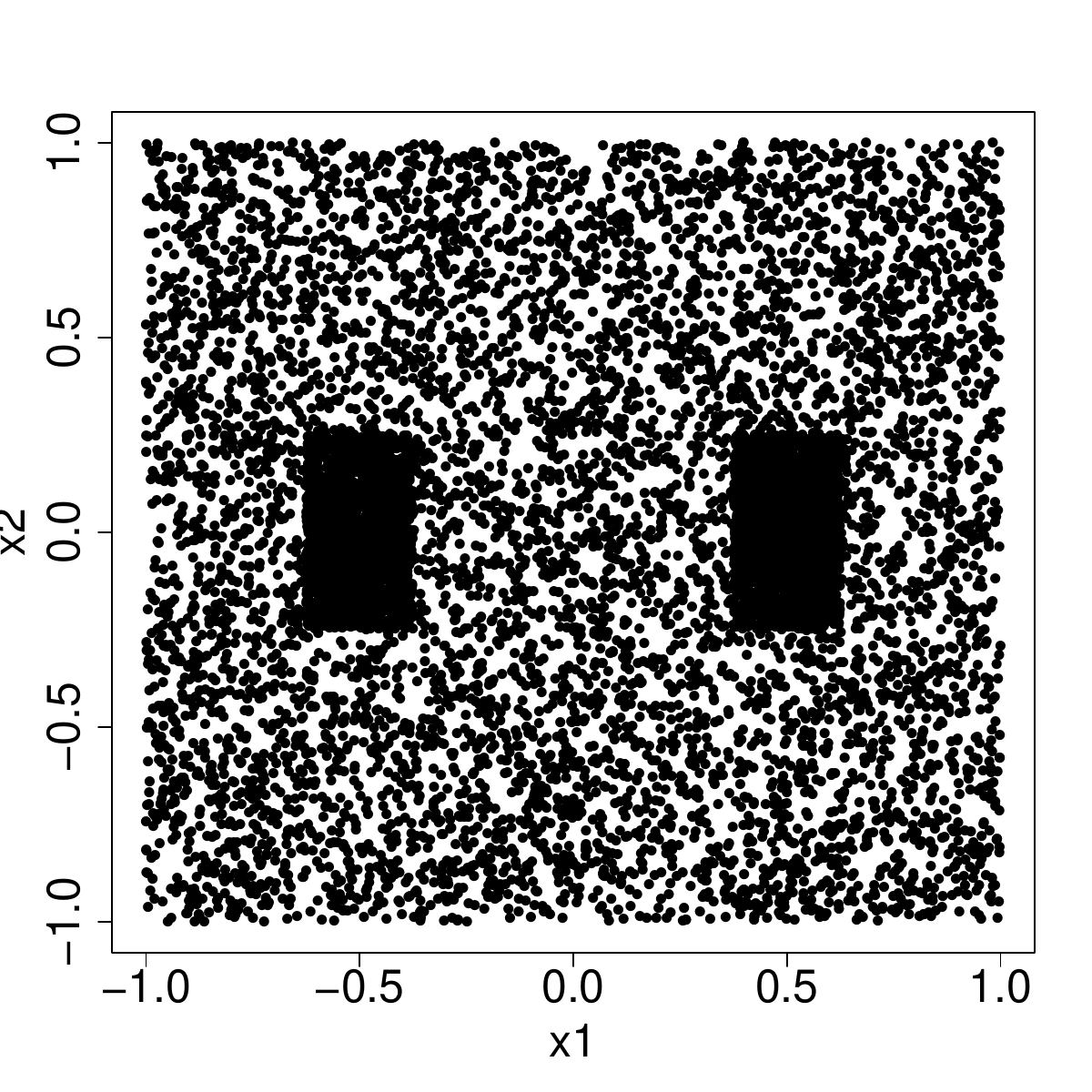}
	\includegraphics[width=0.24\textwidth]{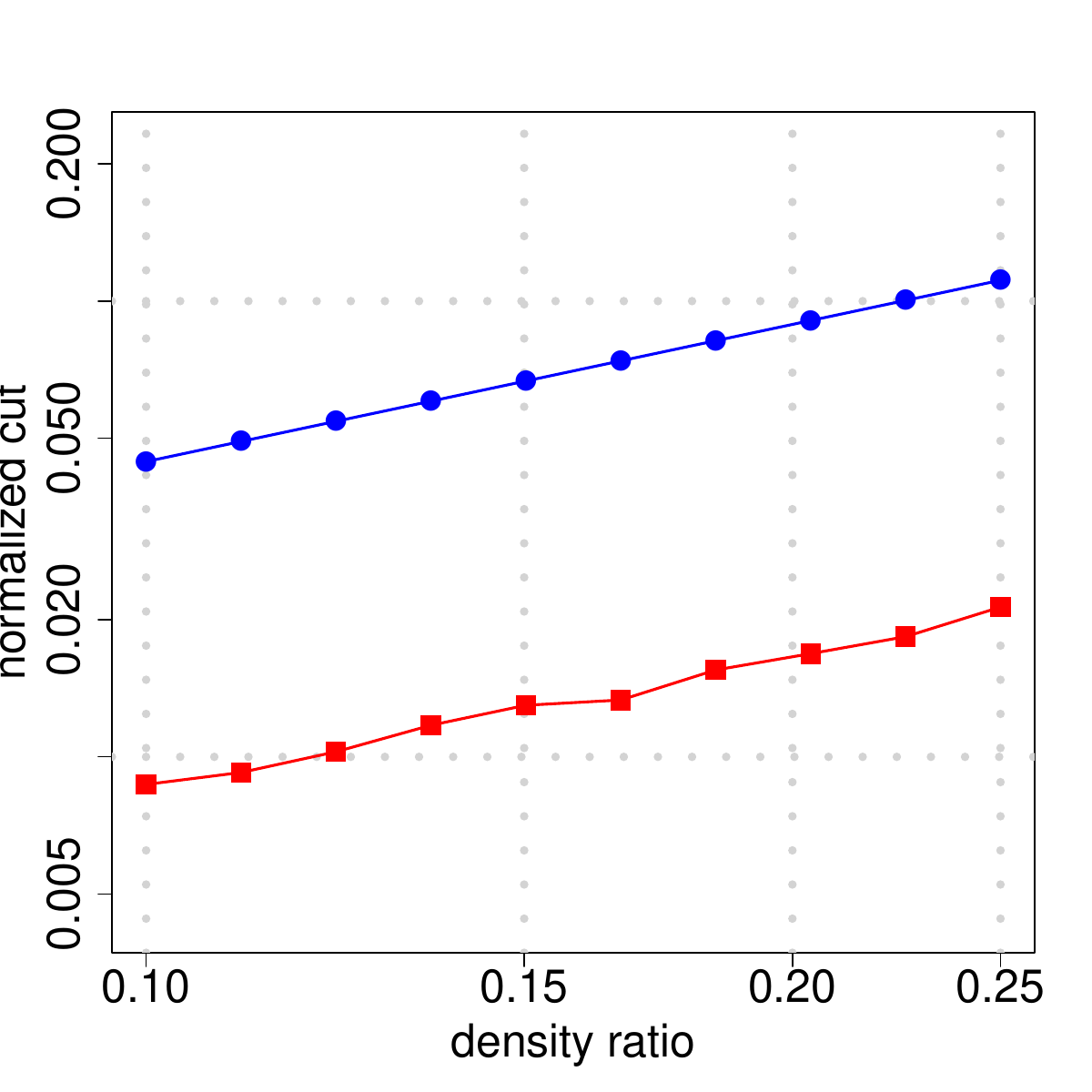}
	\includegraphics[width=0.24\textwidth]{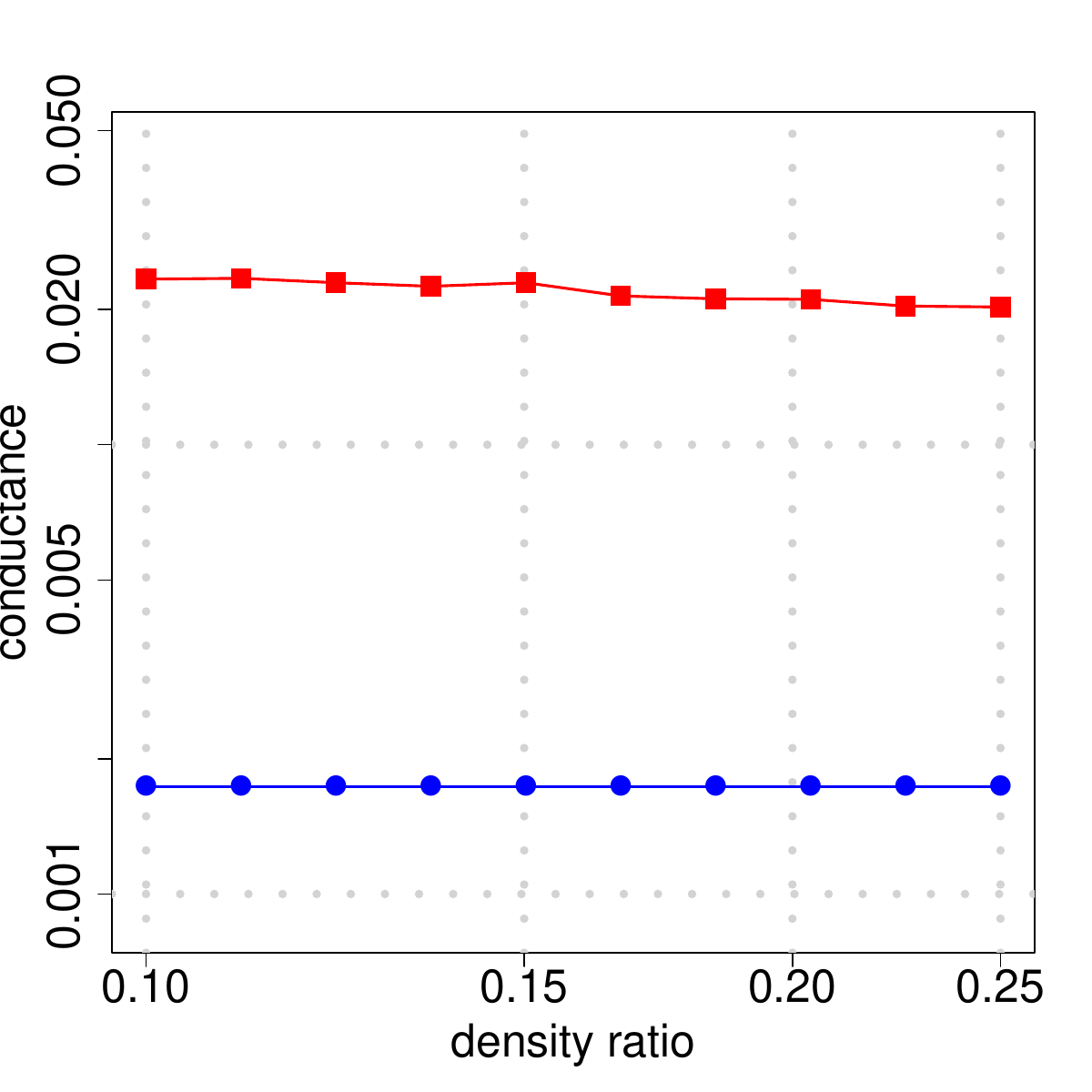}
	\includegraphics[width=0.24\textwidth]{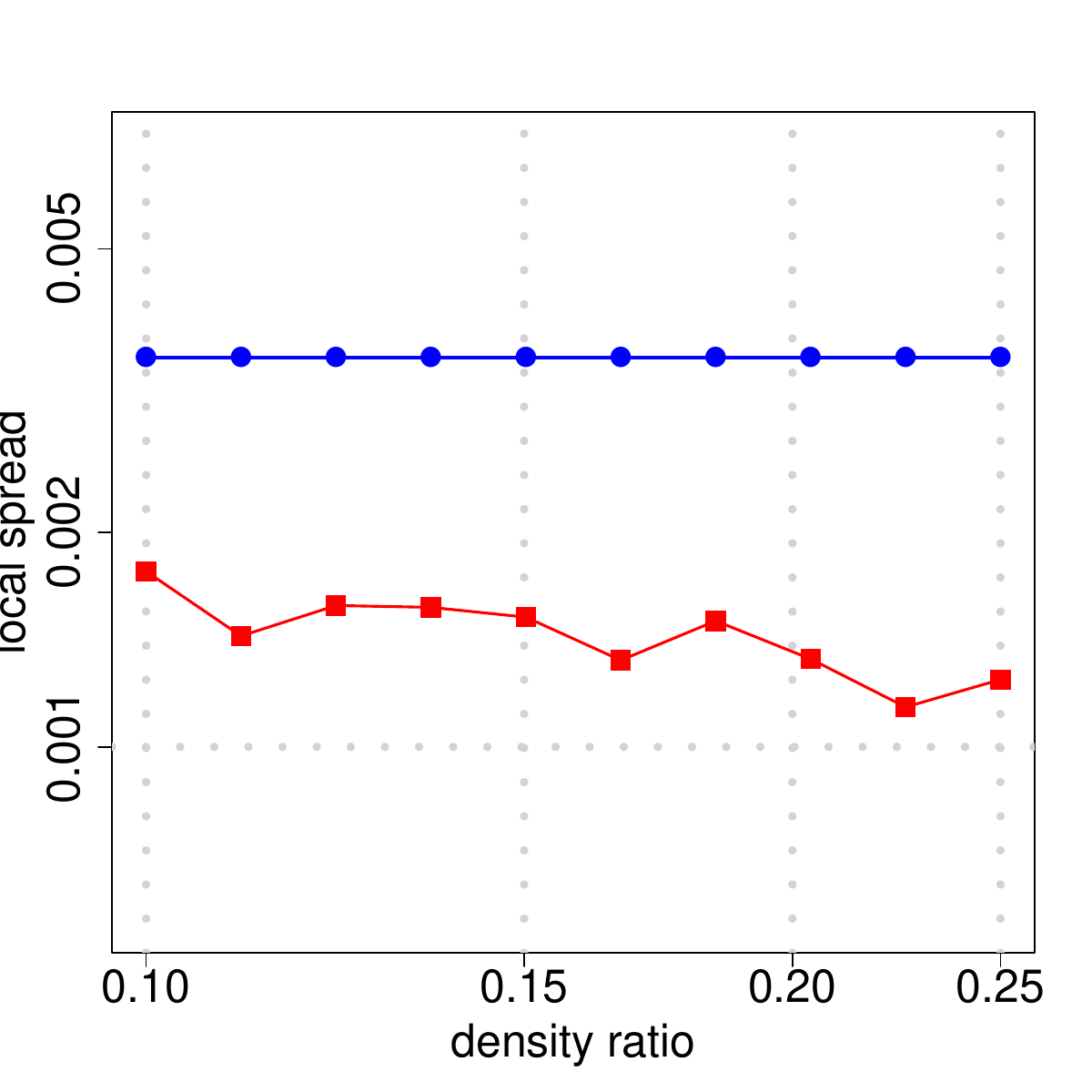}
	\caption{\small Empirical normalized cut, conductance, and local spread (in \textcolor{red}{red}), versus their theoretical bounds (in \textcolor{blue}{blue}). In the first row we vary the diameter $\rho$, in the second row we vary the thickness $\sigma$, and in the third row we vary the density ratio $(\lambda - \theta)/\lambda$. The first column shows $n = 8000$ samples for three different parameter values.}
	\label{fig:bounds}
\end{figure}

We investigate the tightness of Lemma~\ref{lem:density_cluster_local_spread} and Propositions~\ref{prop:density_cluster_normalized_cut} and \ref{prop:density_cluster_conductance}--- i.e., the bounds on population functionals required for the eventual density cluster recovery result in Theorem~\ref{thm:density_cluster_volume_ssd_ub}---via simulation. Figure \ref{fig:bounds} compares our bounds on normalized cut, conductance, and local spread of a density cluster with the actual empirically-computed quantities, when samples are drawn from a mixture of uniform distributions over rectangular clusters. In the first row we vary the diameter $\rho$ of the candidate cluster, in the second row we vary the width $\sigma$, and in the third row we vary the ratio $(\lambda - \theta)/\lambda$ of the density within and outside the cluster. In almost all cases, it is encouraging to see that our bounds track closely with their empirical counterparts, and are loose by roughly an order of magnitude at most. The one exception to this is the dependence of local spread on the width $\sigma$; this theoretical deficiency stems from a loose bound on the volume of sets with large aspect ratio (meaning $\rho/\sigma$ is much greater than $1$), but in any case the local spread contributes only $\log$ factors to the ultimate bound on cluster recovery. On the other hand, the looseness in each of these bounds will propagate to our eventual upper bound on \smash{$\Delta(\wh{C},\mc{C}_{\lambda,\sigma}[X])/\vol_{n,r}(\mc{C}_{\lambda,\sigma}[X])$}, which as a result is loose by several orders of magnitude. 

\subsection{Empirical Behavior of PPR} 
\label{subsec:empirical_behavior_ppr}
\begin{figure}
	\centering
	\includegraphics[width=0.24\textwidth,scale = .5]{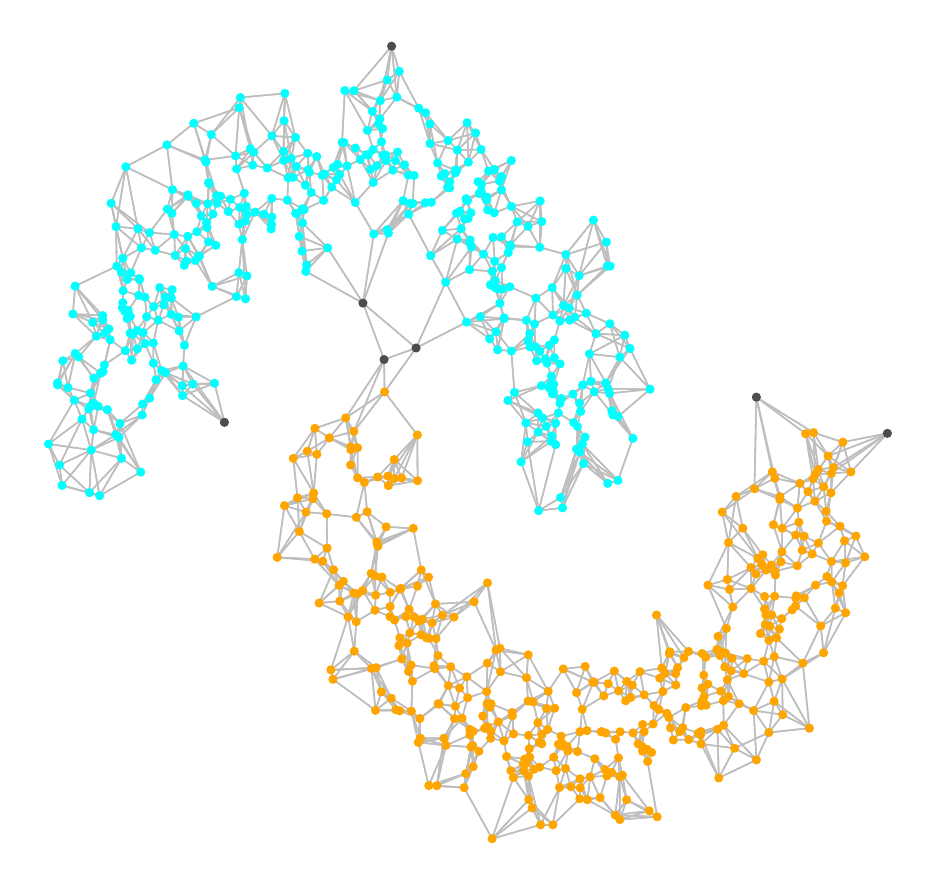}
	\includegraphics[width=0.24\textwidth]{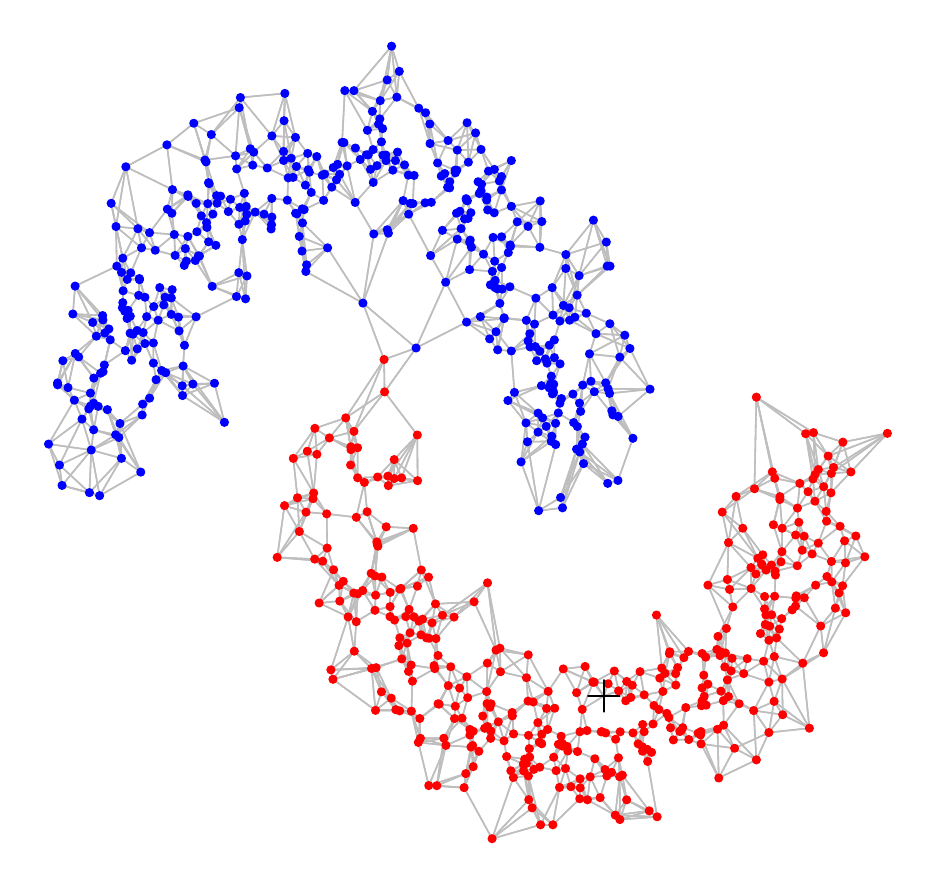}
	\includegraphics[width=0.24\textwidth]{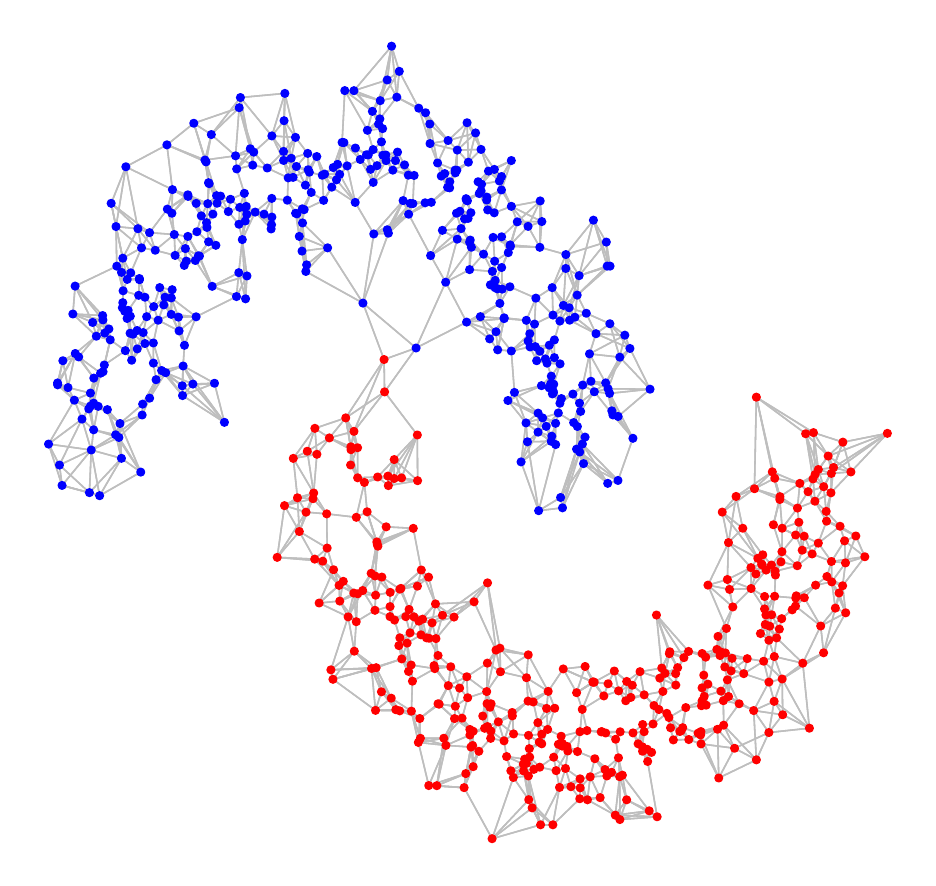}
	\includegraphics[width=0.24\textwidth]{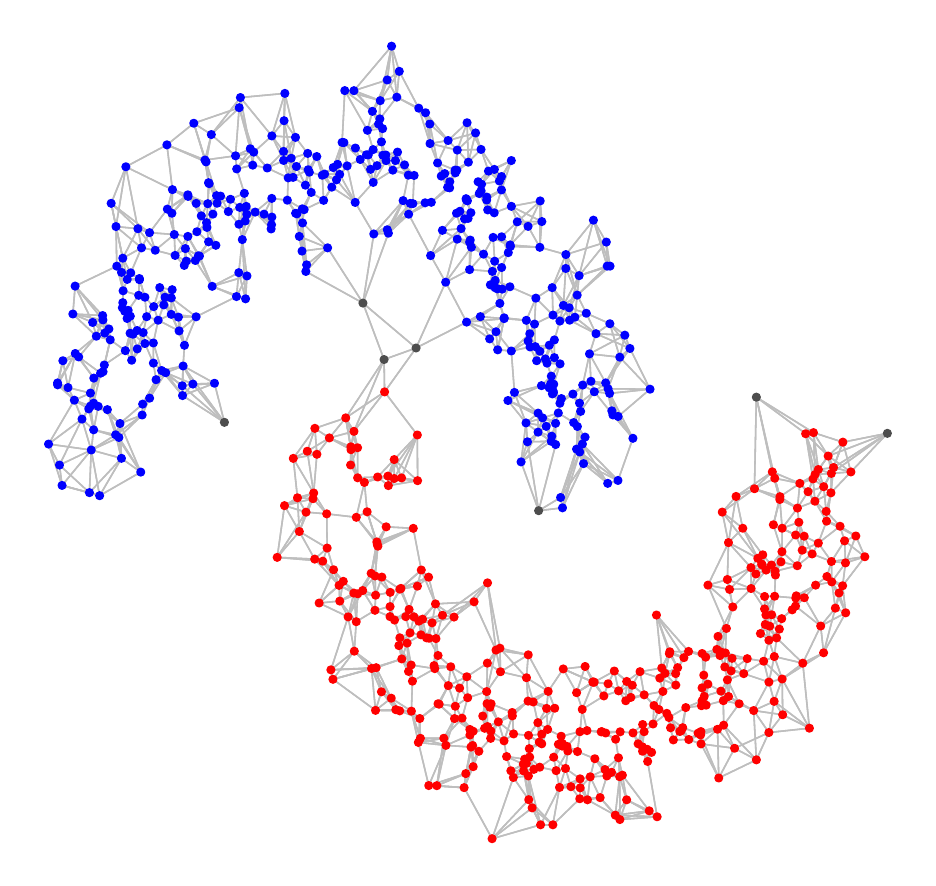}
	\includegraphics[width=0.24\textwidth]{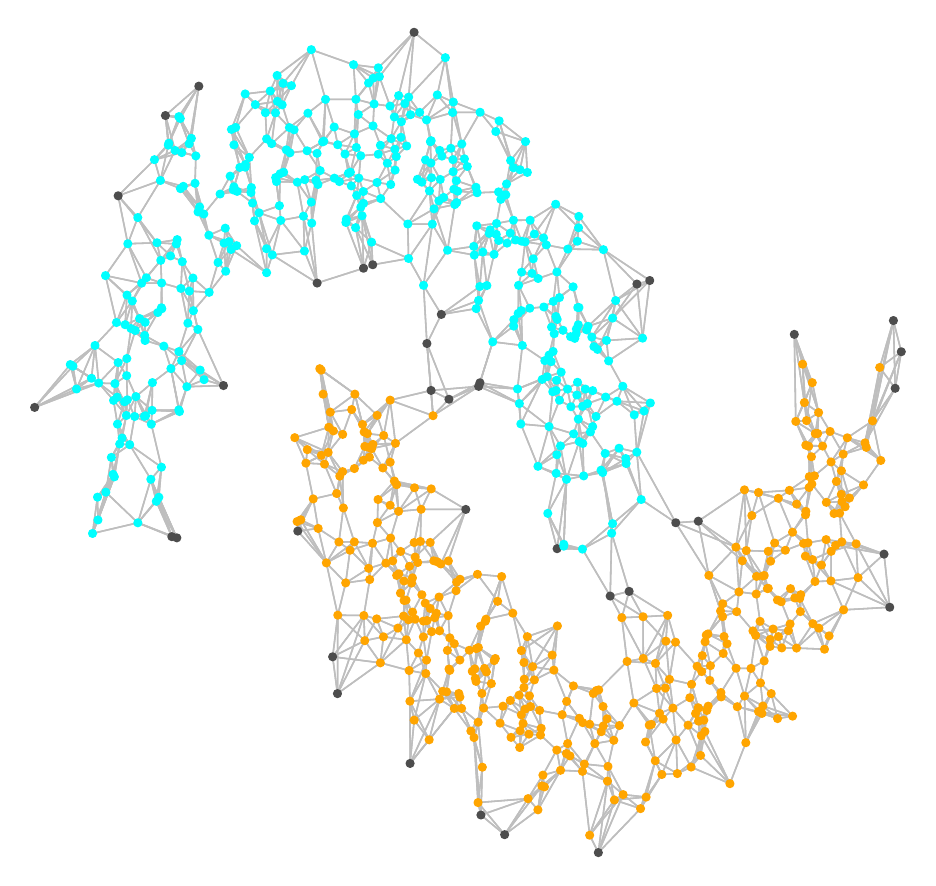}
	\includegraphics[width=0.24\textwidth]{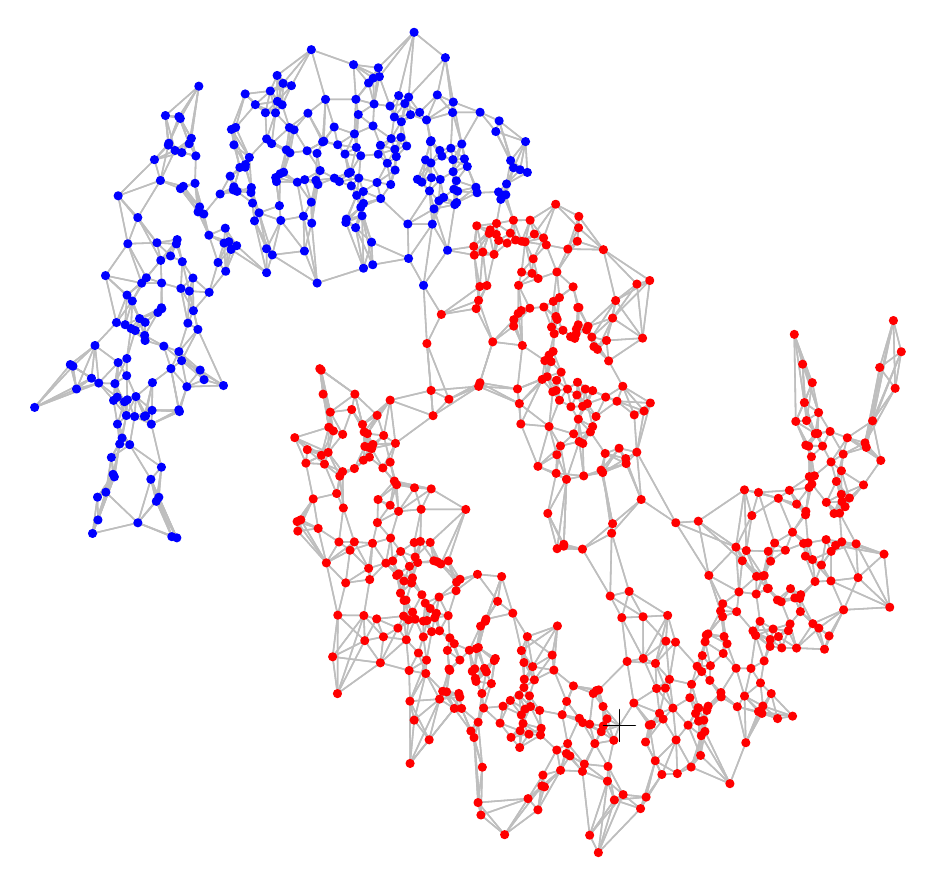}
	\includegraphics[width=0.24\textwidth]{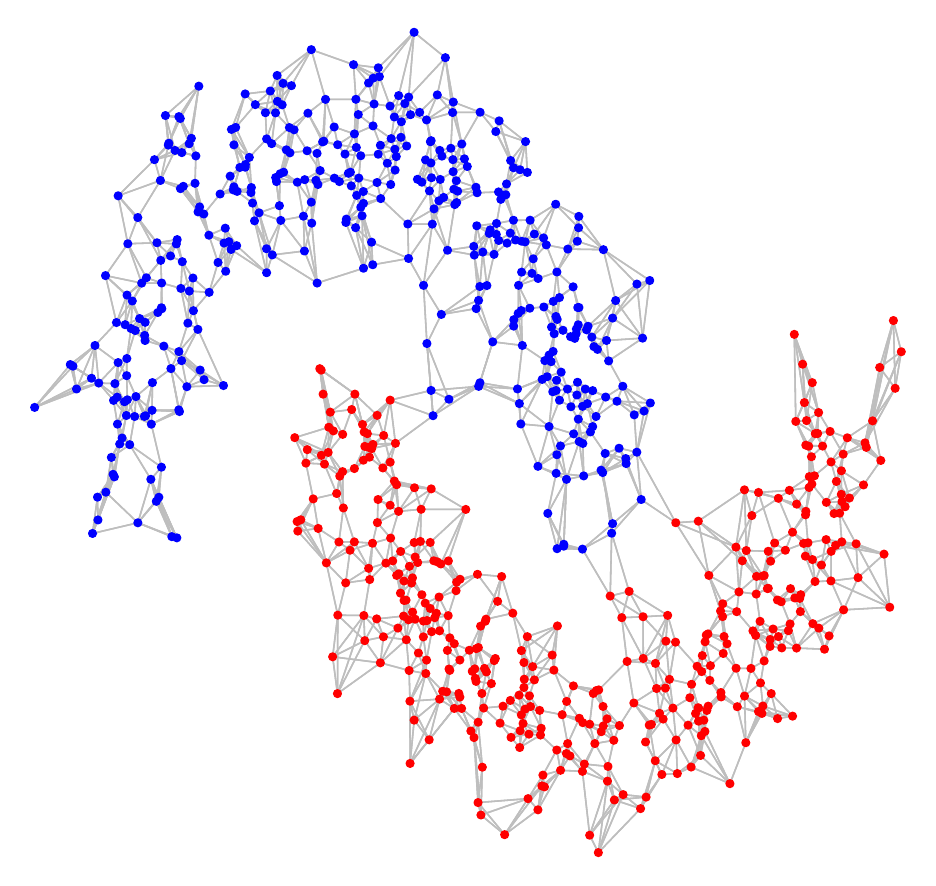}
	\includegraphics[width=0.24\textwidth]{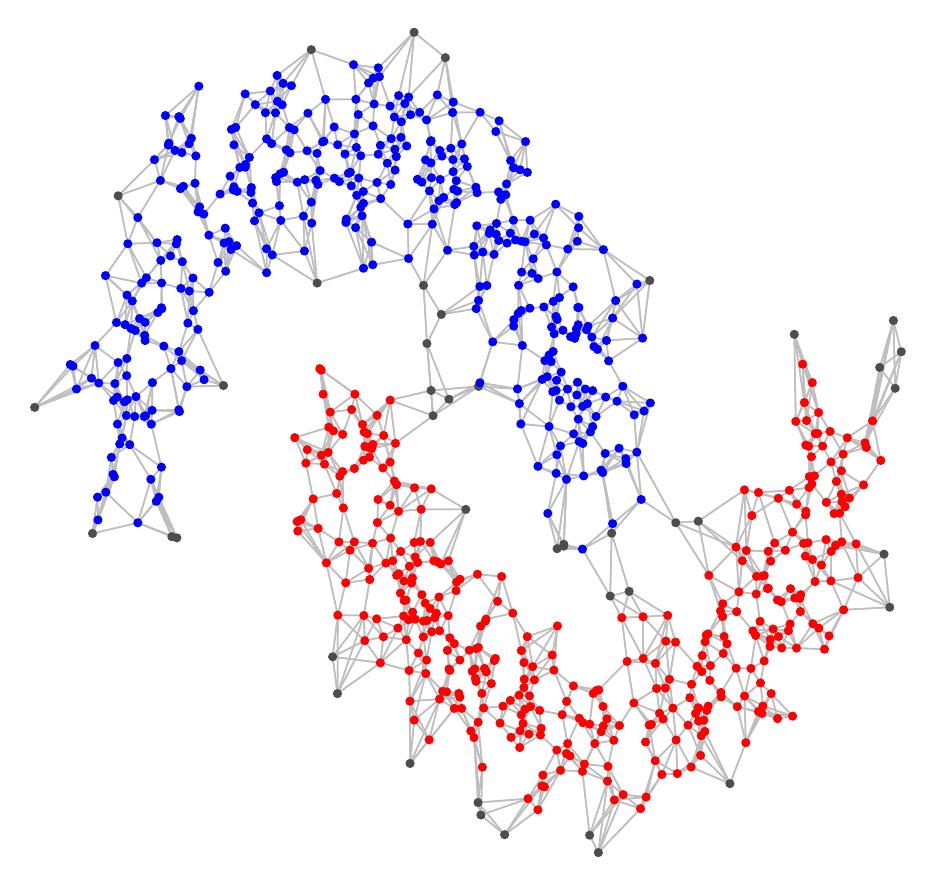}
	\includegraphics[width=0.24\textwidth]{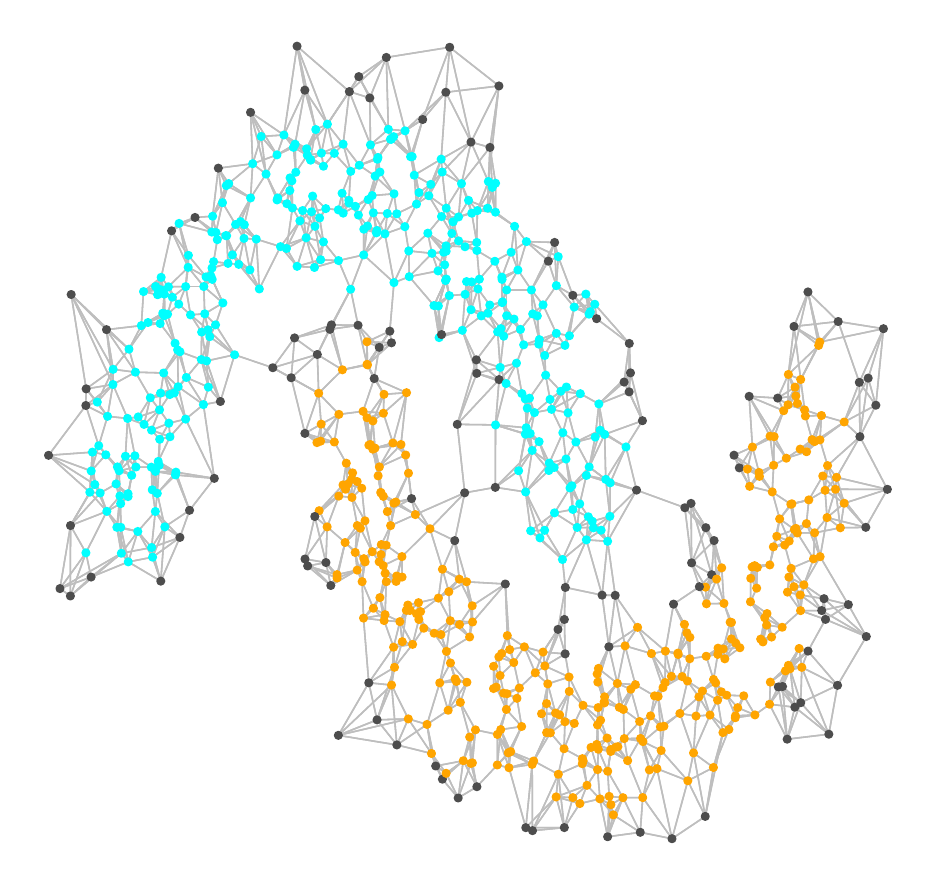}
	\includegraphics[width=0.24\textwidth]{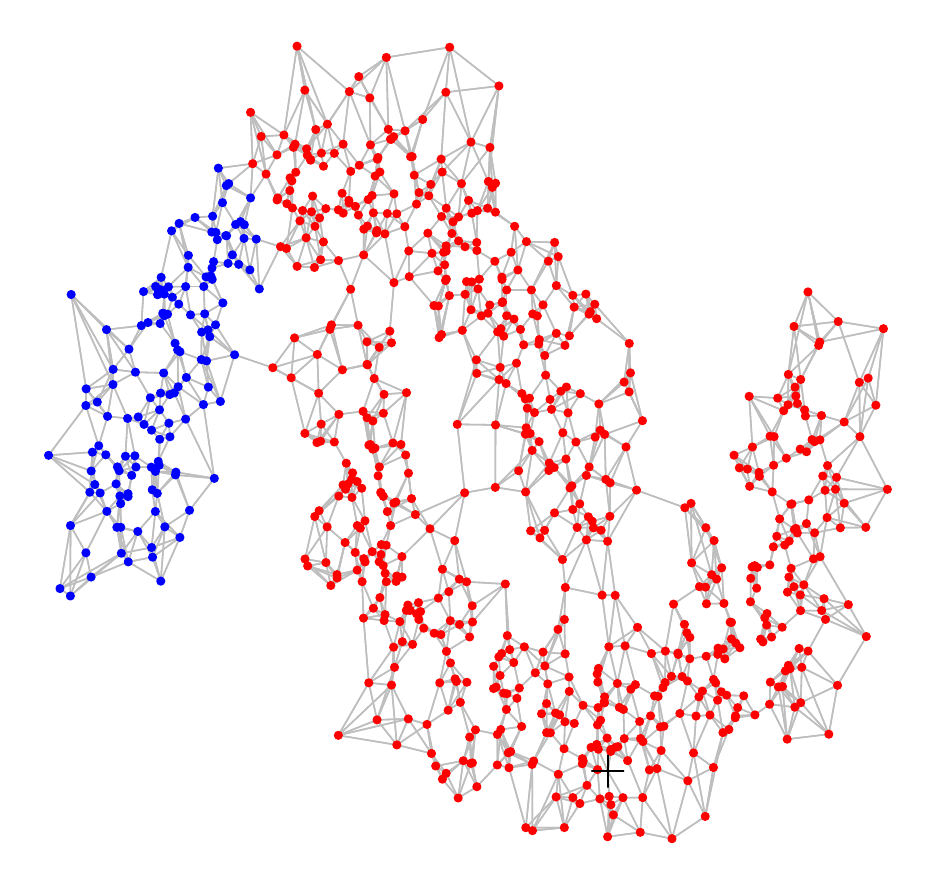}
	\includegraphics[width=0.24\textwidth]{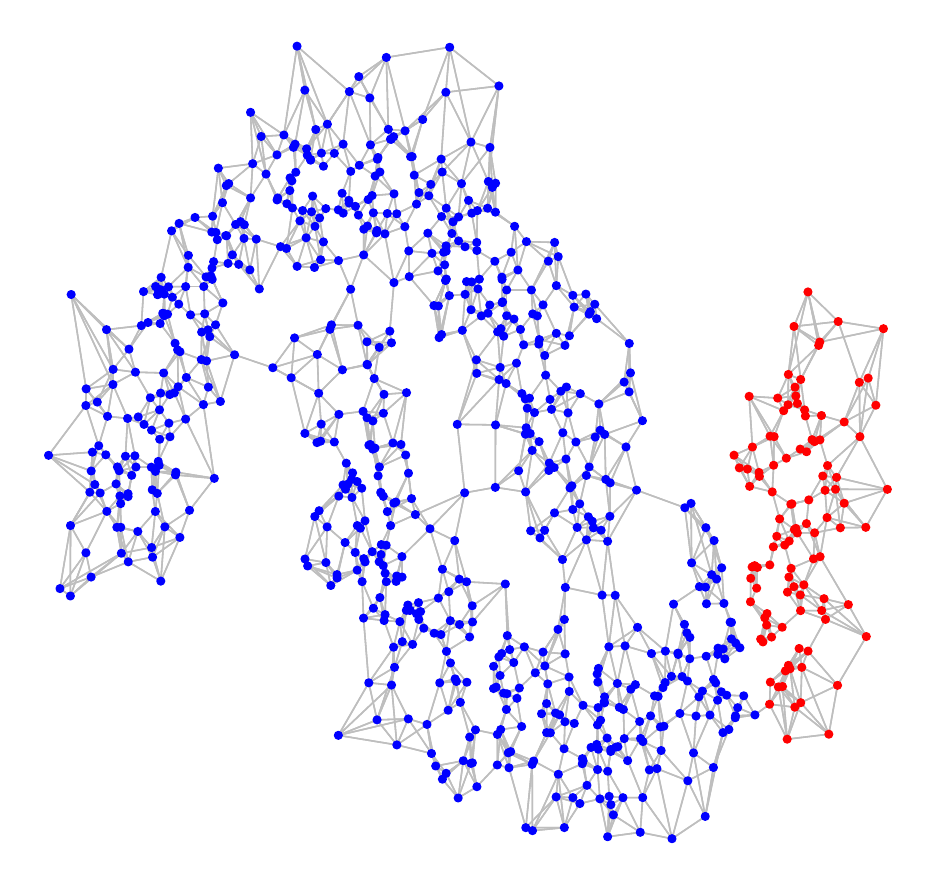}
	\includegraphics[width=0.24\textwidth]{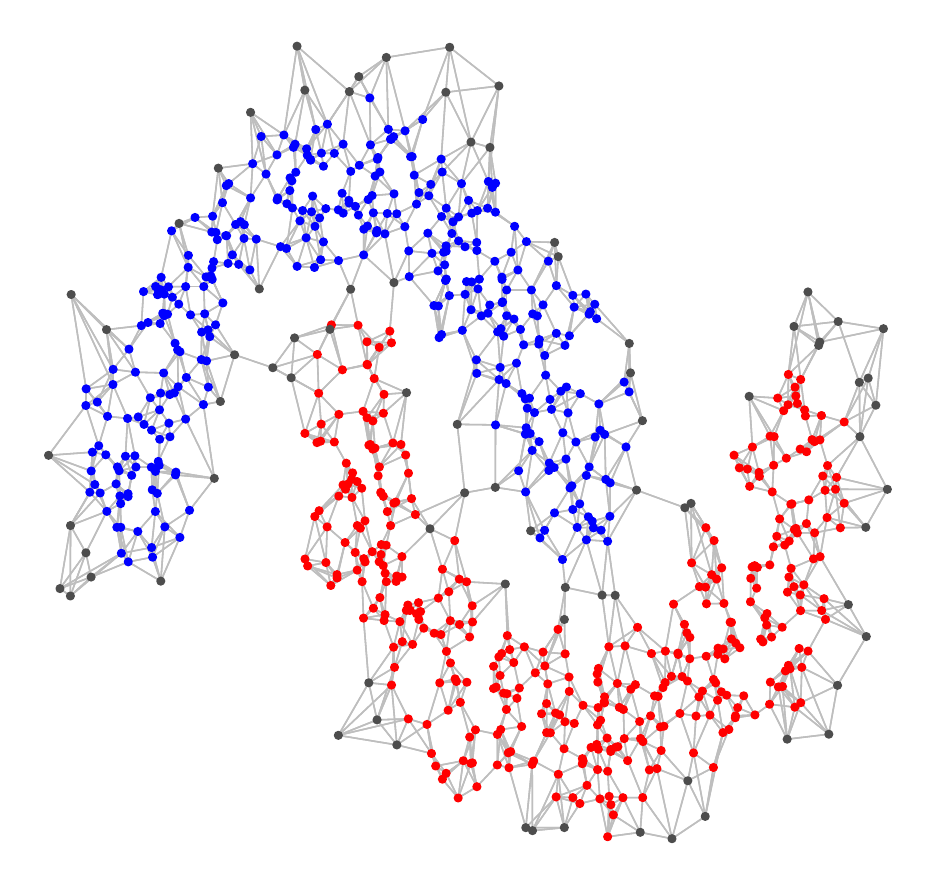}
	\caption{\small True density (column 1), PPR (column 2), minimum normalized
		cut (column 3) and estimated density (column 4) clusters for 3 different 
		simulated data sets. Seed node for PPR denoted by a black cross.} 
	\label{fig:moons}
\end{figure}

In Figure \ref{fig:moons}, to drive home the implications of Sections~\ref{sec:ppr_density_cluster} and~\ref{sec:lower_bound}, we compare the behavior of PPR and the density clustering algorithm of \citet{chaudhuri2010} on the well-known ``two moons'' dataset (with added 2d Gaussian noise), considered a prototypical success story for spectral clustering algorithms. We also examine the cluster which minimizes the normalized cut; as we have discussed previously, this can be seen as a middle ground between the geometric sensitivity of PPR, and the geometric insensitivity of density clustering. The first column shows the empirical density clusters $\mc{C}_{\lambda}[X]$ and $\mc{C}_{\lambda}'[X]$ for a particular threshold $\lambda$ of the density function; the second column shows the cluster recovered by PPR; the third column shows the global minimum normalized cut, computed according to the algorithm of \citet{bresson2013}; and the last column shows a cut of the density cluster tree estimator of \citet{chaudhuri2010}. Each row corresponds to a different separation between the two moons. In the second row, we see that as the two moons become less well-separated, PPR becomes unable to recover the density clusters, but normalized cut still succeeds in doing so. In the third row, we see that the Chaudhuri-Dasgupta algorithm succeeds even when both PPR and normalized cut fail. This supports one of our main messages, which is that PPR recovers only geometrically well-conditioned density clusters.

\section{Discussion}
\label{sec:discussion}
In this work, we have analyzed the behavior of PPR in the classical setup of nonparametric statistics. We have shown how PPR depends on the distribution $\Pbb$ through the population normalized cut, conductance, and local spread, and established upper bounds on the error with which PPR recovers an arbitrary candidate cluster $\mc{C} \subseteq \Rd$.  In the particularly important case where $\mc{C} = \mc{C}_{\lambda}$ is a $\lambda$-density cluster, we have shown that PPR recovers $\mc{C}_{\lambda}$ if and only if both the density cluster and density are well-conditioned. We now conclude by summarizing a couple of interesting directions for future work.

Letting the radius of the neighborhood graph shrink, $r \to 0$ as $n \to \infty$, would be computationally attractive, as it would ensure that the graph $G_{n,r}$ is sparse. However, the bounds~\eqref{eqn:volume_ssd_ub} and~\eqref{eqn:density_cluster_volume_ssd_ub} will blow up as the radius $r$ goes to $0$, preventing us from making claims about the behavior of PPR in this regime. Although the restriction to a kernel function fixed in $n$ is common in spectral clustering theory \citep{vonluxburg2008, schiebinger2015, singer2017}, recent works~\citep{shi2015, calder2019, garciatrillos18, garciatrillos2020, yuan2020} have demonstrated that spectral methods have meaningful continuum limits when $r \to 0$ as $n \to \infty$, and given precise rates of convergence.~\cite{garciatrillos19} have applied these results to analyze global spectral clustering in the nonparametric mixture model, obtaining asymptotic upper bounds that do not depend on $r$; it seems plausible that similar bounds could be obtained for local spectral clustering with PPR, although the arguments would necessarily be quite different.

In another direction, it would be very useful to find reasonable conditions under which the ratio \smash{$\Delta(\wh{C},\mc{C}[X])/\vol_{n,r}(\mc{C}[X])$} would tend to $0$ as $n \to \infty$. It seems likely that such a strong result would entail bounds on the $\Leb^{\infty}$-error of PPR. Though most results thus far derive bounds only on the $\Leb^1$- or $\Leb^2$-error of spectral clustering methods, some recent works~\citep{dunson2020,calder2020} have established $\Leb^{\infty}$-bounds on the error with which the eigenvectors of a graph Laplacian matrix approximate the eigenvectors of a weighted Laplace-Beltrami operator. It is not clear whether the techniques used in these works can be applied to PPR.

\section*{Acknowledgments}
	SB is grateful to Peter Bickel, Martin Wainwright, and Larry Wasserman for early inspiring conversations. This work was supported in part by the NSF grant DMS-1713003.

\bibliographystyle{plainnat}

\clearpage

\appendix
\noindent The proofs of our major theorems largely consist of (at most) three modular parts.
\begin{enumerate}
	\item \textbf{Fixed graph results.} Results which hold with respect to an arbitrary graph $G$, and are stated with respect to functionals (i.e. normalized cut, conductance, and local spread) of $G$;
	\item \textbf{Sample-to-population results.} For the specific choice $G = G_{n,r}$, results relating the aforementioned functionals to their population analogues. 
	\item \textbf{Bounds on population functionals.} (In the case of density clustering only.) When the candidate cluster is a $\lambda$-density cluster, bounds on population functionals as a function of $\lambda$, as well as the other relevant parameters introduced in Section~\ref{sec:ppr_density_cluster}.
\end{enumerate}

Appendices \ref{apdx:fixed_graph}-\ref{apdx:density_cluster_population_functionals} will correspond to each of these three parts. In Appendix~\ref{apdx:pf_major_theorems}, we will combine these parts to prove the major theorems of our main text, Theorems~\ref{thm:volume_ssd_ub} and~\ref{thm:density_cluster_volume_ssd_ub}, as well as our negative result, Theorem~\ref{thm:ppr_lb}. In Appendix~\ref{apdx:appr_misclassification_error} we derive upper bounds for the aPPR vector, and show that under certain conditions the PPR vector can perfectly separate two density clusters. Finally, in Appendix~\ref{apdx:experimental_details} we give relevant details regarding our experiments.

\section{Fixed Graph Results}
\label{apdx:fixed_graph}
In this section, we give all results that hold with respect to an arbitrary graph $G$. For the convenience of the reader, we begin by reviewing some notation from the main text, and also introduce some new notation. \\

\paragraph{Notation.}
The graph $G = (V,E)$ is an undirected and connected but otherwise arbitrary graph, defined over vertices $V = \{1,\ldots,n\}$ with $m = |E|$ total edges. The adjacency matrix of $G$ is $A$, the degree matrix is $D$, and the lazy random walk matrix over $G$ is $W = (I + D^{-1}A)/2$. The lazy random walk originating at node $v \in V$ has distribution $q(v,t;G) = e_vW^t$ after $t$ steps; we use the notational shorthand $q_v^{(t)} := q(v,t;G)$. The stationary distribution of the lazy random walk is $\pi := \pi(G) := \lim_{t \to \infty} q_v^{(t)}$ is given by $\pi(u) = \deg(u;G)/\vol(u;G)$.

For a starting distribution $s$ (by distribution we mean a vector with non-negative entries), the PPR vector $p_s = p(s,\alpha;G)$ is the solution to
\begin{equation}
\label{eqn:ppr}
p_s = \alpha s + (1 - \alpha) p_s W.
\end{equation}
When $s = e_v$, we write $p_v := p_{e_v}$. It is easy to check that $p_s = \alpha \sum_{t = 0}^{\infty} (1 - \alpha)^t q_s^{(t)}$.  Note that $s$ need not be a probability distribution (i.e. its entries need not sum to $1$) to make sense of~\eqref{eqn:ppr}.

Given a distribution $q$ (for instance, $q = q_v^{(t)}$ for $t \in \mathbb{N}$, $q = p_v$, or $q = \pi$) and $\beta \in (0,1)$, the $\beta$-sweep cut of $q$ is
\begin{equation*}
S_{\beta}(q) = \set{u: \frac{q(u)}{\deg(u;G)} > \beta};
\end{equation*} 
in the special case where $q = p_v$ we write $S_{\beta,v}$ for $S_{\beta}(p_v)$. The argument of $S_{\beta}(\cdot)$ will usually be clear from context, in which case we will drop it and simply write $S_{\beta}$. For $j = 1,\ldots,n$, let $\beta_j$ be the smallest value of $\beta \in (0,1)$ such that the sweep cut $S_{\beta_j}$ contains at least $j$ vertices. For notational ease, we will write $S_j := S_{\beta_j}$, and $S_0 = \emptyset$. 

We now introduce the~\emph{Lovasz-Simonovits curve} $h_q(\cdot): [0,2m] \to [0,1]$ to measure the extent to which a distribution $q$ is mixed. To do so, we first define a piecewise linear function $q[\cdot]: [0,2m] \to [0,1]$. Letting $q(S) := \sum_{u \in S} q(u)$, we take $q[\vol(S_j)] = q(S_j)$ for each sweep cut $S_j$, and then extend $q[\cdot]$ by piecewise linear interpolation to be defined everywhere on its domain. Then the mixedness of $q$ is measured by
\begin{equation*}
h_q(k) := q[k] - \frac{k}{2m}.
\end{equation*}
The Lovasz-Simonovits curve is a non-negative function, with $h_q(0) = h_q(2m) = 0$. The stationary distribution $\pi$ is mixed, i.e. $h_{\pi}(k) = 0$ for all $k \in [0,2m]$. Finally, both $q[\cdot]$ and $h_q(\cdot)$ are concave functions, which will be an important fact later on.  

The conductance of $V$ is abbreviated as $\Psi(G) := \Psi(V;G)$, and likewise for the local spread $s(G) := s(V;G)$. Finally, for convenience we introduce the following functionals:
\begin{equation*}
\begin{aligned}
& d_{\max}(C; G) := \max_{u \in C} \deg(u; G), && d_{\min}(C; G) := \min_{u \in C} \deg(u;G) \\
& d_{\max}(G) := d_{\max}(V;G),~~ && d_{\min}(G) := d_{\min}(V; G)
\end{aligned}
\end{equation*}
We note that $d_{\min}(G)^2 \leq d_{\min}(G) \cdot n \leq \vol(G) \leq d_{\max}(G) \cdot n$, and that for any $S \subseteq V$, $|S| \cdot d_{\min}(G) \leq \vol(S;G)$ (where $|S|$ is the cardinality of $S$.) \\

\emph{Organization.} In the following sections we establish: (Section \ref{subsec:pf_lem_zhu}) an upper bound on the misclassification error of PPR in terms of $\alpha$ and $\Phi(C;G)$ (Lemma~\ref{lem:zhu}), and an analogous result for aPPR (Corollary~\ref{cor:zhu}); (Section \ref{subsec:ppr_uniform_bounds}) a uniform bound on the perturbations of the PPR vector, to be used later in the proof of Theorem~\ref{thm:density_cluster_consistent_recovery} (consistency of PPR);  (Section \ref{subsec:lovasz_simonovits_bounds}) upper bounds on the mixedness of $q_v^{(t)}$ (as a function of $t$) and $p_v$ (as a function of $\alpha$), which will be helpful in the proofs of Proposition~\ref{prop:pointwise_mixing_time} and Theorem~\ref{thm:ppr_lb}; (Section \ref{subsec:pf_prop_pointwise_mixing_time}) an upper bound on $\tau_{\infty}(G)$ in terms of $\Psi(G)$ and $s(G)$ (Proposition~\ref{prop:pointwise_mixing_time}); and (Section \ref{subsec:ppr_spectral_partitioning}) an upper bound on the normalized cut $\Phi(\wh{C};G)$ in terms of $\Phi(C;G)$, to be used later in the proof of Theorem~\ref{thm:ppr_lb} (negative example). 

\subsection{Misclassification Error of Clustering with PPR and aPPR}
\label{subsec:pf_lem_zhu}
For a candidate cluster $C \subseteq V$, we use the tilde-notation $\wt{G} = G[C]$ to refer to the subgraph of $G$ induced by $C$. Similarly we write $\wt{q}_v^{(t)} := q(v,t;\wt{G})$ for the $t$-step distribution of the lazy random walk over $\wt{G}$, $\wt{\pi} = \pi(G[C])$ for the stationary distribution of $\wt{q}_v^{(t)}$ (we will always assume $G[C]$ is connected), and $\wt{p}_v := p(v,\alpha;\wt{G})$ for the PPR vector over $\wt{G}$. \\

\emph{Proof of Lemma~\ref{lem:zhu}.}
As mentioned in the main text, Lemma~\ref{lem:zhu} is equivalent, up to constants, to Lemma~3.4 in \cite{zhu2013}, and the proof of Lemma~\ref{lem:zhu} proceeds along very similar lines to the proof of that lemma. In fact, we directly use the following three inequalities, derived in that work:
\begin{itemize}
	\item \textbf{(c.f. Lemma 3.2 of \cite{zhu2013})} For any seed node $v \in C$, the PPR vector is lower bounded,
	\begin{equation}
	\label{pf:zhu1}
	\wt{p}_v(u) \geq \frac{3}{4}\bigl(1 - \alpha \cdot \tau_{\infty}(\wt{G})\bigr) \cdot \wt{\pi}(u),~~\textrm{for every $u \in C$.}
	\end{equation}
	\item \textbf{(c.f. Corollary 3.3 of \cite{zhu2013})} For any seed node $v \in C$, there exists a so-called leakage distribution $\ell = \ell(v)$ such that $\mathrm{supp}(\ell) \subseteq C$, $\|\ell\|_1 \leq 2\Phi(C;G)/\alpha$, and 
	\begin{equation}
	\label{pf:zhu2}
	p_v(u) \geq \wt{p}_v(u) - \wt{p}_{\ell}(u),~~\textrm{for every $u \in C$.}
	\end{equation}
	\item \textbf{(c.f. Lemma 3.1 of \cite{zhu2013})} There exists a set $C^g \subset C$ with $\vol(C^g;G) \geq \frac{1}{2}\vol(C;G)$ such that for any seed node $v \in C^g$, the following inequality holds
	\begin{equation}
	\label{pf:zhu3}
	p_v(C^c) \leq 2\frac{\Phi(C;G)}{\alpha}.
	\end{equation}
\end{itemize}
We use~\eqref{pf:zhu1}-\eqref{pf:zhu3} to separately upper bound $\vol(S_{\beta,v} \setminus C;G)$, $\vol(C^{\mathrm{int}} \setminus S_{\beta,v};G)$ and $\vol(C^{\mathrm{bdry}} \setminus S_{\beta,v};G)$; here $C^{\mathrm{int}} \cup C^{\mathrm{bdry}} = C$ is a partition of $C$, with
\begin{equation*}
C^{\mathrm{int}} := \Bigl\{u \in C: \deg(u;\wt{G}) > \bigl(1 - \alpha \cdot \beta \cdot \vol(C;G)\bigr) \deg(u;G) \Bigr\},
\end{equation*}
consisting of those vertices $u \in C$ with sufficient large degree in $\wt{G}$. 

First we upper bound $\vol(S_{\beta,v} \setminus C;G)$. Observe that for any $u \in S_{\beta,v} \setminus C$, $p_v(u) > \beta \cdot \deg(u;G)$. Summing up over all such vertices, from~\eqref{pf:zhu3} we conclude that
\begin{equation}
\label{pf:zhu3.5}
\vol(S_{\beta,v} \setminus C; G) \leq \frac{p_v(C^c)}{\beta} \leq 2\frac{\Phi(C;G)}{\beta \cdot \alpha}.
\end{equation} 
Next we upper bound $\vol(C^{\mathrm{int}} \setminus S_{\beta,v};G)$. From~\eqref{pf:zhu1} and~\eqref{pf:zhu2} we see that 
\begin{equation*}
p_v(u) \geq \frac{3}{4}\bigl(1 - \alpha \cdot \tau_{\infty}(\wt{G})\bigr) \cdot \wt{\pi}(u) - \wt{p}_{\ell}(u)~~\textrm{for all $u \in C$.}
\end{equation*}
If additionally $u \not\in S_{\beta,v}$ then $p_v(u) \leq \beta \deg(u;G)$, and for all such $u \in C \setminus S_{\beta,v}$,
\begin{equation}
\label{pf:zhu4}
\frac{3}{4}\bigl(1 - \alpha \cdot \tau_{\infty}(\wt{G})\bigr) \cdot \wt{\pi}(u) -  \beta\deg(u;G) \leq \wt{p}_{\ell}(u).
\end{equation}
On the other hand, for any $u \in C^{\mathrm{int}}$ it holds that
\begin{equation*}
\wt{\pi}(u) = \frac{\deg(u;\wt{G})}{\vol(\wt{G})} \geq \frac{\deg(u;\wt{G})}{\vol(G)} \geq \frac{(1 - \alpha \beta \vol(C;G))\deg(u;G)}{\vol(C;G)};
\end{equation*}
by plugging this in to~\eqref{pf:zhu4} we obtain
\begin{equation*}
\biggl(\frac{3(1 - \alpha \beta \vol(C;G))\cdot\bigl(1 - \alpha \tau_{\infty}(\wt{G})\bigr)}{4\vol(C;G)} - \beta\biggr) \cdot \deg(u;G) \leq \wt{p}_{\ell}(u),~~\textrm{for all $u \in C^{\mathrm{int}} \setminus S_{\beta,v}$};
\end{equation*}
and summing over all such $u$ gives
\begin{equation*}
\biggl(\frac{3(1 - \alpha \beta \vol(C;G))\cdot\bigl(1 - \alpha \tau_{\infty}(\wt{G})\bigr)}{4\vol(C;G)} - \beta\biggr) \cdot \vol\bigl(C^{\mathrm{int}} \setminus S_{\beta,v}; G\bigr) \leq \wt{p}_{\ell}\bigl(C^{\mathrm{int}} \setminus S_{\beta,v}\bigr) \leq 2\frac{\Phi(C;G)}{\alpha}.
\end{equation*}
The upper bounds on $\alpha$ and $\beta$ in~\eqref{eqn:zhu_condition} imply
\begin{equation*}
\biggl(\frac{3(1 - \alpha \beta \vol(C;G))\cdot\bigl(1 - \alpha \tau_{\infty}(\wt{G})\bigr)}{4\vol(C;G)} - \beta\biggr) \geq \frac{2}{3}\beta,
\end{equation*}
and we conclude that
\begin{equation}
\label{pf:zhu5}
\vol(C^{\mathrm{int}} \setminus S_{\beta,v}; G) \leq \frac{3\Phi(C;G)}{\alpha\beta}.
\end{equation}
Finally, we upper bound $\vol(C^{\mathrm{bdry}} \setminus S_{\beta,v};G)$. Indeed, for any $u \in C^{\mathrm{bdry}}$,
\begin{equation*}
\frac{1}{\vol(C;G)}\sum_{w \not\in C} \1((u,w) \in E) \geq \alpha \cdot \beta \cdot \deg(u;G)
\end{equation*}
and summing over all such vertices yields
\begin{equation}
\label{pf:zhu6}
\vol(C^{\mathrm{bdry}};G) \leq \frac{1}{\alpha \beta \vol(C;G)}\sum_{\substack{u \in C^{\mathrm{bdry}} \\ w \not\in C}} \1((u,w) \in E) \leq \frac{\Phi(C;G)}{\alpha \cdot \beta}.
\end{equation} 
The claim follows upon summing the upper bounds in~\eqref{pf:zhu3.5}, \eqref{pf:zhu5} and~\eqref{pf:zhu6}. \qed

If the cluster estimate $\wh{C}$ is instead obtained by sweep cutting the aPPR vector $p_v^{(\varepsilon)}$, a similar upper bound on $\vol(\wh{C} \vartriangle C)$ holds, provided that $\varepsilon$ is sufficiently small.
\begin{corollary}
	\label{cor:zhu}
	For a set $C \subseteq V$, suppose that $\alpha,\beta$ satisfy~\eqref{eqn:zhu_condition}, and additionally that 
	\begin{equation}
	\label{eqn:zhu_condition_2}
	\varepsilon \leq \frac{1}{25\vol(C;G)}.
	\end{equation}
	Then there exists a set $C^g \subset C$ with $\vol(C^g;G) \geq \frac{1}{2}\vol(C^g;G)$ such that for any $v \in C^g$, the sweep cut $S_{\beta,v}$ of the aPPR vector $p_v^{(\varepsilon)}$ satisfies
	\begin{equation}
	\label{eqn:zhu_ub_appr}
	\vol(S_{\beta,v} \vartriangle C;G) \leq 6\frac{\Phi(C;G)}{\alpha \beta}.
	\end{equation}
\end{corollary}

\noindent \emph{Proof of Corollary~\ref{cor:zhu}.}
Recall that the upper bound~\eqref{eqn:zhu_ub} on $\vol(\wh{C} \vartriangle C;G)$ comes from combining the upper bounds on $\vol(\wh{C} \setminus C;G)$, $\vol(C^{\mathrm{int}} \setminus \wh{C};G)$ and $\vol(C^{\mathrm{bdry}} \setminus \wh{C};G)$ in~\eqref{pf:zhu3.5},~\eqref{pf:zhu5} and~\eqref{pf:zhu6}. From the upper bound $p_v^{(\varepsilon)}(u) \leq p_v(u)$ for all $u \in V$, it is clear that both~\eqref{pf:zhu3.5} and~\eqref{pf:zhu6} continue to hold when the aPPR vector is used instead of the PPR vector. 

It remains only to establish an upper bound on~$\vol(C^{\mathrm{int}} \setminus \wh{C};G)$. For any $u \in C \setminus S_{\beta,v}$, from inequality~\eqref{pf:zhu2} and the lower bound $p_v^{(\varepsilon)}(u) \geq p_v(u)  - \varepsilon \deg(u;G)$ in~\eqref{eqn:appr_error} we deduce that
\begin{equation}
\frac{3}{4}\bigl(1 - \alpha \cdot \tau_{\infty}(\wt{G})\bigr) \cdot \wt{\pi}(u) -  (\beta + \varepsilon)\deg(u;G) \leq \wt{p}_{\ell}(u).
\end{equation}
Following the same steps as used in the proof of Lemma~\ref{lem:zhu} yields the following inequality:
\begin{equation*}
\biggl(\frac{3(1 - \alpha \beta \vol(C;G))\cdot\bigl(1 - \alpha \tau_{\infty}(\wt{G})\bigr)}{4\vol(C;G)} - \beta\biggr) \cdot \vol\bigl(C^{\mathrm{int}} \setminus S_{\beta,v}; G\bigr) \leq \wt{p}_{\ell}\bigl(C^{\mathrm{int}} \setminus S_{\beta,v}\bigr) \leq 2\frac{\Phi(C;G)}{\alpha}.
\end{equation*}
The upper bounds on $\alpha,\beta$ in~\eqref{eqn:zhu_condition},  and on $\varepsilon$ in~\eqref{eqn:zhu_condition_2}, imply that
\begin{equation*}
\biggl(\frac{3(1 - \alpha \beta \vol(C;G))\cdot\bigl(1 - \alpha \tau_{\infty}(\wt{G})\bigr)}{4\vol(C;G)} - \beta\biggr) \geq \frac{2}{3}\beta,
\end{equation*}
and we conclude that
\begin{equation}
\label{pf:zhu2_1}
\vol(C^{\mathrm{int}} \setminus S_{\beta,v}; G) \leq \frac{3\Phi(C;G)}{\alpha\beta}.
\end{equation}
Summing the right hand sides of~\eqref{pf:zhu2},~\eqref{pf:zhu6}, and~\eqref{pf:zhu2_1} yields the claim. \qed

\subsection{Uniform Bounds on PPR}
\label{subsec:ppr_uniform_bounds}
As mentioned in our main text, in order to prove Theorem~\ref{thm:density_cluster_consistent_recovery}, we require a uniform bound on the PPR vector. Actually, we require two such bounds: for a candidate cluster $C \subseteq V$ and an alternative cluster $C' \subseteq V$, we require a lower bound on $p_v(u)$ for all $u \in C$, and an upper bound on $p_v(u')$ for all $u' \in C'$. In Lemma~\ref{lem:ppr_uniform_bound} we establish an upper bound that holds for all vertices $u$ in the interior $C_{o}$ of $C$, and a lower bound holds for all vertices $u'$ in the interior of $C_{o}'$ of $C'$; here
\begin{equation*}
C_{o} = \Bigl\{u \in C: \deg(u,\wt{G}) =  \deg(u;G)\Bigr\},~~\textrm{and}~~C_{o}'= \Bigl\{u \in C': \deg(u,G[C']) =  \deg(u;G)\Bigr\},
\end{equation*}
and we remind the reader that $\wt{G} = G[C]$. 
\begin{lemma}
	\label{lem:ppr_uniform_bound}
	Let $C$ and $C'$ be disjoint subsets of $V$, and suppose that
	\begin{equation*}
	\alpha \leq \frac{1}{2\tau_{\infty}(\wt{G})}.
	\end{equation*}
	Then there exists a set $C^g \subseteq C$ with $\vol(C^g;G) \geq \vol(C;G)/2$ such that for any $v \in C^g$,
	\begin{equation}
	\label{eqn:ppr_uniform_bound_C}
	p_v(u) \geq \frac{3}{8}\wt{\pi}(u) - \frac{2 \Phi(C;G)}{d_{\min}(\wt{G})\cdot \alpha}~~\textrm{for all $u \in C_{o}$}
	\end{equation}
	and
	\begin{equation}
	\label{eqn:ppr_uniform_bound_Cprime}
	p_v(u') \leq \frac{2\Phi(C;G)}{d_{\min}(C';G) \cdot \alpha}~~\textrm{for all $u \in C_{o}'$.}
	\end{equation}
\end{lemma}

\noindent\emph{``Leakage'' and ``soakage'' vectors.} To prove Lemma~\ref{lem:ppr_uniform_bound}, we will make use of the following explicit representation of the \emph{leakage} distribution $\ell$ from~\eqref{pf:zhu3}, as well as an analogously defined \emph{soakage} distribution $s$:
\begin{equation}
\label{eqn:leakage_soakage}
\begin{aligned}
\ell^{(t)} & := e_v(W \wt{I})^t(I - D^{-1}\wt{D}),~~&& \ell = \sum_{t = 0}^{\infty} (1 - \alpha)^t \ell^{(t)} \\
s^{(t)} & := e_v(W \wt{I})^t W (I - \wt{I}),~~&& s = \sum_{t = 0}^{\infty} (1 - \alpha)^t s^{(t)}.
\end{aligned}
\end{equation}
In the above, $\wt{I} \in \Reals^{n \times n}$ is a diagonal matrix with $I_{uu} = 1$ if $u \in C$ and $0$ otherwise, and $\wt{D}$ is the diagonal matrix with $\wt{D}_{uu} = \deg(u;\wt{G})$ if $u \in C$, and $0$ otherwise. 

These quantities admit a natural interpretation in terms of random walks. For $u \in C$, $\ell^{(t)}(u)$ is the probability that a lazy random walk over $G$ originating at $v$ stays within $C$ for $t$ steps, arriving at $u$ on the $t$th step, and then ``leaks out'' of $C$ on the $(t + 1)$st step. On the other hand, for $u \not\in C$, $s^{(t)}(u)$ is the probability that a lazy random walk over $G$ originating at $v$ stays within $C$ for $t$ steps and is then ``soaked up'' into $u$ on the $(t + 1)$st step. The vectors $\ell$ and $s$ then give the total mass leaked and soaked, respectively, by the PPR vector. 

Three properties of $\ell$ and $s$ are worth pointing out. First, $\mathrm{supp}(\ell) \subseteq C \setminus C_o$ and $\mathrm{supp}(s) \subseteq V \setminus C$. Second, $\|\ell^{(t)}\|_1 = \|s^{(t)}\|_1$ for all $t \in \mathbb{N}$, and so $\|\ell\|_1 = \|s\|_1$. Third, for any $u \in V \setminus C$, $p_v(u) = p_s(u)$. The first two properties are immediate. The third property follows by the law of total probability, which implies that
\begin{equation*}
q_v^{(\tau)}(u) = \sum_{t = 0}^{\tau} q_{s^{(t)}}^{(\tau - t)}(u),~~\textrm{for all $u \in V \setminus C$.}
\end{equation*}
or in terms of the PPR vector,
\begin{equation*}
p_v(u) = \alpha \sum_{\tau = 0}^{\infty} (1 - \alpha)^{\tau} q_v^{(\tau)}(u) = \alpha \sum_{\tau = 0}^{\infty} \sum_{t = 0}^{\tau} (1 - \alpha)^{\tau} q_{s^{(t)}}^{(\tau - t)}(u).
\end{equation*}
Substituting $\Delta = \tau + t$ and rearranging gives the claimed property, as
\begin{equation*}
p_v(u) = \alpha \sum_{\tau = 0}^{\infty} \sum_{t = 0}^{\tau} (1 - \alpha)^{\tau} q_{s^{(t)}}^{(\tau - t)}(u) = \sum_{\Delta = 0}^{\infty} \sum_{t = 0}^{\infty} (1 - \alpha)^{\Delta + t} q_{s^{(t)}}^{(\Delta)}(u) = \alpha \sum_{\Delta = 0}^{\infty} (1 - \alpha)^{\Delta} q_s^{(\Delta)}(u) = p_s(u).
\end{equation*}

\noindent\emph{Proof of Lemma~\ref{lem:ppr_uniform_bound}.}
We first show~\eqref{eqn:ppr_uniform_bound_C}. From~\eqref{pf:zhu3} and~\eqref{pf:zhu2}, we have that
\begin{equation*}
p_v(u) \geq \frac{3}{4}\bigl(1 - \alpha \cdot \tau_{\infty}(\wt{G})\bigr) \cdot \wt{\pi}(u) - \wt{p}_{\ell}(u)~~\textrm{for all $u \in C$,}
\end{equation*}
where $\ell$ has support $\mathrm{supp}(\ell) \subseteq C$ with $\|\ell\|_1 \leq 2\Phi(C;G)/\alpha$. Recalling that $u \in C_{o}$ implies that $u \not\in \mathrm{supp}(\ell)$, as a consequence of~\eqref{pf:interpolator_bound_max_entry}, 
\begin{equation*}
\wt{p}_{\ell}(u) \leq \frac{\|\ell\|_1}{d_{\min}(\wt{G})}~~\textrm{for all $u \in C_{o}$,}
\end{equation*}
establishing~\eqref{eqn:ppr_uniform_bound_C}. The proof of~\eqref{eqn:ppr_uniform_bound_Cprime} follows similarly:
\begin{equation*}
p_v(u) = p_s(u) \overset{(i)}{\leq} \frac{\|s\|_1}{d_{\min}(C';G)} = \frac{\|\ell\|_1}{d_{\min}(C';G)},~~\textrm{for all $u \in C_{o}'$},
\end{equation*}
where the presence of $d_{\min}(C';G)$ on the right hand side of $(i)$ can be verified by inspecting~\eqref{pf:interpolator_bound_max_entry_inductive_step}.  \qed	

\subsection{Mixedness of Lazy Random Walk and PPR Vectors}
\label{subsec:lovasz_simonovits_bounds}
In this subsection, we give upper bounds on $h^{(t)} := h_{q_v^{(t)}}$ and $h^{(\alpha)} := h_{p_v}$. Although similar bounds exist in the literature (see in particular Theorem~1.1 of \citep{lovasz1990} and Theorem~3 of~\citep{andersen2006}), we could not find precisely the results we needed, and so for completeness we state and prove these results ourselves. 

\begin{theorem}
	\label{thm:mixing_time_rw}
	For any $k \in [0, 2m]$, $t_0 \in \mathbb{N}$ and $t \geq t_0$,
	\begin{equation}
	\label{eqn:mixing_time_rw_1}
	h^{(t)}(k) \leq \frac{1}{2^{t_0}} + \frac{d_{\max}(G)}{d_{\min}(G)^2} + \frac{m}{d_{\min}(G)^2} \biggl(1 - \frac{\Psi(G)^2}{8}\biggr)^{t - t_0}.
	\end{equation}
\end{theorem}

\begin{theorem}
	\label{thm:mixing_time_PPR}
	Let $\phi$ be any constant in $[0,1]$. Either the following bound holds for any $t \in \mathbb{N}$ and any $k \in [d_{\max}(G),2m - d_{\max}(G)]$:
	\begin{equation*}
	h^{(\alpha)}(k) \leq \alpha t + \frac{2\alpha}{1 + \alpha} + \frac{d_{\max}(G)}{d_{\min}(G)^2} + \frac{m}{d_{\min}(G)^2} \left(1 - \frac{\phi^2}{8}\right)^{t},
	\end{equation*}
	or there exists some sweep cut $S_j$ of $p_v$ such that $\Phi(S_j;G) < \phi$.
\end{theorem}

The proofs of these upper bounds will be similar to each other (in places word-for-word alike), and will follow a similar approach and use similar notation to that of \citep{lovasz1990,andersen2006}. For $h: [0,2m] \to [0,1]$, $0 \leq K_0 \leq m$ and $k \in [K_0,2m - K_0]$, define
\begin{equation*}
L_{K_0}(k;h) = \frac{2m - K_0 - k}{2m - 2K_0}h(K_0) + \frac{k - K_0}{2m - 2K_0}h(2m - K_0)
\end{equation*}
to be the linear interpolant of $h(K_0)$ and $h(2m - K_0)$, and additionally let
\begin{equation*}
C(K_0;h) := \max\set{\frac{h(k) - L_{K_0}(k;h)}{\sqrt{\wb{k}}}: K_0 \leq k \leq  2m - K_0}.
\end{equation*}
where we use the notation $\wb{k} := \min\{k, 2m - k\}$, and treat $0/0$ as equal to $1$. Our first pair of Lemmas upper bound $h^{(t)}$ and $h^{(\alpha)}$ as a function of $L_{K_0}$ and $C(K_0)$. Lemma~\ref{lem:mixing_random_walk} implies that if $t$ is large relative to $\Psi(G)$, then $h^{(t)}(\cdot)$ must be small.

\begin{lemma}[\textbf{c.f. Theorem 1.2 of~\citep{lovasz1990}}]
	\label{lem:mixing_random_walk}
	For any $K_0 \in [0,m]$, $k \in [K_0, 2m - K_0]$, $t_0 \in \mathbb{N}$ and $t \geq t_0$,
	\begin{equation}
	\label{eqn:mixing_random_walk}
	h^{(t)}(k) \leq L_{K_0}(k;h^{(t_0)}) + C(K_0; h^{(t_0)}) \sqrt{\wb{k}} \cdot \Bigl(1 - \frac{\Psi(G)^2}{8}\Bigr)^{t - t_0}
	\end{equation}
\end{lemma}

Lemma~\ref{lem:mixing_time_PPR} implies that if the PPR random walk is not well mixed, then some sweep cut of $p_v$ must have small normalized cut.

\begin{lemma}[\textbf{c.f Theorem~3 of \citep{andersen2006}}]
	\label{lem:mixing_time_PPR}
	Let $\phi \in [0,1]$. Either the following bound holds for any $t \in \mathbb{N}$, any $K_0 \in [0,m]$, and any $k \in [K_0,2m - K_0]$:
	\begin{equation}
	\label{eqn:mixing_time_PPR}
	h^{(\alpha)}(k) \leq \alpha t + L_{K_0}(k; h^{(\alpha)}) + C(K_0;h^{(\alpha)})\sqrt{\wb{k}}\left(1 - \frac{\phi^2}{8}\right)^t
	\end{equation}
	or else there exists some sweep cut $S_{j}$ of $p_v$ such that $\Phi(S_j;G) < \phi$.
\end{lemma}

In order to make use of these Lemmas, we require upper bounds on $L_{K_0}(\cdot,h)$ and $C(K_0;h)$, for each of $h = h^{(t_0)}$ and $h = h^{(\alpha)}$. Of course, trivially $L_{K_0}(k;h) \leq \max\{h(K_0); h(2m - K_0)\}$ for any $k \in [K_0, 2m - K_0]$. As it happens, this observation will lead to sufficient upper bounds on $L_{K_0}(k,h)$ for both $h = h^{(t_0)}$ (Lemma~\ref{lem:interpolator_bound_rw}) and $h = h^{(\alpha)}$ (Lemma~\ref{lem:interpolator_bound_ppr}).  
\begin{lemma}
	\label{lem:interpolator_bound_rw}
	For any $t_0 \in \mathbb{N}$ and $K_0 \in [0,m]$, the following inequalities hold:
	\begin{equation}
	\label{eqn:interpolator_bound_rw}
	h^{(t_0)}\bigl(2m - K_0\bigr) \leq \frac{K_0}{2m}~~\textrm{and}~~h^{(t_0)}\bigl(K_0\bigr) \leq \frac{K_0}{d_{\min}(G)^2} + \frac{1}{2^{t_0}}.
	\end{equation}
	As a result, for any $k \in [K_0, 2m - K_0]$,
	\begin{equation}
	\label{eqn:interpolator_bound_rw_2}
	L_{K_0}(k;h^{(t_0)}) \leq \max\Bigl\{\frac{K_0}{2m}, \frac{K_0}{d_{\min}(G)^2} + \frac{1}{2^{t_0}}\Bigr\} = \frac{K_0}{d_{\min}(G)^2} + \frac{1}{2^{t_0}}.
	\end{equation}
\end{lemma}

\begin{lemma}
	\label{lem:interpolator_bound_ppr}
	For any $\alpha \in [0,1]$ and $K_0 \in [0,m]$, the following inequalities hold:
	\begin{equation}
	\label{eqn:interpolator_bound_ppr}
	h^{(\alpha)}\bigl(2m - K_0\bigr) \leq \frac{K_0}{2m}~~\textrm{and}~~h^{(\alpha)}\bigl(K_0\bigr) \leq \frac{K_0}{d_{\min}(G)^2} + \frac{2\alpha}{1 + \alpha}.
	\end{equation}
	As a result, for any $k \in [K_0, 2m - K_0]$,
	\begin{equation}
	\label{eqn:interpolator_bound_ppr_2}
	L_{K_0}(k;h^{(\alpha)}) \leq \max\Bigl\{\frac{K_0}{2m}, \frac{K_0}{d_{\min}(G)^2} + \frac{2\alpha}{1 + \alpha}\Bigr\} = \frac{K_0}{d_{\min}(G)^2} + \frac{2\alpha}{1 + \alpha}.
	\end{equation}
\end{lemma}

We next establish an upper bound on $C_{K_0}(k;h)$, which rests on the following key observation: since $h(k)$ is concave and $L_{K_0}(K_0;h) = h(K_0)$, it holds that
\begin{equation}
\label{pf:linearization_bound_1}
\frac{h(k) - L_{K_0}(k)}{\sqrt{\wb{k}}} \leq
\begin{cases}
h'(K_0) \sqrt{k},~& k \leq m \\
-h'(2m - K_0) \sqrt{2m - k},~& k > m.
\end{cases}
\end{equation}
(Since $h$ is not differentiable at $k = k_j$, here $h'$ refers to the right derivative of $h$.)  

Lemma~\ref{lem:linearization_bound} gives good estimates for $h'(K_0)$ and $h'(2m - K_0)$, which hold for both $h = h^{(t_0)}$ and $h = h^{(\alpha)}$, and result in an upper bound on $C(K_0;h)$. Both the statement and proof of this Lemma rely on the following explicit representation of the Lovasz-Simonovits curve $h_q(\cdot)$. Order the vertices $q(u_{(1)})/\deg(u_{(1)};G) \geq q(u_{(2)})/\deg(u_{(2)};G) \geq \cdots \geq q(u_{(n)})/\deg(u_{(n)};G)$. Then for each $j = 0,\ldots,n - 1$, and for all $k \in [\vol(S_j),\vol(S_{j + 1}))$,  the function $h_q(k)$ satisfies
\begin{equation}
\label{eqn:lovasz_simonovits}
h_q(k) = \sum_{i = 0}^{j} \left(q(u_{(i)}) - \pi(u_{(i)})\right) + \frac{\bigl(k - \vol(S_j;G)\bigr)}{\deg(u_{(j + 1)};G)} \left(q(u_{(j+1)}) - \pi(u_{(j+1)})\right). 
\end{equation}
\begin{lemma}
	\label{lem:linearization_bound}
	The following statements hold for both $h = h^{(\alpha)}$ and $h = h^{(t_0)}$. 
	\begin{itemize}
		\item Let $K_0 = k_1 = \deg(v;G)$ if $u_{(1)} = v$, and otherwise $K = 0$. Then  
		\begin{equation}
		\label{eqn:right_derivative_1}
		h'\bigl(K_0\bigr) \leq \frac{1}{d_{\min}(G)^2}.
		\end{equation}
		\item For all $K_0 \in [0,m]$,
		\begin{equation}
		\label{eqn:right_derivative_2}
		h'(2m - K_0) \geq -\frac{d_{\max}(G)}{d_{\min}(G)\cdot \vol(G)}.
		\end{equation}
	\end{itemize}
	As a result, letting $K_0 = \deg(v;G)$ if $u_{(1)} = v$, and otherwise letting $K_0 = 0$, we have
	\begin{equation*}
	C(K_0,h) \leq \frac{\sqrt{m}}{d_{\min}(G)^2}.
	\end{equation*}
\end{lemma}

\subsubsection{Proof of Theorems~\ref{thm:mixing_time_rw} and~\ref{thm:mixing_time_PPR}}

\noindent \emph{Proof of Theorem~\ref{thm:mixing_time_rw}.}
Take $K_0 = 0$ if $u_{(1)} \neq v$, and otherwise take $K_0 = \deg(v;G)$. Combining Lemmas~\ref{lem:mixing_random_walk},~\ref{lem:interpolator_bound_rw} and~\ref{lem:linearization_bound}, we obtain that for any $k \in [K_0,2m - K_0]$, 
\begin{align*}
h^{(t)}(k) & \leq \frac{1}{2^{t_0}} + \frac{K_0}{d_{\min}(G)^2}  + \frac{\sqrt{m}}{d_{\min}(G)^2} \sqrt{\wb{k}} \Bigl(1 - \frac{\Psi^2(G)}{8}\Bigr)^{t-t_0} \\
& \leq \frac{1}{2^{t_0}} + \frac{d_{\max}(G)}{d_{\min}(G)^2}  + \frac{m}{d_{\min}(G)^2} \Bigl(1 - \frac{\Psi^2(G)}{8}\Bigr)^{t-t_0},
\end{align*}
where the second inequality follows since we have chosen $K_0 \leq d_{\max}(G)$, and since $\wb{k} \leq m$. If $K_0 = 0$, we are done. 

Otherwise, we suppose $k \in [0, \deg(v;G)) ~\cup~ (2m - \deg(v;G),2m]$. If $k \in [0, \deg(v;G))$ then
\begin{equation}
\label{pf:mixing_random_walk_1}
h^{(t)}(k) \overset{\eqref{pf:mixing_random_walk_inductive_step}}{\leq} h^{(t_0)}(k) \overset{(\textrm{i})}{\leq} h^{(t_0)}(K_0) \overset{\eqref{eqn:interpolator_bound_rw}}{\leq} \frac{K_0}{d_{\min}(G)^2} + \frac{1}{2^{t_0}},
\end{equation}
where $(\textrm{i})$ follows since $k \in [0,K_0]$, and $h^{(t_0)}$ is linear over $[0,K_0)$ with $h^{(t_0)}(0) = 0$ and $h^{(t_0)}(K_0) \geq 0$. For similar reasons, if $k \in (2m - \deg(v;G),2m]$ then
\begin{equation}
\label{pf:mixing_random_walk_2}
h^{(t)}(k) \leq h^{(t_0)}(k) \leq h^{(t_0)}(2m - K_0) \leq \frac{\deg(v;G)}{2m}.
\end{equation}
Since the ultimate upper bounds in~\eqref{pf:mixing_random_walk_1} and~\eqref{pf:mixing_random_walk_2} are each no greater than that of~\eqref{eqn:mixing_random_walk}, the claim follows. \qed \\

\noindent \emph{Proof of Theorem~\ref{thm:mixing_time_PPR}.}
The proof of Theorem~\ref{thm:mixing_time_PPR} follows immediately from Lemmas~\ref{lem:mixing_time_PPR},~\ref{lem:interpolator_bound_ppr} and~\ref{lem:linearization_bound}, taking $K_0 = 0$ if $u_{(1)} \neq v$ and otherwise $K_0 = \deg(v;G)$. \qed

\subsubsection{Proofs of Lemmas}
In what follows, for a distribution $q$ and vertices $u,w \in V$, we write $q(u,w) := q(u)/d(u) \cdot 1\{(u,w) \in E\}$, and similarly for a collection of dyads $\wt{E} \subseteq V \times V$ we write $q(\wt{E}) := \sum_{(u,w) \in \wt{E}} q(u,w)$. \\

\noindent\emph{Proof of Lemma~\ref{lem:mixing_random_walk}.}
We will prove Lemma~\ref{lem:mixing_random_walk} by induction on $t$. In the base case $t = t_0$, observe that $C(K_0;h^{(t_0)}) \cdot \sqrt{\wb{k}} \geq h^{(t_0)}(k) - L_{K_0}(k;h^{(t_0)})$ for all $k \in [K_0, 2m - K_0]$, which implies 
\begin{equation*}
L_{K_0}(k;h^{(t_0)}) + C(K_0; h^{(t_0)}) \cdot \sqrt{\wb{k}} \geq h^{(t_0)}(k).
\end{equation*}

Now, we proceed with the inductive step, assuming that the inequality holds for $t_0,t_0 + 1,\ldots,t - 1$, and proving that it thus also holds for $t$. By the definition of $L_{K_0}$, the inequality~\eqref{eqn:mixing_random_walk} holds when $k = K_0$ or $k = 2m - K_0$. We will additionally show that~\eqref{eqn:mixing_random_walk} holds for every $k_j = \vol(S_j), j = 1,2,\ldots,n$ such that $k_j \in [K_0, 2m - K_0]$. This suffices to show that the inequality~\eqref{eqn:mixing_random_walk} holds for all $k \in [K_0,2m - K_0]$, since the right hand side of~\eqref{eqn:mixing_random_walk} is a concave function of $k$.

Now, we claim that for each $k_j$, it holds that
\begin{equation}
\label{pf:mixing_random_walk_inductive_step}
q_v^{(t)}[k_j] \leq \frac{1}{2}\Bigl(q_v^{(t - 1)}[k_j - \wb{k}_j \Psi(G)] + q_v^{(t - 1)}[k_j + \wb{k}_j \Psi(G)]\Bigr).
\end{equation}
To establish this claim, we note that for any $u \in V$
\begin{equation*}
q_v^{(t)}(u) = \frac{1}{2}q_v^{(t - 1)}(u) + \frac{1}{2}\sum_{w \in V}q_v^{(t - 1)}(w,u) = \frac{1}{2} \sum_{w \in V} \bigl(q_v^{(t - 1)}(u,w) + q_v^{(t - 1)}(w,u)\bigr),
\end{equation*}
and consequently for any $S \subset V$,
\begin{align*}
q_v^{(t)}(S) & = \frac{1}{2}\bigl\{q_v^{(t - 1)}(\mathrm{in}(S)) +  q_v^{(t - 1)}(\mathrm{out}(S))\bigr\} \\
& = \frac{1}{2}\bigl\{q_v^{(t - 1)}\bigl(\mathrm{in}(S) \cup \mathrm{out}(S)\bigr) +  q_v^{(t - 1)}\bigl(\mathrm{in}(S) \cap \mathrm{out}(S)\bigr)\bigr\}
\end{align*}
where $\mathrm{in}(S) = \{(u,w) \in E: u \in S\}$ and $\mathrm{out}(S) = \{(w,u) \in E: w \in S\}$. We deduce that
\begin{align*}
q_v^{(t)}[k_j] = q_v^{(t)}(S_j) & = \frac{1}{2}\bigl\{q_v^{(t - 1)}\bigl(\mathrm{in}(S_j) \cup \mathrm{out}(S_j)\bigr) +  q_v^{(t - 1)}\bigl(\mathrm{in}(S_j) \cap \mathrm{out}(S_j)\bigr)\bigr\} \\
& \leq \frac{1}{2}\bigl\{q_v^{(t - 1)}\bigl[|\mathrm{in}(S_j) \cup \mathrm{out}(S_j)|\bigr] +  q_v^{(t - 1)}\bigl[|\mathrm{in}(S_j) \cap \mathrm{out}(S_j)|\bigr]\bigr\} \\
& = \frac{1}{2}\bigl\{q_v^{(t - 1)}\bigl[k_j + \mathrm{cut}(S_j;G)\bigr] +  q_v^{(t - 1)}\bigl[k_j - \mathrm{cut}(S_j;G)\bigr]\bigr\} \\
& \leq \frac{1}{2}\bigl\{q_v^{(t - 1)}\bigl[k_j + \wb{k}_j\Phi(S_j;G)\bigr] +  q_v^{(t - 1)}\bigl[k_j - \wb{k}_j\Phi(S_j;G)\bigr]\bigr\} \\
& \leq \frac{1}{2}\bigl\{q_v^{(t - 1)}\bigl[k_j + \wb{k}_j\Psi(G)\bigr] +  q_v^{(t - 1)}\bigl[k_j - \wb{k}_j\Psi(G)\bigr]\bigr\},
\end{align*}
establishing~\eqref{pf:mixing_random_walk_inductive_step}. The final two inequalities both follow from the concavity of $q_v^{(t)}[\cdot]$. 

Subtracting $k_j/2m$ from both sides, we get
\begin{equation}
\label{pf:mixing_random_walk_inductive_step_h}
h^{(t)}(k_j) \leq \frac{1}{2}\bigl\{h^{(t - 1)}\bigl(k_j + \wb{k}_j\Psi(G)\bigr)+  h^{(t - 1)}\bigl(k_j - \wb{k}_j\Psi(G)\bigr)\bigr\}.
\end{equation}
At this point, we divide our analysis into cases.

\textbf{Case 1.}
Assume $k_j - \Psi(G) \wb{k}_j$ and $k_j + 2 \Psi(G) \wb{k}_j$ are both in $[K_0,2m  - K_0]$. We are therefore in a position to apply our inductive hypothesis to both terms on the right hand side of~\eqref{pf:mixing_random_walk_inductive_step_h}, and obtain the following:
\begin{align*}
h^{(t)}(k_j) & \leq \frac{1}{2}\biggl(L_{K_0}\bigl(k_j - \Psi(G) \wb{k}_j; h^{(t_0)}\bigr) + L_{K_0}\bigl(k_j + \Psi(G) \wb{k}_j; h^{(t_0)}\bigr) \biggr) ~~+ \\
& \quad~ \frac{1}{2}C\bigl(K_0; h^{(t_0)}\bigr) \cdot \Bigl(\sqrt{\overline{k_j - \Psi(G) \wb{k}_j}} + \sqrt{\overline{k_j + \Psi(G) \wb{k}_j}}\Bigr)\left(1 - \frac{\Psi(G)^2}{8}\right)^{t-t_0 - 1} \\
& = L_{K_0}(k;h^{(t_0)}) ~~+ \\
& \quad~\frac{1}{2}C(K_0;h^{t_0})\biggl(\sqrt{\overline{k_j - \Psi(G) \wb{k}_j}} + \sqrt{\overline{k_j + \Psi(G) \wb{k}_j}}\biggr)\left(1 - \frac{\Psi(G)^2}{8}\right)^{t-t_0 - 1} \\
& \leq L_{K_0}(k;h^{(t_0)}) ~~+ \\
& \quad~\frac{1}{2}C(K_0;h^{(t_0)})\biggl(\sqrt{\wb{k}_j - \Psi(G) \wb{k}_j} + \sqrt{\wb{k}_j + \Psi(G) \wb{k}_j}\biggr)\left(1 - \frac{\Psi(G)^2}{8}\right)^{t-t_0 - 1}.
\end{align*}
A Taylor expansion of $\sqrt{1 + \Psi(G)}$ around $\Psi(G) = 0$ yields the following bound:
\begin{equation*}
\sqrt{1 + \Psi(G)} + \sqrt{1 - \Psi(G)} \leq 2 - \frac{\Psi(G)^2}{4},
\end{equation*}
and therefore
\begin{align*}
h^{(t)}(k_j) & \leq L_{K_0}(k;h^{(t_0)}) + \frac{C(K_0;h^{(t_0)})}{2}\cdot \sqrt{\wb{k}_j}\cdot\left(2 - \frac{\Psi(G)^2}{4}\right)\left(1 - \frac{\Psi(G)^2}{8}\right)^{t-1} \\
&= L_{K_0}(k_j;h^{(t_0)}) + C(K_0;h^{(t_0)})\sqrt{\wb{k}_j}\left(1 - \frac{\Psi(G)^2}{8}\right)^{t - t_0}.
\end{align*}

\textbf{Case 2.} Otherwise one of $k_j - 2 \Psi(G) \wb{k}_j$ or $k_j + 2 \Psi(G) \wb{k}_j$ is not in $[K_0,2m  - K_0]$. Without loss of generality assume $k_j < m$, so that (i) we have $k_j - 2 \Psi(G) \wb{k}_j < K_0$ and (ii) $k_j + (k_j - K_0) \leq 2m - K_0$. We deduce the following:
\begin{align*}
h^{(t)}(k_j) & \overset{(\textrm{i})}{\leq} \frac{1}{2}\Bigl(h^{(t - 1)}(K_0) + h^{(t - 1)}\bigl(k_j + (k_j - K_0)\bigr)\Bigr) \\
& \overset{(\textrm{ii})}{\leq} \frac{1}{2}\Bigl(h^{(t_0)}(K_0) + h^{(t)}\bigl(k_j + (k_j - K_0)\bigr)\Bigr) \\
& \overset{(\textrm{iii})}{\leq}\frac{1}{2}\Bigl(L_{K_0}(K_0;h^{(t_0)}) + L_{K_0}(2k_j - K_0; h^{(t_0)}\bigr) ~~+ \\
& \quad~ C(K_0;h^{(t_0)})\sqrt{\overline{2k_j - K_0}}\left(1 - \frac{\Psi(G)^2}{8}\right)^{t - t_0 - 1}\Bigr) \\
& \leq L_{K_0}(k_j; h^{(t_0)}) + C(K_0;h^{(t_0)}) \frac{\sqrt{2\wb{k}_j}}{2} \left(1 - \frac{\Psi(G)^2}{8}\right)^{t - t_0 - 1} \\
& \leq L_{K_0}(k_j;h^{(t_0)}) + C(K_0;h^{(t_0)}) \sqrt{\wb{k}_j} \cdot \left(1 - \frac{\Psi(G)^2}{8}\right)^{t - t_0}
\end{align*}
where $\textrm{(i)}$ follows from~\eqref{pf:mixing_random_walk_inductive_step_h} and the concavity of $h^{(t - 1)}$,  we deduce $(\textrm{ii})$ from~\eqref{pf:mixing_random_walk_inductive_step_h}, which implies that $h^{(t)}(k) \leq h^{(t_0)}(k)$, and $(\textrm{iii})$ follows from applying the inductive hypothesis to $h^{(t - 1)}(2k_j - K_0)$. \qed \\

\noindent\emph{Proof (of Lemma~\ref{lem:mixing_time_PPR}).}
We will show that if $\Phi(S_j; g) \geq \phi$ for each $j = 1,\ldots,n$, then \eqref{eqn:mixing_time_PPR} holds for all $t$ and any $k \in [K_0,2m - K_0]$.

We proceed by induction on $t$. Our base case will be $t = 0$. Observe that $C(K_0;h^{(\alpha)}) \cdot \sqrt{\wb{k}} \geq h^{(\alpha)}(k) - L_{K_0}(k;h^{(\alpha)})$ for all $k \in [K_0,2m - K_0]$, which implies
\begin{equation*}
L_{K_0}(k;h^{(\alpha)}) + C(K_0;h^{(\alpha)}) \cdot \sqrt{\wb{k}} \geq h^{(\alpha)}(k).
\end{equation*}

Now, we proceed with the inductive step. By the definition of $L_{K_0}$, the inequality~\eqref{eqn:mixing_time_PPR} holds when $k = K_0$ or $k = 2m - K_0$. We will additionally show that~\eqref{eqn:mixing_time_PPR} holds for every $k_j = \vol(S_j), j = 1,2,\ldots,n$ such that $k_j \in [K_0, 2m - K_0]$. This suffices to show that the inequality~\eqref{eqn:mixing_time_PPR} holds for all $k \in [K_0,2m - K_0]$, since the right hand side of~\eqref{eqn:mixing_time_PPR} is a concave function of $k$.

By Lemma 5 of \citet{andersen2006}, we have that
\begin{align}
p_v[k_j] & \leq \alpha + \frac{1}{2} \bigl(p_v[k_j - \mathrm{cut}(S_j;G)] + p_v[k_j + \mathrm{cut}(S_j;G)] \bigr) \nonumber\\
& \leq \alpha + \frac{1}{2} \bigl(p_v[k_j - \Phi(S_j;G) \wb{k}_j] + p_v[k_j + \Phi(S_j;G) \wb{k}_j]  \bigr) \nonumber \\
& \leq \alpha + \frac{1}{2} \bigl(p_v[k_j - \phi \wb{k}_j] + p_v[k_j + \phi \wb{k}_j]\bigr) \nonumber
\end{align}
and subtracting $k_j/2m$ from both sides, we get
\begin{equation}
\label{eqn:mixing_time_PPR_pf1}
h^{(\alpha)}(k_j) \leq \alpha + \frac{1}{2} \bigl(h^{(\alpha)}(k_j - \phi \wb{k}_j) + h^{(\alpha)}(k_j +  \phi \wb{k}_j) \bigr)
\end{equation}
From this point, we divide our analysis into cases. 

\textbf{Case 1.}
Assume $k_j - 2 \phi \wb{k}_j$ and $k_j + 2 \phi \wb{k}_j$ are both in $[K_0,2m  - K_0]$. We are therefore in a position to apply our inductive hypothesis to \eqref{eqn:mixing_time_PPR_pf1}, yielding
\begin{align*}
h^{(\alpha)}(k_j) & \leq \alpha + \alpha(t-1) \frac{1}{2}\biggl(L_{K_0}(k_j - \phi \wb{k}_j) + L_{K_0}(k_j + \phi \wb{k}_j)\biggr) ~~+ \\
& \quad~ \frac{1}{2}C(K_0;h^{(\alpha)})\bigl(\sqrt{\overline{k_j - \phi \wb{k}_j}} + \sqrt{\overline{k_j + \phi \wb{k}_j}}\bigr)\left(1 - \frac{\phi^2}{8}\right)^{t-1} \\
& \leq \alpha t + L_{K_0}(k;h^{(\alpha)}) + \frac{1}{2}\biggl(C(K_0;h^{(\alpha)})\bigl(\sqrt{\overline{k_j - \phi \wb{k}_j}} + \sqrt{\overline{k_j + \phi \wb{k}_j}}\bigr)\left(1 - \frac{\phi^2}{8}\right)^{t-1} \biggr) \\
& \leq \alpha t + L_{K_0}(k;h^{(\alpha)}) + \frac{1}{2}\biggl(C(K_0;h^{(\alpha)})\bigl(\sqrt{\wb{k}_j - \phi \wb{k}_j} + \sqrt{\wb{k}_j + \phi \wb{k}_j}\bigr)\left(1 - \frac{\phi^2}{8}\right)^{t-1} \biggr).
\end{align*}
and therefore
\begin{align*}
h^{(\alpha)}(k_j) & \leq  \alpha t + L_{K_0}(k;h^{(\alpha)}) + \frac{C(K_0;h^{(\alpha)})}{2}\cdot \sqrt{\wb{k}_j}\cdot\left(2 - \frac{\phi^2}{4}\right)\left(1 - \frac{\phi^2}{8}\right)^{t-1} \\
& = \alpha t + L_{K_0}(k;h^{(\alpha)}) + C(K_0;h^{(\alpha)})\sqrt{\wb{k}_j}\left(1 - \frac{\phi^2}{8}\right)^{t}.
\end{align*}

\textbf{Case 2.} Otherwise one of $k_j - 2 \phi \wb{k}_j$ or $k_j + 2 \phi \wb{k}_j$ is not in $[K_0,2m  - K_0]$. Without loss of generality assume $k_j < m$, so that (i) we have $k_j - 2 \phi \wb{k}_j < K_0$ and (ii) $k_j + (k_j - K_0) \leq 2m - K_0$. By the concavity of $h$, and applying the inductive hypothesis to $h^{(\alpha)}2k_j - K_0)$, we have
\begin{align*}
h^{(\alpha)}(k_j) & \leq \alpha + \frac{1}{2}\Bigl(h^{(\alpha)}(K_0) + h\bigl(k_j + (k_j - K_0)\bigr)\Bigr) \\
& \leq\alpha + \frac{\alpha(t - 1)}{2} + \frac{1}{2}\Bigl(L_{K_0}(K_0;p^{\alpha}) + L_{K_0}(2k_j - K_0\bigr)\Bigr) ~~+ \\
& \quad~ C(K_0;h^{(\alpha)})\sqrt{\wb{2k_j - K_0}}\left(1 - \frac{\phi^2}{8}\right)^{t - 1}\Bigr) \\
& \leq \alpha t + L_{K_0}(k_j) + C(K_0;h^{(\alpha)}) \frac{\sqrt{2\wb{k}_j}}{2} \left(1 - \frac{\phi^2}{8}\right)^{t - 1} \\
& \leq \alpha t + L_{K_0}(k_j) + C(K_0;h^{(\alpha)}) \sqrt{\wb{k}_j} \cdot \left(1 - \frac{\phi^2}{8}\right)^{t}.
\end{align*} \qed

\noindent \emph{Proof of Lemma~\ref{lem:interpolator_bound_rw}.}
We will prove that the inequalities of~\eqref{eqn:interpolator_bound_rw} hold at the knot points of $h^{(t_0)}$, whence they follow for all $K_0 \in [0,m]$. 

We first prove the upper bound on $h^{(t_0)}(2m - K_0)$, when $2m - K_0 = k_j$ for some $j = 0,\ldots,n - 1$. Indeed, the following manipulations show the upper bound holds for $h_q(\cdot)$ regardless of the distribution $q$. Noting that $h_q(2m) = 0$, we have that,
\begin{equation*}
h_q(k_j) = h_q(k_j) - h_q(2m) = \sum_{i = j + 1}^{n} q(u_{(i)}) - \pi(u_{(i)}) \leq \sum_{i = j + 1}^{n} \pi(u_{(i)}) = 1 - \frac{k_j}{2m} = \frac{K_0}{2m}.
\end{equation*}

In contrast, when $K_0 = k_j$ the upper bound on $h^{(t_0)}(\cdot)$ depends on the properties of $q = q_v^{(t_0)}$. In particular, we claim that for any $t \in \mathbb{N}$,
\begin{equation}
\label{pf:interpolator_bound_max_entry}
q_v^{(t)}(u) \leq
\begin{dcases*}
\frac{1}{d_{\min}(G)},& ~~\textrm{if $u \neq v$} \\
\frac{1}{d_{\min}(G)} + \frac{1}{2^t},& ~~\textrm{if $u = v$.}
\end{dcases*}
\end{equation}
This claim follows straightforwardly by induction. In the base case $t = 0$, the claim is obvious. If the claim holds true for a given $t \in \mathbb{N}$, then for $u \neq v$,
\begin{equation}
\label{pf:interpolator_bound_max_entry_inductive_step}
\begin{aligned}
q_v^{(t + 1)}(u) & = \frac{1}{2}\sum_{w \neq u}q_v^{(t)}(w,u)  + \frac{1}{2}q_v^{(t)}(u) \\
& \leq \frac{1}{2d_{\min}(G)}\sum_{w \neq u}q_v^{(t)}(w)  + \frac{1}{2d_{\min}(G)} \\
& \leq \frac{1}{d_{\min}(G)},
\end{aligned}
\end{equation}
where the last inequality holds because $q_v^{(t)}$ is a probability distribution (i.e. the sum of its entries is equal to $1$). Similarly, if $u = v$, then
\begin{align*}
q_v^{(t + 1)}(v) & = \frac{1}{2}\sum_{w \neq v}q_v^{(t)}(w,v)  + \frac{1}{2}q_v^{(t)}(v) \\
& \leq \frac{1}{2d_{\min}(G)}\sum_{w \neq u}q_v^{(t)}(w)  + \frac{1}{2d_{\min(G)}} + \frac{1}{2^{t + 1}} \\
& \leq \frac{1}{d_{\min}(G)} + \frac{1}{2^{t + 1}},
\end{align*}
and the claim~\eqref{pf:interpolator_bound_max_entry} is shown. The upper bound on $h^{(t_0)}(K_0)$ for $K_0 = k_j$ follows straightforwardly:
\begin{equation*}
h^{(t_0)}(K_0) \leq \sum_{i = 0}^{j} q_v^{(t_0)}(u_{(j)}) \leq \frac{j}{d_{\min}(G)} + \frac{1}{2^{t_0}} \leq \frac{K_0}{d_{\min}(G)^2} + \frac{1}{2^{t_0}},
\end{equation*}
where the last inequality follows since $\vol(S) \geq |S| \cdot d_{\min}(G)$ for any set $S \subseteq V$. \qed \\

\noindent \emph{Proof of Lemma~\ref{lem:interpolator_bound_ppr}.}
We have already established the first upper bound in~\eqref{eqn:interpolator_bound_ppr}, in the proof of Lemma~\ref{lem:interpolator_bound_rw}. Then, noting that from~\eqref{pf:interpolator_bound_max_entry}, 
\begin{equation}
\label{eqn:ppr_max_entry}
p_v(u) = \alpha \sum_{t = 0}^{\infty} (1 - \alpha)^t q_v^{(t)}(u) \leq
\begin{dcases*}
\alpha \sum_{t = 0}^{\infty} (1 - \alpha)^t \Bigl(\frac{1}{d_{\min}(G)} + \frac{1}{2^t}\Bigr) = \frac{1}{d_{\min}(G)} + \frac{2\alpha}{1 - \alpha}& ~~\textrm{if $u = v$} \\
\alpha \sum_{t = 0}^{\infty} (1 - \alpha)^t \frac{1}{d_{\min}(G)} = \frac{1}{d_{\min}(G)} & ~~\textrm{if $u \neq v$,}
\end{dcases*}
\end{equation}
the second upper bound in~\eqref{eqn:interpolator_bound_ppr} follows similarly to the proof of the equivalent upper bound in~Lemma~\ref{lem:interpolator_bound_rw}. \qed \\

\noindent \emph{Proof of Lemma~\ref{lem:linearization_bound}.}
The result of the Lemma follows obviously from~\eqref{pf:linearization_bound_1}, once we show \eqref{eqn:right_derivative_1}-\eqref{eqn:right_derivative_2}. We begin by showing~\eqref{eqn:right_derivative_1}. Inspecting the representation~\eqref{eqn:lovasz_simonovits}, we see that for any distribution $q$ and knot point $k_j$, the right derivative of $h_q$ can always be upper bounded,
\begin{equation*}
h_{q}'(k_j) \leq \frac{q(u_{(j + 1)})}{\deg(u_{(j + 1)};G)}.
\end{equation*}
We have chosen $K_0 = k_j$ so that $v \neq u_{(j + 1)}$, and so~\eqref{pf:interpolator_bound_max_entry} implies that $h_{q}'(k_j) \leq 1/(d_{\min}(G)^2)$, for either $q = q_v^{(t)}$ or $q = p_v$.

On the other hand, the inequality \eqref{eqn:right_derivative_2} follows immediately from the representation \eqref{eqn:lovasz_simonovits}, since for any $K_0 \in [0,m]$, taking $j$ so that $2m - K_0 \in [k_j, k_{j + 1})$, 
\begin{equation*}
h'(2m - K_0) \geq -\frac{\pi(u_{(j+1)})}{\mathrm{deg}(u_{(j + 1)};G)} \geq -\frac{d_{\max}(G)}{d_{\min}(G) \cdot \vol(G)}.
\end{equation*} \qed

\subsection{Proof of Proposition~\ref{prop:pointwise_mixing_time}}
\label{subsec:pf_prop_pointwise_mixing_time}
To prove Proposition~\ref{prop:pointwise_mixing_time}, we first give an upper bound on the total variation distance between $q_v^{(t)}$ and its limiting distribution $\pi$, then upgrade to the desired uniform upper bound~\eqref{eqn:pointwise_mixing_time}. The \emph{total variation distance} between distributions $q$ and $p$ is
\begin{equation*}
\mathrm{TV}(q,p) := \frac{1}{2}\sum_{u \in v} \bigl|q(u) - p(u)\bigr|
\end{equation*}
It follows from the representation~\eqref{eqn:lovasz_simonovits} that
\begin{equation*}
\mathrm{TV}(q,\pi) = \max_{S \subseteq V} \Bigl\{q(S) - \pi(S)\Bigr\} = \max_{j = 1,\ldots,n} \Bigl\{q(S_j) - \pi(S_j)\Bigr\} = \max_{k \in [0,2m]} h_q(k),
\end{equation*}
so that Theorem~\ref{thm:mixing_time_rw} gives an upper bound on $\mathrm{TV}(q_v^{(t)},\pi)$. We can then use the following result to upgrade to a uniform upper bound.

\begin{lemma}
	\label{lem:tv_to_pointwise}
	For any $t_{\ast} \in \mathbb{N}$,
	\begin{equation*}
	\max_{u \in V} \Bigl\{\frac{\pi(u) - q_v^{(t_{\ast} + 1)}(u)}{\pi(u)}\Bigr\} \leq \frac{1}{s(G)} \biggl(\sum_{t = 0}^{t_{\ast}} \frac{\mathrm{TV}(q_v^{(t)}, \pi)}{2^{t_{\ast} - t}} + \frac{\mathrm{TV}(q_v^{(0)}, \pi)}{2^{t_{\ast}}}\biggr).
	\end{equation*}
\end{lemma}
The proof of Proposition~\ref{prop:pointwise_mixing_time} is then straightforward. \\

\noindent\emph{Proof of Proposition~\ref{prop:pointwise_mixing_time}.}
Put $T = 8/(\Psi(G)^2) \ln(4/s(G)) + 4$. We will use Theorem~\ref{thm:mixing_time_rw} to show that $\mathrm{TV}(q_v^{(T)},\pi) \leq 1/4$. This will in turn imply (\cite{montenegro2002} pg. 13) that for all $t \geq t_{\ast} := T \log_2(32/s(G))$,
\begin{equation}
\label{pf:pointwise_mixing_time_1}
\mathrm{TV}(q_v^{(t)},\pi) \leq \frac{1}{32}s(G).
\end{equation}
Finally, let $\tau_{\ast} = t_{\ast} + 4\log_2(1/s(G))$. Applying Lemma~\ref{lem:tv_to_pointwise} gives
\begin{align*}
\max_{u \in V} \Bigl\{\frac{\pi(u) - q_v^{(\tau_{\ast} + 1)}(u)}{\pi(u)}\Bigr\} & \leq \frac{1}{s(G)}\biggl(\sum_{t = 0}^{\tau_{\ast}} \frac{\mathrm{TV}(q_v^{(t)}, \pi)}{2^{\tau_{\ast} - t}} + \frac{\mathrm{TV}(q_v^{(0)}, \pi)}{2^{\tau_{\ast}}}\biggr) \\
& = \frac{1}{s(G)}\biggl(\sum_{t = t_{\ast} + 1}^{\tau_{\ast}} \frac{\mathrm{TV}(q_v^{(t)}, \pi)}{2^{\tau_{\ast} - t}} + \sum_{t = 0}^{t_{\ast}} \frac{\mathrm{TV}(q_v^{(t)}, \pi)}{2^{\tau_{\ast} - t}} + \frac{\mathrm{TV}(q_v^{(0)}, \pi)}{2^{\tau_{\ast}}}\biggr) \\
& \leq \frac{1}{4},
\end{align*}
where the final inequality follows from~\eqref{pf:pointwise_mixing_time_1} and the crude upper bound $\mathrm{TV}(q,\pi) \leq 1$, which holds for any distribution $q$. Taking maximum over all $v \in V$, we conclude that $\tau_{\infty}(G) \leq \tau_{\ast} + 1$, which implies the claim of Proposition~\ref{prop:pointwise_mixing_time}. 

It remains to show that $\mathrm{TV}(q_v^{(T)},\pi) \leq 1/4$. Choosing $t_0 = 4$ in the statement of Theorem~\ref{thm:mixing_time_rw}, we have that
\begin{align*}
\mathrm{TV}(q_v^{(T)},\pi) & \leq \frac{1}{16} + \frac{d_{\max}(G)}{d_{\min}(G)^2} + \frac{1}{2s(G)} \Bigl(1 - \frac{\Psi(G)^2}{8}\Bigr)^{T - 4} \\
& \leq \frac{1}{8} + \frac{1}{2s(G)} \Bigl(1 - \frac{\Psi(G)^2}{8}\Bigr)^{T - 4} \\
& \leq \frac{1}{8} + \frac{1}{2s(G)} \exp\Bigl(-\frac{\Psi(G)^2}{8}(T - 4)\Bigr) = \frac{1}{4},
\end{align*}
where the middle inequality follows by assumption. \qed \\ 

\noindent \emph{Proof of Lemma~\ref{lem:tv_to_pointwise}.}
Our goal will be to establish the recurrence relation~\eqref{pf:tv_to_pointwise_2}. 
To derive~\eqref{pf:tv_to_pointwise_2}, the key observation is the following equivalence (see equation (16) of \cite{morris2005}):
\begin{align}
\label{pf:tv_to_pointwise_1}
\frac{\pi(u) - q_v^{(t + 1)}(u)}{\pi(u)} & = \sum_{w \in V} \bigl(\pi(w) - q_v^{(t)}(w) \bigr) \cdot \Bigl(\frac{q_w^{(1)}(u) - \pi(u)}{\pi(u)}\Bigr) \nonumber \\
& = \sum_{w \neq u} \bigl(\pi(w) - q_v^{(t)}(w) \bigr) \cdot \Bigl(\frac{q_w^{(1)}(u) - \pi(u)}{\pi(u)}\Bigr) ~~+ \\
& \quad~\bigl(\pi(u) - q_v^{(t)}(u) \bigr) \cdot \Bigl(\frac{q_u^{(1)}(u) - \pi(u)}{\pi(u)}\Bigr).
\end{align}
We separately upper bound each term on the right hand side of~\eqref{pf:tv_to_pointwise_1}. The sum over all $w \neq u$ can be related to the TV distance between $q_v^{(t)}$ and $\pi$ using H{\"o}lder's inequality,
\begin{align*}
\sum_{w \neq u} \bigl(\pi(w) - q_v^{(t)}(w) \bigr) \cdot \Bigl(\frac{q_w^{(1)}(u) - \pi(u)}{\pi(u)}\Bigr) & \leq 2\mathrm{TV}(q_v^{(t)},\pi) \cdot \max_{w \neq u}\Bigl|\frac{q_w^{(1)}(u) - \pi(u)}{\pi(u)}\Bigr| \\
& \leq 2\mathrm{TV}(q_v^{(t)},\pi) \cdot \max\biggl\{1, \max_{w \neq u} \frac{q_w^{(1)}(u)}{\pi(u)} \biggr\} \\
& \leq 2\mathrm{TV}(q_v^{(t)},\pi) \cdot \frac{m}{d_{\min}(G)^2} = \frac{\mathrm{TV}(q_v^{(t)},\pi) }{s(G)}.
\end{align*}
On the other hand, the second term on the right hand side of~\eqref{pf:tv_to_pointwise_1} satisfies
\begin{align*}
\bigl(\pi(u) - q_v^{(t)}(u) \bigr) \cdot \Bigl(\frac{q_u^{(1)}(u) - \pi(u)}{\pi(u)}\Bigr) \leq \bigl(\pi(u) - q_v^{(t)}(u) \bigr) \cdot \Bigl(\frac{1/2 - \pi(u)}{\pi(u)}\Bigr) \leq \frac{\pi(u) - q_v^{(t)}(u)}{2\pi(u)},
\end{align*}
so that we obtain the recurrence relation
\begin{equation}
\label{pf:tv_to_pointwise_2}
\frac{\pi(u) - q_v^{(t + 1)}(u)}{\pi(u)} \leq \frac{\mathrm{TV}(q_v^{(t)},\pi) }{s(G)} +\frac{\pi(u) - q_v^{(t)}(u)}{2\pi(u)}.
\end{equation}
From~\eqref{pf:tv_to_pointwise_2} along with the initial condition
\begin{equation*}
\Bigl\{\frac{\pi(u) - q_v^{(1)}(u)}{\pi(u)}\Bigr\} \leq 1 \leq 2(1 - \pi(v)) \leq 2\frac{\mathrm{TV}(q_v^{(0)},\pi)}{s(G)},
\end{equation*}
---where the second inequality follows because $\pi(v) \leq d_{\max}(G)/(2m) \leq d_{\max}(G)/d_{\min}(G)^2 \leq 1/16$---we obtain the upper bound
\begin{equation*}
\frac{\pi(u) - q_v^{(t + 1)}(u)}{\pi(u)} \leq \frac{1}{s(G)} \biggl(\sum_{t = 0}^{t_{\ast}} \frac{\mathrm{TV}(q_v^{(t)}, \pi)}{2^{t_{\ast} - t}} + \frac{\mathrm{TV}(q_v^{(0)}, \pi)}{2^{t_{\ast}}}\biggr).
\end{equation*}
This inequality holds for each $u \in V$, and taking the maximum over $u$ completes the proof of Lemma~\ref{lem:tv_to_pointwise}. \qed

\subsection{Spectral Partitioning Properties of PPR}
\label{subsec:ppr_spectral_partitioning}
The following theorem is the main result of Section~\ref{subsec:ppr_spectral_partitioning}. It relates the normalized cut of the sweep sets $\Phi(S_{\beta};G)$ to the normalized cut of a candidate cluster $C \subseteq V$, when $p_v$ is properly initialized within $C$.

\begin{theorem}[\textbf{c.f. Theorem~6 of \cite{andersen2006}}]
	\label{thm:normalized_cut_ppr}
	Suppose that
	\begin{equation}
	\label{eqn:normalized_cut_ppr_vol}
	d_{\max}(G) \leq \vol(C;G) \leq \max\Bigl\{\frac{2}{3}\vol(G); \vol(G) - d_{\max}(G)\Bigr\}
	\end{equation}
	and
	\begin{equation}
	\label{eqn:normalized_cut_ppr_ncut}
	\max\Bigl\{288\Phi(C;G)\cdot \ln\Bigl(\frac{36}{s(G)}\Bigr),72\Phi(C;G) + \frac{d_{\max}(G)}{d_{\min}(G)^2}\Bigr\} < \frac{1}{18}.
	\end{equation}
	Set $\alpha = 36 \cdot \Phi(C;G)$. The following statement holds: there exists a set $C^g \subseteq C$ of large volume, $\vol(C^g;G) \geq 5/6 \cdot \vol(C;G)$, such that for any $v \in C^g$,  the minimum normalized cut of the sweep sets of $p_v$ satisfies 
	\begin{equation}
	\label{eqn:normalized_cut_ppr}
	\min_{\beta \in (0,1)}\Phi(S_{\beta,v};G) < 72\sqrt{\Phi(C;G) \cdot \ln\Bigl(\frac{36}{s(G)}\Bigr)}.
	\end{equation}
\end{theorem}
A few remarks:
\begin{itemize}
	\item Theorem~\ref{thm:normalized_cut_ppr} is similar to Theorem 6 of \citet{andersen2006}, but crucially the above bound depends on $\log\bigl(1/s(G)\bigr)$ rather than $\log m$. In the case where $d_{\min}(G)^2 \asymp \vol(G)$ and thus $s(G) \asymp 1$, this amounts to replacing a factor of $O(\log m)$ by a factor of ${O}(1)$, and therefore allows us to obtain meaningful results in the limit as $m \to \infty$. 
	\item For simplicity, we have chosen to state Theorem~\ref{thm:normalized_cut_ppr} with respect to a specific choice of $\alpha = 36 \cdot \Phi(C;G)$, but if $\alpha \approx 36 \cdot \Phi(C;G)$ then the Theorem will still hold up to constant factors.
\end{itemize}

It follows from Markov's inequality (see Theorem~4 of \cite{andersen2006}) that there exists a set $C^g \subseteq C$ of volume $\vol(C^g;G) \geq 5/6 \cdot \vol(C;G)$ such that for any $v \in C^g$,
\begin{equation}
\label{eqn:ppr_leakage}
p_v(C) \geq 1 - \frac{6\Phi(C;G)}{\alpha}.
\end{equation}
The claim of Theorem~\ref{thm:normalized_cut_ppr} is a consequence of~\eqref{eqn:ppr_leakage} along with Theorem~\ref{thm:mixing_time_PPR}, as we now demonstrate. \\

\noindent \emph{Proof of Theorem~\ref{thm:normalized_cut_ppr}.}
From~\eqref{eqn:ppr_leakage}, the upper bound in~\eqref{eqn:normalized_cut_ppr_vol}, and the choice of $\alpha = 36 \cdot \Phi(C;G)$, 
\begin{equation}
\label{pf:normalized_cut_ppr}
p_v(C) - \pi(C) \geq \frac{1}{3} - \frac{6 \Phi(C;G)}{\alpha} = \frac{1}{6}.
\end{equation}
Now, put 
\begin{equation*}
t_{\ast} = \frac{1}{648 \Phi(C;G)},~~\phi_{\ast}^2 = \frac{8}{t_{\ast}} \cdot \ln\Bigl(\frac{36}{s(G)}\Bigr),
\end{equation*}
and note that by~\eqref{eqn:normalized_cut_ppr_ncut} $\phi_{\ast}^2 \in [0,1]$. It therefore follows from~\eqref{pf:normalized_cut_ppr} and Theorem~\ref{thm:mixing_time_PPR} that either
\begin{equation}
\label{pf:normalized_cut_ppr_2}
\frac{1}{6} \leq p_v(C) - \pi(C) \leq \frac{1}{18} + 72 \Phi(C;G) + \frac{d_{\max}(G)}{d_{\min}(G)^2} + \frac{1}{2s(G)} \cdot \left(1 - \frac{\phi_{\ast}^2}{8}\right)^{t_{\ast}},
\end{equation}
or $\min_{\beta \in (0,1)} \Phi(S_{\beta,v};G) \leq \phi_{\ast}^2$. But by~\eqref{eqn:normalized_cut_ppr_ncut}
\begin{equation*}
72 \Phi(C;G) + \frac{d_{\max}(G)}{d_{\min}(G)^2} < \frac{1}{18},
\end{equation*}
and we have chosen $\phi_{\ast}$ precisely so that
\begin{equation*}
\frac{1}{2s(G)} \cdot \left(1 - \frac{\phi_{\ast}^2}{8}\right)^{t_{\ast}} \leq \frac{1}{2s(G)} \exp\Bigl(-\frac{\phi_{\ast}^2 t_{\ast}}{8}\Bigr) \leq \frac{1}{18}.
\end{equation*}
Thus the inequality~\eqref{pf:normalized_cut_ppr_2} cannot hold, and so it must be that $\min_{\beta \in (0,1)} \Phi(S_{\beta,v};G) \leq \phi_{\ast}^2$. This is exactly the claim of the theorem. \qed

\section{Sample-to-Population Bounds}
\label{apdx:sample_to_population}
In this appendix, we prove Propositions~\ref{prop:sample_to_population_1} and~\ref{prop:sample_to_population_2}, by establishing high-probability finite-sample bounds on various functionals of the random graph $G_{n,r}$: cut, volume, and normalized cut~(\ref{subsec:sample_to_population_ncut}), minimum and maximum degree, and local spread~(\ref{subsec:sample_to_population_local_spread}), and conductance~(\ref{subsec:sample_to_population_conductance}). To establish these results, we will use several different concentration inequalities, and we begin by reviewing these in~(\ref{subsec:concentration}). Throughout, we denote the empirical probability of a set $\mc{S} \subseteq \Rd$ as $\Pbb_n(\mc{S}) = \sum_{i = 1}^{n} \1\{x_i \in \mc{S}\}/n$, and the conditional (on being in $\mc{C}$) empirical probability as $\wt{\Pbb}_n = \sum_{i = 1}^{n} \1\{x_i \in (\mc{S} \cap \mc{C})\}/\wt{n}$, where $\wt{n} = |\mc{C}[X]|$ is the number of sample points that are in $\mc{C}$. For a probability measure $\Qbb$, we also write
\begin{equation}
\label{eqn:population_minmax_degree}
d_{\min}(\Qbb) := \inf_{x \in \mathrm{supp}(\mbb{Q})} \deg_{\Pbb,r}(x),~~\textrm{and}~~d_{\max}(\Qbb) := \sup_{x \in \mathrm{supp}(\mbb{Q})} \deg_{\Pbb,r}(x).
\end{equation}

\subsection{Review: Concentration Inequalities}
\label{subsec:concentration}
We use Bernstein's inequality to control the deviations of the empirical probability of $\mc{S}$.
\begin{lemma}[Bernstein's Inequality.]
	\label{lem:hoeffding_2}
	Fix $\delta \in (0,1)$. For any measurable $\mc{S} \subseteq \Rd$, each of the inequalities,
	\begin{equation*}
	(1 - \delta) \Pbb(\mathcal{S}) \leq \Pbb_n(\mathcal{S})~~\textrm{and}~~\Pbb_n(\mathcal{S}) \leq (1 + \delta)\Pbb(\mathcal{S}),
	\end{equation*}
	hold with probability at least $1 - \exp\set{-n\delta^2\Pbb(\mathcal{S})/(2 + 2\delta)} \geq 1 - \exp\set{-n\delta^2\Pbb(\mathcal{S})/4}$. 
\end{lemma}

Many graph functionals are order-2 U-statistics, and we use Bernstein's inequality to control the deviations of these functionals from their expectations. Recall that $U_n$ is an order-2 U-statistic with kernel $\varphi: \Rd \times \Rd \to \Reals$ if 
\begin{equation*}
U_n = \frac{1}{n(n-1)}\sum_{i = 1}^{n} \sum_{j \neq i} \varphi(x_i,x_j).
\end{equation*}
We write $\norm{\varphi}_{\infty} = \sup_{x,y} |\varphi(x,y)|$. 
\begin{lemma}[Bernstein's Inequality for Order-2 U-statistics.]
	\label{lem:hoeffding}
	Fix $\delta \in (0,1)$. Assume $\norm{\varphi}_{\infty} \leq 1$. Then each of the inequalities,
	\begin{equation*}
	(1 - \delta) \mathbb{E}U_n \leq U_n~~\textrm{and}~~ U_n\leq (1 + \delta) \mathbb{E}U_n,
	\end{equation*}
	hold with probability at least $1 - \exp\{-n\delta^2 \Ebb U_n/(4 + 4\delta/3)\} \geq 1 - \exp\{-n\delta^2 \Ebb U_n/6\}$.
\end{lemma}

Finally, we use Lemma~\ref{lem:bernstein_union}---a combination of Bernstein's inequality and a union bound---to upper and lower bound $d_{\max}(G_{n,r})$ and $d_{\min}(G_{n,r})$. For measurable sets $\mc{S}_1,\ldots,\mc{S}_M$, we denote $p_{\min} := \min_{m = 1,\ldots,M} \Pbb(\mathcal{A}_m)$, and likewise let $p_{\max} := \max_{m = 1,\ldots,M} \Pbb(\mathcal{A}_m)$
\begin{lemma}[Bernstein's inequality + union bound.]
	\label{lem:bernstein_union}
	Fix $\delta \in (0,1)$. For any measurable $\mathcal{S}_1,\ldots,\mathcal{S}_M \subseteq \Rd$, each of the inequalities
	\begin{equation*}
	(1 - \delta) p_{\min} \leq \min_{m = 1,\ldots,M} \Pbb_n(\mathcal{A}_m),~~\textrm{and}~~\max_{m = 1,\ldots,M}  \Pbb_n(\mathcal{A}_m) \leq (1 + \delta) p_{\max}
	\end{equation*}
	hold with probability at least $1 - M \exp \{-n\delta^2p_{\min}/(2 + 2\delta) \} \geq 1 - M \exp\{-n\delta^2p_{\min}/4\}$.  
\end{lemma}

\subsection{Sample-to-Population: Normalized Cut}
\label{subsec:sample_to_population_ncut}
In this subsection we establish~\eqref{eqn:sample_to_population_normalized_cut}. For a set $\mc{S} \subseteq \Rd$, both $\mathrm{cut}_{n,r}(\mc{S}[X])$ and $\mathrm{vol}_{n,r}(\mc{S}[X])$ are order-$2$ U-statistics:
\begin{align*}
\mathrm{cut}_{n,r}(\mc{S}[X]) & = \sum_{i = 1}^{n} \sum_{j \neq i} \1\{\|x_i - x_j\| \leq r\} \cdot \1\{x_i \in \mc{S}\} \cdot \1\{x_j \not\in \mc{S}\},
\intertext{and}
\mathrm{vol}_{n,r}(\mc{S}[X]) & = \sum_{i = 1}^{n} \sum_{j \neq i} \1\{\|x_i - x_j\| \leq r\} \cdot \1\{x_i \in \mc{S}\}.
\end{align*}
Therefore with probability at least $1 - \exp\{-n\delta^2\cut_{\Pbb,r}(\mc{S})/4\}$,
\begin{equation*}
\frac{1}{n(n - 1)} \mathrm{cut}_{n,r}(\mc{S}[X]) \leq (1 + \delta) \cut_{\Pbb,r}(\mc{S}),
\end{equation*}
and likewise with probability at least $1 - \exp\{-n\delta^2 \vol_{\Pbb,r}(\mc{S})/4\} - \exp\{-n\delta^2 \vol_{\Pbb,r}(\mc{S}^c)/4\}$,
\begin{equation*}
(1 - \delta) \vol_{\Pbb,r}(\mc{S}) \leq \frac{1}{n(n - 1)} \vol_{n,r}(\mc{S}[X]),~~\textrm{and}~~(1 - \delta) \vol_{\Pbb,r}(\mc{S}^c) \leq \frac{1}{n(n - 1)} \vol_{n,r}(\mc{S}^c[X])
\end{equation*}
Consequently, for any $\delta \in (0,1/3)$,  
\begin{equation*}
\Phi_{n,r}(\mc{C}[X]) \leq \frac{1 + \delta}{1 - \delta} \cdot  \frac{\cut_{\Pbb,r}(\mc{C})}{\min\{\vol_{\Pbb,r}(\mc{C}), \vol_{\Pbb,r}(\mc{C}^c)\}} = \frac{1 + \delta}{1 - \delta} \cdot \Phi_{\Pbb,r}(\mc{C}) \leq (1 + 3\delta) \cdot \Phi_{\Pbb,r}(\mc{C})  
\end{equation*}
with probability at least $1 - 3\exp\{-n \delta^2 \cut_{\Pbb,r}(\mc{C})/4\}$. This establishes~\eqref{eqn:sample_to_population_normalized_cut} upon taking $b_1 := 3\cut_{\Pbb,r}(\mc{C})/4$. \qed

\subsection{Sample-to-Population: Local Spread}
\label{subsec:sample_to_population_local_spread}
In this subsection we establish~\eqref{eqn:sample_to_population_local_spread}. To ease the notational burden, let $\wt{G}_{n,r} := G_{n,r}\bigl[\mc{C}[X]\bigr]$. 
Conditional on~$\wt{n}$, it follows from Lemma~\ref{lem:bernstein_union} that with probability at least $1 - \wt{n}\exp\{-(\wt{n} - 1) \delta^2d_{\min}(\wt{\Pbb})/4\}$,
\begin{equation}
\label{pf:sample_to_population_local_spread_2}
(1 - \delta) \cdot d_{\min}(\wt{\Pbb}) \leq \frac{1}{\wt{n} - 1} d_{\min}(\wt{G}_{n,r}),
\end{equation}
Likewise it follows from Lemma~\ref{lem:hoeffding} that with probability at least $1 - \exp\{-\wt{n}\delta^2 \vol_{\wt{\Pbb},r}(\mc{C})/6\}$,
\begin{equation*}
(1 - \delta) \cdot \vol_{\wt{\Pbb},r}(\mc{C}) \leq \frac{1}{\wt{n}(\wt{n} - 1)} \vol(\wt{G}_{n,r}).
\end{equation*}
Finally, it follows from Lemma~\ref{lem:hoeffding_2} that with probability at least $1 - \exp\{-n\delta^2 \Pbb(\mc{C})/4\}$
\begin{equation}
\label{pf:sample_to_population_local_spread_1}
\wt{n} \geq (1 - \delta) \cdot n \cdot \Pbb(\mc{C}),
\end{equation}
and therefore by~\eqref{eqn:sample_to_population_local_spread_sample_complexity}, $(\wt{n} - 1)/\wt{n} \geq 1 - \delta$. Consequently for any $\delta \in (0,1/3)$,
\begin{equation*}
s_{n,r}(\mc{C}[X]) = \frac{d_{\min}(\wt{G}_{n,r})^2}{\vol(\wt{G}_{n,r})} = \frac{\wt{n} - 1}{\wt{n}} \cdot \frac{\frac{1}{(\wt{n} - 1)^2}d_{\min}(\wt{G}_{n,r})^2}{\frac{1}{\wt{n}(\wt{n} - 1)}\vol(\wt{G}_{n,r})} \geq \frac{(1 - \delta)^3}{(1 + \delta)} \cdot  \frac{d_{\min}(\wt{\Pbb})^2}{\vol_{\wt{\Pbb},r}(\mc{C})} \geq (1 - 4\delta) \cdot s_{\Pbb,r}(\mc{C}).
\end{equation*}
with probability at least $1 - n\exp\{-n\Pbb(\mc{C}) \cdot \delta^2d_{\min}(\wt{\Pbb})/9\} - \exp\{-n\Pbb(\mc{C})\delta^2 \cdot \vol_{\wt{\Pbb},r}(\mc{C})/14\} - \exp\{-n\delta^2 \Pbb(\mc{C})/4\}$. This establishes~\eqref{eqn:sample_to_population_local_spread} upon taking $b_2 := \Pbb(\mc{C}) \cdot d_{\min}(\wt{\Pbb})/14$. \qed


\subsection{Sample-to-Population: Conductance}
\label{subsec:sample_to_population_conductance}

In this section we establish~\eqref{eqn:sample_to_population_conductance}. As mentioned in our main text, the proof of~\eqref{eqn:sample_to_population_conductance} relies on a high-probability upper bound of the $\infty$-transportation distance between $\mbb{P}$ and $\mbb{P}_n$, from~\citep{garciatrillos16b}. We begin by reviewing this upper bound, which we restate in Theorem~\ref{thm:garciatrillos16}. Subsequently in Proposition~\ref{prop:conductance_lb_transportation_distance}, we relate the $\infty$-transportation distance between two measures $\mbb{Q}_1$ and $\mbb{Q}_2$ to the difference of their conductances. Together these results will imply~\eqref{eqn:sample_to_population_conductance}.\\

\noindent \emph{Review: $\infty$-transportation distance and transportation maps.}
We give a brief review of some of the main ideas regarding $\infty$-transportation distance, and transportation maps. This discussion is largely taken from~\citep{garciatrillos16b,garciatrillos16}, and the reader should consult these works for more detail. 

For two measures $\mbb{Q}_1$ and $\mbb{Q}_2$ on a domain $D$, the \emph{$\infty$-transportation distance} $\Delta_{\infty}(\mbb{Q}_1,\mbb{Q}_2)$ is
\begin{equation*}
\Delta_{\infty}(\mbb{Q}_1,\mbb{Q}_2) := \inf_{\gamma}\Bigl\{\mathrm{ess sup}_{\gamma} \bigl\{|x - y|: (x,y) \in D \times D\bigr\}: \gamma \in \Gamma(\mbb{Q}_1,\mbb{Q}_2) \Bigr\}
\end{equation*}
where $\Gamma(\mbb{Q}_1,\mbb{Q}_2)$ is the set of all couplings of $\mbb{Q}_1$ and $\mbb{Q}_2$, that is the set of all probability measures on $D \times D$ for which the marginal distribution in the first variable is $\mbb{Q}_1$, and the marginal distribution in the second variable is $\mbb{Q}_2$. 

Suppose $\mbb{Q}_1$ is absolutely continuous with respect to the Lebesgue measure. Then $\Delta_{\infty}(\mbb{Q}_1,\mbb{Q}_2)$ can be more simply defined in terms of push-forward measures and transportation maps. For a Borel map $T: D \to D$, the \emph{push-forward} of $\mbb{Q}_1$ by $T$ is $T_{\sharp}\mbb{Q}_1$, defined for Borel sets $U$ as 
\begin{equation*}
T_{\sharp}\mbb{Q}_1(U) = \mbb{Q}_1(T^{-1}(U)).
\end{equation*}
A \emph{transportation map} from $\mbb{Q}_1$ to $\mbb{Q}_2$ is a Borel map $T$ for which $T_{\sharp}\mbb{Q}_1 = \mbb{Q}_2$. Transportation maps satisfy two important properties. First, the transportation distance can be formulated in terms of transportation maps:
\begin{equation*}
\Delta_{\infty}(\mbb{Q}_1,\mbb{Q}_2) = \inf_{T} \|\mathrm{Id} - T\|_{\Leb^{\infty}(\mbb{Q}_1)}
\end{equation*}
where $\mathrm{Id}: D \to D$ is the identity mapping, and the infimum is over transportation maps $T$ from $\mbb{Q}_1$ to $\mbb{Q}_2$. Second, they result in the following change of variables formula; if $T_{\sharp}\mbb{Q}_1 = \mbb{Q}_2$, then for any $g \in L^1(\mbb{Q}_2)$,
\begin{equation}
\label{eqn:transportation_map_change_of_variable}
\int g(y) \,d\mbb{Q}_2(y) = \int g(T(x)) \,d\mbb{Q}_1(x).
\end{equation} 

\noindent \emph{$\infty$-transportation distance between empirical and population measures.}
We now review the relevant upper bound on $\Delta_{\infty}(\Pbb,\Pbb_n)$, which holds under the following mild regularity conditions.
\begin{enumerate}[label=(A\arabic*)]
	\setcounter{enumi}{10}
	\item 
	\label{asmp:bounded_density_full} 
	The distribution $\Pbb$ has density $g: D \to (0,\infty)$ such that there exist $g_{\min} \leq 1 \leq g_{\max}$ for which
	\begin{equation*}
	(\forall x \in D)~~ g_{\min} \leq g(x) \leq g_{\max}.
	\end{equation*}
	\item 
	\label{asmp:domain_full} 
	The distribution $\Pbb$ is defined on a bounded, connected, open domain $D \subseteq \Rd$. If $d \geq 2$ then additionally $D$ has Lipschitz boundary.
\end{enumerate}
When $d = 1$, it follows from Proposition~6.2 of \cite{dudley1968} that $\Delta_{\infty}(\Pbb,\Pbb_n) \leq B_5 \|F - F_n\|_{\infty}$ for some positive constant $B_5$, and in turn from the DKW inequality that
\begin{equation}
\label{eqn:ot_distance_1d}
\Delta_{\infty}(\Pbb,\Pbb_n) \leq B_5 \sqrt{\frac{\ln(2n/B_2)}{n}}
\end{equation}
with probability at least $1 - B_2/n$.

When $d \geq 2$, \cite{garciatrillos16b} derive an upper bound on the transportation distance $\Delta_{\infty}(\mbb{P},\mbb{P}_n)$.
\begin{theorem}[Theorem~1.1 of \cite{garciatrillos16b}]
	\label{thm:garciatrillos16}
	Suppose $\Pbb$ satisfies~\ref{asmp:bounded_density_full} and~\ref{asmp:domain_full}. Then, there exists positive constants $B_{2}$ and $B_5$ that do not depend on $n$, such that with probability at least $1 - B_2/n$:
	\begin{equation*}
	\Delta_{\infty}(\mbb{P},\mbb{P}_n) \leq B_5 \cdot 
	\begin{dcases*}
	\frac{\ln(n)^{3/4}}{n^{1/2}},&~~\textrm{if $d = 2$,} \\
	\frac{\ln(n)^{1/d}}{n^{1/d}},&~~\textrm{if $d \geq 3$.}
	\end{dcases*}
	\end{equation*}
\end{theorem}
Assuming the candidate cluster $\mc{C}$ and conditional distribution $\wt{\Pbb}$ satisfy~\ref{asmp:bounded_density} and~\ref{asmp:domain}, then~\eqref{eqn:ot_distance_1d} ($d = 1)$ or Theorem~\ref{thm:garciatrillos16} ($d \geq 2$) apply to $\Delta_{\infty}(\wt{\Pbb},\wt{\Pbb}_n)$; we will use these upper bounds on $\Delta_{\infty}(\wt{\Pbb},\wt{\Pbb}_n)$ to show~\eqref{eqn:sample_to_population_conductance}.\\

\noindent \emph{Lower bound on conductance using transportation maps.}
Let $\mbb{Q}_1$ and $\mbb{Q}_2$ be probability measures, with $\mbb{Q}_1$ absolutely continuous with respect to Lebesgue measure, and let $T$ be a transportation map from $\mbb{Q}_1$ to $\mbb{Q}_2$. We write $\Delta_T(\mbb{Q}_1,\mbb{Q}_2) := \|\mathrm{Id} - T\|_{\Leb^{\infty}(\mbb{Q}_1)}$. To facilitate easy comparison between the conductances of two arbitrary distributions, let $\Psi_r(\mbb{Q}) := \Psi_{\mbb{Q},r}(\mathrm{supp}(\mbb{Q}))$ for a distribution $\mbb{Q}$. In the following Proposition, we lower bound $\Psi_r(\mbb{Q}_2)$ by $\Psi_r(\mbb{Q}_1)$, plus an error term that depends on $\Delta(\mbb{Q}_1,\mbb{Q}_2)$.
\begin{proposition}
	\label{prop:conductance_lb_transportation_distance}
	Let $\mbb{Q}_1$ be a probability measure that admits a density $g$ with respect to $\nu(\cdot)$, let $\mbb{Q}_2$ be an arbitrary probability measure, and let $T$ be a transportation map from $\mbb{Q}_1$ to $\mbb{Q}_2$. Suppose $\Delta_T(\mbb{Q}_1,\mbb{Q}_2) \leq r/(4(d - 1))$. It follows that
	\begin{equation}
	\label{eqn:conductance_lb_transportation_distance}
	\Psi_r(\mbb{Q}_2) \geq \Psi_r(\mbb{Q}_1) \cdot \biggl(1 - \frac{B_6\Delta_T(\mbb{Q}_1,\mbb{Q}_2)}{\bigl(1 - \Psi_r(\mbb{Q}_1)\bigr) \cdot \bigl(d_{\min}(\mbb{Q}_2)\bigr)^2}\biggr) - \frac{B_6 \Delta_T(\mbb{Q}_1,\mbb{Q}_2)}{\bigl(1 - \Psi_r(\mbb{Q}_1)\bigr) \cdot \bigl(d_{\min}(\mbb{Q}_2)\bigr)^2},
	\end{equation}
	where $B_6 := 4d \nu_d r^{d - 1} \cdot \max_{x \in \Rd}\{g(x)\}$ is a positive constant that does not depend on $\mbb{Q}_2$. 
\end{proposition}
We note that the lower bound can also be stated with respect to the $\infty$-optimal transport distance $\Delta_{\infty}(\mbb{Q}_1,\mbb{Q}_2)$.\\

\noindent \emph{Proof of Proposition~\ref{prop:conductance_lb_transportation_distance}.}
Throughout this proof, we will write $\Delta_{12} = \Delta_T(\mbb{Q}_1,\mbb{Q}_2)$, and $\wb{\vol}_{\mbb{Q},r}(\mc{R}) = \min\bigl\{\vol_{\mbb{Q},r}\bigl(\mc{R}\bigr),\vol_{\mbb{Q},r}\bigl(\mc{R}^c\bigr)\bigr\}$ for conciseness. Naturally, the proof of Proposition~\ref{prop:conductance_lb_transportation_distance} involves using the transportation map $T$ to relate $\cut_{\mbb{Q}_2,r}(\cdot)$ to $\cut_{\mbb{Q}_1,r}(\cdot)$, and likewise $\vol_{\mbb{Q}_2,r}(\cdot)$ to $\vol_{\mbb{Q}_1,r}(\cdot)$. Define the remainder term $R_{\epsilon,\mbb{Q}_1}^{(\Delta)}(x) = \int \1\{\epsilon \leq \|x - y\| \leq \epsilon + \Delta\} \,d\mbb{Q}_1(y)$ for any $\epsilon, \Delta > 0$. Then for any set $\mc{S} \subseteq \mathrm{supp}(\mbb{Q}_2)$, we have that
\begin{align}
\cut_{\mbb{Q}_2,r}(\mc{S}) & = \iint \1\{\|x - y\|\leq r\} \cdot \1\{x \in \mc{S} \} \cdot \1\{y \in \mc{S}^c \} \,d\mbb{Q}_2(y) \,d\mbb{Q}_2(x) \nonumber \\ 
& \overset{\mathrm{(i)}}{=} \iint \1\{\|T(x) - T(y)\|\leq r\} \cdot \1\{x \in T^{-1}(\mc{S}) \} \cdot \1\{y \in T^{-1}(\mc{S})^c \} \,d\mbb{Q}_1(y) \,d\mbb{Q}_1(x) \nonumber \\
& \overset{\mathrm{(ii)}}{\geq} \iint \1\{\|x - y\|\leq r - 2\Delta_{12}\} \cdot \1\{x \in T^{-1}(\mc{S}) \} \cdot \1\{y \in T^{-1}(\mc{S}^c) \} \,d\mbb{Q}_1(y) \,d\mbb{Q}_1(x) \nonumber \\
& = \cut_{\mbb{Q}_1,r}\bigl(T^{-1}(\mc{S})\bigr) - \int R_{r - 2\Delta_{12} ,\mbb{Q}_1}^{(2\Delta_{12})}(x) \,d\mbb{Q}_1(x) \label{pf:conductance_lb_transportation_distance_1}
\end{align} 
where $\mathrm{(i)}$ follows from the change of variables formula~\eqref{eqn:transportation_map_change_of_variable}, and $\mathrm{(ii)}$ follows from the triangle inequality. Similar reasoning implies that
\begin{equation}
\label{pf:conductance_lb_transportation_distance_2}
\vol_{\mbb{Q}_2,r}(\mc{S}) \leq \vol_{\mbb{Q}_1,r}\bigl(T^{-1}(\mc{S})\bigr) + \int R_{r,\mbb{Q}_1}^{(2\Delta_{12})}(x)  \,d\mbb{Q}_1(x).
\end{equation}
For any $x \in \Rd$, since $0 \leq \Delta_{12} \leq r/(4(d - 1))$, the remainder terms can be upper bounded as follows:
\begin{equation*}
R_{r - 2\Delta_{12} ,\mbb{Q}_1}^{(2\Delta_{12})}(x) \leq \nu_dr^d\Bigl\{1 - \Bigl(1 - \frac{2\Delta_{12}}{r}\Bigr)^d \Bigr\} \cdot \max_{x \in \Rd}\{g(x)\} \leq \underbrace{4 d \nu_d r^{d - 1} \cdot \max_{x \in \Rd}\{g(x)\}}_{= B_6} \cdot \Delta_{12},
\end{equation*}
and
\begin{equation*}
R_{r,\mbb{Q}_1}^{(2\Delta_{12})}(x) \leq \nu_dr^d\Bigl\{\Bigl(1 + \frac{2\Delta_{12}}{r}\Bigr)^d - 1\Bigr\} \cdot \max_{x \in \Rd}\{g(x)\} \leq B_6 \cdot \Delta_{12}.
\end{equation*} 
Plugging these bounds on the remainder terms back into~\eqref{pf:conductance_lb_transportation_distance_1} and~\eqref{pf:conductance_lb_transportation_distance_2} respectively, we see that
\begin{align*}
\Phi_{\mbb{Q}_2,r}(\mc{S}) & \geq \frac{\cut_{\mbb{Q}_1,r}\bigl(T^{-1}(\mc{S})\bigr) - B_6 \Delta_{12}}{\wb{\vol}_{\mbb{Q}_1,r}(T^{-1}(\mc{S})) + B_6 \Delta_{12}} \\ & = \Phi_{\mbb{Q}_1,r}(T^{-1}(\mc{S})) \cdot \biggl(\frac{\wb{\vol}_{\mbb{Q}_1,r}(T^{-1}(\mc{S}))}{\wb{\vol}_{\mbb{Q}_1,r}(T^{-1}(\mc{S})) + B_6\Delta_{12}}\biggr) - \frac{B_6\Delta_{12}}{\wb{\vol}_{\mbb{Q}_1,r}(T^{-1}(\mc{S})) + B_6\Delta_{12}} \\
& \overset{\eqref{pf:conductance_lb_transportation_distance_2}}{\geq} \Phi_{\mbb{Q}_1,r}(T^{-1}(\mc{S})) \cdot \biggl(\frac{\wb{\vol}_{\mbb{Q}_2,r}(\mc{S}) - B_6\Delta_{12}}{\wb{\vol}_{\mbb{Q}_2,r}(\mc{S})}\biggr) - \frac{B_6\Delta_{12}}{\wb{\vol}_{\mbb{Q}_2,r}(\mc{S})}.
\end{align*}
We would like to conclude by taking an infimum over $\mc{S}$ on both sides, but in order to ensure that the remainder term is small we must specially handle the case where $\wb{\vol}_{\mbb{Q}_2,r}(\mc{S})$ is small. Let
\begin{equation*}
\mathfrak{L}_{r}(\mbb{Q}_1,\mbb{Q}_2) = \bigl\{\mc{S} \subseteq \mathrm{supp}(\mbb{Q}_2): \wb{\vol}_{\mbb{Q}_2,r}(\mc{S}) \geq (1 - \Psi_r(\mbb{Q}_1)) \cdot d_{\min}(\mbb{Q}_2)^2 \bigr\}.
\end{equation*}
On the one hand, taking an infimum over all sets $\mc{S} \in \mathfrak{L}_{r}(\mbb{Q}_1,\mbb{Q}_2)$, we have that
\begin{equation*}
\inf_{\mc{S}: \mc{S} \in \mathfrak{L}_{r}(\mbb{Q}_1,\mbb{Q}_2)} \Phi_{\mbb{Q}_2,r}(\mc{S}) \geq \Psi_{r}(\mbb{Q}_1) \cdot \biggl(1 - \frac{B_6 \Delta_{12}}{(1 - \Psi_r(\mbb{Q}_1)) \cdot d_{\min}(\mbb{Q}_2)^2}\biggr) - \frac{B_6\Delta_{12}}{(1 - \Psi_r(\mbb{Q}_1)) \cdot d_{\min}(\mbb{Q}_2)^2}
\end{equation*}
On the other hand, we claim that 
\begin{equation}
\label{pf:conductance_lb_transportation_distance_3}
\Phi_{r,\mbb{Q}_2}(\mc{R}) \geq \Psi_{r}(\mbb{Q}_1),~~\textrm{for any $\mc{R} \not\in \mathfrak{L}(\mbb{Q}_1,\mbb{Q}_2)$}.
\end{equation}
To derive~\eqref{pf:conductance_lb_transportation_distance_3}, suppose that $\mc{R} \subseteq \mathrm{supp}(\mbb{Q}_2)$ and $\mc{R} \not\in \mathfrak{L}(\mbb{Q}_1,\mbb{Q}_2)$. Without loss of generality, we shall assume that $\vol_{\mbb{Q}_2,r}(\mc{R}) \leq (1 - \Psi_r(\mbb{Q}_1)) \cdot d_{\min}(\mbb{Q}_2)^2$ (otherwise we can work with respect to $\mc{R}^c$.)  Then, for all $x \in \mc{R}$,
\begin{align*}
\int \1\{\|x - y\| \leq r\}\cdot \1\{y \in \mc{R}^c\} \,d\mbb{Q}_2(y) & \geq \deg_{\Qbb_2,r}(x) -  \mbb{Q}_2(\mc{R}) \\
& \geq \deg_{\Qbb_2,r}(x) - \frac{\vol_{\mbb{Q}_2,r}(\mc{R})}{d_{\min}(\mbb{Q}_2)} \\
& \geq d_{\min}(\mbb{Q}_2) \cdot \Psi_{r}(\mbb{Q}_2),
\end{align*}
whence integrating over all $x \in \mc{R}$ and dividing by~$\vol_{\mbb{Q}_2,r}(\mc{R})$ yields~\eqref{pf:conductance_lb_transportation_distance_3}. This completes the proof of Proposition~\ref{prop:conductance_lb_transportation_distance}. \\

\noindent \emph{Putting the pieces together.} 
First, we note that
\begin{equation*}
\Psi_r(\wt{\mbb{P}}) = \Psi_{\wt{\mbb{P}},r}(\mathrm{supp}(\wt{\mbb{P}})) = \Psi_{\Pbb,r}(\mc{C}), ~~\textrm{and}~~\Psi_{r}(\wt{\mbb{P}}_n) = \Psi_{\wt{\mbb{P}}_n, r}(\mathrm{supp}(\wt{\mbb{P}}_n)) = \Psi_{n,r}(\mc{C}[X]),
\end{equation*}
so that we may apply Proposition~\ref{prop:conductance_lb_transportation_distance} to get a lower bound on $\Psi_{n,r}(\mc{C}[X])$ in terms of $\Psi_{\Pbb,r}(\mc{C})$, $\Delta_{\infty}(\wt{\Pbb}, \wt{\Pbb}_n)$, and $d_{\min}(\wt{\Pbb}_{n})$. To begin, we recall from Section~\ref{subsec:sample_to_population_local_spread} that the lower bound on minimum degree,
\begin{equation*}
d_{\min}(\wt{\Pbb}_n) = \frac{1}{\wt{n}} \bigl(d_{\min}(\wt{G}_{n,r}) + 1\bigr) \geq \frac{1}{\sqrt{2}} d_{\min}(\wt{\Pbb}),
\end{equation*}
is satisfied with probability at least $1 - (n + 1)\exp\{-nb_2/16\}$. On the other hand, taking
\begin{equation*}
b_6 := \frac{1}{B_6}\Psi_r(\wt{\Pbb}) \cdot (1 - \Psi_r(\wt{\Pbb})) \cdot d_{\min}(\wt{\Pbb})^2,~~\textrm{and}~~ B_1 := B_5 \Bigl(\min\Bigl\{b_6, \frac{r}{4(d - 1)} \Bigr\}\Bigr)^{-1},
\end{equation*}
by~\eqref{eqn:ot_distance_1d} (if $d = 1$) or Theorem~\ref{thm:garciatrillos16} (if $d \geq 2$) along with~\eqref{eqn:sample_to_population_conductance_sample_complexity}, we have that
\begin{equation*}
\Delta_{\infty}(\wt{\Pbb},\wt{\Pbb}_n) \leq B_5 \frac{(\log n)^{p_d}}{\min\{n^{1/2},n^{1/d}\}} \leq \min\Bigl\{b_6, \frac{r}{4(d - 1)} \Bigr\} \cdot \delta
\end{equation*}
with probability at least $1 - B_2/n$. Finally, appealing first to  Proposition~\ref{prop:conductance_lb_transportation_distance} and then to the bounds we have just established on $d_{\min}(\wt{\Pbb}_n)$ and $\Delta_{\infty}(\wt{\Pbb},\wt{\Pbb}_n)$, we conclude that the sample conductance is lower bounded,
\begin{align*}
\Psi_r(\wt{\Pbb}_n) & \geq \Psi_r(\wt{\Pbb}) \cdot \biggl(1 - \frac{B_6\Delta_{\infty}(\wt{\Pbb}_n,\wt{\Pbb})}{\bigl(1 - \Psi_r(\wt{\Pbb})\bigr) \cdot \bigl(d_{\min}(\wt{\Pbb}_n)\bigr)^2}\biggr) - \frac{B_6 \Delta_{\infty}(\wt{\Pbb}_n,\wt{\Pbb})}{\bigl(1 - \Psi_r(\wt{\Pbb})\bigr) \cdot \bigl(d_{\min}(\wt{\Pbb}_n)\bigr)^2} \\
& \geq \Psi_r(\wt{\Pbb}) (1 - 2\delta),
\end{align*}
with probability at least $1 - B_2/n - (n + 1)\exp\{-nb_2/16\}$, establishing~\eqref{eqn:sample_to_population_conductance} upon taking $b_3 := b_2/16$. \qed

\section{Population Functionals for Density Clusters}
\label{apdx:density_cluster_population_functionals}
In this appendix, we prove Lemma~\ref{lem:density_cluster_local_spread} (in Section~\ref{subsec:density_cluster_local_spread}), Proposition~\ref{prop:density_cluster_normalized_cut} (in Section~\ref{subsec:density_cluster_ncut}), and Proposition~\ref{prop:density_cluster_conductance} (in Section~\ref{subsec:density_cluster_conductance}), by establishing bounds on the population-level local spread, normalized cut, and conductance of a thickened density cluster $\mc{C}_{\lambda,\sigma}$. In these proofs, we make use of some estimates on the volume of spherical caps (given in Section~\ref{subsec:spherical_caps}); some isoperimetric inequalities (Section~\ref{subsec:isoperimetric_inequalities}), and some reverse isoperimetric inequalities (Section~\ref{subsec:reverse_isoperimetric_inequalities}). Finally, in Section~\ref{subsec:density_cluster_hard_case}, for the hard case distribution $\Pbb$ defined in~\eqref{eqn:lb_density} and $\mc{L}$ defined in~\eqref{eqn:lower_set}, we establish bounds on the population-level normalized cut $\Phi_{\Pbb,r}(\mc{L})$ and local spread $s_{\Pbb,r}(\mc{X})$; these will be useful in the proof of Theorem~\ref{thm:ppr_lb}. Throughout, we write $\nu_d := \nu(B(0,1))$ for the Lebesgue measure of a $d$-dimensional unit ball.

\subsection{Balls, Spherical Caps, and Associated Estimates}
\label{subsec:spherical_caps}
In this section, we derive lower bounds on the volume of the intersection between two balls in $\Rd$, and the volume of a spherical cap. Results of this type are well-known, but since we could not find exactly the statements we desire, for completeness we also supply proofs. We use the notation $B(x,r)$ for a ball of radius $r$ centered at $x \in \Rd$, and $\mathrm{cap}_{r}(h)$ for a spherical cap of height $h$ and radius $r$. Recall that the Lebesgue measure of a spherical cap is
\begin{equation*}
\nu\bigl(\mathrm{cap}_r(h)\bigr) = \frac{1}{2} \nu_d r^d I_{1 - a}\left(\frac{d + 1}{2}; \frac{1}{2}\right),
\end{equation*}
where $a = (r - h)^2/r^2$, and
\begin{equation*}
I_{1 - a}(z,w) = \frac{\Gamma(z + w)}{\Gamma(z) \Gamma(w)} \int_{0}^{1 - a} u^{z - 1} (1 - u)^{w - 1} du,
\end{equation*}
is the cumulative distribution function of a $\mathrm{Beta}(z,w)$ distribution, evaluated at $1 - a$. (Here $\Gamma(\cdot)$ is the gamma function).
\begin{lemma}
	\label{lem:overlap_balls}
	For any $x,y \in \Rd$ and $r > 0$, it holds that
	\begin{equation}
	\label{eqn:overlap_balls_1}
	\nu\bigl(B(x,r) \cap B(y,r)\bigr) \geq \nu_d r^d\biggl(1 - \frac{\|x - y\|}{r} \sqrt{\frac{d + 2}{2\pi}}\biggr).
	\end{equation}
	For any $x,y \in \Rd$ and $r,\sigma > 0$ such that $\|x - y\| \leq \sigma$, it holds that,
	\begin{equation}
	\label{eqn:overlap_balls_2}
	\nu\bigl(B(x,r) \cap B(y,\sigma)\bigr) \geq \frac{1}{2} \nu_d r^d\biggl(1 - \frac{r}{\sigma}\sqrt{\frac{d + 2}{2\pi}}\biggr).
	\end{equation}
\end{lemma}
\begin{lemma}
	\label{lem:volume_of_spherical_cap}
	For any $0 < h \leq r$, and $a = 1 - (2 r h - h^2)/r^2$,
	\begin{equation*}
	\nu\bigl(\mathrm{cap}_r(h)\bigr) \geq \frac{1}{2}\nu_dr^d\biggl(1 - 2\sqrt{a} \cdot \sqrt{\frac{d + 2}{2\pi}}\biggr).
	\end{equation*}
\end{lemma}
An immediate implication of~\eqref{eqn:overlap_balls_2} is that for any $x \in \mc{C}_{\lambda,\sigma}$,
\begin{equation}
\label{eqn:uniform_local_conductance}
\nu\bigl(B(x,r) \cap \mc{C}_{\lambda,\sigma}\bigr) \geq \frac{1}{2} \nu_d r^d\biggl(1 - \frac{r}{\sigma}\sqrt{\frac{d + 2}{2\pi}}\biggr).
\end{equation}  
\emph{Proof of Lemma~\ref{lem:overlap_balls}.}
First, we prove~\eqref{eqn:overlap_balls_1}. The intersection $B(x,r) \cap B(y,r)$ consists of two symmetric spherical caps, each of height $h = r - \frac{\|x - y\|}{2}$. 
As a result, by Lemma~\ref{lem:volume_of_spherical_cap} we have
\begin{equation*}
\nu\bigl(B(x,r) \cap B(y,r)\bigr) \geq \nu_d r^d \bigl(1 - 2\sqrt{a} \cdot \sqrt{\frac{d + 2}{2\pi}}\bigr)
\end{equation*}
where $a = \|x - y\|^2/(4r^2)$, and the claim follows.

Next we prove~\eqref{eqn:overlap_balls_2}. Assume that $\|x - y\| = \sigma$, as otherwise if $0 \leq \|x - y\| < \sigma$ the volume of the overlap will only be larger. Then $B(x,r) \cap B(y,\sigma)$ contains a spherical cap of radius $r$ and height $h = r - \frac{r^2}{2\sigma}$, from Lemma~\ref{lem:volume_of_spherical_cap} we deduce
\begin{equation*}
\nu\bigl(B(x,r) \cap B(y,\sigma)\bigr) \geq \frac{1}{2}\nu_dr^d\biggl(1 - 2\sqrt{a}\cdot\sqrt{\frac{d + 2}{2\pi}}\biggr)
\end{equation*}
for $a = (r - h)^2/r^2 = r^2/(4\sigma^2)$, and the claim follows. \qed \\

\noindent\emph{Proof of Lemma~\ref{lem:volume_of_spherical_cap}.}
For any $0 \leq a \leq 1$, we have that
\begin{equation*}
\int_{0}^{1 - a}u^{(d-1)/2}(1 - u)^{-1/2}du = \int_{0}^{1}u^{(d-1)/2}(1 - u)^{-1/2}du - \int_{1 - a}^{1}u^{(d-1)/2}(1 - u)^{-1/2}du. 
\end{equation*}
The first integral is simply
\begin{equation*}
\int_{0}^{1}u^{(d-1)/2}(1 - u)^{-1/2}du = \frac{\Gamma\bigl(\frac{d + 1}{2}\bigr)\Gamma\bigl(\frac{1}{2}\bigr)}{ \Gamma\bigl(\frac{d}{2}+ 1\bigr)},
\end{equation*}
whereas for all $u \in [0,1]$ and $d \geq 1$, the second integral can be upper bounded as follows:
\begin{equation*}
\int_{1 - a}^{1}u^{(d-1)/2}(1 - u)^{-1/2}du \leq \int_{1 - a}^{1}(1 - u)^{-1/2}du = \int_{0}^{a} u^{-1/2}du = 2\sqrt{a}.
\end{equation*}
As a result, 
\begin{equation*}
\nu\bigl(\mathrm{cap}_r(h)\bigr) \geq \frac{1}{2}\nu_dr^d \biggl(1 - 2\sqrt{a}\frac{\Gamma(\frac{d}{2} + 1)}{\Gamma(\frac{d + 1}{2})\Gamma(\frac{1}{2})}\biggr)  \geq \frac{1}{2}\nu_dr^d \biggl(1 - 2\sqrt{a} \cdot \sqrt{\frac{d + 2}{2\pi}}\biggr).
\end{equation*} \qed

\subsection{Isoperimetric Inequalities}
\label{subsec:isoperimetric_inequalities}
\citet{dyer1991b} establish the following isoperimetric inequality for convex sets.
\begin{lemma}[Isoperimetry of a convex set.]
	\label{lem:convex_isoperimetric_inequality}
	For any partition $(\mc{R}_1,\mc{R}_2,\mc{R}_3)$ of a convex set $\mc{K} \subseteq \Rd$, it holds that
	\begin{equation*}
	\nu(\mc{R}_3) \geq 2\frac{\mathrm{dist}(\mc{R}_1, \mc{R}_2)}{\mathrm{diam}(\mc{K})} \min(\nu(\mc{R}_1), \nu(\mc{R}_2)).
	\end{equation*}
\end{lemma}
\citet{abbasi-yadkori2016a} points out that if $\mc{S}$ is the image of a convex set under a Lipschitz measure-preserving mapping $g: \Rd \to \Rd$, a similar inequality can be obtained.
\begin{corollary}[Isoperimetry of Lipschitz embeddings of convex sets.]
	\label{cor:nonconvex_isoperimetric_inequality}
	Suppose $\mc{S}$ is the image of a convex set $\mathcal{K}$ under a mapping $g:\Rd \to \Rd$ such that
	\begin{equation*}
	\|g(x) - g(y)\| \leq M \cdot \|x - y\|,~~\textrm{for all $x,y \in \mc{K}$, and}~~\det(\nabla g(x)) = 1~~\textrm{for all $x \in \mc{K}$.}
	\end{equation*}
	Then for any partition $(\Omega_1,\Omega_2,\Omega_3)$ of $\mc{S}$, 
	\begin{equation*}
	\nu(\Omega_3) \geq 2\frac{\mathrm{dist}(\Omega_1, \Omega_2)}{\diam(\mc{K}) M} \min(\nu(\Omega_1), \nu(\Omega_2)).
	\end{equation*}
\end{corollary}

\subsection{Reverse Isoperimetric Inequalities}
\label{subsec:reverse_isoperimetric_inequalities}

For any set $\mc{C} \subseteq \Rd$ and $\sigma > 0$, let $\mc{C}_{\sigma} := \{x: \mathrm{dist}(x,\mc{C}) \leq \sigma\}$. We begin with an upper bound on the volume of $\mc{C}_{\sigma + \delta}$ as compared to $\mc{C}_{\sigma}$. 
\begin{lemma}
	\label{lem:reverse_isoperimetric_inequality}
	For any bounded set $\mc{C} \subseteq \Rd$ and  $\sigma, \delta > 0$, it holds that
	\begin{equation}
	\label{eqn:reverse_isoperimetric_inequality}
	\nu(\mc{C}_{\sigma + \delta}) \leq \nu(\mc{C}_{\sigma}) \cdot \Bigl(1 + \frac{\delta}{\sigma}\Bigr)^d.
	\end{equation}
\end{lemma}
Lemma~\ref{lem:reverse_isoperimetric_inequality} is a reverse isoperimetric inequality. To see this, note that if $\delta \leq \sigma/d$ then $\bigl(1 + \delta/\sigma\bigr)^d \leq 1 + d \cdot \delta/(\sigma - d\delta)$, and we deduce from~\eqref{eqn:reverse_isoperimetric_inequality} that
\begin{equation}
\label{eqn:reverse_isoperimetric_inequality_2}
\nu(\mc{C}_{\sigma + \delta} \setminus \mc{C}_{\sigma}) = \nu(\mc{C}_{\sigma + \delta}) - \nu(\mc{C}_{\delta}) \leq \frac{d\delta}{\sigma - d\delta} \cdot \nu(\mc{C}_{\sigma}).
\end{equation}
We use~\eqref{eqn:reverse_isoperimetric_inequality_2} along with Assumption~\ref{asmp:low_noise_density} to derive a density-weighted reverse isoperimetric inequality. 
\begin{lemma}
	\label{lem:reverse_isoperimetric_inequality_density_weighted}
	Let $\mc{C}_{\lambda,\sigma}$ satisfy Assumption~\ref{asmp:bounded_density} and~\ref{asmp:low_noise_density} for some $\theta, \gamma$ and $\lambda_{\sigma}$. Then for any $0 < r \leq \sigma/d$, it holds that
	\begin{equation}
	\label{eqn:reverse_isoperimetric_inequality_density_weighted}
	\Pbb\bigl(\mc{C}_{\lambda,\sigma + r} \setminus \mc{C}_{\lambda,\sigma}\bigr) \leq \Bigl(1 + \frac{dr}{\sigma - dr}\Bigr) \cdot \frac{dr}{\sigma} \cdot \left(\lambda_{\sigma} - \theta\frac{r^{\gamma}}{\gamma + 1}\right) \cdot \nu(\mc{C}_{\lambda,\sigma}).
	\end{equation}
\end{lemma}
\emph{Proof of Lemma~\ref{lem:reverse_isoperimetric_inequality}.}
Fix $\delta' > 0$, and take $\epsilon = \delta + \delta'$. We will show that
\begin{equation}
\label{pf:reverse_isoperimetric_inequality}
\nu(\mc{C}_{\sigma + \epsilon}) \leq \nu(\mc{C}_{\sigma}) \cdot \Bigl(1 + \frac{\epsilon}{\sigma}\Bigr)^d,
\end{equation}
whence taking a limit as $\delta' \to 0$ yields the claim.

To show~\eqref{pf:reverse_isoperimetric_inequality}, we need to construct a particular disjoint covering $\mc{A}_1(\sigma + \epsilon), \ldots, \mc{A}_N(\sigma + \epsilon)$ of $\mc{C}_{\sigma + \delta}$. To do so, we first take a finite set of points $x_1,\ldots,x_N$ such that the net $B(x_1,\sigma + \epsilon),\ldots,B(x_N,\sigma + \epsilon)$ covers $\mc{C}_{\sigma + \delta}$.  Note that such a covering exists for some finite $N = N(\epsilon)$ because $\mc{C}_{\sigma + \delta}$ is bounded, and the closure of $\mc{C}_{\sigma + \delta}$ is thus a compact subset of $\cup_{x \in \mc{C}} B(x,\sigma + \epsilon)$. Defining $\mc{A}_1(s), \ldots, \mc{A}_N(s)$ for a given $s > 0$ to be
\begin{equation*}
\mc{A}_1(s) := B(x_1,s),~~\textrm{and}~~ \mc{A}_{j + 1}(s) := B(x_{j + 1},s) \setminus \bigcup_{i = 1}^{j} B(x_i,s) ~~\textrm{for $j = 1,\ldots,N - 1$},
\end{equation*}
we have that $\mc{A}_1(\sigma + \epsilon),\ldots,\mc{A}_N(\sigma + \epsilon)$ is a disjoint covering of $\mc{C}_{\sigma + \delta}$, and so $\nu(\mc{C}_{\sigma + \delta}) \leq \sum_{j = 1}^{N} \nu\bigl(\mc{A}_j(\sigma + \epsilon)\bigr)$. 

We claim that for all $j = 1,\ldots,N$, the function $s \mapsto \nu\bigl(\mc{A}_j(s)\bigr)/\nu\bigl(B(x_j,s)\bigr)$ is monotonically non-increasing in $s$. Once this claim is verified, it follows that
\begin{align*}
\nu(\mc{A}_j(\sigma + \epsilon)) & = \nu\bigl(B(x_j,\sigma + \epsilon)\bigr) \cdot \frac{\nu\bigl(\mc{A}_j(\sigma + \epsilon)\bigr)}{\nu\bigl(B(x_j,\sigma + \epsilon)\bigr)} \\
& \leq \Bigl(1 + \frac{\epsilon}{\sigma}\Bigr)^d \cdot  \nu\bigl(B(x_j,\sigma)\bigr) \cdot \frac{\nu\bigl(\mc{A}_j(\sigma)\bigr)}{\nu\bigl(B(x_j,\sigma)\bigr)} \\
& = \Bigl(1 + \frac{\epsilon}{\sigma}\Bigr)^d \cdot \nu\bigl(\mc{A}_j(\sigma),
\end{align*}
and summing over $j$, we see that
\begin{equation*}
\nu(\mc{C}_{\sigma + \delta}) \leq \sum_{j = 1}^{N} \nu\bigl(\mc{A}_j(\sigma + \epsilon)\bigr) \leq \Bigl(1 + \frac{\epsilon}{\sigma}\Bigr)^d \cdot \sum_{j = 1}^{N} \nu\bigl(\mc{A}_j(\sigma)\bigr) \leq \Bigl(1 + \frac{\epsilon}{\sigma}\Bigr)^d \nu(\mc{C}_{\sigma}).
\end{equation*}
The last inequality follows since $\mc{A}_1(\sigma),\ldots,\mc{A}_N(\sigma)$ are disjoint subsets of the closure of $\mc{C}_{\sigma}$.

It remains to verify that $s \mapsto \nu\bigl(\mc{A}_j(s)\bigr)/\nu\bigl(B(x_j,s)\bigr)$ is monotonically non-increasing. For any $0 < s < t$ and $j = 1,\ldots,N$, suppose $x \in \mc{A}_j(T) - x_j$, meaning $x \in B(0,t)$ and $x \not\in B(x_i - x_j,t)$ for any $i = 1,\ldots,j - 1$. Thus $(s/t)x \in B(0,s)$, and
\begin{equation*}
\|(s/t)x - (x_i - x_j)\| \geq \|x - (x_i - x_j)\| - \|x - (s/t)x\| > t - (1 - s/t)\|x\| \geq t - (1 - s/t)t = s,
\end{equation*}
or in other words $(s/t)x \not\in B(x_i - x_j,s)$ for any $i = 1,\ldots,j - 1$. Consequently,
\begin{equation*}
\Bigl(\mc{A}_j(t) - x_j\Bigr) \subset \frac{t}{s} \cdot \Bigl(\mc{A}_{j}(s)- x_j\Bigr),
\end{equation*}
and applying $\nu(\cdot)$ to both sides yields the claim. \qed \\

\noindent \emph{Proof of Lemma~\ref{lem:reverse_isoperimetric_inequality_density_weighted}.}
Fix $k \in \mathbb{N}$. To establish~\eqref{eqn:reverse_isoperimetric_inequality_density_weighted}, we partition $\mc{C}_{\lambda,\sigma + r} \setminus \mc{C}_{\lambda,\sigma}$ into thin tubes $\mc{T}_{1},\ldots,\mc{T}_{k}$, with the $j$th tube $\mc{T}_j$ defined as $\mc{T}_j := \mc{C}_{\lambda,\sigma + jr/k} \setminus \mc{C}_{\lambda,\sigma + (j - 1)r/k}.$ We upper bound the Lebesgue measure of each tube $\mc{T}_j$ using~\eqref{eqn:reverse_isoperimetric_inequality_2}:\footnote{Note that $\mc{C}$ must be bounded, since the density $f(x) \geq \lambda_{\sigma} > 0$ for all $x \in \mc{C}$.}
\begin{equation*}
\nu(\mc{T}_j) \leq \frac{dr/k}{\sigma - dr/k} \nu(\mc{C}_{\lambda,\sigma + (j - 1)r/k}) \leq \frac{dr/k}{\sigma - dr/k} \nu(\mc{C}_{\lambda,\sigma + r}) \leq \Bigl(1 + \frac{dr}{\sigma - dr}\Bigr) \cdot \frac{dr/k}{\sigma - dr/k}  \cdot \nu(\mc{C}_{\lambda,\sigma}),
\end{equation*}
and the maximum density within each tube using~\ref{asmp:low_noise_density}:
\begin{equation*}
\max_{x \in \mc{T}_j} f(x) \leq \lambda_{\sigma} - \theta \Bigl(\frac{j - 1}{k}r\Bigr)^{\gamma};
\end{equation*}
combining these upper bounds, we see that
\begin{equation}
\label{pf:reverse_isoperimetric_inequality_density_weighted_1}
\Pbb\bigl(\mc{C}_{\lambda,\sigma + r} \setminus \mc{C}_{\lambda,\sigma}\bigr) = \sum_{j = 1}^{k} \Pbb(\mc{T}_j) \leq \Bigl(1 + \frac{dr}{\sigma - dr}\Bigr) \cdot  \frac{dr/k}{\sigma - dr/k} \cdot \nu(\mc{C}_{\lambda,\sigma}) \cdot \biggl(\sum_{j = 0}^{k - 1}\lambda_{\sigma} -  \theta r^{\gamma} \Bigl(\frac{j}{k}\Bigr)^{\gamma} \biggr).
\end{equation}
Treating the sum in the previous expression as a Riemann sum of a non-increasing function evaluated at $0,\ldots,k -1$ gives the upper bound
\begin{equation*}
\sum_{j = 0}^{k - 1}\lambda_{\sigma} -  \theta r^{\gamma} \Bigl(\frac{j}{k}\Bigr)^{\gamma } \leq \lambda_{\sigma} + \int_{0}^{k - 1} \Bigl(\lambda_{\sigma} -  \theta r^{\gamma} \Bigl(\frac{x}{k}\Bigr)^{\gamma}\Bigr) \,dx \leq k\lambda_{\sigma} + (k - 1)\frac{\theta r^{\gamma}}{\gamma + 1} \Bigl(\frac{k - 1}{k}\Bigr)^{\gamma},
\end{equation*}
and plugging back in to~\eqref{pf:reverse_isoperimetric_inequality_density_weighted_1}, we obtain
\begin{equation*}
\Pbb\bigl(\mc{C}_{\lambda,\sigma + r} \setminus \mc{C}_{\lambda,\sigma}\bigr) \leq  \Bigl(1 + \frac{dr}{\sigma - dr}\Bigr) \cdot \frac{dr}{\sigma - dr/k} \nu(\mc{C}_{\lambda,\sigma}) \cdot \biggl(\lambda - \frac{\theta r^{\gamma}}{\gamma + 1} \cdot \Bigl(\frac{k - 1}{k}\Bigr)^{\gamma + 1} \biggl).
\end{equation*}
The above inequality holds for any $k \in \mathbb{N}$, and taking the limit of the right hand side as $k \to \infty$ yields the claim. \qed

\subsection{Proof of Lemma~\ref{lem:density_cluster_local_spread}}
\label{subsec:density_cluster_local_spread}
The population-level local spread of $\mc{C}_{\lambda,\sigma}$ is
\begin{equation*}
s_{\Pbb,r}(\mc{C}_{\lambda,\sigma}) = \frac{\bigl(d_{\min}(\wt{\Pbb})\bigr)^2}{\vol_{\wt{\Pbb},r}(\mc{C}_{\lambda},\sigma)},
\end{equation*}
where we recall that $\wt{\Pbb}(\mc{S}) = \frac{\Pbb(\mc{S} \cap \mc{C}_{\lambda,\sigma})}{\Pbb(\mc{C}_{\lambda,\sigma})}$ for Borel sets $\mc{S}$, and $d_{\min}(\wt{\Pbb}) := \min_{x \in \mc{C}_{\lambda,\sigma}}\{\deg_{\wt{\Pbb},r}(x)\}^2$. To lower bound $s_{\Pbb,r}(\mc{C}_{\lambda,\sigma})$, we first lower bound $d_{\min}(\wt{\Pbb})$, and then upper bound $\vol_{\wt{\Pbb},r}(\mc{C}_{\lambda},\sigma)$. Using the lower bound $f(x) \geq \lambda_{\sigma}$ for all $x \in \mc{C}_{\lambda,\sigma}$ stipulated in~\ref{asmp:lambda_bounded_density}, we deduce that
\begin{align*}
d_{\min}(\wt{\Pbb}) & = \min_{x \in \mc{C}_{\lambda,\sigma}} \Bigl\{ \int \1\{\|x - y\| \leq r\} \,d\wt{\Pbb}(y) \Bigr\} \\
& \geq \frac{\lambda_{\sigma}}{\Pbb(\mc{C}_{\lambda,\sigma})} \cdot \min_{x \in \mc{C}_{\lambda,\sigma}} \Bigl\{ \int_{\mc{C}_{\lambda,\sigma}} \1\{\|x - y\| \leq r\} \,dy \Bigr\} \\
& \geq \frac{\lambda_{\sigma}}{\Pbb(\mc{C}_{\lambda,\sigma})} \cdot \frac{1}{2}\nu_dr^d \cdot \Bigl(1 - \frac{r}{\sigma} \sqrt{\frac{d + 2}{2\pi}}\Bigr),
\end{align*}
where the final inequality follows from Lemma~\ref{lem:overlap_balls}. 

On the other hand, using the upper bound $f(x) \leq \Lambda_{\sigma}$ for all $x \in \mc{C}_{\lambda,\sigma}$, we deduce that
\begin{align*}
\vol_{\wt{\Pbb},r}(\mc{C}_{\lambda,\sigma}) & = \iint \1\{\|x - y\|\leq r\} \,d\wt{\Pbb}(y) \,d\wt{\Pbb}(x) \\
& \leq \frac{\Lambda_{\sigma}^2}{\Pbb(\mc{C}_{\lambda,\sigma})^2} \cdot \int_{\mc{C}_{\lambda,\sigma}} \int_{\mc{C}_{\lambda,\sigma}} \1\{\|x - y\|\leq r\} \,dy \,dx \\
& \leq \frac{\Lambda_{\sigma}^2}{\Pbb(\mc{C}_{\lambda,\sigma})^2} \cdot \nu_d r^d \cdot \nu(\mc{C}_{\lambda,\sigma}) \\
& \leq \frac{\Lambda_{\sigma}^2}{\Pbb(\mc{C}_{\lambda,\sigma})^2} \cdot \nu_d^2 r^d  \cdot \Bigl(\frac{\rho}{2}\Bigr)^d;
\end{align*}
the final inequality follows from~\ref{asmp:embedding}, which along with the isodiametric inequality for convex sets~\citep{gruber2007} implies $\nu(\mc{C}_{\lambda,\sigma}) = \nu(\mc{K}) \leq \nu_d(\rho/2)^d$. The claim of Lemma~\ref{lem:density_cluster_local_spread} follows. \qed

\subsection{Proof of Proposition~\ref{prop:density_cluster_normalized_cut}}
\label{subsec:density_cluster_ncut}
By Assumption~\ref{asmp:bounded_volume}, we have that $\Phi_{\Pbb,r}(\mc{C}_{\lambda,\sigma}) = \cut_{\Pbb,r}(\mc{C}_{\lambda,\sigma})/\vol_{\Pbb,r}(\mc{C}_{\lambda,\sigma})$, and to prove Proposition~\ref{prop:density_cluster_normalized_cut} we must therefore upper bound $\cut_{\Pbb,r}(\mc{C}_{\lambda,\sigma})$ and lower bound $\vol_{\Pbb,r}(\mc{C}_{\lambda,\sigma})$. 

Let $\mc{C}_{\lambda,\sigma + r} = \{x: \mathrm{dist}(x,\mc{C}_{\lambda}) \leq \sigma + r\}$. We upper bound $\cut_{\Pbb,r}(\mc{C}_{\lambda,\sigma})$ in terms of the probability mass of $\mc{C}_{\lambda,\sigma + r} \setminus \mc{C}_{\lambda,\sigma}$:
\begin{align*}
\cut_{\Pbb,r}(\mc{C}_{\lambda,\sigma}) & = \iint \1\{\|x - y\| \leq r \} \cdot \1\{x \in \mc{C}_{\lambda,\sigma}\} \cdot \1\{y \not\in \mc{C}_{\lambda,\sigma} \} \,d\Pbb(y) \,d\Pbb(x) \\
& \leq  \iint \1\{\|x - y\| \leq r \} \cdot \1\{x \in \mc{C}_{\lambda,\sigma}\} \cdot \1\{y \in \mc{C}_{\lambda,\sigma + r} \setminus \mc{C}_{\lambda,\sigma}  \} \,d\Pbb(y) \,d\Pbb(x) \\
& \leq \lambda \nu_d r^d \cdot \Pbb\bigl(\mc{C}_{\lambda,\sigma + r} \setminus \mc{C}_{\lambda,\sigma}\bigr).
\end{align*}
On the other hand, using the lower bound $f(x) \geq \lambda_{\sigma} $ for all $x \in \mc{C}_{\lambda,\sigma}$, we lower bound $\vol_{\Pbb,r}(\mc{C}_{\lambda,\sigma})$ in terms of the Lebesgue measure of $\mc{C}_{\lambda,\sigma}$:
\begin{align*}
\vol_{\Pbb,r}(\mc{C}_{\lambda,\sigma}) & = \iint \1\{\|x - y\| \leq r \} \cdot \1\{x \in \mc{C}_{\lambda,\sigma}\} \,d\Pbb(y) \,d\Pbb(x) \\
& \geq \lambda_{\sigma}^2 \cdot \iint \1\{\|x - y\| \leq r \} \cdot \1\{x,y \in \mc{C}_{\lambda,\sigma}\} \,dy \,dy \\
& \geq \lambda_{\sigma}^2 \cdot \frac{1}{2} \nu_d r^d \cdot \Bigl(1 - \frac{r}{\sigma} \sqrt{\frac{d + 2}{2\pi}}\Bigr) \cdot \nu(\mc{C}_{\lambda,\sigma}).
\end{align*}
The claim of Proposition~\ref{prop:density_cluster_normalized_cut} follows upon using Lemma~\ref{lem:reverse_isoperimetric_inequality_density_weighted} to upper bound $\Pbb\bigl(\mc{C}_{\lambda,\sigma + r} \setminus \mc{C}_{\lambda,\sigma}\bigr)$. \qed

\subsection{Proof of Proposition~\ref{prop:density_cluster_conductance}}
\label{subsec:density_cluster_conductance}
The following Lemma lower bounds the population-level uniform conductance $\Psi_{\nu,r}(\mc{C}_{\lambda,\sigma})$. We note that for convex sets, results of this type are well known (see e.g.~\cite{vempala2005} and references therein).  
\begin{lemma}
	\label{lem:uniform_density_cluster_conductance}
	Suppose $\mc{C}_{\lambda,\sigma}$ satisfies Assumption~\ref{asmp:embedding} with respect to some $\rho \in (0,\infty)$ and $M \in [1,\infty)$. For any $0 < r \leq \sigma \cdot \sqrt{2\pi/(d + 2)}$, it holds that
	\begin{equation}
	\label{eqn:uniform_density_cluster_conductance}
	\Psi_{\nu,r}(\mc{C}_{\lambda,\sigma}) \geq \Bigl(1 - \frac{r}{4\rho M}\Bigr) \cdot \Bigl(1 - \frac{r}{\sigma}\sqrt{\frac{d + 2}{2\pi}}\Bigr)^2 \cdot \frac{\sqrt{2\pi}}{36} \cdot \frac{r}{\rho M \sqrt{d + 2}}.
	\end{equation}
\end{lemma}
Noting that $\Psi_{\Pbb,r}(\mc{C}_{\lambda,\sigma}) \geq \Psi_{\nu,r}(\mc{C}_{\lambda,\sigma}) \cdot \lambda_{\sigma}^2/\Lambda_{\sigma}^2$, Proposition~\ref{prop:density_cluster_conductance} follows from~\eqref{eqn:uniform_density_cluster_conductance}. \\

\noindent \emph{Proof of Lemma~\ref{lem:uniform_density_cluster_conductance}.}
For ease of notation, throughout this proof we write $\wt{\nu}$ for the uniform probability measure over $\mc{C}_{\lambda,\sigma}$, put $\ell \nu_d r^d := \min_{x \in \mc{C}_{\lambda,\sigma}}\nu(B(x,r) \cap \mc{C}_{\sigma})$ and $a := r/(2\rho M)$. 

Let $\mc{S}$ be an arbitrary measurable subset of $\mc{C}_{\lambda,\sigma}$, and let $\mc{R} = \mc{C}_{\lambda,\sigma} \setminus \mc{S}$. For a given $\delta \in (0,1)$, let the $\delta$-interior of $\mc{S}$ be 
\begin{equation*}
\mc{S}^{\delta} := \{x \in \mc{S}: \nu\bigl(B(x,r) \cap \mc{R}\bigr) \leq \ell \delta \nu_d r^d\};
\end{equation*}
define $\mc{R}^{\delta}$ likewise, and let $\mc{B}^{\delta} = \mc{C}_{\lambda,\sigma} \setminus (\mc{S}^{\delta} \cup \mc{R}^{\delta})$ consist of the remaining boundary points.
As is standard (see for example~\cite{dyer1991b,lovasz1990}), the proof of Lemma~\ref{lem:uniform_density_cluster_conductance} uses several inequalities to lower bound the normalized cut $\Phi_{\wt{\nu},r}(\mc{S})$.
\begin{itemize}
	\item \textbf{Bounds on cut and volume.}
	We can lower bound $\cut_{\wt{\nu},r}(\mc{S})$ as follows:
	\begin{align*}
	\nu(\mc{C}_{\lambda,\sigma})^2 \cdot \cut_{\wt{\nu},r}(\mc{S}) & = \int_{\mc{S}} \int_{\mc{R}} \1(\|x - y\| \leq r) \,dy \,dx \\
	& = \frac{1}{2}\Bigl(\int_{\mc{S}} \int_{\mc{R}} \1(\|x - y\| \leq r) \,dy \,dx + \int_{\mc{R}} \int_{\mc{S}} \1(\|x - y\| \leq r) \,dy \,dx\Bigr) \\
	& \geq \frac{1}{2} \delta \ell \nu_d r^d \cdot \nu(\mc{B}^{\delta}).
	\end{align*}
	We can upper bound $\vol_{\wt{\nu},r}(\mc{S})$ as follows:
	\begin{equation*}
	\nu(\mc{C}_{\lambda,\sigma})^2 \cdot \vol_{\wt{\nu},r}(\mc{S}) = \int_{\mc{C}_{\lambda,\sigma}} \int_{\mc{S}} \1(\|x - y\| \leq r) \,dy \,dx \leq \nu_d r^d \nu(\mc{S})
	\end{equation*}
	and likewise for $\vol_{\wt{\nu},r}(\mc{R})$. Therefore,
	\begin{equation}
	\label{pf:uniform_density_cluster_conductance_1}
	\Phi_{\wt{\nu},r}(\mc{S}) \geq \frac{\delta \ell \cdot \nu(\mc{B}^{\delta})}{2 \cdot \min\{\nu(\mc{S}), \nu(\mc{R})\}}.
	\end{equation}
	\item \textbf{Isoperimetric inequality.}
	Applying Corollary~\ref{cor:nonconvex_isoperimetric_inequality}, we have that
	\begin{equation}
	\label{pf:uniform_density_cluster_conductance_2}
	\nu(\mc{B}^{\delta}) \geq \frac{2 \cdot \mathrm{dist}(\mc{S}^{\delta},\mc{R}^{\delta})}{\rho M} \cdot \min\bigl\{\nu(\mc{S}^{\delta}),\nu(\mc{R}^{\delta})\bigr\}.
	\end{equation}
	\item \textbf{Lebesgue measure of $\delta$-interiors.} Suppose $\nu(\mc{S}^{\delta}) \leq (1 - a) \cdot \nu(\mc{S})$ or $\nu(\mc{R}^{\delta}) \leq (1 - a) \cdot \nu(\mc{R})$. Then 
	$\nu(\mc{B}^{\delta}) \geq a \cdot \min\{\nu(\mc{S}), \nu(\mc{R})\}$, and combined with~\eqref{pf:uniform_density_cluster_conductance_1} we have that $\Phi_{\wt{\nu},r}(\mc{S}) \geq \delta a \ell/2$. Otherwise,
	\begin{equation}
	\label{pf:uniform_density_cluster_conductance_3}
	\min\bigl\{\nu(\mc{S}^{\delta}),\nu(\mc{R}^{\delta})\bigr\} \geq (1- a) \cdot \min\bigl\{\nu(\mc{S}),\nu(\mc{R})\bigr\}.
	\end{equation}
	\item \textbf{Distance between $\delta$-interiors.} For any $x \in \mc{S}^{\delta}$ and $y \in \mc{R}^{\delta}$, we have that 
	\begin{align*}
	\nu\bigl(B(x,r) \cap B(y,r)\bigr) & = \nu\bigl(B(x,r) \cap B(y,r) \cap \mc{R}\bigr) + \nu\bigl(B(x,r) \cap B(y,r) \cap \mc{S}\bigr) ~~+ \\
	& \quad~\nu\bigl(B(x,r) \cap B(y,r) \cap \mc{C}_{\lambda,\sigma}^c\bigr) \\
	& \leq \nu\bigl(B(x,r) \cap \mc{R}\bigr) + \nu\bigl(B(y,r) \cap \mc{S}\bigr) + \nu\bigl(B(x,r) \cap \mc{C}_{\lambda,\sigma}^c\bigr) \\
	& \leq \bigl(2\ell\delta + (1 - \ell)\bigr) \cdot \nu_dr^d.
	\end{align*}
	It follows from~\eqref{eqn:overlap_balls_1} that
	\begin{equation*}
	\|x - y\| \geq \frac{r}{\nu_d r^d} \cdot \Bigl(\nu_d r^d - \nu\bigl(B(x,r) \cap B(y,r)\bigr)\Bigr) \cdot \sqrt{\frac{2\pi}{d + 2}} \geq r \cdot \ell \cdot (1 - 2\delta) \cdot \sqrt{\frac{2\pi}{d + 2}},
	\end{equation*}
	and taking the infimum over all $x \in \mc{S}^{\delta}$ and $y \in \mc{R}^{\delta}$, we have
	\begin{equation}
	\label{pf:uniform_density_cluster_conductance_4}
	\mathrm{dist}(\mc{S}^{\delta},\mc{R}^{\delta}) \geq r \cdot \ell \cdot (1 - 2\delta) \cdot \sqrt{\frac{2\pi}{d + 2}}.
	\end{equation}
\end{itemize}	
Combining~\eqref{pf:uniform_density_cluster_conductance_1}-\eqref{pf:uniform_density_cluster_conductance_4} and taking $\delta = 1/3$ implies that
\begin{equation*}
\Phi_{\wt{\nu},r}(\mc{S}) \geq \min\biggl\{(1 - a) \cdot \frac{r}{\rho M} \cdot \frac{\ell^2}{9} \cdot \sqrt{\frac{2\pi}{d + 2}}, \frac{a\ell}{6}\biggr\}
\end{equation*}
and the claim follows from~\eqref{eqn:uniform_local_conductance}, which implies that $\ell \geq 1/2 \cdot (1 - r/\sigma)\sqrt{2\pi/(d + 2)}$. \qed

\subsection{Population Functionals, Hard Case}
\label{subsec:density_cluster_hard_case}
Let $\Pbb$ be the hard case distribution over rectangular domain $\mc{X}$, defined as in~\eqref{eqn:lb_density}, and $\mc{L}$ the lower half of $\mc{X}$. Suppose $r \in (0,\sigma/2)$. Then the population normalized cut $\Phi_{\Pbb,r}(\mathcal{L})$ is upper bounded,
\begin{equation}
\label{eqn:ncut_lower_set}
\Phi_{\Pbb,r}(\mc{L}) \leq \frac{8}{3} \cdot \frac{r}{\rho}.
\end{equation}
and the population local spread $s_{\Pbb,r}(\mc{X})$ is lower bounded,
\begin{equation}
\label{eqn:local_spread_ub}
s_{\Pbb,r}(\mc{X}) \geq \frac{\pi r^2 \epsilon^2}{2 \rho\sigma}
\end{equation}
\emph{Proof of~\eqref{eqn:ncut_lower_set}.}
Noting that $\vol_{\Pbb,r}(\mathcal{L}) = \vol_{\Pbb,r}(\mathcal{X} \setminus \mathcal{L})$, it suffices to upper bound $\cut_{\Pbb,r}(\mathcal{L})$ and lower bound $\vol_{\Pbb,r}(\mathcal{L})$. Note that for any $x = (x_1,x_2) \in \mathcal{L}$, if $x_2 \leq -r$ the ball $B(x,r)$ and the set $\mathcal{X}\setminus\mathcal{L}$ are disjoint. As a result,
\begin{equation*}
\cut_{\Pbb,r}(\mathcal{L}) \leq \Pbb\Bigl(\bigl\{x \in \mathcal{X}: x_2 \in (-r,0)\bigr\}\Bigr) \cdot d_{\max}(\Pbb) \leq \frac{r}{2 \rho} \cdot \frac{\pi r^2}{2 \sigma \rho}.
\end{equation*}
On the other hand, noting that $\deg_{\Pbb,r}(x) \geq \frac{\pi r^2}{2 \sigma \rho}$ for all $x \in \mc{C}^{(1)}$ such that $\mathrm{dist}(x, \partial \mc{C}^{(1)}) > r$, we have
\begin{align*}
\vol_{\Pbb,r}(\mathcal{L}) & \geq \Pbb\Bigl(\bigl\{x \in \mc{C}^{(1)} \cap \mathcal{L}: \mathrm{dist}(x, \partial \mc{C}^{(1)}) > r\bigr\}\Bigr) \cdot \frac{\pi r^2}{2 \sigma \rho} \\
& = \frac{(\sigma - 2r)(\rho - r)}{2 \sigma \rho} \cdot \frac{\pi r^2}{2 \sigma \rho}  \\
& \geq \frac{3}{16} \cdot \frac{\pi r^2}{2 \sigma \rho}
\end{align*}
where the last inequality follows since $r \leq \frac{1}{4}\sigma \leq \frac{1}{4}\rho$. \qed \\ 

\noindent \emph{Proof of~\eqref{eqn:local_spread_ub}.}
The statement follows since 
\begin{equation*}
d_{\min}(\Pbb) \geq \frac{\pi r^2}{2} \cdot \min_{x \in \mc{X}} f(x) = \frac{\pi r^2}{2} \cdot \frac{\epsilon}{\rho \sigma},
\end{equation*}
and
\begin{equation*}
\vol_{\Pbb,r}(\mc{X}) \leq d_{\max}(\Pbb) \leq \frac{\pi r^2}{2 \sigma \rho}.
\end{equation*} \qed

\section{Proof of Major Theorems}
\label{apdx:pf_major_theorems}
We now prove the three major theorems of our paper: Theorem~\ref{thm:volume_ssd_ub} (in Section~\ref{subsec:pf_volume_ssd_ub}), Theorem~\ref{thm:density_cluster_volume_ssd_ub}, and Theorem~\ref{thm:ppr_lb}.  Throughout, we use the notation $\wt{n} = |\mc{C}[X]|$ and $\wt{G}_{n,r} = G_{n,r}\bigl[\mc{C}[X]\bigr]$ as defined above.

\subsection{Proof of Theorem~\ref{thm:volume_ssd_ub}}
\label{subsec:pf_volume_ssd_ub}
We begin by recalling some probabilistic estimates needed for the proof of Theorem~\ref{thm:volume_ssd_ub}, along with the probability with which they hold.\\

\noindent\emph{Probabilistic estimates.}
Throughout the proof of Theorem~\ref{thm:volume_ssd_ub}, we will assume (i) that the inequalities~\eqref{eqn:sample_to_population_normalized_cut}-\eqref{eqn:sample_to_population_conductance} are satisfied; (ii) that the volume of $\mc{C}[X]$ is upper and lower bounded,
\begin{equation}
\label{pf:volume_ssd_ub_1}
(1 - \delta) \cdot \vol_{\Pbb,r}(\mc{C}) \leq \frac{1}{n(n-1)}\vol_{n,r}(\mc{C}[X]) \leq (1 + \delta) \cdot \vol_{\Pbb,r}(\mc{C});
\end{equation}
(iii) that the number of sample points in $\mc{C}$ is lower bounded,
\begin{equation}
\label{pf:volume_ssd_ub_2}
\wt{n} \geq (1 - \delta) \cdot n \cdot \mbb{P}(\mc{C}) \overset{\eqref{eqn:sample_to_population_local_spread_sample_complexity}}{\Longrightarrow} \wt{n} - 1 \geq (1 - \delta)^2 \cdot n \cdot \mbb{P}(\mc{C});
\end{equation}
and finally (iv) that the minimum and maximum degree of $\wt{G}_{n,r}$ are lower and upper bounded respectively,
\begin{equation}
\label{pf:volume_ssd_ub_3}
~~\frac{1}{\wt{n} - 1} d_{\min}(\wt{G}_{n,r}) \geq (1 - \delta) \cdot d_{\min}(\wt{\Pbb}),~~\textrm{and}~~\frac{1}{\wt{n} - 1} d_{\max}(\wt{G}_{n,r}) \leq (1 + \delta) \cdot d_{\max}(\wt{\Pbb}).
\end{equation}
By Propositions~\ref{prop:sample_to_population_1} and \ref{prop:sample_to_population_2}, and Lemmas~\ref{lem:hoeffding_2}-\ref{lem:bernstein_union}, these inequalities are satisfied with probability at least $1 - B_2/n - 4\exp\{-b_1\delta^2n\} - (2n + 2)\exp\{-b_2\delta^2n\} - (n + 1)\exp\{-b_3n\}$.\\


\noindent \emph{Proof of Theorem~\ref{thm:volume_ssd_ub}.}
We use Lemma~\ref{lem:zhu} to upper bound $\Delta(\wh{C},\mc{C}[X])$. In order to do so, we must verify that the tuning parameters $\alpha$ and $(L,U)$ satisfy the condition \eqref{eqn:zhu_condition} of this lemma, i.e. that $\alpha \leq 1/\bigl(2\tau_{\infty}(\wt{G}_{n,r})\bigr)$ and $U \leq 1/\bigl(5 \vol_{n,r}(\mc{C}[X])\bigr)$. In order to verify the upper bound on $\alpha$, we will use Proposition~\ref{prop:pointwise_mixing_time} to upper bound $\tau_{\infty}(\wt{G}_{n,r})$, which we may validly apply because
\begin{equation*}
\frac{d_{\max}(\wt{G}_{n,r})}{\bigl(d_{\min}(\wt{G}_{n,r})\bigr)^2} \leq \frac{(1 + \delta)}{(1 - \delta)^2} \cdot \frac{d_{\max}(\wt{\Pbb})}{(\wt{n} - 1) \cdot \bigl(d_{\min}(\wt{\Pbb})\bigr)^2} \leq \frac{(1 + \delta)}{(1 - \delta)^4} \cdot \frac{d_{\max}(\wt{\Pbb})}{n\Pbb(\mc{C})\bigl(d_{\min}(\wt{\Pbb})\bigr)^2} \leq \frac{1}{16}.
\end{equation*}
The last inequality in the above follows by taking $B_3 := 16 \cdot d_{\max}(\wt{\Pbb})/\bigl(d_{\min}(\wt{\Pbb})\bigr)^2$ in~\eqref{eqn:volume_ssd_ub_sample_complexity}.

Therefore by Proposition~\ref{prop:pointwise_mixing_time}, along with inequalities~\eqref{eqn:sample_to_population_local_spread} and~\eqref{eqn:sample_to_population_conductance} and the initialization conditions~\eqref{eqn:initialization} and \eqref{eqn:alpha_initialization}, we have that $\alpha \leq 1/45 \wedge 1/\bigl(2\tau_{\infty}(\wt{G}_{n,r})\bigr)$.  
On the other hand, by the upper bound on $\vol_{n,r}(\mc{C}[X])$ given in~\eqref{pf:volume_ssd_ub_1} and the initialization condition~\eqref{eqn:initialization}, we have that $U \leq 1/\bigl(5 \vol_{n,r}(\mc{C}[X])\bigr)$. In summary, we have confirmed that the condition~\eqref{eqn:zhu_condition} is satisfied.

Invoking Lemma~\ref{lem:zhu}, we conclude that there exists a set $\mc{C}[X]^g \subset \mc{C}[X]$ of volume at least $\vol_{n,r}(\mc{C}[X]^g) \geq \vol_{n,r}(\mc{C}[X])/2$, such that for any $\beta \in (L,U)$,
\begin{equation*}
\vol_{n,r}(S_{\beta,v} \vartriangle \mc{C}[X]) \leq 60 \cdot \frac{\Phi_{n,r}(\mc{C}[X])}{\alpha L} \leq 60 \frac{(1 + 2\delta)}{(1 - 4\delta)^2} \cdot \frac{\Phi_{n,r}(\mc{C}[X])}{\alpha_{\Pbb,r}(\mc{C},\delta)}\cdot n(n-1)\vol_{\Pbb,r}(\mc{C})
\end{equation*}
Noting that $\wh{C} = S_{\beta,v}$ for some $\beta \in (L,U)$, the claimed upper bound~\eqref{eqn:volume_ssd_ub} on $\Delta(\wh{C},\mc{C}[X])$ then follows from the upper bound~\eqref{eqn:sample_to_population_normalized_cut} on $\Phi_{n,r}(\mc{C}[X])$ and the upper bound on $\vol_{\Pbb,r}(\mc{C})$ in~\eqref{pf:volume_ssd_ub_1}. \qed

\subsection{Proof of Theorem~\ref{thm:density_cluster_volume_ssd_ub}}
From Theorem~\ref{thm:volume_ssd_ub}, we have that with probability $1 - B_2/n - 4\exp\{-b_1\delta^2n\} - (2n + 2)\exp\{-b_2\delta^2n\} - (n + 1)\exp\{-b_3n\}$, there exists a set $\mc{C}_{\lambda,\sigma}[X]^g \subset \mc{C}_{\lambda,\sigma}[X]$ of volume at least $\vol_{n,r}(\mc{C}_{\lambda,\sigma}[X]^g) \geq \vol_{n,r}(\mc{C}_{\lambda,\sigma}[X])/2$, such that
\begin{align*}
\frac{\Delta(\wh{C},\mc{C}_{\lambda,\sigma}[X])}{\vol_{n,r}(\mc{C}_{\lambda,\sigma}[X])} & \leq 60 \cdot \frac{(1 + 3\delta)(1 + 2\delta)}{(1 - 4\delta)^2(1 - 2\delta)} \cdot \frac{\Phi_{\Pbb,r}(\mc{C}_{\lambda,\sigma})}{\alpha_{\Pbb,r}(\mc{C}_{\lambda,\sigma},\delta)} \\
& \leq \frac{1020}{\ln(2)} \cdot \frac{(1 + 3\delta)(1 + 2\delta)}{(1 - 4\delta)^2(1 - 2\delta)} \cdot \frac{\Phi_{\Pbb,r}(\mc{C}_{\lambda,\sigma})}{\Psi_{\Pbb,r}(\mc{C}_{\lambda,\sigma})^2} \cdot \ln^2\Bigl(\frac{32}{(1 - 3\delta) s_{\Pbb,r}(\mc{C})}\Bigr)
\end{align*}
The claimed upper bound~\eqref{eqn:density_cluster_volume_ssd_ub} on $\Delta(\wh{C},\mc{C}_{\lambda,\sigma}[X])$ then follows from the bounds~\eqref{eqn:density_cluster_local_spread}-\eqref{eqn:density_cluster_conductance} on the population-level local spread, normalized cut, and conductance of $\mc{C}_{\lambda,\sigma}$, noting that the condition $r \leq \sigma/(4d)$ implies that $(1 - r/(4\rho L)) \geq 1 - 1/16$ and $1 - r/\sigma \cdot \sqrt{(d + 2)/(2\pi)} \geq 1 - 1/\sqrt{32}$, and taking
\begin{align*}
C_{1,\delta} & := \frac{1175040}{\pi \cdot \ln(2) \cdot (1 - 1/16)^2 \cdot (1 - 1/\sqrt{32})^4} \cdot \frac{(1 + 3\delta)(1 + 2\delta)}{(1 - 4\delta)^2(1 - 2\delta)}\\
C_{2,\delta} & := \frac{144}{(1 - 1/\sqrt{32})} \cdot \frac{1}{1-3\delta}.
\end{align*} \qed

\subsection{Proof of Theorem~\ref{thm:ppr_lb}}
We start by defining some constants, to make our proof statements easier to digest. Put
\begin{align*}
C_{3,\delta} & := \frac{288(1 + \delta)}{(1 - \delta)} \sqrt{8/3 + 8 \delta}
,&&C_{4,\delta} := \frac{72}{(1 - 3\delta)\pi}, \\
B_{1,\delta} & := 768 \cdot (1 + 3\delta) \cdot \ln\Bigl(C_{4,\delta} \frac{\rho\sigma}{r^2 \epsilon^2} \Bigr),&&B_{2,\delta} := \frac{(1 + \delta)^2}{(1 - \delta)^2} \cdot \frac{\rho \sigma}{r^2 \epsilon^2} \\
B_4 & := 1 + \frac{48\sigma\rho}{r^2} + \frac{4\rho}{r},&&b_4 :=  b_8 \wedge \cut_{\Pbb,r}(\mc{L}) \wedge d_{\min}(\Pbb)/14 \wedge \vol_{\Pbb,r}(\mc{L} \cap \mc{C}^{(1)}), \\
b_8 &:= \vol_{\Pbb,r}(\mc{X})/4 \wedge \frac{\epsilon r^2}{4\rho\sigma} \wedge \frac{\pi r^3}{8\sigma \rho^2}.
\end{align*}
To prove Theorem~\ref{thm:ppr_lb}, we use Theorem~\ref{thm:normalized_cut_ppr}, Proposition~\ref{prop:sample_to_population_1} and~\eqref{eqn:ncut_lower_set} to show that the cluster estimate $\wh{C}$ must have a small normalized cut. On the other hand, in Lemma~\ref{lem:normalized_cut_lb} we establish that any set $Z \subseteq X$ which is close to $\mc{C}^{(1)}[X]$---meaning $\vol_{n,r}(Z \vartriangle \mc{C}^{(1)}[X])$ is small---has a large normalized cut.
\begin{lemma}
	\label{lem:normalized_cut_lb}
	Fix $\delta \in (0,1)$. With probability at least $1 - B_4\exp\{-n\delta^2b_8\}$, the following statement holds:
	\begin{equation}
	\label{eqn:normalized_cut_lb}
	\Phi_{n,r}(Z) \geq \frac{(1 - \delta)^2}{4(1 + \delta)\pi} \left(1 - 2 \frac{\sigma \rho}{(1 - \delta) r^2 n^2} \vol_{n,r}(Z \vartriangle \mc{C}^{(1)}[X]) \right) \frac{\epsilon^2 r}{\sigma},~~\textrm{for all $Z \subseteq X$}.
	\end{equation}
\end{lemma}
We therefore conclude that $\vol_{n,r}(\wh{C} \vartriangle \mc{C}^{(1)}[X])$ must be large. In the remainder of this, we detail the probabilistic estimates used in the proof of Theorem~\ref{thm:ppr_lb}, and then give a formal proof of Theorem~\ref{thm:ppr_lb} and then of Lemma~\ref{lem:normalized_cut_lb}.\\

\noindent \emph{Probabilistic estimates.}
In addition to~\eqref{eqn:normalized_cut_lb}, we will assume (i) that the graph normalized cut of $\mc{L}$ and local spread of $\mc{X}$ are respectively upper and lower bounded,
\begin{equation*}
\Phi_{n,r}(\mc{L}[X]) \leq (1 + 3\delta) \cdot \Phi_{\Pbb,r}(\mc{L}),~~\textrm{and}~~s_{n,r}(X) \geq (1 - 3\delta) \cdot s_{\Pbb,r}(\mc{X});
\end{equation*}
(ii) that the graph volume of $\mc{L}$ is upper and lower bounded,
\begin{equation*}
(1 - \delta) \vol_{\Pbb,r}(\mc{L}) \leq \frac{1}{n(n - 1)}\vol_{n,r}(\mc{L}[X]) \leq (1 + \delta) \vol_{\Pbb,r}(\mc{L});
\end{equation*}
(iii) that the graph volumes of $\mc{L} \cap \mc{C}^{(1)}$ and $\mc{C}^{(1)}$ are respectively lower and upper bounded,
\begin{align*}
(1 - \delta) \vol_{\Pbb,r}(\mc{L} \cap \mc{C}^{(1)}) & \leq \frac{1}{n(n - 1)}\vol_{n,r}({\mc{L}[X] \cap \mc{C}^{(1)}}[X]) \\
\frac{1}{n(n - 1)}\vol_{n,r}(\mc{C}^{(1)}[X]) & \leq (1 + \delta) \vol_{\Pbb,r}(\mc{C}^{(1)});
\end{align*}
(iv) that the graph volume of $\mc{X}$ is lower bounded,
\begin{equation*}
\frac{1}{n(n - 1)}\vol_{n,r}(X) \geq (1 - \delta) \vol_{\Pbb,r}(\mc{X});
\end{equation*}
and finally (v) that the maximum degree of $G_{n,r}$ is upper bounded,
\begin{equation*}
\frac{1}{n - 1}d_{\max}(G_{n,r}) \leq (1 + \delta) d_{\max}(\Pbb).
\end{equation*}
It follows from Lemma~\ref{lem:normalized_cut_lb}, Propositions~\ref{prop:sample_to_population_1} and \ref{prop:sample_to_population_2}, and Lemmas~\ref{lem:hoeffding_2} and~\ref{lem:bernstein_union} that these estimates are together satisfied with probability at least $1 - B_4\exp\{-n\delta^2b_8\} - 3\exp\{-n\delta^2 \cut_{\Pbb,r}(\mc{L})\} - (2n + 2)\exp\{-n\delta^2 \cdot d_{\min}(\Pbb)/14\} - 5\exp\{-n\delta^2 \vol_{\Pbb,r}(\mc{L} \cap \mc{C}^{(1)})\} \geq 1 - (B_4 + 2n + 10)\exp\{-n\delta^2 b_4\}$. \\

\noindent \emph{Proof of Theorem~\ref{thm:ppr_lb}.}
As mentioned, we would like to use Theorem~\ref{thm:normalized_cut_ppr} to upper bound $\Phi_{n,r}(\wh{C})$, and so we first verify that the conditions of Theorem~\ref{thm:normalized_cut_ppr} are met. In particular, we have each of the following.
\begin{itemize}
	\item Recall that $n \geq 8 \cdot (1 + \delta)/(1 - \delta)$ \eqref{eqn:ppr_lb_condition} and that $\vol_{\Pbb,r}(\mc{L}) \geq 3/16 \cdot \pi r^2/(2\sigma\rho)$ (as shown in the proof of~\eqref{eqn:ncut_lower_set}). It is additionally clear that $d_{\max}(\Pbb) \leq \pi r^2 /(2\rho\sigma)$, and consequently,
	\begin{equation}
	\label{pf:ppr_lb_1}
	d_{\max}(G_{n,r}) \leq (n - 1) \cdot (1 + \delta) d_{\max}(\Pbb) \leq \frac{1}{3}n^2 (1 - \delta) \vol_{\Pbb,r}(\mc{L}) \leq \frac{1}{3}\vol_{n,r}(\mc{L}[X]).
	\end{equation}
	Therefore the lower bound in condition~\eqref{eqn:normalized_cut_ppr_vol} is satisfied.
	\item Note that $\delta \in (0,1/7)$ implies $(1 - \delta)/(1 + \delta) > 3/4$, and additionally that $\vol_{\Pbb,r}(\mc{L}) \leq \vol_{\Pbb,r}(\mc{X})/2$. It follows that
	\begin{equation}
	\begin{aligned}
	\label{pf:ppr_lb_2}
	\vol_{n,r}(X) & \geq n(n - 1)(1 - \delta) \vol_{\Pbb,r}(\mc{X}) \\
	& \geq 2n(n - 1) (1 - \delta) \vol_{\Pbb,r}(\mc{L}) \geq 2\frac{(1 - \delta)}{(1 + \delta)} \vol_{n,r}(\mc{L}[X]) \\
	& \geq \frac{3}{2} \vol_{n,r}(\mc{L}[X]).
	\end{aligned}
	\end{equation}
	Therefore the upper bound in condition~\eqref{eqn:normalized_cut_ppr_vol} is satisfied.
	\item By~\eqref{eqn:ncut_lower_set}, the normalized cut of $\mc{L}$ satisfies the following upper bound,
	\begin{equation}
	\Phi_{n,r}(\mc{L}[X]) \leq (1 + 3\delta) \cdot \Phi_{\Pbb,r}(\mc{L}) \leq (8/3 + 8\delta) \cdot  \frac{r}{\rho}, 
	\end{equation}
	and by~\eqref{eqn:local_spread_ub} the local spread of $\mc{X}$ satisfies the following lower bound,
	\begin{equation}
	s_{n,r}(X) \geq (1 - 3\delta) \cdot s_{\Pbb,r}(\mc{X}) \geq (1 - 3\delta) \cdot \frac{\pi r^2}{2\rho\sigma}.
	\end{equation}
	The constants $B_{1,\delta}$ and $B_{2,\delta}$ in assumption~\eqref{eqn:ppr_lb_condition} are chosen so that condition~\eqref{eqn:normalized_cut_ppr_ncut} is satisfied.
\end{itemize}
As a result, we may apply Theorem~\ref{thm:normalized_cut_ppr}, and deduce the following: there exists a set $\mc{L}[X]^g \subset \mc{L}$ of large volume, $\vol_{n,r}(\mc{L}[X]^g) \geq 5/6 \cdot \vol_{n,r}(\mc{L}[X])$, such that for any seed node $v \in \mc{L}[X]^g$, the normalized cut of the PPR cluster estimate $\Phi_{n,r}(\wh{C})$ satisfies the following upper bound:
\begin{equation*}
\Phi_{n,r}(\wh{C}) < 72\sqrt{\Phi_{n,r}(\mc{L}[X]) \cdot \ln\Bigl(\frac{36}{s_{n,r}(X)}\Bigr)} \leq  72\sqrt{(8/3 + 8\delta) \cdot \frac{r}{\rho} \cdot \ln\Bigl(\frac{72 \rho \sigma}{(1 - 3\delta)\pi r^2 \epsilon^2}\Bigr)}.
\end{equation*}	
Combined with Lemma~\ref{lem:normalized_cut_lb}, this implies
\begin{equation}
\begin{aligned}
\label{pf:ppr_lb_3}
&\frac{(1 - \delta)^2}{4(1 + \delta)\pi} \left(1 - 2 \frac{\sigma \rho}{(1 - \delta) r^2 n^2} \vol_{n,r}(\wh{C} \vartriangle \mc{C}^{(1)}[X]) \right) \frac{\epsilon^2 r}{\sigma} \\
& \leq 72\sqrt{(8/3 + 8\delta) \cdot \frac{r}{\rho} \cdot \ln\Bigl(\frac{72 \rho \sigma}{(1 - 3\delta)\pi r^2 \epsilon^2}\Bigr)},
\end{aligned}
\end{equation}
and solving for $\vol_{n,r}(\wh{C} \vartriangle \mc{C}^{(1)}[X])$ yields~\eqref{eqn:ppr_lb}. 

We conclude by observing that the set $\mc{L}[X]^g$ must have significant overlap with $\mc{C}^{(1)}[X]$. In particular,
\begin{align*}
\vol_{n,r}(\mathcal{L}[X]^g \cap \mc{C}^{(1)}[X]) & \geq \vol_{n,r}\bigl((\mathcal{L} \cap \mc{C}^{(1)})[X]\bigr) - \frac{1}{6}\vol_{n,r}\bigl(\mc{L}[X]\bigr) \\
& \overset{\mathrm{(i)}}{\geq} n(n-1) \cdot \left((1 - \delta) - \frac{1}{2}(1 + \delta)\right)\vol_{\Pbb,r}( \mathcal{L} \cap \mc{C}^{(1)}) \\
& \overset{\mathrm(ii)}{\geq} n(n - 1) \cdot \frac{1}{7} \vol_{\Pbb,r}( \mathcal{L} \cap \mc{C}^{(1)}) \\
& \overset{\mathrm(iii)}{\geq} \frac{1}{8} \vol_{n,r}( \mathcal{L} \cap \mc{C}^{(1)})
\end{align*}
where in $\mathrm{(i)}$ we have used $\vol_{\Pbb,r}(\mc{L}) \leq 3 \vol_{\Pbb,r}(\mc{C}^{(1)})$, and in $\mathrm{(ii)}$ and $\mathrm{(iii)}$ we have used $\delta \in (0,1/7)$. \qed \\

\noindent \emph{Proof of Lemma~\ref{lem:normalized_cut_lb}.}
To lower bound the normalized cut $\Phi_{n,r}(Z)$, it suffices to lower bound $\cut_{n,r}(Z)$ and upper bound $\vol_{n,r}(Z)$. A crude upper bound on the volume is simply 
\begin{equation}
\label{eqn:normalized_cut_lb_pf4}
\vol_{n,r}(Z) \leq \vol_{n,r}(G_{n,r}) \overset{(i)}{\leq} (1 + \delta)  \vol_{\Pbb,r}(\mathcal{X}) n (n - 1) \leq (1 + \delta)\frac{\pi r^2}{\rho \sigma} n^2
\end{equation}
where by Lemma~\ref{lem:hoeffding}, inequality $(i)$ holds with probability at least $1 - \exp\set{-n\delta^2 \vol_{\Pbb,r}(\mathcal{X})/4}$. This crude upper bound will suffice for our purposes.

We turn to lower bounding $\cut_{n,r}(Z)$. Establishing this lower bound is considerably more involved, and we start by giving an outline of our approach to build intuition.
\begin{itemize}
	\item \emph{Discretization}. Intuitively, we will approximate the cut of $Z$ by discretizing the space $\mathcal{X}$ into a collection of rectangular bins $\wb{Q}$. We will consider a subset of bins $\partial \wb{Z} \subseteq \wb{Q}$, which we call the discretized boundary of $Z$, and which consists of those bins $Q \in \wb{Q}$ containing many points $x_i \in Z$ adjacent, in the graph $G_{n,r}$, to many $x_j \in X \setminus Z$. This discretization allows us to relate $\cut_{n,r}(Z)$ to the number of bins in the discretized boundary $\partial \wb{Z}$, as in~\eqref{eqn:normalized_cut_lb_pf1}.
	\item \emph{Slicing}. In order to lower bound the number of bins in the discretized boundary $|\partial \wb{Z}|$, we consider a second partition of $\mc{X}$ into horizontal slices $R$. We argue that for each slice $R$, one of two things must be true: either $R$ contains a bin belonging to the discretized boundary of $Z$, $Q \in \partial \wb{Z}$, or $(Z \vartriangle \mc{C}^{(1)}[X]) \cap R$ has a ``substantial'' volume,
	\begin{equation*}
	\vol_{n,r}\Bigl((Z \vartriangle \mc{C}^{(1)}[X]) \cap R\Bigr) \geq R_{\min},
	\end{equation*}
	where $R_{\min}$ is a random quantity defined in~\eqref{eqn:normalized_cut_lb_pf2}.
	Summing over all slices $R$ then gives a lower bound on $|\partial \wb{Z}|$, as in~\eqref{eqn:normalized_cut_lb_pf2}. 
	\item \emph{Probabilistic Guarantees}. The lower bound~\eqref{eqn:normalized_cut_lb_pf2} depends on two random quantities: the minimum number of points in any bin (denoted by $Q_{\min}$), and $R_{\min}$, which was discussed in the previous bullet. We give suitable lower bounds on each quantity, which hold with high probability, and which suffice to complete the proof of Lemma~\ref{lem:normalized_cut_lb}.
\end{itemize}

Now, we fill in the details, working step-by-step.\\

\noindent \emph{Step 1: Discretization.}
We begin by discretizing $\mc{X}$ into a collection of bins. To that end, to each $k_1 \in \bigl[\frac{12\sigma}{r}\bigr], k_2 \in \bigl[\frac{4\rho}{r}\bigr]$ associate the bin\footnote{Here we are assuming, without loss of generality, that $12\sigma/r \in \mathbb{N}$ and $4\rho/r \in \mathbb{N}$, and using the notation $[x] = \{1,\ldots,x\}$ for $x \in \mathbb{N}$.}
\begin{equation*}
Q_{(k_1,k_2)} := \biggl[-\frac{3\sigma}{2} + \frac{(k_1 - 1)}{4}r, -\frac{3\sigma}{2} + \frac{k_1}{4}r\biggr] \times \biggl[-\frac{\rho}{2} + \frac{(k_2 - 1)}{4}r, -\frac{\rho}{2} + \frac{k_2}{4}r\biggr];
\end{equation*}
and let $\overline{Q} = \set{Q_{(k_1,k_2)}: k_1 \in \left[\frac{12\sigma}{r}\right], k_2 \in \bigl[\frac{4\rho}{r}\bigr]}$ be the collection of such bins. Next, we define the binned set $\overline{Z} \subset \overline{Q}$ to contain all bins $Q \in \wb{Q}$ for which $Z \cap Q$ is larger than $X\setminus Z \cap Q$; in mathematical notation,
\begin{equation*}
\overline{Z} := \set{Q \in \overline{Q}: \Pbb_n(Z \cap Q) \geq \frac{1}{2}\Pbb_n(Q)}.
\end{equation*}
Then the discretized boundary $\partial \wb{Z}$ is 
\begin{equation*}
\partial \overline{Z} := \set{Q_{(k_1,k_2)} \in \overline{Z}: \exists (\ell_1,\ell_2) \in \Bigl[\frac{12\sigma}{r}\Bigr] \times \Bigl[\frac{4\rho}{r}\Bigr]~~\textrm{s.t.}~~Q_{(\ell_1,\ell_2)} \not\in \wb{Z}, \norm{(k_1,k_2) - (\ell_1,\ell_2)}_1 = 1}.
\end{equation*}
Intuitively, every point $x_i \in Z$ which belongs to a cube in the boundary set $\partial\overline{Z}$ will have many edges to $X\setminus Z$. Formally, letting $Q_{\min} := \min_{Q \in \overline{Q}} \Pbb_n(Q)$, we have
\begin{equation}
\label{eqn:normalized_cut_lb_pf1}
\cut_{n,r}(Z) \geq \cut_{n,r}(Z \cap \set{x_i \in \overline{Z}}) \geq \frac{1}{4} \abs{\partial \overline{Z}} Q_{\min}^2.
\end{equation}
To establish the last inequality, we reason as follows: first, for every cube $Q_{(k_1,k_2)} \in \partial\overline{Z}$, there exists a cube $Q_{(\ell_1,\ell_2)} \not\in \overline{Z}$ such that $\norm{(k_1,k_2) - (\ell_1,\ell_2)}_1 \leq 1$; second, since each cube has side length $r/4$, this implies that for every $x_i \in Q_{(k_1,k_2)}$ and $x_j \in Q_{(\ell_1,\ell_2)}$ the edge $(x_i,x_j)$ belongs to $G_{n,r}$.\\ 

\noindent \emph{Step 2: Slicing.}
Now we move on lower bounding the size of the discretized boundary $\abs{\partial\wb{Z}}$. To do so, we divide $\mathcal{X}$ into slices horizontally. Let $$R_k = \set{(x_1,x_2) \in \mathcal{X}: x_2 \in \bigl[-\frac{\rho}{2} + \frac{(k - 1)}{4}r, -\frac{\rho}{2} + \frac{k}{4}r\bigr]}$$ 
be the $k$th horizontal slice, and $\wb{R}_k = \set{Q_{(k_1,k)} \in \overline{Q}:k_1 \in [\frac{12\sigma}{r}]}$ be the binned version of $R_k$. For every $k$, we are in one of three cases:
\begin{enumerate}
	\item $\wb{R}_k \cap \wb{Z} = \emptyset$, in which case
	\begin{equation*}
	\vol_{n,r}\Bigl( \bigl(Z \vartriangle \mc{C}^{(1)}[X]\bigr) \cap R_k \Bigr) \geq \frac{1}{2}\vol_{n,r}(\mc{C}^{(1)}[X] \cap R_k), ~~ \textrm{or}
	\end{equation*}
	\item $\wb{R}_k \cap \wb{Z} = \wb{R}_k$, in which case
	\begin{equation*}
	\vol_{n,r}\Bigl( \bigl(Z \vartriangle \mc{C}^{(1)}[X]\bigr) \cap R_k \Bigr) \geq \frac{1}{2}\vol_{n,r}(\mc{C}^{(2)}[X] \cap R_k), ~~\textrm{or}
	\end{equation*}
	\item $\wb{R}_k \cap \partial \wb{Z} \neq \emptyset$.
\end{enumerate}
In words, these cases describe, for each slice $R$, the situation in which the set $Z$ has a ``substantial'' boundary (Case 3), or the situation in which $Z$ does not have a substantial boundary, and instead has a substantial difference with $\mc{C}^{(1)}$ (Cases 1 and 2). 

Now, let $N(R)$ be the number of slices which contain a boundary cell,
\begin{equation*}
N(R) := \#\biggl\{k \in \biggl[\frac{2\rho}{r}\biggr]: \wb{R}_k \cap \partial \wb{Z} \neq \emptyset\biggr\}.
\end{equation*}
and let $R_{\min}$ be the minimum volume of the intersection of any slice $R_k$ with the density clusters $\mc{C}^{(1)}, \mc{C}^{(2)}$; formally,
\begin{equation}
\label{eqn:normalized_cut_lb_pf1.5}
R_{\min} := \min_{k}\Bigl\{\vol_{n,r}(\mc{C}^{(1)}[X] \cap R_k) \wedge \vol_{n,r}(\mc{C}^{(2)}[X] \cap R_k)\Bigr\}.
\end{equation} 
Summing over all slices $R_k$, we have from our previous analysis that,
\begin{equation*}
\vol_{n,r}\biggl(Z \vartriangle \mc{C}^{(1)}[X]\biggr) = \sum_{k \in [2\rho/r]} \vol_{n,r}\biggl( (Z \vartriangle \mc{C}^{(1)}[X]) \cap R_k\biggr) \geq R_{\min}  \left[\frac{2\rho}{r} - \frac{N(R)}{2}\right].
\end{equation*}
Rearranging this expression and noting that $\abs{\partial\wb{Z}} \geq N(R)$,  we obtain the following lower bound on the cardinality of the discretized boundary,
\begin{equation}
\label{eqn:normalized_cut_lb_pf2}
\abs{\partial\wb{Z}} \geq N(R) \geq 2\Bigl(\frac{2\rho}{r} - \frac{\vol_{n,r}(Z \vartriangle \mc{C}^{(1)}[X])}{R_{\min}}\Bigr).
\end{equation}
Combining \eqref{eqn:normalized_cut_lb_pf1} and \eqref{eqn:normalized_cut_lb_pf2}, we have that
\begin{align}
\cut_{n,r}(Z) & \geq \frac{1}{4}N(R) Q_{\min}^2 \nonumber \\
& \geq \left(\frac{\rho}{r} - \frac{\vol_{n,r}(Z \vartriangle \mc{C}^{(1)}[X])}{2R_{\min}}\right) Q_{\min}^2 \label{eqn:normalized_cut_lb_pf3}
\end{align}
for all $Z \subset X$.\\

\noindent \emph{Step 3: Probabilistic guarantees.}
It remains to lower bound the random quantities $R_{\min}$ and $Q_{\min}$.
To do so, we first lower bound the expected probability of any cell $Q$,
\begin{equation}
\label{eqn:normalized_cut_lb_pf4.5}
\min_{Q \in \overline{Q}} \Pbb(Q) \geq \frac{\epsilon r^2}{4 \rho \sigma},
\end{equation}
and the expected volume of $\mc{C}^{(1)}[X] \cap R_k$ and $\mc{C}^{(2)}[X] \cap R_k$,
\begin{equation}
\label{eqn:normalized_cut_lb_pf5}
\vol_{\Pbb,r}(\mc{C}^{(1)} \cap R_k) = \vol_{\Pbb,r}(\mc{C}^{(2)} \cap R_k) \geq \frac{\pi r^3}{4 \sigma \rho^2}~~\textrm{for all $k$.} 
\end{equation}

Since $Q_{\min}$ and $R_{\min}$ are obtained by taking the minimum of functionals over a fixed number of sets in $n$, they concentrate tightly around their means. Specifically, note that the total number of cubes is $\abs{\overline{Q}} = \frac{48 \sigma \rho}{r^2}$, and the total number of horizontal slices is $\frac{4\rho}{r}$. Along with \eqref{eqn:normalized_cut_lb_pf4.5} and \eqref{eqn:normalized_cut_lb_pf5}, by Lemma~\ref{lem:bernstein_union} 
\begin{equation*}
Q_{\min} \geq (1 - \delta)\frac{\epsilon r^2}{4 \rho \sigma}n ~~\textrm{and}~~ R_{\min} \geq (1 - \delta)\frac{\pi r^3}{4\sigma \rho^2 }n(n - 1),
\end{equation*}
with probability at least $1 - \frac{48 \sigma \rho}{r^2}\exp\set{-\frac{n\delta^2\epsilon r^2}{4 \rho \sigma}} - \frac{4\rho}{r}\exp\set{-\frac{n(n - 1)\delta^2\pi r^3}{8\sigma\rho^2}}$. Combining these lower bounds with \eqref{eqn:normalized_cut_lb_pf4} and \eqref{eqn:normalized_cut_lb_pf3}, we obtain 
\begin{equation*}	
\Phi_{n,r}(Z) \geq \frac{(1 - \delta)^2}{4(1 + \delta)\pi} \left(1 - 2 \frac{\sigma \rho}{(1 - \delta) r^2 n^2} \vol_{n,r}(Z \vartriangle \mc{C}^{(1)}[X]) \right) \frac{\epsilon^2 r}{\sigma},
\end{equation*}
for all $Z \subseteq X$. \qed 

\section{Additional Results: aPPR and Separation of Clusters via PPR}
\label{apdx:appr_misclassification_error}
In this appendix, we prove two additional results regarding PPR remarked upon in our main text. In Section~\ref{subsec:appr_volume_ssd_ub}, we show that clustering using the aPPR vector satisfies an equivalent guarantee to Theorem~\ref{thm:volume_ssd_ub}. In Section~\ref{subsec:consistent_recovery_density_clusters}, we show that the PPR vector can perfectly distinguish two distinct density clusters $\mc{C}_{\lambda},\mc{C}_{\lambda}'$. 

\subsection{Generic Cluster Recovery with aPPR}
\label{subsec:appr_volume_ssd_ub}
Our formal claim regarding cluster recovery with aPPR is contained in Corollary~\ref{cor:appr}.
\begin{corollary}
	\label{cor:appr}
	Consider instead of
	Algorithm \ref{alg:ppr} using the approximate PPR vector from
	\citet{andersen2006} satisfying \eqref{eqn:appr_error}, and forming the 
	corresponding cluster estimate \smash{$\wh{C}$} in the same manner.  Then 
	provided we take 
	\begin{equation}
	\label{eqn:appr_parameter}
	\varepsilon = \frac{1}{25(1 + \delta)n(n - 1)\vol_{\Pbb,r}(\mc{C})} ,
	\end{equation}
	under the assumptions of Theorem~\ref{thm:volume_ssd_ub} the upper bound on symmetric set difference in \eqref{eqn:volume_ssd_ub} still
	holds.
\end{corollary}	
\emph{Proof of Corollary~\ref{cor:appr}.}
Note that the choice of $\varepsilon$ in~\eqref{eqn:appr_parameter} implies $\varepsilon \leq 1/\bigl(25\vol_{n,r}(\mc{C}[X])\bigr)$ with probability at least $1 - \exp\bigl\{-n\delta^2\vol_{\Pbb,r}(\mc{C})\bigr\}$. The proof of Corollary~\ref{cor:appr} is then identical to that of Theorem~\ref{thm:volume_ssd_ub}, except one uses Corollary~\ref{cor:zhu} rather than Lemma~\ref{lem:zhu} to relate the symmetric set difference to the graph normalized cut and mixing time. \qed 

\subsection{Perfectly Distinguishing Two Density Clusters}
\label{subsec:consistent_recovery_density_clusters}

As mentioned in our main text, the symmetric set difference does not measure whether \smash{$\wh{C}$}
can (perfectly) distinguish any two distinct clusters \smash{$\mc{C}_{\lambda},\mc{C}_{\lambda}' \in 
	\mathbb{C}_f(\lambda)$}. We therefore also study a second notion of cluster 
estimation, first introduced by \citet{hartigan1981}.

\begin{definition}
	\label{def:density_cluster_consistency}
	For an estimator \smash{$\wh{C} \subseteq X$} and distinct clusters \smash{$\mc{C}_{\lambda}, \mc{C}_{\lambda}' \in \mathbb{C}_f(\lambda)$}, we say \smash{$\wh{C}$} \emph{separates} $\mc{C}_{\lambda}$ from $\mc{C}_{\lambda}'$ if 
	\begin{equation}
	\label{eqn:density_cluster_consistency}
	\mc{C}_{\lambda}[X] \subseteq \wh{C} \quad \text{and} \quad
	\wh{C} \cap \mc{C}_{\lambda}'[X] = \emptyset.
	\end{equation}
\end{definition}

The bound on symmetric set difference \eqref{eqn:density_cluster_volume_ssd_ub} does not imply \eqref{eqn:density_cluster_consistency}, which requires a
uniform bound over the PPR vector $p_v$. As an example, suppose that we were
able to show that for all \smash{$\mc{C}' \in \mathbb{C}_f(\lambda), \mc{C}' \neq \mc{C}$}, and each $u \in \mc{C}, u' \in \mc{C}'$,  
\begin{equation}
\label{eqn:ppr_gap}
\frac{p_v(u')}{\deg(u';G)} \leq \frac{1}{10} \cdot \frac{1}{n (n - 1) \vol_{\Pbb,r}(\mc{C}_{\lambda,\sigma})} <
\frac{1}{5} \cdot \frac{1}{n (n-1) \vol_{\Pbb,r}(\mc{C}_{\lambda,\sigma})} \leq \frac{p_v(u)}{\deg(u;G)}. 
\end{equation}
Then, any $(L,U)$ satisfying~\eqref{eqn:initialization} and any sweep cut
$S_{\beta}$ for $\beta \in (L,U)$ would result in a cluster estimate $\wh{C}$ fulfilling both conditions laid out in
\eqref{eqn:density_cluster_consistency}. In Theorem 
\ref{thm:density_cluster_consistent_recovery}, we show that a sufficiently 
small upper bound on $\Delta(\wh{C},\mc{C}_{\lambda,\sigma}[X])$ ensures that with high probability the uniform bound~\eqref{eqn:ppr_gap} is satisfied, and hence implies
that \smash{$\wh{C}$} will separate $\mc{C}_{\lambda}$ from $\mc{C}_{\lambda}'$. In what follows, put
\begin{equation*}
c_{1,\delta} := \frac{(1 - \delta)^5}{4} \cdot \min\Bigl\{\Bigl(\frac{3}{8(1 + \delta)} - \frac{1}{5}\Bigr),\frac{1}{10}\Bigr\}
\end{equation*}
and note that if $\delta \in (0,7/8)$ then $c_{1,\delta} > 0$. (In fact, we will have to take $\delta \in (0,1/4)$ in order to use Propositions~\ref{prop:sample_to_population_1} and~\ref{prop:sample_to_population_2}). Additionally, denote $\Pbb'$ for the conditional distribution of a sample point given that it falls in $\mc{C}_{\lambda,\sigma}'$, i.e. $\Pbb'(\mc{S}) := \Pbb(\mc{S} \cap \mc{C}_{\lambda,\sigma}')/\Pbb(\mc{C}_{\lambda,\sigma}')$, and $G_{n,r}' := G_{n,r}[\mc{C}_{\lambda,\sigma}']$ for the subgraph of $G_{n,r}$ induced by $\mc{C}_{\lambda,\sigma}'$. 
\begin{theorem}
	\label{thm:density_cluster_consistent_recovery}
	For any $\delta \in (0,1/4)$ any $n \in \mathbb{N}$ such that
	\begin{equation}
	\label{eqn:density_cluster_consistent_recovery_sample_complexity}
	\frac{1}{n} \leq \delta \cdot \frac{4\Pbb(\mc{C}_{\lambda,\sigma}')}{3}
	\end{equation} 
	and otherwise under the same conditions as Theorem~\ref{thm:volume_ssd_ub}, the following statement holds with probability at least $1 - B_2/n - 4\exp\{-b_1\delta^2n\} - (n + 2)\exp\{-b_3n\} - 3(n + 3)\exp\{-b_7\delta^2n\}$: there exists a set $\mc{C}_{\lambda,\sigma}[X]^g \subseteq \mc{C}_{\lambda,\sigma}[X]$ of large volume, $\vol_{n,r}(\mc{C}_{\lambda,\sigma}[X]^g) \geq \vol_{n,r}(\mc{C}_{\lambda,\sigma}[X])/2$, such that if Algorithm~\ref{alg:ppr} is $\delta$-well-initialized and run with any seed node $v \in \mc{C}_{\lambda,\sigma}[X]^g$, and moreover
	\begin{equation}
	\label{eqn:density_cluster_consistent_recovery_condition}
	\kappa_{\Pbb,r}(\mc{C}_{\lambda,\sigma},\delta) \leq c_{1,\delta} \cdot \frac{\min\bigl\{\Pbb(\mc{C}_{\lambda,\sigma})^2 \cdot d_{\min}(\wt{\Pbb})^2, \Pbb(\mc{C}_{\lambda,\sigma}')^2 \cdot d_{\min}(\Pbb')^2\bigl\}}{\vol_{\Pbb,r}(\mc{C}_{\lambda,\sigma})}
	\end{equation}
	then the PPR estimated cluster \smash{$\wh{C}$} satisfies \eqref{eqn:density_cluster_consistency}.
\end{theorem}

Before we prove Theorem~\ref{thm:density_cluster_consistent_recovery}, we make a few brief remarks:
\begin{itemize}
	\item In one sense, Theroem~\ref{thm:density_cluster_consistent_recovery} is a strong result: if the density clusters $\mc{C}_{\lambda},\mc{C}_{\lambda}'$ satisfies the requirement~\eqref{eqn:density_cluster_consistent_recovery_condition}, and we are willing to ignore the behavior of the algorithm in low-density regions, Theorem~\ref{thm:density_cluster_consistent_recovery} guarantees that PPR will \emph{perfectly distinguish} the candidate cluster $\mc{C}_{\lambda}$ from $\mc{C}_{\lambda}'$.
	
	\item On the other hand, unfortunately the requirement~\eqref{eqn:density_cluster_consistent_recovery_condition} is rather restrictive. Suppose the density cluster $\mc{C}_{\lambda,\sigma}$ satisfies~\ref{asmp:lambda_bounded_density}. Then from the following chain of inequalities,
	\begin{equation*}
	\frac{\Delta(\wh{C}, \mc{C}_{\lambda,\sigma}[X])}{\vol_{n,r}(\mc{C}_{\lambda,\sigma}[X])} \overset{(\textrm{Thm.~\ref{thm:volume_ssd_ub}})}{\leq} \kappa_{\Pbb,r}(\mc{C},\delta) \overset{\eqref{eqn:density_cluster_consistent_recovery_condition}}{\leq} c_{1,\delta} \cdot \frac{\Pbb(\mc{C}_{\lambda,\sigma})^2 \cdot d_{\min}(\wt{\Pbb})^2}{\vol_{\Pbb,r}(\mc{C}_{\lambda,\sigma})} \overset{\ref{asmp:lambda_bounded_density}}{\leq} c_{1,\delta} \cdot \frac{\Lambda_{\sigma}}{\lambda_{\sigma}} \nu_d r^d,
	\end{equation*}
	we see that in order for~\eqref{eqn:density_cluster_consistent_recovery_condition} to be met, it is necessary that \smash{$\Delta(\wh{C}, \mc{C}_{\lambda,\sigma}[X])/\vol_{n,r}(\mc{C}_{\lambda,\sigma}[X])$} be on the order of $r^d$. In plain terms, we 
	are able to recover a density cluster $\mc{C}_{\lambda}$ in the strong sense of
	\eqref{eqn:density_cluster_consistency} only when we can guarantee the volume of the symmetric set difference will be very small. This strong condition is
	the price we pay in order to obtain the uniform bound in~\eqref{eqn:ppr_gap}. 
	\item The proof of Theorem~\ref{thm:density_cluster_consistent_recovery} relies heavily on Lemma~\ref{lem:ppr_uniform_bound}. This lemma---or more accurately, the equation~\eqref{pf:interpolator_bound_max_entry} used in the proof of the lemma---can be thought of as a smoothness result for the PPR vector, showing that the mass of $p_v(\cdot)$ cannot be overly concentrated at any one vertex $u \in V$. However,~\eqref{pf:interpolator_bound_max_entry} is a somewhat crude bound. By plugging a stronger result on the smoothness of $p_v(\cdot)$ in to the proof of Lemma~\ref{lem:ppr_uniform_bound}, we could improve the uniform bounds of the lemma, and in turn show that the conclusion of Theorem~\ref{thm:density_cluster_consistent_recovery} holds under weaker conditions than~\eqref{eqn:density_cluster_consistent_recovery_condition}. 
\end{itemize}  

Now we recall some probabilistic estimates before proceeding to the proof of Theorem~\ref{thm:density_cluster_consistent_recovery}.\\

\noindent \emph{Probabilistic estimates.}
As in the proof of Theorem~\ref{thm:volume_ssd_ub}, we will assume that the inequalities~\eqref{eqn:sample_to_population_normalized_cut}-\eqref{eqn:sample_to_population_conductance} and~\eqref{pf:volume_ssd_ub_1}-\eqref{pf:volume_ssd_ub_3} are satisfied. We will additionally assume that
\begin{equation}
n' \geq (1 - \delta) \cdot n \cdot \Pbb(\mc{C}_{\lambda,\sigma}') \overset{\eqref{eqn:density_cluster_consistent_recovery_sample_complexity}}{\Longrightarrow} n' - 1 \geq (1 - \delta)^2 \cdot n \cdot \Pbb(\mc{C}_{\lambda,\sigma}')
\end{equation}
and that 
\begin{equation}
\frac{1}{n' - 1} d_{\min}(G_{n,r}') \geq (1 - \delta) \cdot d_{\min}(\Pbb').
\end{equation}
By Propositions~\ref{prop:sample_to_population_1}-\ref{prop:sample_to_population_2} and Lemmas~\ref{lem:hoeffding_2}-\ref{lem:bernstein_union}, these inequalities hold with probability at least $1 - B_2/n - 4\exp\{-b_1\delta^2n\} - (n + 1)\exp\{-b_3n\} - (3n + 3)\exp\{-b_7\delta^2n\}$, taking $b_7 := b_2 \wedge \Pbb(\mc{C}_{\lambda,\sigma}') \cdot d_{\min}(\Pbb')/9$.\\

\noindent \emph{Proof of Theorem~\ref{thm:density_cluster_consistent_recovery}.}
We have already verified in the proof of Theorem~\ref{thm:volume_ssd_ub} that $\alpha \leq 1/(2\tau_{\infty}(\wt{G}_{n,r}))$, and we may therefore apply Lemma~\ref{lem:ppr_uniform_bound}, which gives a lower bound on $p_v(u)$ for all $u \in \mc{C}_{\lambda,\sigma}[X]_o$ and an upper bound on $p_v(u')$ for all $p_v(u')$ for all $u' \in \mc{C}_{\lambda,\sigma}'[X]_o$. These bounds are useful because $r \leq \sigma$, which implies that $\mc{C}_{\lambda}[X] \subseteq\mc{C}_{\lambda,\sigma}[X]_o$ and likewise that $C_{\lambda}'[X] \subseteq\mc{C}_{\lambda,\sigma}'[X]_o$. We will show that these bounds in turn imply~\eqref{eqn:ppr_gap}, from which the claim of the theorem follows.

We begin with the lower bound in~\eqref{eqn:ppr_gap}. From (in order) Lemma~\ref{lem:ppr_uniform_bound}, our assorted probabilistic estimates, and the assumed lower bound~\eqref{eqn:density_cluster_consistent_recovery_condition} on $\kappa_{\Pbb,r}(\mc{C},\delta)$, we have that for all $u \in \mc{C}_{\lambda}[X]$,
\begin{align*}
\frac{p_v(u)}{\deg_{n,r}(u)} & \geq \frac{3}{8\vol_{n,r}(\mc{C}[X])} - 2\frac{\Phi_{n,r}(\mc{C}_{\lambda,\sigma}[X])}{d_{\min}(\wt{G}_{n,r})^2 \alpha} \\
& \geq \frac{1}{n(n - 1)}\biggl(\frac{3}{8(1 + \delta)\vol_{\Pbb,r}(\mc{C}_{\lambda,\sigma})} - 4 \cdot \frac{n(n - 1)}{(\wt{n} - 1)^2} \cdot \frac{\kappa_{\Pbb,r}(\mc{C}_{\lambda,\sigma},\delta)}{(1 - \delta) d_{\min}(\wt{\Pbb})^2}\biggr) \\
& \geq \frac{1}{n(n - 1)}\biggl(\frac{3}{8(1 + \delta)\vol_{\Pbb,r}(\mc{C}_{\lambda,\sigma})} - 4 \cdot \frac{\kappa_{\Pbb,r}(\mc{C}_{\lambda,\sigma},\delta)}{(1 - \delta)^5 \Pbb(\mc{C})^2 d_{\min}(\wt{\Pbb})^2}\biggr) \\
& \geq \frac{1}{5 n^2\vol_{\Pbb,r}(\mc{C}_{\lambda,\sigma})}.
\end{align*}
An equivalent derivation implies the upper bound in~\eqref{eqn:ppr_gap}: for all $u' \in \mc{C}_{\lambda}'[X]$,
\begin{align*}
\frac{p_v(u')}{\deg_{n,r}(u')} & \leq 2 \frac{\Phi_{n,r}(\mc{C}_{\lambda,\sigma}[X])}{d_{\min}(G_{n,r}')^2 \alpha} \\
& \leq 4 \cdot \frac{1}{(n' - 1)^2} \cdot \frac{ \kappa(\mc{C}_{\lambda,\sigma},\delta)}{(1 - \delta) d_{\min}(\Pbb')^2} \\
& \leq 4\frac{ \kappa(\mc{C}_{\lambda,\sigma},\delta)}{n^2 (1 - \delta)^5 \Pbb(\mc{C}_{\lambda,\sigma}')^2 d_{\min}(\Pbb')^2} \leq \frac{1}{10 n^2\vol_{\Pbb,r}(\mc{C}_{\lambda,\sigma})},
\end{align*}
completing the proof of Theorem~\ref{thm:density_cluster_consistent_recovery}. \qed

\section{Experimental Details}
\label{apdx:experimental_details}
Finally, we detail the settings of our experiments, and include an additional figure.

\subsection{Experimental Settings for Figure~\ref{fig:bounds}}
Let $\mc{R}_{\sigma,\rho} = [-\sigma/2,\sigma/2] \times [-\rho/2,\rho/2]$ be the two-dimensional rectangle of width $\sigma$ and height $\rho$, centered at the origin. We sample $n = 8000$ points according to the density function $f_{\rho,\sigma,\lambda}$, defined over domain $\mc{X} = [-1,1]^2$ and parameterized by $\rho, \sigma$ and $\lambda$ as follows:
\begin{equation}
f_{\rho,\sigma,\lambda}(x) :=
\begin{dcases}
\lambda,~~& \textrm{if $x \in R_{\sigma,\rho} - (-.5,0)$ or $x \in R_{\sigma,\rho} + (-.5,0)$} \\
\frac{4 - 2 \lambda \rho \sigma}{1 - 2\rho\sigma},~~& \textrm{if $x \in \mc{X}$, $x \not\in R_{\sigma,\rho} - (-.5,0)$ and $x \not\in R_{\sigma,\rho} + (-.5,0)$.}
\end{dcases}
\end{equation}
Then $\theta := \lambda - \frac{4 - 2 \lambda \rho \sigma}{1 - 2\rho\sigma}$ measures the difference in density between the density clusters and the rest of the domain. The first column displays $n = 8000$ points sampled from three different parameterizations of $f_{\rho,\sigma,\lambda}$:
\begin{align*}
\rho &= .913, &&\sigma = .25, &&&(\lambda - \theta)/\lambda = .25 \tag{top panel} \\
\rho &= .25, &&\sigma = , &&&(\lambda - \theta)/\lambda = .05 \tag{middle panel} \\
\rho &= .5, &&\sigma = .25, &&&(\lambda - \theta)/\lambda = .12 \tag{bottom panel.}
\end{align*}
In each of the first, second, and third rows, we fix two parameters and vary the third. In the first row, we fix $\sigma = .25$, $(\lambda - \theta)/\lambda = .25$, and vary $\rho$ from $.25$ to $2$. In the second row, we fix $\rho = 1.8$, $(\lambda - \theta)/\lambda = .05$, and vary $\sigma$ from $.1$ to $.2$ In the third row, we fix $\rho = .5$, $\sigma = .25$ and vary $(\lambda - \theta)/\lambda$ from $.1$ to $.25$. In the first and third rows, we take $r = \sigma/8$; in the second row, where we vary $\sigma$, we take $r = .1/8$. 

\subsection{Experimental Settings for Figure \ref{fig:moons}}

To form each of the three rows in Figure \ref{fig:moons}, $n = 800$ points are independently sampled following a 'two moons plus Gaussian noise model'. Formally, the (respective) generative models for the data are
\begin{align}
Z & \sim \textrm{Bern}(1/2), \theta \sim \textrm{Unif}(0, \pi) \\
X(Z,\theta) & = 
\begin{cases}
\mu_1 + (r \cos(\theta), r \sin(\theta)) + \sigma \epsilon,~ & \text{if}~ Z = 1 \\
\mu_2 + (r \cos(\theta), - r \sin(\theta)) + \sigma \epsilon,~ & \text{if}~ Z = 0
\end{cases}
\end{align}
where 
\begin{align*}
\mu_1 & = (-.5, 0),~ \mu_2 = (0,0),~ \epsilon \sim N(0, I_2) \tag{row 1} \\
\mu_1 & = (-.5, -.07),~ \mu_2 = (0,.07),~ \epsilon \sim N(0, I_2) \tag{row 2} \\
\mu_1 & = (-.5, -.125),~ \mu_2 = (0,.125),~ \epsilon \sim N(0, I_2) \tag{row 3} 
\end{align*}
for $I_d$ the $d \times d$ identity matrix. In all cases $\sigma = .07$. In each case $\lambda$ is taken as small as possible such that there exist exactly two distinct density clusters, which we call $\mc{C}_{\lambda}$ and $\mc{C}_{\lambda}'$; $r$ is taken as small as possible so that each vertex has at least $2$ neighbors. The first column consists of the empirical density clusters $\mc{C}_{\lambda}[X]$ and $\mc{C}_{\lambda}'[X]$ for a particular threshold $\lambda$ of the density function; the second column shows the PPR plus minimum normalized sweep cut cluster, with hyperparameter $\alpha$ and all sweep cuts considered; the third column shows the global minimum normalized cut, computed according to the algorithm of \cite{bresson2013}; and the last column shows a cut of the density cluster tree estimator of \cite{chaudhuri2010}. \\

\begin{figure}
	\centering
	\begin{adjustbox}{minipage=\linewidth,scale=0.8}
		\begin{subfigure}{.24\linewidth}
			\includegraphics[width=\linewidth]{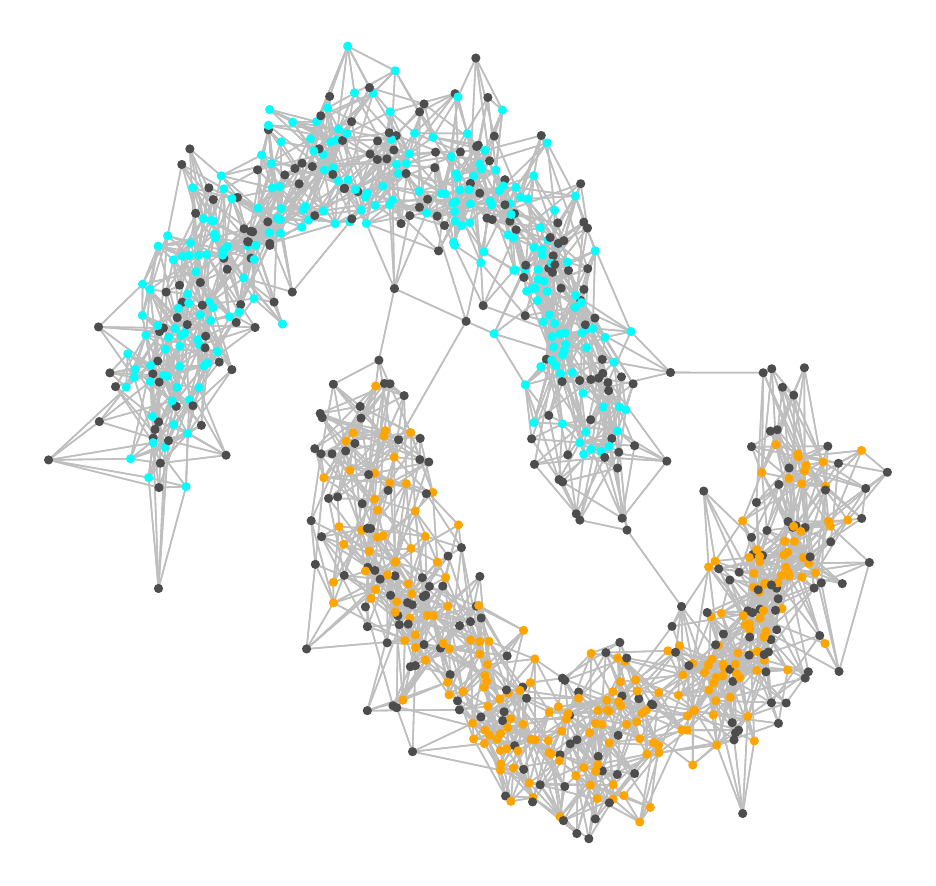}
			\caption{}
		\end{subfigure}
		\begin{subfigure}{.24\linewidth}
			\includegraphics[width=\linewidth]{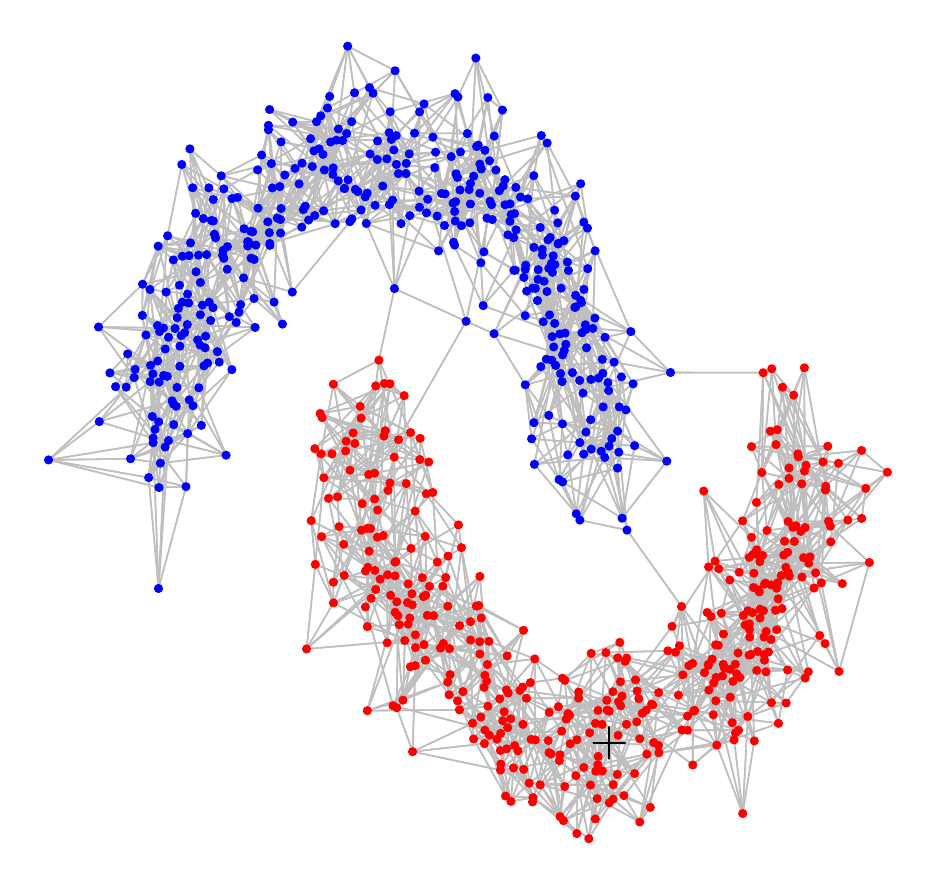}
			\caption{}
		\end{subfigure}
		\begin{subfigure}{.24\linewidth}
			\includegraphics[width=\linewidth]{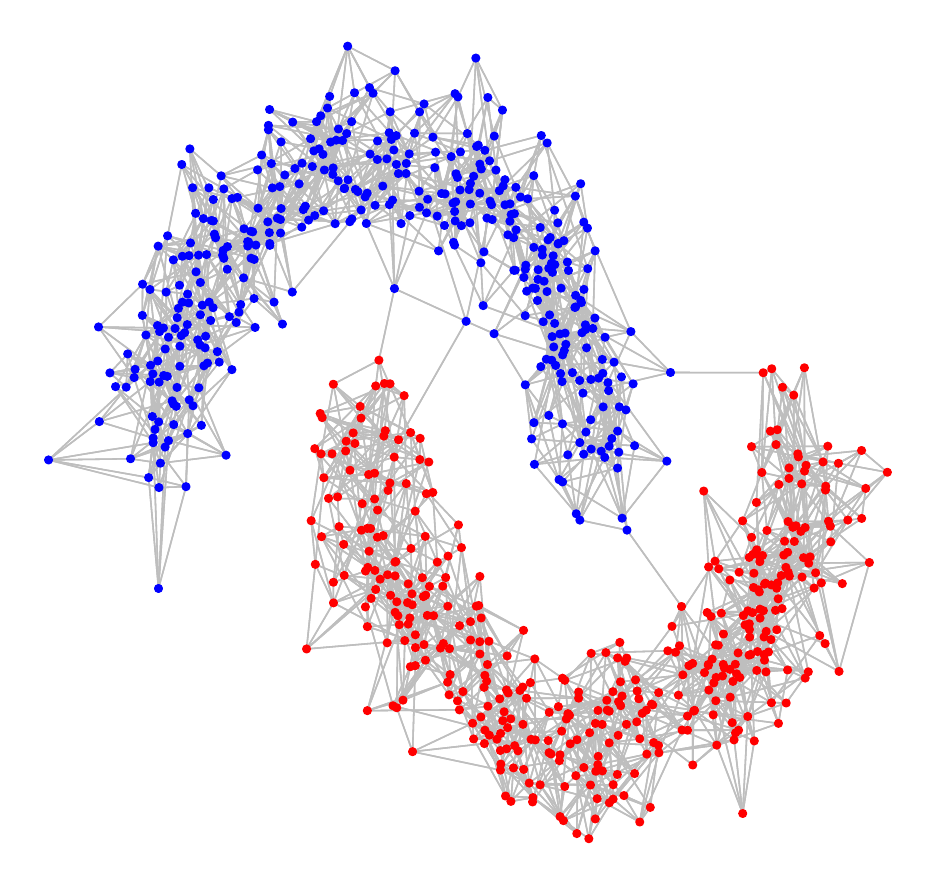}
			\caption{}
		\end{subfigure}
		\begin{subfigure}{.24\linewidth}
			\includegraphics[width=\linewidth]{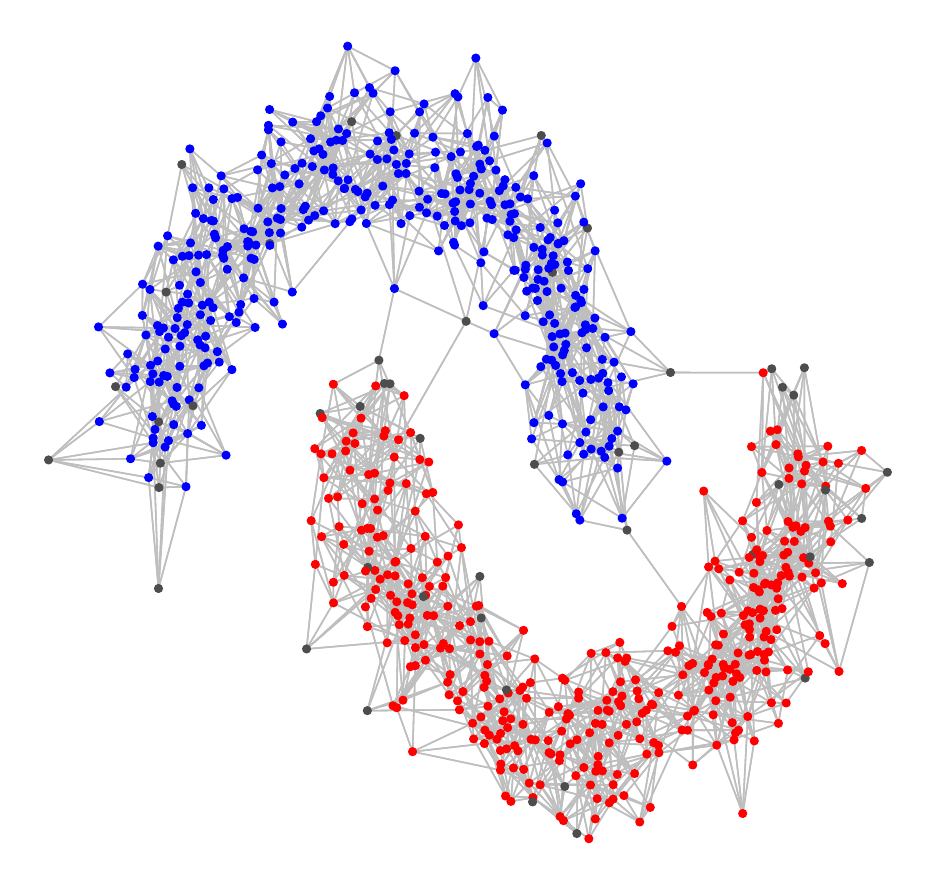}
			\caption{}
		\end{subfigure}
		\caption{\small True density (column 1), PPR (column 2), minimum normalized
			cut (column 3) and estimated density (column 4) clusters for two-moons with 10 dimensional noise. Seed node for PPR denoted by a black cross.} 
		\label{fig:moons_hd}
	\end{adjustbox}
\end{figure}

\noindent \emph{Performance of PPR with high-dimensional noise.} Figure \ref{fig:moons_hd} is similar to Figure \ref{fig:moons} of the main text, but with parameters
\begin{equation*}
\mu_1 = (-.5, -.025),~ \mu_2 = (0,.025),~ \epsilon \sim N(0, I_{10}).
\end{equation*}
The gray dots in $(a)$ (as in the left-hand column of Figure \ref{fig:moons} in the main text) represent observations in low-density regions. While the PPR sweep cut $(b)$ has relatively high symmetric set difference with the chosen density cut, it still recovers separates $C_{\lambda}[X]$ and $\mc{C}_{\lambda}'[X]$, in the sense of Definition \ref{def:density_cluster_consistency}.

\begin{thebibliography}{70}
\providecommand{\natexlab}[1]{#1}
\providecommand{\url}[1]{\texttt{#1}}
\expandafter\ifx\csname urlstyle\endcsname\relax
  \providecommand{\doi}[1]{doi: #1}\else
  \providecommand{\doi}{doi: \begingroup \urlstyle{rm}\Url}\fi

\bibitem[Abbasi-Yadkori(2016)]{abbasi-yadkori2016a}
Yasin Abbasi-Yadkori.
\newblock Fast mixing random walks and regularity of incompressible vector
  fields.
\newblock \emph{arXiv preprint arXiv:1611.09252}, 2016.

\bibitem[Abbasi-Yadkori et~al.(2017)Abbasi-Yadkori, Bartlett, Gabillon, and
  Malek]{abbasi-yadkori2016}
Yasin Abbasi-Yadkori, Peter Bartlett, Victor Gabillon, and Alan Malek.
\newblock {Hit-and-Run for Sampling and Planning in Non-Convex Spaces}.
\newblock In \emph{Proceedings of the 20th International Conference on
  Artificial Intelligence and Statistics}, pages 888--895, 2017.

\bibitem[Allen-Zhu et~al.(2013)Allen-Zhu, Lattanzi, and Mirrokni]{zhu2013}
Zeyuan Allen-Zhu, Silvio Lattanzi, and Vahab~S Mirrokni.
\newblock A local algorithm for finding well-connected clusters.
\newblock In \emph{Proceedings of the 30th International Conference on
  International Conference on Machine Learning}, pages 396--404, 2013.

\bibitem[Andersen and Peres(2009)]{andersen2009}
Reid Andersen and Yuval Peres.
\newblock Finding sparse cuts locally using evolving sets.
\newblock In \emph{Proceedings of the 41st Annual ACM Symposium on Theory of
  Computing}, pages 235--244, 2009.

\bibitem[Andersen et~al.(2006)Andersen, Chung, and Lang]{andersen2006}
Reid Andersen, Fan Chung, and Kevin Lang.
\newblock Local graph partitioning using pagerank vectors.
\newblock In \emph{Proceedings of the 47th Annual IEEE Symposium on Foundations
  of Computer Science}, pages 475--486, 2006.

\bibitem[Andersen et~al.(2012)Andersen, Gleich, and Mirrokni]{andersen2012}
Reid Andersen, David~F Gleich, and Vahab Mirrokni.
\newblock Overlapping clusters for distributed computation.
\newblock In \emph{Proceedings of the Fifth ACM International Conference on Web
  Search and Data Mining}, pages 273--282, 2012.

\bibitem[Arias-Castro(2009)]{arias-castro2009}
Ery Arias-Castro.
\newblock Clustering based on pairwise distances when the data is of mixed
  dimensions.
\newblock \emph{arXiv preprint arXiv:0909.2353}, 2009.

\bibitem[Balakrishnan et~al.(2011)Balakrishnan, Xu, Krishnamurthy, and
  Singh]{balakrishnan2011}
Sivaraman Balakrishnan, Min Xu, Akshay Krishnamurthy, and Aarti Singh.
\newblock Noise thresholds for spectral clustering.
\newblock In \emph{Advances in Neural Information Processing Systems 24}, pages
  954--962, 2011.

\bibitem[Balakrishnan et~al.(2013)Balakrishnan, Narayanan, Rinaldo, Singh, and
  Wasserman]{balakrishnan2013}
Sivaraman Balakrishnan, Srivatsan Narayanan, Alessandro Rinaldo, Aarti Singh,
  and Larry Wasserman.
\newblock Cluster trees on manifolds.
\newblock In \emph{Advances in Neural Information Processing Systems 26}, pages
  2679--2687, USA, 2013.

\bibitem[Belkin and Niyogi(2007)]{belkin07}
Mikhail Belkin and Partha Niyogi.
\newblock Convergence of laplacian eigenmaps.
\newblock In \emph{Advances in Neural Information Processing Systems 19}, pages
  129--136. 2007.

\bibitem[Bresson et~al.(2012)Bresson, Laurent, Uminsky, and
  Brecht]{bresson2013}
Xavier Bresson, Thomas Laurent, David Uminsky, and James Brecht.
\newblock Convergence and energy landscape for cheeger cut clustering.
\newblock In \emph{Advances in Neural Information Processing Systems 25}, pages
  1385--1393, 2012.

\bibitem[Calder and Garc\'{\i}a~Trillos(2019)]{calder2019}
Jeff Calder and Nicolas Garc\'{\i}a~Trillos.
\newblock Improved spectral convergence rates for graph laplacians on
  epsilon-graphs and k-nn graphs.
\newblock \emph{arXiv preprint arXiv:1910.13476}, 2019.

\bibitem[Calder et~al.(2022)Calder, Garc\'{\i}a~Trillos, and
  Lewicka]{calder2020}
Jeff Calder, Nicolas Garc\'{\i}a~Trillos, and Marta Lewicka.
\newblock Lipschitz regularity of graph laplacians on random data clouds.
\newblock \emph{To appear, SIAM Journal on Mathematical Analysis}, 2022.

\bibitem[Chaudhuri and Dasgupta(2010)]{chaudhuri2010}
Kamalika Chaudhuri and Sanjoy Dasgupta.
\newblock Rates of convergence for the cluster tree.
\newblock In \emph{Advances in Neural Information Processing Systems 23}, pages
  343--351. 2010.

\bibitem[Chung(1997)]{chung1997}
Fan~RK Chung.
\newblock \emph{Spectral graph theory}.
\newblock American Mathematical Society, 1997.

\bibitem[Clauset et~al.(2008)Clauset, Moore, and Newman]{clauset08}
Aaron Clauset, Cristopher Moore, and MEJ Newman.
\newblock Hierarchical structure and the prediction of missing links in
  networks.
\newblock \emph{Nature}, 453\penalty0 (7191):\penalty0 98--102, 2008.

\bibitem[Dudley(1968)]{dudley1968}
R.~M. Dudley.
\newblock Distances of probability measures and random variables.
\newblock \emph{Ann. Math. Statist.}, 39\penalty0 (5):\penalty0 1563--1572,
  1968.

\bibitem[Dunson et~al.(2021)Dunson, Wu, and Wu]{dunson2020}
David~B. Dunson, Hau-Tieng Wu, and Nan Wu.
\newblock Spectral convergence of graph laplacian and heat kernel
  reconstruction in l$\infty$ from random samples.
\newblock \emph{Applied and Computational Harmonic Analysis}, 55:\penalty0
  282--336, 2021.

\bibitem[Dyer and Frieze(1991)]{dyer1991b}
Martin Dyer and Alan Frieze.
\newblock Computing the volume of convex bodies: a case where randomness
  provably helps.
\newblock Technical Report 91-104, Carnegie Mellon University, 1991.

\bibitem[Garc\'{\i}a~Trillos and Slep\v{c}ev(2015)]{garciatrillos16b}
Nicolas Garc\'{\i}a~Trillos and Dejan Slep\v{c}ev.
\newblock On the rate of convergence of empirical measures in
  infinity-transportation distance.
\newblock \emph{Canadian Journal of Mathematics}, 67\penalty0 (6):\penalty0
  1358--1383, 2015.

\bibitem[Garc\'{\i}a~Trillos and Slep\v{c}ev(2018)]{garciatrillos18}
Nicol\'{a}s Garc\'{\i}a~Trillos and Dejan Slep\v{c}ev.
\newblock A variational approach to the consistency of spectral clustering.
\newblock \emph{Applied and Computational Harmonic Analysis}, 45\penalty0
  (2):\penalty0 239--281, 2018.

\bibitem[Garc\'{\i}a~Trillos et~al.(2016)Garc\'{\i}a~Trillos, Slep\v{c}ev,
  Von~Brecht, Laurent, and Bresson]{garciatrillos16}
Nicol\'{a}s Garc\'{\i}a~Trillos, Dejan Slep\v{c}ev, James Von~Brecht, Thomas
  Laurent, and Xavier Bresson.
\newblock Consistency of cheeger and ratio graph cuts.
\newblock \emph{Journal of Machine Learning Research}, 17\penalty0
  (1):\penalty0 6268--6313, 2016.

\bibitem[Garc\'{\i}a~Trillos et~al.(2020)Garc\'{\i}a~Trillos, Gerlach, Hein,
  and Slep{\v{c}}ev]{garciatrillos2020}
Nicol{\'a}s Garc\'{\i}a~Trillos, Moritz Gerlach, Matthias Hein, and Dejan
  Slep{\v{c}}ev.
\newblock Error estimates for spectral convergence of the graph laplacian on
  random geometric graphs toward the laplace--beltrami operator.
\newblock \emph{Foundations of Computational Mathematics}, 20\penalty0
  (4):\penalty0 827--887, 2020.

\bibitem[Garc\'{\i}a~Trillos et~al.(2021)Garc\'{\i}a~Trillos, Hoffmann, and
  Hosseini]{garciatrillos19}
Nicolas Garc\'{\i}a~Trillos, Franca Hoffmann, and Bamdad Hosseini.
\newblock Geometric structure of graph laplacian embeddings.
\newblock \emph{Journal of Machine Learning Research}, 22\penalty0
  (63):\penalty0 1--55, 2021.

\bibitem[Gharan and Trevisan(2012)]{gharan2012}
Shayan~Oveis Gharan and Luca Trevisan.
\newblock Approximating the expansion profile and almost optimal local graph
  clustering.
\newblock In \emph{Proceedings of the 2012 IEEE 53rd Annual Symposium on
  Foundations of Computer Science}, pages 187--196, 2012.

\bibitem[Gleich and Seshadhri(2012)]{gleich2012}
David~F Gleich and C~Seshadhri.
\newblock Vertex neighborhoods, low conductance cuts, and good seeds for local
  community methods.
\newblock In \emph{Proceedings of the 18th ACM SIGKDD International Conference
  on Knowledge Discovery and Data Mining}, pages 597--605, 2012.

\bibitem[Gruber(2007)]{gruber2007}
Peter~M Gruber.
\newblock \emph{Convex and discrete geometry}, volume 336.
\newblock Springer, 2007.

\bibitem[Guattery and Miller(1995)]{guattery1995}
Stephen Guattery and Gary~L Miller.
\newblock On the performance of spectral graph partitioning methods.
\newblock In \emph{Proceedings of the Sixth Annual ACM-SIAM Symposium on
  Discrete Algorithms}, volume~95, pages 233--242, 1995.

\bibitem[Hartigan(1975)]{hartigan1975}
John~A Hartigan.
\newblock \emph{Clustering algorithms}.
\newblock John Wiley \& Sons, Inc., 1975.

\bibitem[Hartigan(1981)]{hartigan1981}
John~A. Hartigan.
\newblock Consistency of single-linkage for high-density clusters.
\newblock \emph{Journal of the American Statistical Association}, 76\penalty0
  (374):\penalty0 388--394, 1981.

\bibitem[Haveliwala(2003)]{haveliwala2003}
Taher~H Haveliwala.
\newblock Topic-sensitive pagerank: A context-sensitive ranking algorithm for
  web search.
\newblock \emph{IEEE Transactions on Knowledge and Data Engineering},
  15\penalty0 (4):\penalty0 784--796, 2003.

\bibitem[Hein and B{\"u}hler(2010)]{hein2010}
Matthias Hein and Thomas B{\"u}hler.
\newblock An inverse power method for nonlinear eigenproblems with applications
  in 1-spectral clustering and sparse pca.
\newblock In \emph{Advances in Neural Information Processing Systems 23}, pages
  847--855, 2010.

\bibitem[Hoffmann et~al.(2019)Hoffmann, Hosseini, Oberai, and
  Stuart]{hoffmann2019}
Franca Hoffmann, Bamdad Hosseini, Assad~A Oberai, and Andrew~M Stuart.
\newblock Spectral analysis of weighted laplacians arising in data clustering.
\newblock \emph{arXiv preprint arXiv:1909.06389}, 2019.

\bibitem[Jiang(2017)]{jiang2017}
Heinrich Jiang.
\newblock Density level set estimation on manifolds with {DBSCAN}.
\newblock In \emph{Proceedings of the 34th International Conference on Machine
  Learning}, pages 1684--1693, 2017.

\bibitem[Koltchinskii and Gine(2000)]{koltchinskii2000}
Vladimir Koltchinskii and Evarist Gine.
\newblock Random matrix approximation of spectra of integral operators.
\newblock \emph{Bernoulli}, 6\penalty0 (1):\penalty0 113--167, 2000.

\bibitem[Korostelev and Tsybakov(1993)]{korostelev1993}
Aleksandr~P. Korostelev and Alexandre~B. Tsybakov.
\newblock \emph{Minimax theory of image reconstruction}.
\newblock Springer, 1993.

\bibitem[Kpotufe and von Luxburg(2011)]{kpotufe11}
Samory Kpotufe and Ulrike von Luxburg.
\newblock Pruning nearest neighbor cluster trees.
\newblock In \emph{Proceedings of the 28th International Conference on Machine
  Learning}, pages 225--232, 2011.

\bibitem[Lei and Rinaldo(2015)]{lei2015}
Jing Lei and Alessandro Rinaldo.
\newblock Consistency of spectral clustering in stochastic block models.
\newblock \emph{Ann. Statist.}, 43\penalty0 (1):\penalty0 215--237, 2015.

\bibitem[Leoni(2017)]{leoni2017}
Giovanni Leoni.
\newblock \emph{A First Course in Sobolev Spaces}.
\newblock American Mathematical Society, 2017.

\bibitem[Leskovec et~al.(2010)Leskovec, Lang, and Mahoney]{leskovec2010}
Jure Leskovec, Kevin~J. Lang, and Michael Mahoney.
\newblock Empirical comparison of algorithms for network community detection.
\newblock In \emph{Proceedings of the 19th International Conference on World
  Wide Web}, page 631–640, 2010.

\bibitem[Li et~al.(2020)Li, Lei, Bhattacharyya, den Berge, Sarkar, Bickel, and
  Levina]{li2018}
Tianxi Li, Lihua Lei, Sharmodeep Bhattacharyya, Koen~Van den Berge, Purnamrita
  Sarkar, Peter~J. Bickel, and Elizaveta Levina.
\newblock Hierarchical community detection by recursive partitioning.
\newblock \emph{Journal of the American Statistical Association}, pages 1--18,
  2020.

\bibitem[Little et~al.(2020)Little, Maggioni, and Murphy]{little2020}
Anna~V Little, Mauro Maggioni, and James~M Murphy.
\newblock Path-based spectral clustering: Guarantees, robustness to outliers,
  and fast algorithms.
\newblock \emph{Journal of Machine Learning Research}, 21\penalty0
  (6):\penalty0 1--66, 2020.

\bibitem[Lov{\'a}sz and Simonovits(1990)]{lovasz1990}
L{\'a}szl{\'o} Lov{\'a}sz and Mikl{\'o}s Simonovits.
\newblock The mixing rate of markov chains, an isoperimetric inequality, and
  computing the volume.
\newblock In \emph{Proceedings of the 31st Annual Symposium on Foundations of
  Computer Science}, pages 346--354, 1990.

\bibitem[Mahoney et~al.(2012)Mahoney, Orecchia, and Vishnoi]{mahoney2012}
Michael~W. Mahoney, Lorenzo Orecchia, and Nisheeth~K. Vishnoi.
\newblock A local spectral method for graphs: with applications to improving
  graph partitions and exploring data graphs locally.
\newblock \emph{Journal of Machine Learning Research}, 13:\penalty0 2339--2365,
  2012.

\bibitem[McSherry(2001)]{mcsherry2001}
Frank McSherry.
\newblock Spectral partitioning of random graphs.
\newblock In \emph{Proceedings of the 42nd IEEE Symposium on Foundations of
  Computer Science}, pages 529--537, 2001.

\bibitem[Montenegro(2002)]{montenegro2002}
Ravi Montenegro.
\newblock \emph{Faster mixing by isoperimetric inequalities}.
\newblock PhD thesis, Yale University, 2002.

\bibitem[Morris and Peres(2005)]{morris2005}
Ben Morris and Yuval Peres.
\newblock Evolving sets, mixing and heat kernel bounds.
\newblock \emph{Probability Theory and Related Fields}, 133\penalty0
  (2):\penalty0 245--266, 2005.

\bibitem[Pelletier and Pudlo(2011)]{pelletier2011}
Bruno Pelletier and Pierre Pudlo.
\newblock Operator norm convergence of spectral clustering on level sets.
\newblock \emph{Journal of Machine Learning Research}, 12\penalty0
  (12):\penalty0 385--416, 2011.

\bibitem[Polonik(1995)]{polonik1995}
Wolfgang Polonik.
\newblock Measuring mass concentrations and estimating density contour
  clusters-an excess mass approach.
\newblock \emph{Ann. Statist.}, 23\penalty0 (3):\penalty0 855--881, 1995.

\bibitem[Rigollet and Vert(2009)]{rigollet2009}
Philippe Rigollet and R{\'e}gis Vert.
\newblock Optimal rates for plug-in estimators of density level sets.
\newblock \emph{Bernoulli}, 15\penalty0 (4):\penalty0 1154--1178, 2009.

\bibitem[Rinaldo and Wasserman(2010)]{rinaldo2010}
Alessandro Rinaldo and Larry Wasserman.
\newblock Generalized density clustering.
\newblock \emph{Ann. Statist.}, 38\penalty0 (5):\penalty0 2678--2722, 2010.

\bibitem[Rohe et~al.(2011)Rohe, Chatterjee, and Yu]{rohe2011}
Karl Rohe, Sourav Chatterjee, and Bin Yu.
\newblock Spectral clustering and the high-dimensional stochastic blockmodel.
\newblock \emph{Ann. Statist.}, 39\penalty0 (4):\penalty0 1878--1915, 08 2011.

\bibitem[Rosasco et~al.(2010)Rosasco, Belkin, and Vito]{rosasco10}
Lorenzo Rosasco, Mikhail Belkin, and Ernesto~De Vito.
\newblock On learning with integral operators.
\newblock \emph{Journal of Machine Learning Research}, 11:\penalty0 905--934,
  2010.

\bibitem[Schiebinger et~al.(2015)Schiebinger, Wainwright, and
  Yu]{schiebinger2015}
Geoffrey Schiebinger, Martin~J. Wainwright, and Bin Yu.
\newblock The geometry of kernelized spectral clustering.
\newblock \emph{Ann. Statist.}, 43\penalty0 (2):\penalty0 819--846, 04 2015.

\bibitem[Shi and Malik(2000)]{shi00}
Jianbo Shi and Jitendra Malik.
\newblock Normalized cuts and image segmentation.
\newblock \emph{IEEE Transactions on Pattern Analysis and Machine
  Intelligence}, 22\penalty0 (8), 2000.

\bibitem[Shi et~al.(2009)Shi, Belkin, and Yu]{shi2009}
Tao Shi, Mikhail Belkin, and Bin Yu.
\newblock Data spectroscopy: Eigenspaces of convolution operators and
  clustering.
\newblock \emph{Ann. Statist.}, 37\penalty0 (6B):\penalty0 3960--3984, 12 2009.

\bibitem[Shi(2015)]{shi2015}
Zuoqiang Shi.
\newblock Convergence of laplacian spectra from random samples.
\newblock \emph{arXiv preprint arXiv:1507.00151}, 2015.

\bibitem[Singer and Wu(2017)]{singer2017}
Amit Singer and Hau-Tieng Wu.
\newblock Spectral convergence of the connection laplacian from random samples.
\newblock \emph{Information and Inference: A Journal of the IMA}, 6\penalty0
  (1):\penalty0 58--123, 2017.

\bibitem[Singh et~al.(2009)Singh, Scott, and Nowak]{singh2009}
Aarti Singh, Clayton Scott, and Robert Nowak.
\newblock Adaptive hausdorff estimation of density level sets.
\newblock \emph{Ann. Statist.}, 37\penalty0 (5B):\penalty0 2760--2782, 10 2009.

\bibitem[Spielman and Teng(2011)]{spielman2011}
Daniel~A Spielman and Shang-Hua Teng.
\newblock Spectral sparsification of graphs.
\newblock \emph{SIAM Journal on Computing}, 40\penalty0 (4):\penalty0
  981--1025, 2011.

\bibitem[Spielman and Teng(2013)]{spielman2013}
Daniel~A Spielman and Shang-Hua Teng.
\newblock A local clustering algorithm for massive graphs and its application
  to nearly linear time graph partitioning.
\newblock \emph{SIAM Journal on Computing}, 42\penalty0 (1):\penalty0 1--26,
  2013.

\bibitem[Spielman and Teng(2014)]{spielman2014}
Daniel~A Spielman and Shang-Hua Teng.
\newblock Nearly linear time algorithms for preconditioning and solving
  symmetric, diagonally dominant linear systems.
\newblock \emph{SIAM Journal on Matrix Analysis and Applications}, 35\penalty0
  (3):\penalty0 835--885, 2014.

\bibitem[Steinwart(2015)]{steinwart2015}
Ingo Steinwart.
\newblock Fully adaptive density-based clustering.
\newblock \emph{Ann. Statist.}, 43\penalty0 (5):\penalty0 2132--2167, 2015.

\bibitem[Steinwart et~al.(2017)Steinwart, Sriperumbudur, and
  Thomann]{steinwart2017}
Ingo Steinwart, Bharath~K Sriperumbudur, and Philipp Thomann.
\newblock Adaptive clustering using kernel density estimators.
\newblock \emph{arXiv preprint arXiv:1708.05254}, 2017.

\bibitem[Tsybakov(1997)]{tsybakov1997}
Alexandre~B Tsybakov.
\newblock On nonparametric estimation of density level sets.
\newblock \emph{Ann. Statist.}, 25\penalty0 (3):\penalty0 948--969, 1997.

\bibitem[Vempala(2005)]{vempala2005}
Santosh Vempala.
\newblock Geometric random walks: a survey.
\newblock \emph{Combinatorial and computational geometry}, 52\penalty0 (2),
  2005.

\bibitem[von Luxburg et~al.(2008)von Luxburg, Belkin, and
  Bousquet]{vonluxburg2008}
Ulrike von Luxburg, Mikhail Belkin, and Olivier Bousquet.
\newblock Consistency of spectral clustering.
\newblock \emph{Ann. Statist.}, 36\penalty0 (2):\penalty0 555--586, 04 2008.

\bibitem[Wang et~al.(2019)Wang, Lu, and Rinaldo]{wang2019}
Daren Wang, Xinyang Lu, and Alessandro Rinaldo.
\newblock Dbscan: Optimal rates for density-based cluster estimation.
\newblock \emph{Journal of Machine Learning Research}, 20\penalty0
  (170):\penalty0 1--50, 2019.

\bibitem[Wu et~al.(2012)Wu, Li, So, Wright, and fu~Chang]{wu2012}
{Xiao-Ming} Wu, Zhenguo Li, Anthony~M. So, John Wright, and Shih fu~Chang.
\newblock Learning with partially absorbing random walks.
\newblock In \emph{Advances in Neural Information Processing Systems 25}, pages
  3077--3085. 2012.

\bibitem[Yuan et~al.(2021)Yuan, Calder, and Osting]{yuan2020}
Amber Yuan, Jeff Calder, and Braxton Osting.
\newblock A continuum limit for the pagerank algorithm.
\newblock \emph{European Journal of Applied Mathematics}, pages 1--33, 2021.

\end{thebibliography}
\end{document}